\documentclass[12pt,twoside,leqno]{article}
\usepackage{amsmath}
\usepackage{amssymb}
\usepackage{amsxtra}
\usepackage{amscd}
\usepackage{amsthm}
\usepackage[mathscr]{eucal}
\usepackage{color}

\theoremstyle{plain}

\newtheorem{sbthm}[subsubsection]{Theorem}
\newtheorem{sbprop}[subsubsection]{Proposition}
\newtheorem{sbcor}[subsubsection]{Corollary}
\newtheorem{sblem}[subsubsection]{Lemma}

\theoremstyle{definition}

\newtheorem{para}[subsection]{}

\newtheorem{sbrem}[subsubsection]{Remark}

\newtheorem{sbpara}[subsubsection]{}

\newenvironment{pf}{\proof[\proofname]}{\endproof}

\usepackage{fancyhdr}
\usepackage{amssymb}
\usepackage{graphicx}

\begin{document}
\title{Classifying spaces of degenerating mixed Hodge structures, VI:
Log real analytic functions and log $C^{\infty}$ functions}
\author
{Kazuya Kato, Chikara Nakayama, Sampei Usui}

\maketitle

\renewcommand{\mathbb}{\bold}

\newcommand\Cal{\mathcal}
\newcommand\define{\newcommand}

\define\gp{\mathrm{gp}}%
\define\fs{\mathrm{fs}}%
\define\an{\mathrm{an}}%
\define\mult{\mathrm{mult}}%
\define\Ker{\mathrm{Ker}\,}%
\define\Coker{\mathrm{Coker}\,}%
\define\Hom{\mathrm{Hom}\,}%
\define\Ext{\mathrm{Ext}\,}%
\define\rank{\mathrm{rank}\,}%
\define\gr{\mathrm{gr}}%
\define\cHom{\Cal{Hom}}
\define\cExt{\Cal Ext\,}%
\define\cA{\Cal A}
\define\cB{\Cal B}
\define\cC{\Cal C}
\define\cD{\Cal D}
\define\cI{\Cal I}
\define\cO{\Cal O}
\define\cS{\Cal S}
\define\cT{\Cal T}
\define\cM{\Cal M}
\define\cG{\Cal G}
\define\cH{\Cal H}
\define\cE{\Cal E}
\define\cF{\Cal F}
\define\cN{\Cal N}
\define\fF{\frak F}
\define\Dc{\check{D}}
\define\Ec{\check{E}}
\define\Cc{\check{C}}
\define\Hc{\check{H}}

\newcommand{\N}{{\mathbb{N}}}
\newcommand{\Q}{{\mathbb{Q}}}
\newcommand{\Z}{{\mathbb{Z}}}
\newcommand{\R}{{\mathbb{R}}}
\newcommand{\C}{{\mathbb{C}}}
\newcommand{\bN}{{\mathbb{N}}}
\newcommand{\bQ}{{\mathbb{Q}}}
\newcommand{\bF}{{\mathbb{F}}}
\newcommand{\bZ}{{\mathbb{Z}}}
\newcommand{\bP}{{\mathbb{P}}}
\newcommand{\bR}{{\mathbb{R}}}
\newcommand{\bC}{{\mathbb{C}}}
\newcommand{\bbQ}{{\bar \mathbb{Q}}}
\newcommand{\ol}[1]{\overline{#1}}
\newcommand{\too}{\longrightarrow}
\newcommand{\respect}{\rightsquigarrow}
\newcommand{\compatible}{\leftrightsquigarrow}
\newcommand{\upc}[1]{\overset {\lower 0.3ex \hbox{${\;}_{\circ}$}}{#1}}
\newcommand{\Gmlog}{\bG_{m, \log}}
\newcommand{\Gm}{\bG_m}
\newcommand{\ep}{\varepsilon}
\newcommand{\Spec}{\operatorname{Spec}}
\newcommand{\val}{{\mathrm{val}}} 
\newcommand{\n}{\operatorname{naive}}
\newcommand{\bs}{\operatorname{\backslash}}
\newcommand{\Gal}{\operatorname{{Gal}}}
\newcommand{\gal}{{\rm {Gal}}({\bar \Q}/{\Q})}
\newcommand{\galp}{{\rm {Gal}}({\bar \Q}_p/{\Q}_p)}
\newcommand{\gall}{{\rm{Gal}}({\bar \Q}_\ell/\Q_\ell)}
\newcommand{\wep}{W({\bar \Q}_p/\Q_p)}
\newcommand{\wel}{W({\bar \Q}_\ell/\Q_\ell)}
\newcommand{\Ad}{{\rm{Ad}}}
\newcommand{\BS}{{\rm {BS}}}
\newcommand{\even}{\operatorname{even}}
\newcommand{\End}{{\rm {End}}}
\newcommand{\odd}{\operatorname{odd}}
\newcommand{\GL}{\operatorname{GL}}
\newcommand{\np}{\text{non-$p$}}
\newcommand{\g}{{\gamma}}
\newcommand{\G}{{\Gamma}}
\newcommand{\Lam}{{\Lambda}}
\newcommand{\La}{{\Lambda}}
\newcommand{\lam}{{\lambda}}
\newcommand{\la}{{\lambda}}
\newcommand{\uL}{{{\hat {L}}^{\rm {ur}}}}
\newcommand{\uQp}{{{\hat \Q}_p}^{\text{ur}}}
\newcommand{\sel}{\operatorname{Sel}}
\newcommand{\dt}{{\rm{Det}}}
\newcommand{\Sig}{\Sigma}
\newcommand{\fil}{{\rm{fil}}}
\newcommand{\SL}{{\rm{SL}}}
\newcommand{\spl}{{\rm{spl}}}
\newcommand{\st}{{\rm{st}}}
\newcommand{\Isom}{{\rm {Isom}}}
\newcommand{\Mor}{{\rm {Mor}}}
\newcommand{\bg}{\bar{g}}
\newcommand{\id}{{\rm {id}}}
\newcommand{\cone}{{\rm {cone}}}
\newcommand{\al}{a}
\newcommand{\ChL}{{\cal{C}}(\La)}
\newcommand{\Image}{{\rm {Image}}}
\newcommand{\toric}{{\operatorname{toric}}}
\newcommand{\torus}{{\operatorname{torus}}}
\newcommand{\Aut}{{\rm {Aut}}}
\newcommand{\Qp}{{\mathbb{Q}}_p}
\newcommand{\barQp}{{\mathbb{Q}}_p}
\newcommand{\Qpur}{{\mathbb{Q}}_p^{\rm {ur}}}
\newcommand{\Zp}{{\mathbb{Z}}_p}
\newcommand{\Zl}{{\mathbb{Z}}_l}
\newcommand{\Ql}{{\mathbb{Q}}_l}
\newcommand{\Qlur}{{\mathbb{Q}}_l^{\rm {ur}}}
\newcommand{\F}{{\mathbb{F}}}
\newcommand{\eps}{{\epsilon}}
\newcommand{\epsLa}{{\epsilon}_{\La}}
\newcommand{\epsLaVxi}{{\epsilon}_{\La}(V, \xi)}
\newcommand{\epsOLaVxi}{{\epsilon}_{0,\La}(V, \xi)}
\newcommand{\Qplin}{{\mathbb{Q}}_p(\mu_{l^{\infty}})}
\newcommand{\otimesQplin}{\otimes_{\Qp}{\mathbb{Q}}_p(\mu_{l^{\infty}})}
\newcommand{\galFl}{{\rm{Gal}}({\bar {\Bbb F}}_\ell/{\Bbb F}_\ell)}
\newcommand{\gallur}{{\rm{Gal}}({\bar \Q}_\ell/\Q_\ell^{\rm {ur}})}
\newcommand{\galFF}{{\rm {Gal}}(F_{\infty}/F)}
\newcommand{\galFv}{{\rm {Gal}}(\bar{F}_v/F_v)}
\newcommand{\galF}{{\rm {Gal}}(\bar{F}/F)}
\newcommand{\epsVxi}{{\epsilon}(V, \xi)}
\newcommand{\epsOVxi}{{\epsilon}_0(V, \xi)}
\newcommand{\plim}{\lim_
{\scriptstyle 
\longleftarrow \atop \scriptstyle n}}
\newcommand{\sig}{{\sigma}}
\newcommand{\ga}{{\gamma}}
\newcommand{\del}{{\delta}}
\newcommand{\Vss}{V^{\rm {ss}}}
\newcommand{\Bst}{B_{\rm {st}}}
\newcommand{\Dpst}{D_{\rm {pst}}}
\newcommand{\Dcrys}{D_{\rm {crys}}}
\newcommand{\DdR}{D_{\rm {dR}}}
\newcommand{\Fin}{F_{\infty}}
\newcommand{\Kla}{K_{\lambda}}
\newcommand{\Ola}{O_{\lambda}}
\newcommand{\Mla}{M_{\lambda}}
\newcommand{\Det}{{\rm{Det}}}
\newcommand{\Sym}{{\rm{Sym}}}
\newcommand{\LaSa}{{\La_{S^*}}}
\newcommand{\cX}{{\cal {X}}}
\newcommand{\MHG}{{\frak {M}}_H(G)}
\newcommand{\tauMla}{\tau(M_{\lambda})}
\newcommand{\Fvur}{{F_v^{\rm {ur}}}}
\newcommand{\Lie}{{\rm {Lie}}}
\newcommand{\cL}{{\cal {L}}}
\newcommand{\cW}{{\cal {W}}}
\newcommand{\fq}{{\frak {q}}}
\newcommand{\cont}{{\rm {cont}}}
\newcommand{\SC}{{SC}}
\newcommand{\Om}{{\Omega}}
\newcommand{\dR}{{\rm {dR}}}
\newcommand{\crys}{{\rm {crys}}}
\newcommand{\hatSig}{{\hat{\Sigma}}}
\newcommand{\rdet}{{{\rm {det}}}}
\newcommand{\ord}{{{\rm {ord}}}}
\newcommand{\BdR}{{B_{\rm {dR}}}}
\newcommand{\BdRO}{{B^0_{\rm {dR}}}}
\newcommand{\Bcrys}{{B_{\rm {crys}}}}
\newcommand{\Qw}{{\mathbb{Q}}_w}
\newcommand{\barkappa}{{\bar{\kappa}}}
\newcommand{\cP}{{\Cal {P}}}
\newcommand{\cZ}{{\Cal {Z}}}
\newcommand{\br}{{{\rm {\bold r}}}}
\newcommand{\oppLa}{{\Lambda^{\circ}}}
\newcommand{\add}{{{\rm {add}}}}
\newcommand{\red}{{{\rm {red}}}}
\newcommand{\new}{{{\rm {new}}}}
\newcommand\twoheaddownarrow{\mathrel{\rotatebox[origin=c]{90}{$\twoheadleftarrow$}}}
\renewcommand{\Re}{{\mathrm{Re}\,}}
\renewcommand{\Im}{{\mathrm{Im}\,}}
\newcommand{\triv}{{\rm {triv}}}
\newcommand{\Supp}{{\rm {Supp}}}
\newcommand{\Log}{{\rm {Log}}}
\newcommand{\spe}{{\rm {sp}}}
\newcommand{\BM}{{\rm {BM}}}
\newcommand{\pts}{{\rm {pts}}}
\newcommand {\pe}{\frak p}
\newcommand{\LMH}{{\rm {LMH}}}
\newcommand{\adm}{{\rm {adm}}}

\begin{abstract} 
  To advance our log Hodge theory, we introduce log real analytic functions and log $C^{\infty}$ functions, define how to integrate them, and 
prove the log Poincar\'e lemma. 
  We give better understandings of the degeneration of Hodge structure, including a geometric interpretation of the theory of $\SL(2)$-orbits.
\end{abstract}


\section*{Contents}

\noindent \S\ref{s:intro}. Introduction

\noindent \S\ref{s:function}. 
Log real analytic functions and log $C^{\infty}$ functions

\noindent \S\ref{s:ls}. 
Local coordinates of log smooth morphisms in log $C^{\infty}$ geometry

\noindent \S\ref{s:integration}. 
Log integration

\noindent \S\ref{s:newa}. 
Log real analytic functions and $\SL(2)$-orbits

\noindent \S\ref{s:HDI}.
Higher direct images (plan and a special case)


\renewcommand{\thefootnote}{\fnsymbol{footnote}}
\footnote[0]{MSC2020: Primary 14A21; 
Secondary 14D07, 32G20} 

\footnote[0]{Keywords: variation of Hodge structure, degenerations, period integrals, smooth functions}
\renewcommand{\thefootnote}{\arabic{footnote}}

\setcounter{section}{-1}
\section{Introduction}\label{s:intro}

\begin{para}\label{Intro1}
In this Part VI of our series of papers, we start geometric applications of the methods and materials considered  in our  book \cite{KU} and in our series of papers  Part I--Part V (\cite{KNU1}--\cite{KNU5}). 
Here we define log real analytic functions and log $C^{\infty}$ functions and  consider the following  subjects (1) and (2).

(1) The theory of integration in log geometry by using log $C^{\infty}$ functions.  Here the divergence of the integral is treated nicely. See \ref{Intro2}. This subject is discussed in Section \ref{s:integration}. 

(2) Good understandings of the degeneration of Hodge structure by using log real analytic functions, including a geometric understanding of the theory of $\SL(2)$-orbits. See \ref{Intro3}. This subject is discussed in Section \ref{s:newa}. 

We plan that by using these (1) and (2), we study in Part VII of this series of papers the higher direct images of variations of Hodge structure with degeneration. See \ref{Intro4}. 

Real analytic functions and $C^{\infty}$ functions are useful in Hodge theory. For example, the Hodge decomposition  of a variation of Hodge structure associated to its Hodge filtration can be described by using real analytic functions,  and the classical Hodge theory is constructed  by using harmonic forms which are $C^{\infty}$ functions. However, to treat   degenerating Hodge structures, usual real analytic functions and usual $C^{\infty}$ functions do not have enough powers. In this Part VI and in Part VII, we will show that the  classical theories can  revive in degeneration with log real analytic functions and log $C^{\infty}$ functions.

In this paper, we use the space of ratios defined in \cite{KNU4}. See \ref{Intro6}. 

The definitions of log real analytic functions and log $C^{\infty}$ functions and their properties are explained in Section \ref{s:function} of this paper. 

\end{para}

\begin{para}\label{Intro2}  Concerning the subject (1) of \ref{Intro1}. By using our log $C^{\infty}$ functions, in Section \ref{s:integration}, we develop the theory of integrations on log complex analytic spaces and prove log Poincar\'e lemma (Theorem \ref{LPL}). 

Our theory of integration using log geometry has the following shape. 
In degeneration, various objects diverge and become infinite. Log geometry is good to treat divergence.

For the elliptic curve $E(\tau)= \C/(\Z+\Z\tau)$ with   $\tau=x+iy$, $x,y\in \R$, 
$y>0$, we have 
$$ \int_{\gamma} dz=\tau,$$
where $z$ is the coordinate function of $\C$ and  $dz$ is regarded as a differential form on $E(\tau)$, and $\gamma$ is the loop on $E(\tau)$ given by a path from $0$ to $\tau$ on $\C$.
In degeneration, $\Im(\tau)\to \infty$, $E(\tau)$ degenerates to $E(\infty)={\bf P}^1(\C)/(0\sim \infty)$ ($0$ and $\infty$ of ${\bf P}^1(\C)$ are identified), 
and $dz$ degenerates to $(2\pi i)^{-1} d\log(t)$ on $E(\infty)$, where $t$ is the coordinate function of ${\bf P}^1(\C)$.
We can regard $E(\infty)$ as a log complex analytic space over an fs log point $S$.  In log geometry, the loop $\gamma$ degenerates to a loop $\gamma(\infty)$  in a space $E(\infty)_{[:]}^{\log}$ associated to $E(\infty)$, whose image on $E(\infty)$ is homotope to the loop $[0, \infty]/(0\sim \infty)$ on $E(\infty)$. 
Our theory of integration shows
$$\int_{\gamma(\infty)}(2\pi i)^{-1} d\log(t)= \tau,$$
 where $\tau$ is a log $C^{\infty}$ function on the base    $S$.  The integral $\int_{[0, \infty]/(0\sim \infty)} (2\pi i)^{-1}d\log(t)$ diverges and is just understood as $\infty$  in the usual geometry, but in our log geometry, the above integral along $\gamma(\infty)$  is not just $\infty$ but is a log $C^{\infty}$ function $\tau$. This   $\tau$ belongs to the ring $A_3\otimes_\R \C$ of log $C^{\infty}$ functions of the fs log point $S$ in the notation in Section \ref{s:function}. It has value $\infty$ on $S$ but is not determined by this value. 
  For the details, see Section \ref{ss:intex}.

  \end{para}

\begin{para}\label{Intro3} Concerning the subject (2) in \ref{Intro1}. 

Let $X$ be a smooth  complex analytic space, let $D$ be a divisor on $X$ with normal crossings, and let $H$ be a variation of Hodge structure  on the complement $U=X\smallsetminus D$ with logarithmic degeneration along $D$. Then an important fact  in the theory of degeneration is that a nilpotent orbit 
appears at each point $s$ of $D$. By our log Hodge theory in \cite{KU} and \cite{KNU3}, we can understand this as follows. We endow $X$ with the log structure associated to $D$. Then $H$ extends to a variation of log Hodge structure $\tilde H$ on $X$ and the nilpotent orbit which appears at $s$ is regarded as the fiber of $\tilde H$ at $s$. Now another important fact in the theory of degeneration is that an  $\SL(2)$-orbit appears at each $s$ when we fix some ratios as follows.
The nilpotent orbit which appears at $s\in D$ has the shape $(F_0, N_1, \dots, N_n)$, where $F_0$ is a fixed Hodge filtration and $N_1, \dots, N_n$ are monodromy logarithm and we consider $\exp(\sum_{j=1}^n iy_jN_j)F_0$, where $y_j\in \R$ tends to $\infty$. An $\SL(2)$-orbit appears when furthermore the ratios $y_j:y_k$ converges in the interval $[0,\infty]$ for each $j, k$.  In Section \ref{ss:HA}, we can give the following geometric understanding of this $\SL(2)$-orbit  by using the space of ratios $X_{[:]}$ with a canonical map $X_{[:]}\to X$. 
 When  $y_j\to \infty$ and  $y_j/y_k$ converges in $[0, \infty]$ for each $j,k$, a point $p$ of the space of ratios $X_{[:]}$ lying over $s$ is determined.  We show that we can understand  this $\SL(2)$-orbit  as the log real analytic fiber at $p$ of the pullback of  $\tilde H$ to $X_{[:]}$. 

Usually, the Hodge filtration and the Hodge metric diverge when a Hodge structure degenerates. But as in Section \ref{ss:HA}, by using log real analytic functions, we have a good theory of degeneration in which the Hodge filtration and the Hodge metric converge. The above understanding of $\SL(2)$-orbit bases on this good theory.

In Section \ref{ss:CKS}, by using our log real analytic functions,
we prove the log real analyticity of the important map $D^{\sharp}_{\Sig, [:]}\to D_{\SL(2)}$ (the CKS map) between extended period domains, and give a new proof of the continuity of this map   proved in \cite{KU}, \cite{KNU3}, \cite{KNU4}.

\end{para}

\begin{para}\label{Intro4} 

 In Section \ref{ss:curve}, we give an  application of (1) in \ref{Intro1} to  a special case of the higher direct image, that is, the higher direct image of  the constant Hodge structure in the case of a relative curve.  
 We plan to consider more general cases in Part VII (\cite{KNU7}) by using (1) and (2) in \ref{Intro1} and by discussing log harmonic forms. 

\end{para}

\begin{para}\label{Intro6} In this work, we use  the space of ratios $X_{[:]}$ defined in \cite{KNU4} intensively. We will show that for a variation of log Hodge structure $\tilde H$ on a log complex analytic space $X$, it is nice to  pull back $\tilde H$ to $X_{[:]}$. 
  The space $X_{[:]}$ is topologically very near to $X$ (see \ref{contrac}, \ref{contrac2}), but it is $X_{[:]}$ that has 
  good  relations with the theory of $\SL(2)$-orbits, which 
  appear locally on $X_{[:]}$, though do not necessarily appear locally on $X$. %
  Also, the map $X_{[:]}\to Y_{[:]}$, induced from a morphism $X\to Y$ of log complex analytic spaces, 
  has better   properties  than the original map $X\to Y$, which match well with log $C^{\infty}$ functions (Theorem \ref{SEisS}).  It is important to work on the space of ratios $X_{[:]}$ especially because the theory of $\SL(2)$-orbits is understood in a geometric way on $X_{[:]}$ (Section \ref{ss:HA}).

\end{para}

\begin{para}\label{old} For a smooth complex analytic space $X$ endowed with a divisor with normal crossings, some kinds of log real analytic functions were already considered in Section 5 of Cattani--Kaplan--Schmid \cite{CKS1}, and some kinds 
of log $C^{\infty}$ functions were already used by several authors (for example, \cite{HP}, \cite{H}, \cite{HZ1}, \cite{HZ2}, \cite{KMN}). 
  Compared with these previous works, our theory works on singular spaces like ${\bf P}^1(\C)/(0\sim \infty)$ in \ref{Intro2}. 
\end{para}

\medskip

  K.\ Kato was 
partially supported by NFS grants DMS 1303421, DMS 1601861, and DMS 2001182.
C.\ Nakayama was 
partially supported by JSPS Grants-in-Aid for Scientific Research (B) 23340008, (C) 16K05093, and (C) 21K03199.
S.\ Usui was 
partially supported by JSPS Grants-in-Aid for Scientific Research (B) 23340008, (C) 17K05200, and  (C) 22K03247.

\bigskip

\noindent{\bf Convention.}
  {\it Concerning topology}: A continuous map $f: X \to Y$ of topological spaces is said to be {\it proper} if it is universally closed and separated. 

  A topological space $X$ is said to be {\it compact} if the map from $X$ to a point is proper. 
  That is, we impose the Hausdorffness in the definition of the compactness.

  We impose the Hausdorffness also in the definitions of the local compactness and the paracompactness. 

  {\it Concerning commutative monoids}: Let $P$ be an fs monoid.
  For a multiplicative set $F$ of $P$, we denote $F^{-1}P$ by $P_F$.
  We denote $P/P^{\times}$ by $\overline P$. 
  This $P\mapsto \overline P$ %
  gives a functor from the category of fs monoids to that of sharp fs monoids. 
  Generally, our convention on monoids is compatible with that of \cite{O}.

  {\it Concerning the inverse image of a sheaf}: The inverse image of a sheaf by a continuous map of topological spaces is often denoted by the same letter. 
  In particular, the inverse image of the log structure $M_X$ on an fs log analytic space $X$ to the associated space of ratios $X_{[:]}$ (cf.\ \ref{ss:srsa}) is denoted by the same symbol $M_X$. 
  Hence for a point $p$ of $X_{[:]}$, $M_{X,p}$ denotes the stalk of this inverse image at $p$, which is isomorphic to the stalk of the original $M_X$ on $X$ at the image of $p$. 

  {\it Concerning the direct image with proper supports}: 
  Recall that for a sheaf $F$ of abelian groups on a topological space $T$ and for a section $s\in F(T)$, the support of $s$ is the closed set of $T$ consisting of all points $t$ such that the image of $s$ in the stalk $F_t$ is not zero.
  Let $f:X \to Y$ be a continuous map of topological spaces, and let $F$ be a sheaf of abelian group on $X$. 
  Then $f_!F$ is a subsheaf of $f_*F$ whose set of sections over an open set $V \subset Y$ consists of the sections whose supports are proper over $V$.

\section{Log real analytic functions and log $C^{\infty}$ functions}
\label{s:function}

In  Section \ref{s:function}, we define 
 sheaves of log real analytic functions and sheaves of log $C^{\infty}$ functions, and study their properties.

\subsection{Introduction to Section \ref{s:function}}\label{ss:IntroS1}

\begin{sbpara}\label{remCKS} In this Section \ref{ss:IntroS1}, as an introduction to Section \ref{s:function}, we describe roughly how the sheaves of log real analytic functions and the sheaves of log $C^{\infty}$ functions  appear in this Section \ref{s:function}. 

Many sheaves appear, and so, to have some idea, we first give an introduction to the three especially important sheaves among them, sheaves of log real analytic functions $A^{\an}_{(2)}$, $A^{\an}_{(2+)}$ and a sheaf of  log $C^{\infty}$ functions $A_2$.

Let $X$ be a smooth complex analytic space endowed with a divisor $D$ with normal crossings. 
Let $x\in X$ and fix a family  $(q_j)_{1\leq j\leq n+m}$ of 
local coordinate functions of $X$ at $x$ such that $q_j$ for $1\leq j\leq n$ define the  components of $D$ at $x$. 
Then  in (5.17) of \cite{CKS1}, Cattani, Kaplan and Schmid define a  ring  ${\cal R}_K^b$ 
which is useful to treat  the degeneration of Hodge metric at $x$ of a variation of Hodge structure.

Our sheaf $A^{\an}_{(2)}$  for the log structure of $X$ associated to $D$ is a sheaf version of ${\cal R}^b_K$ and has a similar  definition to that of ${\cal R}^b_K$. 
It is a sheaf on the space of ratios $X_{[:]}$ over $X$. For this $X$ with the above log structure, $A^{\an}_{(2+)}=A^{\an}_{(2)}$.  
Though ${\cal R}^b_K$ depends on the choice of $(q_j)_{1\leq j\leq n+m}$ as above, these sheaves  do not depend on the choice, and is defined  more generally on $X_{[:]}$ for a general object $X$ of the category $\cB(\log)$ in \ref{ABlog} below, and $A^{\an}_{(2+)}$ is useful to treat the degeneration of Hodge metric on such general $X$.

The relation between this ring ${\cal R}_K^b$ and these sheaves is explained in Section \ref{ss:CKS5}. Though the space $X_{[:]}$ is not used in \cite{CKS1}, some atmosphere of $X_{[:]}$ appears in \cite{CKS1} as is explained in Section \ref{ss:CKS5}. 

At $x$, 
if we write $q_j= \exp(2\pi i(x_j+iy_j))$ with $x_j$, $y_j$ being real, the  sheaf  $A_2$ on $X_{[:]}$ around  $x$ is, roughly speaking,  the sheaf of $C^{\infty}$  functions on $X\smallsetminus D\subset X_{[:]}$ whose all iterated derivatives by  $y_j \partial/\partial x_j$ ($1\leq j\leq n$), $y_j\partial/\partial y_j$ ($1\leq j\leq n$), $\partial/\partial \Re(q_{n+j})$ ($1\leq j\leq m$), $\partial/\partial \Im(q_{n+j})$ ($1\leq j\leq m$), extend to $X$ as continuous functions.

\end{sbpara}
\begin{sbpara}\label{ABlog}

Let the categories $\cA$,  $\cB$, $\cA(\log)$, $\cB(\log)$ be as in \cite{KU} Section 3.2. That is, $\cA$ is the category of complex analytic spaces (in the sense of Grothendieck, that is, we allow non-zero nilpotent sections in the structure sheaf), $\cB$ is the category of locally ringed spaces $X$ over $\C$ such that locally on $X$, there is an object $Z$ of $\cA$ and a morphism $i:X\to Z$ of locally ringed spaces over $\C$ satisfying the conditions (i)--(iii) below, and $\cA(\log)$ (resp.\ $\cB(\log)$) is the category of objects of $\cA$ (resp.\ $\cB$) endowed with an fs log structure. (i) The map $i: X\to Z$ is injective. (ii) The topology of $X$ is the strong topology in $Z$ (\cite{KU} Section 3.1). (iii) The map $i^{-1}(\cO_Z)\to \cO_X$ is an isomorphism. 

An object of $\cA(\log)$ is called an {\it fs log analytic space}.
The $X$ in \ref{remCKS} with the log structure associated to $D$ is an example of a log smooth fs log analytic space. %

The category $\cB(\log)$ is important because  the toroidal partial compactifications of period domains constructed in \cite{KU} and \cite{KNU3}  belong to $\cB(\log)$ but not necessarily to $\cA(\log)$. 

\end{sbpara}

\begin{sbpara}\label{sheaves}
Let $X$ be an object of the category $\cB(\log)$.

On the space of ratios $X_{[:]}$ of $X$ (cf.\ \ref{ss:srsa}), we define the sheaves of rings in the following commutative diagram: $$\begin{matrix}    
A_1^{\an} & \subset & A_2^{\an}& \subset &A_{2+}^{\an} &\subset & A_3^{\an} \\ 
&& \downarrow&&\downarrow && \downarrow\\
\downarrow && A_{(2)}^{\an}& \subset & A^{\an}_{(2+)} &\subset &A^{\an}_{(3)} \\
&& \downarrow && & &\downarrow \\
A_1 &\subset &A_2 & \subset & A_{2*} &\subset &A_3.\end{matrix}$$

On the space $X^{\log}_{[:]}=X^{\log} \times_X X_{[:]}$,  we define the  sheaves of rings 
$$\tau_X^{-1}(A^{\an}_k) \subset A_k^{\an,\log} \subset A_k^{\an,\log*},\quad \tau_X^{-1}(A_k) \subset A_k^{\log} \subset A_k^{\log*},$$
where $k=1,2$, etc. and $\tau_X$ is the canonical map $X^{\log}_{[:]}\to X_{[:]}$. We have a commutative diagram as above for $A^{\an,\log}_k$ and $A_k^{\log}$, and also a commutative diagram as above for  $A^{\an,\log*}_k$ and $A_k^{\log*}$.

 The sheaves $A_k^{\an}$ (for various $k$), $A_k^{\an,\log}$  and $A_k^{\an,\log*}$  are called {\it the sheaves of log real analytic functions}. The sheaves $A_k$, $A_k^{\log}$  and $A_k^{\log*}$  are called {\it the sheaves of log $C^{\infty}$ functions}.
 
 We denote the sheaves $A^{\an}_k$, $A^{\log}_k$, etc. also as $A^{\an}_{k,X}$, $A^{\log}_{k,X}$, etc.\  (not as $A^{\an}_{k, X_{[:]}}$, $A^{\log}_{k, X_{[:]}^{\log}}$, etc.).

 If the log structure of $X$ is trivial, then $X=X_{[:]}=X^{\log}=X^{\log}_{[:]}$. If furthermore $X$ is a smooth complex analytic space, then 
all the  above sheaves of log real analytic functions coincide with the sheaf of usual real analytic functions on $X$, and all the above sheaves of log $C^{\infty}$ functions coincide with the sheaf of usual $\R$-valued 
$C^{\infty}$-functions on $X$.

\end{sbpara}

\begin{sbpara}
In the case where $X$ is a log smooth fs log  analytic space, if $U=X_{\triv}$ denotes the dense open set of $X$ consisting of all points at which the log structure is trivial, $A^{\an}_k$ (resp.\ $A_k^{\an,\log}$ and 
$A_k^{\an,\log*}$) are subsheaves of the direct  image on $X_{[:]}$ (resp.\ $X_{[:]}^{\log}$) of the sheaf of real analytic functions on $U$, and $A_k$ (resp.\ $A_k^{\log}$)  are subsheaves of the direct  image on $X_{[:]}$ (resp.\ $X_{[:]}^{\log}$) of the sheaf of $\R$-valued $C^{\infty}$  functions on $U$.

\end{sbpara}
\begin{sbpara}

For a section $f$ of  $A_k^{\an}$ for $k=1,2, 2+, (2), (2+)$ or of  $A_k$ for $k=1,2$ on an open set $U$ of $X_{[:]}$, the value of $f$ at each point of $U$ is defined. For a section $f$ of  $A_k^{\an,\log}$ or $A_k^{\an,\log*}$ for $k=1, 2, 2+, (2), (2+)$ or of  $A_k^{\log}$ or $A_k^{\log*}$ for $k=1,2$
on an open set $U$ of $X^{\log}_{[:]}$, the value of $f$ at each point of $U$ is defined. However, sections of the other sheaves in the above can be unbounded at the boundary.

\end{sbpara}

\begin{sbpara}\label{imp} Among the above sheaves, the sheaves $A^{\an}_{(2)}$, $A^{\an}_{(2+)}$, $A_2$ and $A_2^{\log}$ are especially important (other sheaves play auxiliary roles). The sheaf $A^{\an}_{(2+)}$ is  especially important  to treat degeneration of Hodge structure nicely as in Section \ref{s:newa}. The sheaves $A_2$ and $A_2^{\log}$ are especially important to have integration and iterated derivatives in degeneration as in Section \ref{s:integration}. The sheaf $A_{(2)}^{\an}$ is important to connect $A^{\an}_{(2+)}$ with $A_2$ and $A_2^{\log}$.

\end{sbpara}

\begin{sbpara}\label{shdelta}  As an example, let $X=\Delta=\{q\in \C\;|\; |q|<1\}$ endowed with the log structure generated by the coordinate function $q$. Then the  sheaves in \ref{sheaves} are as follows. On $\Delta\smallsetminus \{0\}$, write  $q=e^{2\pi i (x+iy)}$ with being $x,y$ real.

Then on $\Delta=\Delta_{[:]}$, $A_1^{\an}$ is the sheaf of real analytic functions of $\Re(q), \Im(q), y^{-1/2}$. We have $A^{\an}_2=A^{\an}_{2+}$ and $A^{\an}_{(2)}=A^{\an}_{(2+)}$. As sheaves of rings,  $A_2^{\an}$ (resp.\ $A_3^{\an}$) is  generated over $A^{\an}_1$ algebraically by $y^e\Re(q)$ and $y^e\Im(q)$ for integers $e\geq 0$ (resp.\ by $y$), 
 and $A^{\an}_{(2)}$ (resp.\ $A^{\an}_{(3)}$) is generated over $A_2^{\an}$ (resp.\ $A^{\an}_3$) algebraically by $1/f$ for $f\in A^{\an}_2$ such that the value of $f$ is non-zero at every point. For example, we have the section $y$ of $A_3^{\an}$ which does not  belong to $A^{\an}_{(2+)}$, and we have the section  $(1+\Re(q)y)^{-1}$ of $A^{\an}_{(2)}$ which does not belong to $A^{\an}_3$.

On $\Delta$, $A_1$ is the sheaf of $C^{\infty}$ functions of $\Re(q), \Im(q), y^{-1/2}$. $A_2$ is the sheaf of $C^{\infty}$ functions on $\Delta\smallsetminus \{0\}$ whose all iterated derivatives under $y\partial/\partial x$ and $y\partial /\partial y$ are extended over the origin as continuous functions. 
$A_{2*}$ (resp.\ $A_3$)  is the sheaf of rings generated over $A_2$ algebraically by $\log(y)$ (resp.\ by $y$). 

The space $\Delta_{[:]}^{\log}=\Delta^{\log}$ is understood as $\{|q|\in \R\;|\;0\leq |q|< 1\}\times \R/\Z$ with the canonical map $\Delta_{[:]}^{\log}\to \Delta\;;\; (|q|,x)\mapsto q= |q|e^{2\pi ix}$ which is an isomorphism outside $0\in \Delta$. On  $\Delta_{[:]}^{\log}$,  as sheaves of rings, $A_k^{\an,\log}$  (resp.\  $A_k^{\log}$) is generated over the inverse image of $A^{\an}_k$ (resp.\ $A_k$) algebraically by $x/y$, and $A_k^{\an,\log*}$  (resp.\ $A^{\log*}_k$) is generated over  the inverse image of $A_k^{\an}$ (resp.\ $A_k$) algebraically by $x$. 
\end{sbpara}

\begin{sbpara}\label{1ma1}

  The following (1)--(3) are main results about properties of log real analytic functions, and are proved in Section \ref{ss:Aanp}. 

\medskip

(1)  (Theorem \ref{A2an}) If $X$ is a log smooth fs log analytic space, we have 
\begin{alignat*}{2}
A_2^{\an} &= A_{2+}^{\an}&&= A_3^{\an}\cap \{\text{continuous functions}\},\text{ and} \\
A_{(2)}^{\an}&= A_{(2+)}^{\an}&&= A_{(3)}^{\an}\cap \{\text{continuous functions}\}.
\end{alignat*}

(2) (Theorem \ref{sharp=})  If $X$ is an ideally log smooth (\ref{ILS}) fs log analytic space and is reduced,  then for $R=A^{\an}_2$ (resp.\ $R=A^{\an}_{(2)}$) and $R_+= A^{\an}_{2+}$ (resp.\ $R_+= A^{\an}_{(2+)}$) and for an $R$-module $E$ on $X_{[:]}$ locally free of finite rank, we have
 $$\sig_*(E) \overset{\cong}\to \sig_*(E\otimes_R R_+),$$ 
where $\sig$ is the canonical map $X_{[:]}\to X$ (\ref{tau}). 

(3) (Theorem \ref{A*cont}) Let $X$ be an object of $\cB(\log)$ and let $U$ be an open set of $X_{[:]}$. Then an element $f$ of $A^{\an}_{(2+)}(U)$ gives a continuous map $U\to \R\;;\;p\mapsto f(p)$. 
 
\medskip
 
The definition of ideally log smoothness in (2) is as in \cite{IKN} Definition (1.5) as reviewed in \ref{ILS}. The word ``reduced'' in (2) means that the structure sheaf $\cO_X$ has no non-zero nilpotent local sections.

An important example to which (2) can be applied is as follows. 
  If  $X\to Y$ is a log smooth saturated morphism of fs log analytic spaces, all fibers are ideally log smooth and reduced.  
\end{sbpara}

\begin{sbpara}\label{1ma2} The following is the main result 
about properties of log $C^{\infty}$ functions, and is 
proved in Section \ref{ss:soft}.  

Let $X$ be an fs log analytic space and assume that the underlying complex analytic space of $X$ is Hausdorff and satisfies the second axiom of countability. Then 
for $k=1,2,2*, 3$ and for any open set $U$ of $X_{[:]}$, the sheaf $A_k$ on $U$  is soft (Theorem \ref{softhm}).

\end{sbpara}
 
\subsection{The space of ratios and the space of arguments}
\label{ss:srsa}

Because the space $X_{[:]}$ of ratios (\cite{KNU4} 4.2) and the space $X^{\log}$ of arguments (\cite{KtNk}, \cite{KU} 2.2) associated to an object $X$ of $\cB(\log)$ are very  important in this paper, we briefly review them.

\begin{sbpara}
\label{tau}
These $X_{[:]}$ and $X^{\log}$ are topological spaces endowed with  proper surjective maps $X_{[:]}\to X$ and $X^{\log}\to X$, which are denoted by $\sigma=\sigma_X$ and $\tau=\tau_X$, respectively.

If the log structure of $X$ is trivial, then $X_{[:]}=X$ and $X^{\log}=X$. 

\end{sbpara}

\begin{sbpara}\label{chratio} The following (1)--(3) characterize $X_{[:]}$ for all $X\in \cB(\log)$. 

(1) If $X$ is log smooth over $\bC$, let $U$ be the largest open subset of $X$ on which the log structure is trivial. 
  Then for $f, g\in M_X(X)$ such that $|f|<1, |g|<1$, the ratio $\log|f|/\log|g|: U\to \R_{>0}$ extends uniquely to a continuous map $X_{[:]}\to [0, \infty]$. 

(2) In (1), if  furthermore $\cS\to M_X$ is a chart with $\cS$ being an fs monoid  and if $f_1, \dots, f_m\in \cS$  generate $\cS$ and  $|f_j|<1$ for $1\leq j\leq m$, 
then $X_{[:]}$ is the closure in $X\times [0,\infty]^{m^2}$ 
of the set $\{(p,\varphi(p))\;|\; p\in U\} \subset U \times \R_{>0}^{m^2}$ with $\varphi(p)= (\log|f_j(p)|/\log|f_k(p)|)_{1\leq j\leq m, 1\leq k\leq m}$.

(3) If $Y\to X$ is strict, then $Y_{[:]}= X_{[:]} \times_X Y$. 

\end{sbpara}

\begin{sbpara}\label{Delta2:}

If $X=\Delta^2$ with the log structure generated by the coordinate functions $q_1, q_2$, then the map $X_{[:]}\to X$ is  a homeomorphism outside $(0,0)\in X$, and the fiber of $(0,0)$ in $X_{[:]}$ is homeomorphic to $[0, \infty]$ by $\log|q_2|/\log|q_1|$. 

\end{sbpara}

\begin{sbpara}\label{charg} The following (1)--(3) characterize $X^{\log}$ for all $X\in \cB(\log)$. 

(1) If $X$ is log smooth over $\bC$, let $U$ be the largest open subset of $X$ on which the log structure is trivial. 
  Then for $f\in M^{\gp}_X(X)$ such that $|f|<1$, the argument $\arg(f): U\to \R/2\pi \Z$ extends uniquely to a continuous map $X^{\log}\to \R/2\pi \Z$. 

(2) In (1), if  furthermore $\cS\to M_X$ is a chart with $\cS$ being an fs monoid and if $f_1, \dots, f_m\in \cS^{\gp}$  generate $\cS^{\gp}$, 
then $X^{\log}$ is the closure in $X\times (\R/\Z)^m$
of the set $\{(p, \varphi(p))\;|\; p\in U\}\subset U \times (\R/\Z)^m$ with $\varphi(p)=((2\pi)^{-1}\arg(f_j(p)))_{1\leq j\leq m}$.  

(3) If $Y\to X$ is strict, then $Y^{\log}= X^{\log} \times_X Y$.

\end{sbpara}

\begin{sbpara}

Let $X=\Delta^n$ with coordinate functions $q_1,\dots, q_n$ and with the log structure given by $q_1, \dots, q_r$ (we fix an $r$ such that $1\leq r\leq n$). Then we have a homeomorphism $X^{\log} \cong \R_{\geq 0}^r\times (\R/\Z)^r \times \Delta^{n-r}$ given by $((|q_j|)_{1\leq j\leq r}, ((2\pi)^{-1}\arg(q_j))_{1\leq j\leq r}, (q_j)_{r+1\leq j\leq n})$. 

\end{sbpara}

\begin{sbpara} We define $X^{\log}_{[:]}=X_{[:]} \times_X X^{\log}$. 
  Denote the projection $X^{\log}_{[:]} \to X_{[:]}$ by $\tau=\tau_X$.
  (It is the same notation as the commonly used $\tau: X^{\log} \to X$.) 
Denote the projection $X^{\log}_{[:]} \to X^{\log}$ by $\sig=\sig_X$.
\end{sbpara}
 
 \begin{sbpara}\label{fgcont}  For $X\in\cB(\log)$, for an open set $U$ of $X_{[:]}$  and $f,g\in M_X(U)$ such that $|f|<1, |g|<1$, we have a continuous map $\log|f|/\log|g|: U \to [0, \infty]$ which is as in \ref{chratio}  in the case where $X$ is log smooth and which has the value $r(\overline f,\overline g)$ of  the ratio $r$ (\cite{KNU4} 4.2.1)   at $p\in X_{[:]}$  in the case where at least one of $f,g$ does not belong $\cO^\times_p$ at $p$. 
  Here $\overline f$, $\overline g$ are the images of $f$, $g$ in $(M_X/\cO_X^{\times})_p$, respectively. 
For an open set $U$ of $X^{\log}$ or of $X_{[:]}^{\log}$ and for $f\in M_X^{\gp}(U)$, we have a continuous map $\arg(f): U \to \R/2\pi \Z$  which is as explained in \ref{charg} in the case where $X$ is log smooth.

 \end{sbpara}

\begin{sbpara}\label{contrac}  All fibers of $X_{[:]}\to X$ are weakly contractible and locally contractible. 
This is proved in \cite{KNUc}. (Though \cite{KNUc} is of the style of an announcement paper, it contains a complete 
proof of this fact.) 
  Hence, all fibers of $X_{[:]}^{\log}\to X^{\log}$ are also weakly contractible and locally contractible. 
  
  In general,  if $a: S'\to S$ is a proper continuous map whose all fibers are weakly contractible and locally contractible (as in the  case $S=X$ and $S'=X_{[:]}$  and also as  in the case $S=X^{\log}$ and $S'=X^{\log}_{[:]}$), we have the  following (1) and (2).
    
  (1) The map $H^m(U, \cF)\to H^m(a^{-1}(U), a^{-1}(\cF))$ is bijective for every open set $U$ of $S$ in each of the following cases (i)--(iii).  (i) $\cF$ is a locally constant sheaf of abelian groups, and $m$ is any integer. (ii) $\cF$ is a locally constant sheaf of groups and $m=1$. (iii) $\cF$ is a locally constant sheaf of sets and $m=0$. 
  
  (2) If $R$ is a ring, the category of local systems of finitely generated $R$-modules on $S$ is equivalent to that on $S'$ via the pullback. 
  
  (1) follows from the proper base change theorem. (2) is proved as follows. For a locally constant sheaf $\cF$ on $S$, we have $\cF\overset{\cong}\to a_*a^{-1}(\cF)$ by (1). Next let $\cF$ be a local system of finitely generated $R$-modules on $S'$ and we prove that $a_*\cF$ is locally constant and $\cF=a^{-1}a_*\cF$. 
  Let $s \in S$ and we work around $s$. 
  The inverse image of $\cF$ on the fiber of $s$ is a constant sheaf associated to a finitely generated $R$-module $F$.  By the properness of $a$, we may assume that $\cF$ is isomorphic to $F$ locally on $S'$. 
  Let $G=\Aut(F)$ and consider the sheaf $\cG=\cA ut(F)$ of automorphisms of $F$. Since $F$ is finitely generated, $\cG$ is the constant sheaf associated to  $G$, and 
the sheaf of isomorphisms between $F$ and $\cF$ is a $G$-torsor on $S'$.  By the vanishing of $R^1a_*G$, after shrinking $S$ if necessary, we have an isomorphism $\cF\cong F$ globally on $S'$. 
  Hence $a_*\cF$ is locally constant and  $\cF=a^{-1}a_*\cF$. 

 \end{sbpara}
 
\begin{sbpara}\label{contrac2}
 
In log Hodge theory, for a subring $R$ of $\R$, the $R$-structure $H_R$ of an $R$-log Hodge structure $H$ is a local system of finitely generated $R$-modules   on $X^{\log}$, and we have the Hodge bundle $H_{\cO}$ of $H$ which is a vector bundle on $X$ (see Section \ref{ss:LMH}). By the above (1) and (2), it is harmless to consider $H_R$  on $X_{[:]}^{\log}$ and consider $H_{\cO}$  on $X_{[:]}$ by taking their inverse image sheaves. 
  
  \end{sbpara}

\begin{sbpara} The importance of the space of ratios in the theory of degeneration of Hodge structure is seen, for example, as follows. Let $X=\Delta^2$ with the log structure generated by the coordinate functions $q_1, q_2$. In the study of degenerations of variations of Hodge structure on $(\Delta^*)^2$ at $(0,0)$ (here $\Delta^*:=\Delta\smallsetminus \{0\}$), we often consider the behavior of

(1) $\exp(iy_1N_1+iy_2N_2)F_0\quad  (y_j= -(2\pi)^{-1}\log|q_j|$ for $j=1,2$),  

\noindent 
where $N_j$ is the logarithm of the action of the $j$-th generator of $\pi_1((\Delta^*)^2)=\Z^2$ and $F_0$ is a fixed Hodge filtration such that $(N_1, N_2, F_0)$ generates a nilpotent orbit. It is important that for $a\in [0, \infty]$, we have an $\SL(2)$-orbit as the ``limit'' of (1) for $y_j \to \infty$, $y_1:y_2\to a$. As in Theorem \ref{ourSL2}, we understand that this SL(2)-orbit appears at the point $p\in X_{[:]}$ lying over $(0,0)\in X$ corresponding to $a$  (\ref{Delta2:}). That $y_j$ tends to $\infty$ and $y_1:y_2$ tends to $a$ is equal to that  $(q_1, q_2)\in (\Delta^*)^2$ tends to $p$.

\end{sbpara}

\subsection{The sheaves $A_1^{\an}$ and $A_1$ (log smooth case)}
\label{ss:A1ls}

 \begin{sbpara}
\label{A_1top}
 Assume that $X$ is a log smooth fs log  analytic space. Let $\cO_X$ be the sheaf of holomorphic functions on $X$, 
let $M_X$ be the log structure of $X$, and we denote their inverse image sheaves on $X_{[:]}$ and on $X^{\log}_{[:]}$ simply by $\cO_X$ and $M_X$, respectively. 

Let $U=X_{\triv}$ be the  dense open set of $X$ consisting of all points at which the log structure is trivial, let $j:U\to X_{[:]}$  be the inclusion map,
 and let $A_U^{\an}$ (resp.\ $A_U$) be the sheaf of real analytic (resp.\ $\R$-valued $C^{\infty}$) functions on $U$. 
We define the sheaf $A_1^{\an}$ (resp.\ $A_1$) on $X_{[:]}$  as a sheaf of subrings of $j_*(A^{\an}_U)$ (resp.\ $j_*(A_U))$.

 \end{sbpara}
 
\begin{sbpara}
\label{A_1} 
Let $A_1^{\an}$ (resp.\ $A_1$) on $X_{[:]}$  be the sheaf 
of real analytic (resp.\ $\R$-valued $C^{\infty}$) functions of
$$\left(\frac{\log|f|}{\log|g|}\right)^{1/2}$$
 for $f, g\in M_X$ such that $|f|<1$, $|g|<1$ and such that $\log|f|/\log|g|$ do not have value $\infty$. 
 
 More precisely, for an open set $V$ of $X_{[:]}$, $A_1^{\an}(V)$ (resp.\ $A_1(V))$ is the subring of $j_*(A_U^{\an})(V)=A_U^{\an}(V\cap U)$ (resp.\ $j_*(A_U)(V)=A_U(V\cap U)$) consisting of all elements $\varphi$ such that for some open covering $(V_{\lam})_{\lam}$ of $V$, the restriction $\varphi_{\la}$ of $\varphi$ to each $V_{\lam}$ is presented in the following way:

 \medskip
  
 There are   $f_j, g_j\in M(V_{\la}) $ ($1\leq j\leq n$ for some $n$)  such that $|f_j|<1$, $|g_j|<1$, such that  $\log|f_j|/\log|g_j|$ do not have the value $\infty$ on $(V_{\la})_{[:]}$ (\ref{chratio} (1)), and such that if  $a$ denotes the induced map $((\log|f_j|/\log|g_j|)^{1/2})_{1\leq j\leq n}: V_{\la}\to \R^n$, there are an open set  $W$ of $\R^n$ which contains the image of $V_{\la}$ and a real analytic (resp.\ $\R$-valued $C^{\infty}$) function $\psi$ on $W$ such that  $\varphi_{\la}= \psi\circ a$. 
 
 \medskip

 Of course, we have $A^{\an}_1\subset A_1$. 
 
 The reason why we use the square root $(\log|f|/\log|g|)^{1/2}$ in the above definition is explained in Remark \ref{A1rem}.
 \end{sbpara}
 
 \begin{sbpara}\label{1.3.3} 
 We define  the sheaf of rings $A_0^{\an}$ (resp.\ $A_0$) on $X$ to be the sheaf of $\R$-valued functions on $X$  which are real analytic (resp.\  $C^{\infty}$) functions of $\Re(f)$ for $f\in \cO_X$. We can say ``$\Re(f)$ and $\Im(f)$'' here instead of ``$\Re(f)$'', for $\Im(f)= -\Re(if)$.

We define the sheaf of rings  $A_0^{\an}$ (resp.\ $A_0$) on $X_{[:]}$ to be the inverse image sheaf of the sheaf $A_0^{\an}$ (resp.\ $A_0$) on $X$. 
 
 (As in these cases of $A_0^{\an}$ and $A_0$ and as in the case of $M_X$ in \ref{A_1top}, in this paper, we often denote the inverse image of a sheaf $\cF$ by the same notation $\cF$. 
  Cf.\ Convention after Introduction.)

 An understanding of real analytic functions on objects of  $\cB(\log)$ is given in \ref{XX'}.

 \end{sbpara}
 
 \begin{sbprop}\label{01}  $$A_0^{\an}\subset A_1^{\an},\quad A_0\subset A_1.$$
  \end{sbprop}

 \begin{pf} It is sufficient to prove that for $f\in \cO_X$, $\Re(f)$ belongs to $A_1^{\an}$. Let $p\in X_{[:]}$ with the image $s$ in $X$, and take $c\in \R$ such that $c<0$ and $\Re(f)(s)+c <0$.  Let $g=\exp(f+c)$ and $h= \exp(c)$. Then $|g|<1$ for some neighborhood of $s$ and $|h|<1$. There we have $\log|g|/\log|h|= (\Re(f)+c)/c=c^{-1}\Re(f)+1$, and hence $c^{-1}\Re(f)+1\in A_1^{\an}$ there.  Hence $\Re(f)\in A_1^{\an}$ there. \end{pf}

 \begin{sbcor}
\label{c:cO}
We have $\cO_X\subset A^{\an}_{1,\C}:= A_1^{\an}\otimes_{\R}\C$.  
 
 \end{sbcor}

\begin{pf}
 Since $\cO_X\subset A^{\an}_{0,\C}:= A_0^{\an}\otimes_\R \C$, this follows from Proposition \ref{01}.
\end{pf}

 \begin{sbpara}\label{Ex1}  {\it Example.} 
 Let $X=\Delta$ with the log structure at the origin.
 
  Let $q$ be the coordinate function of $\Delta$ and write $q=e^{2\pi i(x+iy)}$ on $\Delta\smallsetminus \{0\}$ with $x,y$ real. Then $A_1^{\an}$ (resp.\ $A_1$) is the sheaf of real analytic (resp.\ $\R$-valued $C^{\infty}$) functions of $y^{-1/2}$, $\Re(q)$, and $\Im(q)$. 

For the proof, see Proposition \ref{Ex3} below. 
 \end{sbpara}

 \begin{sbpara}\label{Ex2} {\it Example.} 
 Let $X=\Delta^n$ with the coordinate functions $q_1, \dots, q_n$ with the log structure given by $q_1, \dots, q_r$ ($1\leq r\leq n$). Let  $x_j$ and $y_j$ ($1\leq j\leq r$) be as in \ref{Ex1} ($q$ there  is replaced by  $q_j$). By the {\it standard open set} of $X_{[:]}$, we  mean the set of all $p\in X_{[:]}$ such that $y_{j+1}/y_j$ do not have value $\infty$ for $1\leq j\leq r-1$. Then on the standard open set of $X_{[:]}$, $A_1^{\an}$ (resp.\ $A_1$) is the sheaf of real analytic (resp.\ $\R$-valued $C^{\infty}$) functions of $(y_j/y_{j+1})^{-1/2}$ ($1\leq j\leq r-1$), $y_r^{-1/2}$, and $\Re(q_j)$, $\Im(q_j)$ ($1\leq j\leq n$). 
  
 For the proof, see Proposition \ref{Ex3} below. 
 \end{sbpara}
 
  \begin{sbrem}\label{A1rem}  A reason  why we put the square root $(\log|f|/\log|g|)^{1/2}$ of $\log|f|/\log|g|$ in the definition of $A_1^{\an}$ (\ref{A_1}), not only $\log|f|/\log|g|$, is that 
  these square roots (whose examples are $y^{-1/2}$ in \ref{Ex1} and $(y_j/y_{j+1})^{-1/2}$ in \ref{Ex2})  appear in Hodge theory, for example, as  in Proposition \ref{HA}. 
  As in there, it is often convenient to think that  $y^{m/2}$ on $\Delta$ for an integer $m$ has the Hodge theoretic weight $m$ for the relative monodromy filtration in degeneration.

 \end{sbrem}

\begin{sbpara}
\label{m,f_jk}
Only in this paragraph, we cancel the assumption that $X$ is log smooth and we assume only that $X$ is an object of $\cB(\log)$. (We will use the notation introduced in this \ref{m,f_jk} not only for log smooth $X$ as in this Section \ref{ss:A1ls}, but also for more general $X$ later in this paper.)
  Let $p \in X_{[:]}$.

  Let
$M_{X,p}=M^{(0)}\supsetneq M^{(1)}\supsetneq \cdots \supsetneq M^{(m)}=\cO_{X,p}^\times$ be the sequence of faces of $M_{X,p}$ 
which is the pullback
of a linearly ordered sequence of faces of $\cS:=(M_X/\cO_X^\times)_p$ 
in \cite{KNU4} 4.1.6. 
  Here and hereafter, $M_X$ and $\cO_X^{\times}$ are the pullbacks of those on $X$ to $X_{[:]}$ (cf.\ \ref{A_1top}) so that 
$M_{X,p}$ and $\cO_{X,p}^\times$ are the same as 
$M_{X,s}$ and $\cO_{X,s}^\times$, respectively, where $s$ is the image of $p$ in $X$.

  We can take $f_{j,k}\in M^{(j-1)}$  ($1\leq j\leq m$, $1\leq k\leq r(j)$ with $r(j)$ being the rank of $M^{(j-1),\gp}/M^{(j),\gp}$) such that $|f_{j,k}|<1$ and such that for each $j$, $f_{j,k}$ ($1\leq k\leq r(j)$) form a base of the $\Q$-vector space  $\Q\otimes (M^{(j-1),\gp}/M^{(j),\gp})$. Note that $r:=\sum_{j=1}^m r(j)$ is equal to $\rank_{\Z}((M^{\gp}_X/\cO^\times_X)_s)$.
\end{sbpara}

More generally than \ref{Ex1} and \ref{Ex2}, we have the following.  (Again $X$ is log smooth here.) 

\begin{sbprop}
  \label{Ex3} 
  Let $X$ be a log smooth fs log analytic space. 
  Let $p \in X_{[:]}$.
  Let $m$ and $f_{j,k}$ be as in $\ref{m,f_jk}$. 
  Then, $A_1^{\an}$ and $A_1$ near $p$ are described as follows. 
  $A_1^{\an}$ (resp.\ $A_1$) is the sheaf of real analytic functions (resp.\ $\R$-valued $C^{\infty}$ functions) of $\Re(h)$ for $h\in \cO_X$,  $(\log|f_{j,1}|/\log|f_{j+1, 1}|)^{-1/2}$ ($1\leq j\leq m-1$), $(-\log|f_{m, 1}|)^{-1/2}$, and $(\log|f_{j,k}|/\log|f_{j,1}|)^{-1/2}$ ($1\leq j\leq m$, $2\leq k\leq r(j)$).  
\end{sbprop}

\begin{pf}\label{Ex4} 
Let $f,g\in M_{X, p}$, $|f|<1, |g|<1$,  and assume that $\log|f|/\log|g|$ does not have value $\infty$ at $p$. We prove that $(\log|f|/\log|g|)^{1/2}$ is described as stated in 
the proposition.

Step 1. We prove the case $g\in \cO^\times_{X,p}$. In this case,  $f\in \cO_{X,p}^\times$. 
Hence $f=\exp(h_1)$, $g=\exp(h_2)$ for some $h_1,h_2\in \cO_{X,p}$, and $\log|f|/\log|g|= \Re(h_1)/\Re(h_2)$. Since $\Re(h_1)<0$ and $\Re(h_2)<0$, $\Re(h_2)^{-1}$ is a real analytic function of $\Re(h_2)$  and $\Re(h_1)/\Re(h_2)>0$ around $p$, and hence $(\Re(h_1)/\Re(h_2))^{1/2} $ is a real analytic function of $\Re(h_1)$ and $\Re(h_2)$ there. 

Step 2. We prove the case  where $f\in \cO_{X,p}^\times$ and $g= f_{j,1}$ for some $1\leq j\leq m$.
In this case, $f=\exp(h)$ for some $h\in \cO_{X,p}$ and $$(\log|f|/\log|f_{j,1}|)^{1/2}=(- \Re(h))^{1/2}(\log|f_{m,1}|/\log|f_{j,1}|)^{1/2}(-\log|f_{m,1}|)^{-1/2},$$
which is described as stated in the proposition. 

Step 3. We prove the case $f\in M^{(j-1)}$, $f\notin M^{(j)}$ and $g=f_{j,1}$ for some $1\leq j\leq m$. 
Let $a$ be the value  of $\log|f|/\log|f_{j,1}|$ at $p$. Then $0<a<\infty$. Let $\varphi:=(a^{-1}\log|f|/\log|f_{j,1}|)-1$. Then $-1/2\leq \varphi\leq 1/2$ around $p$. We have $(\log|f|/\log|f_{j,1}|)^{1/2}= a^{1/2}(1+\varphi)^{1/2} $ and  this is a real analytic function of $\log|f|/\log|f_{j,1}|$. But $\log|f|$ is a $\Q$-linear combination of $\log|f_{j',k}|$ and $\log|u|$, where $j\leq j'\leq m$, $1\leq k\leq r(j')$, $u\in \cO_{X,p}^\times$. Each  $\log|f_{j',k}|/\log|f_{j,1}|$ is described as stated in the proposition and $\log|u|/\log|f_{j,1}|$ also by Step 2. 

Step 4. We prove the general case. 
By Step 1, we may assume $g\notin \cO_{X,p}^\times$. Let $j$ be  the smallest integer such that $1\leq j \leq m$ and $g\notin M^{(j)}$. Then $f,g\in M^{(j-1)}$. By Step 3,  $(\log|g|/\log|f_{j,1}|)^{1/2}$ is described as stated in the proposition. Since $(\log|g|/\log|f_{j,1}|)^{1/2}$ has a non-zero value at $p$, $(\log|f_{j,1}|/\log|g|)^{1/2}$ is a real analytic function of $(\log|g|/\log|f_{j,1}|)^{1/2}$ around $p$, and 
$(\log|f|/\log|g|)^{1/2}= (\log|f|/\log|f_{j,1}|)^{1/2}(\log|f_{j,1}|/\log|g|)^{1/2}$. It remains to prove that $(\log|f|/\log|f_{j,1}|)^{1/2}$ is described as stated in the proposition.  
  If $f\in \cO_{X,p}^\times$, this is already proved in Step 2. Assume $f\notin \cO_{X,p}^\times$ and let $k$ be the smallest integer such that $j\leq k\leq m$ and $f\notin M^{(k)}$. Then $(\log|f|/\log|f_{j,1}|)^{1/2}= (\log|f|/\log|f_{k,1}|)^{1/2}(\log|f_{k,1}|/\log|f_{j,1}|)^{1/2}$ and we are done by Step 3. 
\end{pf}

\subsection{Log real analytic functions on $X_{[:]}$}\label{ss:Aan}

Let $X$ be an object of $\cB(\log)$.

We define sheaves on $X_{[:]}$ of rings over $\R$ such that 
$$\begin{matrix} A_0^{\an} & \subset & A_1^{\an} & \subset &A_2^{\an}& \subset & A_{2+}^{\an}& \subset&  A_3^{\an}\\
&&&& \downarrow && \downarrow && \downarrow\\
&&&& A^{\an}_{(2)} &\subset & A^{\an}_{(2+)} &\subset& A^{\an}_{(3)}.
\end{matrix}$$  
The sheaf $A_0^{\an}$ is called the {\it sheaf of real analytic functions on $X_{[:]}$}, and the other sheaves are called {\it sheaves of log real analytic functions on $X_{[:]}$}.

\begin{sbpara}\label{strim}
  Let $X$ be an  object of $\cB(\log)$ (cf.\ \ref{ABlog}). 
    
  By a {\it strict immersion of $X$  into log smooth}, we mean a morphism $X\to Z$ in $\cB(\log)$ with $Z$ being a log smooth fs log analytic space satisfying the following conditions (i)--(iii). 
  
  (i) As a map of topological spaces, $X\to Z$ is injective and the topology of $X$ is the strong topology (\cite{KU} 3.1) in $Z$. 
  
  (ii) $\cO_Z\to \cO_X$ is surjective.
  
  (iii) The log structure of $X$ is induced from that of $Z$. 
  
  \end{sbpara}

  \begin{sbpara}\label{Aangeneral}
  
  Using a strict immersion $X\subset Z$ into log smooth  (which always exists locally on $X$), we define $A_k^{\an}$ for $k=0,1$,
  as the quotient of $A^{\an}_k$ of $Z$   by  the ideal generated by $\Re(f)$ for $f\in \Ker(\cO_Z \to \cO_X)$.   

Since $\Im(f)= - \Re(if)$, we have $\Im(f)=0$ in $A_k^{\an}$ of $X$  for $f\in \Ker(\cO_Z \to \cO_X)$. 

We will show in \ref{anindep} below that $A_k^{\an}$ does not depend on the choice of $X\to Z$ and hence it is defined globally on $X$. 
\end{sbpara}

\begin{sbpara}\label{XX'}
  The sheaf $A_0^{\an}$ is also understood as follows, not using the log structure. 

It is the inverse image of the following sheaf $A^{\an}_0$ on $X$.

Let $X'$ be the following locally ringed space over $\C$. As a locally ringed space, $X'$ is the same as $X$. But we regard $\cO_X$ as a sheaf of rings over $\C$ via the map $\C\to \cO_X$ which is the composition of the complex conjugation $\C\to \C$ followed by the inclusion map $\C\overset{\subset}\to \cO_X$. 
Since the category $\cB$ (\ref{ABlog}) has fiber products, we have the object $X\times X'$ of $\cB$ whose underlying set is identified with the product $X\times X$. Consider the continuous map $X\to X \times X'$, the diagonal map. Let the sheaf $A^{\an}_{0, \C}$ of rings over $\R$ on  $X$ be the sheaf-theoretic inverse image of  $\cO_{X\times X'}$ under this diagonal map. 

We have the complex conjugation $\iota: A^{\an}_{0, \C}\to A^{\an}_{0,\C}$ defined as follows. Consider the map $(x,y) \mapsto (y, x)\;;\; X\times X'\overset{\cong}\to X'\times X$. Composing this map with the canonical identification $(X\times X')'= X' \times X$ of locally ringed spaces over  $\C$, we have an isomorphism $X\times X' \overset{\cong}\to (X \times X')'$ of locally ringed spaces over $\R$  which is compatible with the diagonal maps from $X$. This gives an automorphism $\iota$ of $A^{\an}_{0,\C}$ over $\R$ such that $\iota(zf)= {\bar z}\iota(f)$ for $z\in \C$ and $f\in A^{\an}_{0,\C}$. Define $A^{\an}_0\subset A^{\an}_{0,\C}$  to be the fixed part of $\iota$. Then $A^{\an}_{0,\C}=  A^{\an}_0\otimes_\R \C$.

\end{sbpara}

The coincidence of the two definitions of $A_0^{\an}$ will be shown in \ref{an0=an0} below.

\begin{sblem}\label{red} If $S$ is a complex analytic space and reduced, its $A_0^{\an}$ defined as in $\ref{XX'}$ is reduced (that is, it has no non-zero nilpotents) and the map from this $A^{\an}_0$ to the sheaf of $\R$-valued functions is injective. 
\end{sblem}
\begin{pf} This follows from the definition in  \ref{XX'}. With the notation there, $S\times S'$ is reduced because $S$ is reduced. \end{pf}

\begin{sbcor}\label{red2} If $S$ is a complex analytic space and reduced, its $A_0^{\an}$ defined as in $\ref{XX'}$ coincides with the sheaf of subrings of the sheaf of $\R$-valued functions consisting of functions which are written as real analytic functions of $\Re(f)$ for $f\in \cO_X$. 
\end{sbcor}

\begin{sbpara}\label{an0=an0} We prove the coincidence of the two definitions of $A_0^{\an}$. It is enough to prove this in the case where $X$ is log smooth. Then $X$ is reduced as a complex analytic space, and hence Corollary \ref{red2} shows the coincidence for $X$.

\end{sbpara}

\begin{sbpara}\label{c:cO2}  By the understanding of $A^{\an}_0$ as in \ref{XX'}, we have $\cO_X\subset A^{\an}_{0,\C}=A^{\an}_0 \otimes_\R \C\;;\; f\mapsto f \otimes 1$. 

\end{sbpara}

\begin{sbpara}\label{anindep} We prove that the sheaf $A^{\an}_1$  is independent of the choice of $Z\supset X$, which gives the well-defined sheaf globally. 

  In the following, a {\it partial log modification} means an open set of a log modification (\cite{KU} 3.6) of an open set. For a partial log modification $X\to Y$, the map $\Mor(\cdot, X)\to \Mor(\cdot,Y)$ is injective. 

If we have two strict immersions $i_1: X\to Z_1$ and $i_2: X\to Z_2$ into log smooth, locally on $X$, we have a partial log modification $Z_3\to Z_1\times Z_2$ such that the diagonal morphism $X\to Z_1\times Z_2$ factors through a strict immersion  $i_3: X\to Z_3$ into log smooth. 
  Denote the ring $A^{\an}_1$  given by $X\to Z_j$ by $A_{X,Z_j}$ for $j=1,2,3$. 
  On $X$, we have $A_{X, Z_j}\to A_{X, Z_3}$ for $j=1,2$. We prove that these $A_{X,Z_j}\to A_{X, Z_3}$ are isomorphisms for $j=1,2$. 
  We may assume $j=1$. 
  Locally, $Z_3$ is strict smooth over $Z_1$ so that it is enough to show that $A_{X,Z_1} \to A_{X,Z_3}$ is an isomorphism in the case where 
$Z_3= Z_1\times \Delta^n$ with the trivial log structure on $\Delta^n$ and $Z_3\to Z_1$ is the projection. The immersion $X\to Z_3$ is written as $(i_1, f)$ with $f: X\to \Delta^n \subset \C^n$. Locally this  $f$ comes from $g: Z_1\to\C^n$ and hence $X\to Z_3$ factors as $X\to Z_1\to Z_1\times \C^n$, where the second morphism is $z \mapsto (z, g(z))$. 
Thus we are reduced to the statement that $A^{\an}_1$ of $Z_1$  coincides with $A_{Z_1, Z_1 \times \C^n}$ obtained by the immersion $Z_1\to Z_1\times \C^n$. But the isomorphism $Z_1\times \C^n \to Z_1\times \C^n\;;\; (a, b)\mapsto (a, b+g(a))$ (the inverse is $(a,b)\mapsto (a,b-g(a))$) transforms the immersion $Z_1\to Z_1\times \C^n$ to the standard immersion $Z_1\to Z_1\times \C^n\;;\;s\mapsto (s,0)$. 
  Thus we are reduced to the following claim. 

\medskip

{\bf Claim.}  For  a log smooth $X$ and for $X \to Z=X \times \C^n;\;s\mapsto (s, 0)$, the map $A_{X,Z}\to A^{\an}_{1,X}$  is an isomorphism. 

\medskip

This is clear.
\end{sbpara}

\begin{sbpara}\label{fgan} For an open set $U$ of $X_{[:]}$ and for $f,g\in M_X(U)$ such that $|f|<1$ and $|g|<1$, we 
 have $(\log|f|/\log|g|)^{1/2}\in A_1^{\an}(U)$ in the case where its values are  always not $\infty$. This is  clear from definition  in the case where $X$ is log smooth, and in general, this is given by the pullback from the log smooth case by a strict immersion into log smooth.

\end{sbpara}

\begin{sbpara} Let $A_3^{\an}$  on $X_{[:]}$   be the ring theoretic localization of $A_1^{\an}$ by inverting $(\log|f|)^{-1}$ for $f$ in $M_X$ such that $|f|<1$. 

If $X$ is a log smooth fs log analytic space, $A_3^{\an}$ is  the sheaf of subrings of $j_*(A^{\an}_U)$  generated locally over $A_1^{\an}$  algebraically by $\log|f|$ for $f$ in $M_X$ such that $|f|<1$. Here $U$ and $j$ are as in \ref{A_1top}.
 \end{sbpara}

  \begin{sbpara}\label{A2andef} 
Let  $A^{\an}_2$  on $X_{[:]}$ be the sheaf of subrings of $A^{\an}_3$  generated over the image of $A_1^{\an}$ algebraically by 
$$\Re(\alpha(f))(\log|f|)^e$$ 
for $f\in M_X$ such that $|f|<1$ and for $e\in \Z_{\geq 1}$.  Here $\alpha:M_X\to \cO_X$ is the structure map of the log structure. Note that $\Re(\alpha(f))(\log|f|)^e$ for $e\in \Z_{\leq 0}$ are already in $A_1^{\an}$. 

\smallskip

\noindent {\it Remarks.} 
(1) In the place  of $\Re(\alpha(f))(\log|f|)^e$ in the definition of $A^{\an}_2$, we can put   $\Re(\alpha(f))(\log|f|)^e$ and $\Im(\alpha(f))(\log|f|)^e$ to have the same definition because $\Im(h)=-\Re(ih)$ for $h\in \cO_X$.

 (2) We will show in \ref{Poin}--Proposition \ref{d(A2an)}
that $A_2^{\an}$ matches the Poincar\'e base of differential forms. 
\end{sbpara}

\begin{sbpara}
 {\it Example.} 
  Let 
$X=\Delta^n$ and let the notation be as in \ref{Ex2}. 
  Then on the standard open set (\ref{Ex2}) of $\Delta_{[:]}^n$, we have the following: 

As a sheaf of rings, $A_2^{\an}$ is generated over $A_1^{\an}$ algebraically by  $\Re(q_j)y_j^e$ and $\Im(q_j)y_j^e$
for $1\leq j\leq r$ and $e\in \Z_{\geq 1}$.

\end{sbpara}

In the following, for a sharp fs monoid $\La$, $\C\{\La\}$ denotes the ring of convergent series $\sum_{\la\in \La} c_{\la}\la$ ($c_{\la}\in \C$). For an fs monoid $\La$ and variables $T_1, \dots, T_n$, $\C\{\La, T_1, \dots, T_n\}$ denotes $\C\{\La\times \N^n\}$. 

  The stalks of $A_{k,\C}^{\an}$ ($k=0,1,2,3$) are described as follows.
\begin{sbprop}\label{anst} 
  Let $X$ be an object of $\cB(\log)$. 
Let $p\in X_{[:]}$. 
  Let the notation be as in $\ref{m,f_jk}$.
Take a chart $\cS\to M_X$ which induces $\cS\overset{\cong}\to (M_X/\cO_X^\times)_p$. Locally on $X$ at the image of $p$ in $X$, we have a strict immersion into log smooth $X\to \Spec(\C[\cS])^{\an} \times \C^s$ which sends the image of $p$ in $X$ to the origin in the target space, that is, the pair of the unique point of $\Spec(\C[\cS])^{\an}$ lying over the maximal ideal of $\cS$ and the origin of $\bC^s$.
  We fix such an immersion. 

$(1)$ We have $$A^{\an}_{0,\C,p}=\C\{\cS\times \overline{\cS}, S_1, \dots, S_s,\overline{S}_1, \dots, \overline{S}_s \}/(I,\overline{I}).$$
  Here $\overline{\cS}=\{\overline{q}\;|\; q\in \cS\}$ is a copy of $\cS$, $q\in \cS$ gives $q\in \cO_{X,p}$, $\overline{q}\in \overline{\cS}$ gives the complex conjugate of it,  $S_j$ gives the $j$-th coordinate function of $\C^s$, $\overline{S}_j$ gives the complex conjugate of it, and $I$ is the kernel of the surjection $\C\{\cS, S_1, \dots, S_s\}\to \cO_{X,p}$. 

$(2)$ We have $$A^{\an}_{1,\C,p}=\C\{\cS\times \overline{\cS}, S_j, \overline{S}_j\ (1\leq j \leq s), T_{j,k}\ (1\leq j\leq m, 1\leq k \leq r(j))\}/(I, \overline{I}),$$ 
which is compatible with the equality in $(1)$. 
Here  $T_{j,1}$ corresponds to $(\log|f_{j+1,1}|/\log|f_{j,1}|)^{1/2}$ ($f_{m+1}$ denotes $e^{-1}$) and $T_{j,k}$ for $k\neq 1$ corresponds to $(\log|f_{j,k}|/\log|f_{j,1}|)^{1/2}-a_{j,k}$, where $a_{j,k}\in \R$ denotes the value of $(\log|f_{j,k}|/\log|f_{j,1}|)^{1/2}$ at $p$. 

$(3)$ $A^{\an}_{3,\C,p}$ is the localization of  the right-hand-side of $(2)$ 
by inverting $T_{j,1}$ for $1\leq j\leq m$.

$(4)$ $A^{\an}_{2,\C,p}$ is the subring of $A^{\an}_{3,\C,p}$ generated over the right-hand-side of $(2)$ by $q(\prod_{j=a}^m T_{j,1})^{-e}$ and $\overline{q}(\prod_{j=a}^m T_{j,1})^{-e}$ ($1\leq a \leq m$, $q\in \cS$, $q\notin \cS^{(a)}$,  $e\in  \Z_{\geq 1}$). 
  Here $\cS^{(a)}$ is the pullback of $M^{(a)}$ by $\cS \to M_{X,p}$. 
\end{sbprop}

\smallskip

\noindent {\it Remark.}
  Another formulation of the propositions of this kind is as follows$:$ 
  Let $\cC$ be the category of $\R$-algebras of the form $\R\{T_1, \dots, T_n\}/I$ for some ideal $I$. Then for $R\in \cC$, we can define $R\{T_1,\dots, T_n\}$ by the functor $\cC\ni R'\mapsto \text{Mor}_{\cC}(R, R') \times m_{R'}^n$ represented by it. (Here $m_{R'}$ denotes the maximal ideal of $R'$.) If $R=\R\{S_1, \dots, S_m\}/I$, we have $R\{T_1, \dots, T_n\}=\R\{S_1, \dots, S_m, T_1, \dots, T_n\}/(I)$.
  Then we have 
$$A^{\an}_{1,p}=A^{\an}_{0,p}\{T_{jk}\;;\; (1\leq j\leq m, 1\leq k\leq r(j))\},$$ 
which is proved similarly to the above (2). 
  This gives another proof of  the fact that $A_1^{\an}$ is  independent of the choice of the immersion into log smooth (cf.\ \ref{anindep}). 

\begin{pf}
  Let $Z=\Spec(\C[\cS])^{\an} \times \C^s$ and let $p' \in Z_{[:]}$ be the image of $p$. 

(1) We have $$A^{\an}_{0,Z,\C,p'}=\C\{\cS\times \overline{\cS}, S_1, \dots, S_s,\overline{S}_1, \dots, \overline{S}_s \}.$$
  Since $I$ is identified with the stalk of $\Ker(\cO_Z \to \cO_X)$ via $\cO_{Z,p'} \cong \bC\{\cS,S_1,\ldots,S_s\}$, the above equality implies the desired equality. 

(2) Similarly to (1), the equality in (2) is reduced to 
the equality 
$$A^{\an}_{1,Z,\C,p'}=\C\{\cS\times \overline{\cS}, S_j, \overline{S_j}\ (1\leq j \leq s), T_{j,k}\ (1\leq j\leq m, 1\leq k \leq r(j))\}.$$ 
  We prove this. 
  Since the natural map
$$\C\{\cS\times \overline{\cS}, S_j, \overline{S}_j\ (1\leq j \leq s), T_{j,k}\ (1\leq j\leq m, 1\leq k \leq r(j))\} \to A^{\an}_{1,Z,\C,p'}$$ 
is surjective by Proposition \ref{Ex3}, it is enough to see its injectivity. 

It is proved in the following way. By induction on $m$, by fixing non-zero  values of $f_{j,k}$ for $j\geq 2$ and fixing a point of $\C^s$, we are reduced to the case $m=1$ and $s=0$. 
  By fixing values of $\log|f_{1,k}|/\log|f_{1,1}|$ ($2\leq k\leq r(1)$)  and fixing $\arg(f_{1,k})$ ($1\leq k\leq r(1)$), 
we are reduced to the following elementary fact  by taking  $t=(-\log|f_{1,1}|)^{1/2}$ (so $\exp(-t^{-2})= |f_{1,1}|$) and by taking as $B$ the additive submonoid of $\R$ generated by 
$\log|f_{1,k}|/\log|f_{1,1}|$ for $1\leq k\leq r(1)$\ :

Let $F(t_1, t_2)=\sum_{m,n} a_{m,n} t_1^mt_2^n$ ($a_{m,n}\in \C$) be a series, where $(m,n)$ ranges over elements of $\N \times B$ for some  finitely generated additive submonoid $B$ of $\R_{\geq 0}$. Assume that $F$ absolutely converges when $|t_1|$ and $|t_2|$ are sufficiently small. Assume also that  there is a subset $T$ of $\R_{>0}$ whose closure contains $0$ such that  $F(t, \exp(-t^{-2}))=0$ for all sufficiently small $t\in T$. Then $F=0$. 

 (This set $\N \times B$ appears here because we are treating power series in $(-\log|f_{1,1}|)^{1/2}$ and in $|f_{1,k}|$ ($1\leq k\leq r(1)$), and we have $|f_{1,k}|=|f_{1,1}|^b$ with $b=\log|f_{1,k}|/\log|f_{1,1}|\in B$.)

(3) and (4) follow from (2).
\end{pf}

\begin{sbcor}\label{01loc}  For $p\in X_{[:]}$, the stalks  $A^{\an}_{0,p}$ and $A^{\an}_{1,p}$ are local rings. For an element  $f$ of these stalks, $f$ is invertible if and only if the value $f(p)\in \R$ is not zero.

\end{sbcor}

\begin{sbcor}\label{lgpt} Let $X$ be an fs log point and let $p\in X_{[:]}$. Let the notation be as in $\ref{m,f_jk}$. 
  Then the stalk $A^{\an}_{1,p}$  at $p$ is isomorphic to the ring of convergent  power series over $\R$ in $r$ variables $T_{j,k}$ ($1\leq j\leq m$, $1\leq k\leq r(j)$), where $T_{j,1}$ corresponds to $(\log |f_{j+1,1}|/\log |f_{j,1}|)^{1/2}$ for $1\leq j\leq m$ ($f_{m+1,1}$ denotes $e^{-1}$), and $T_{j,k}$ with $k\neq 1$ corresponds to $(\log|f_{j,k}|/\log|f_{j,1}|)^{1/2}-a_{j,k}$, where $a_{j,k}$ is the value of $(\log|f_{j,k}|/\log|f_{j,1}|)^{1/2}$ at $p$.

\end{sbcor}

\begin{sbpara}

The inclusions $A^{\an}_0\subset A^{\an}_1\subset A^{\an}_2\subset A_3^{\an}$ follow from 
 the description of stalks in Proposition \ref{anst}.

\end{sbpara}

 \begin{sbpara}\label{A**} Let $X$ be an object of $\cB(\log)$.
 We define the subsheaf of rings $A^{\an}_{2+}$ of $A^{\an}_3$ as follows. For an open set $V$ of $X_{[:]}$, $A^{\an}_{2+}(V)$ is the set of elements $f$ of $A^{\an}_3(V)$ such that for every $s\in X$, if $V(s)$ denotes the inverse image of $V$ in $s_{[:]}$, the pullback of $f$ to $A^{\an}_3(V(s))$ belongs to $A^{\an}_1(V(s))$. 
 
  \end{sbpara}

\begin{sbpara}\label{A()} Let $X$ be an object of $\cB(\log)$. Define the sheaf of rings $A^{\an}_{(2)}$, $A^{\an}_{(2+)}$, $A^{\an}_{(3)}$ as follows. Let $A^{\an}_{(3)}$ be the localization of $A^{\an}_3$   by inverting local sections $f\in A^{\an}_{2+}$ whose all values $f(p)$  are non-zero. Let
 $A^{\an}_{(2)}$ (resp.\ $A^{\an}_{(2+)}$) be the sheaf of subrings of $A^{\an}_{(3)}$ consisting of local sections of the form $g/f$,
where $f, g\in A^{\an}_2$ (resp.\ $A^{\an}_{2+}$) such that all values $f(p)$ are non-zero. As is easily seen, $A^{\an}_{(2+)}$ is the localization of 
$A^{\an}_{2+}$ by inverting  local sections $f$ whose all values $f(p)$  are non-zero. 

For example, if $X=\Delta$ with the log structure given by the coordinate function $q$,  for an integer $e>0$, $1+\Re(q)y^e$ is not invertible in $A^{\an}_2$ but is invertible in $A^{\an}_{(2)}$. 
\end{sbpara}

 \begin{sbpara}
\label{sheaves_on_logpt}
 If $X$ is an fs log point, we have $A^{\an}_1=A^{\an}_2=A^{\an}_{2+}=A^{\an}_{(2)}=A^{\an}_{(2+)}$. 

  By this, we have $A^{\an}_2\subset A^{\an}_{2+}$ for a general object $X$ of $\cB(\log)$.

\end{sbpara}

\begin{sblem}\label{anloc} For $p\in X_{[:]}$, the stalks 
$A^{\an}_{(2), p}$ and $A^{\an}_{(2+), p}$ at $p$ are local rings. For an element $f$ of these stalks, $f$ is invertible if and only if the value $f(p)\in \R$ is not zero.
\end{sblem}

This is evident.

\begin{sbpara}\label{pts} For $k=1,3$, let $A^{\an}_{k,\mathrm{pts}}$ be the sheaf on $X_{[:]}$ defined as follows. For an open set $U$ of $X_{[:]}$, $A^{\an}_{k,\mathrm{pts}}(U)= \prod_{x\in X} A^{\an}_{k,x}(U(x))$, where $A^{\an}_{k,x}$ is the $A^{\an}_k$ of the fs log point $x$ and $U(x)$ is the inverse image of $U$ in $x_{[:]}$.

\end{sbpara}

\begin{sblem}\label{()pt}

$(1)$ The canonical map from $A^{\an}_2$ (resp.\ $A^{\an}_{2+}$) to $A^{\an}_{1,\pts}$ uniquely factors through $A^{\an}_{(2)}$ (resp.\ $A^{\an}_{(2+)}$). 

$(2)$ The canonical map $A^{\an}_3\to A^{\an}_{3,\pts}$ uniquely factors through $A^{\an}_{(3)}$.

$(3)$ In $A^{\an}_{(3)}$, $A^{\an}_{(2+)}$ coincides with the inverse image of $A^{\an}_{1,\pts}$ under the map $A^{\an}_{(3)}\to A^{\an}_{3,\pts}$.

\end{sblem}

\begin{pf} (1) and (2) follow from Corollary \ref{01loc} applied to $A^{\an}_1$ of fs log points.  We prove (3). A local section of $A^{\an}_{(3)}$ is written locally 
in the form $f/g$ with 
$f\in A_3^{\an}$ and $g\in A^{\an}_{2+}$ such that all values of $g$ are non-zero.
 If the image of $f/g$ in $A^{\an}_{3,\pts}$ belongs 
to $A^{\an}_{	1,\pts}$, then the image of $f$ in $A^{\an}_{3,\pts}$ 
belongs to $A^{\an}_{1,\pts}$ and hence $f$ belongs to $A^{\an}_{2+}$. This proves $f/g\in A^{\an}_{(2+)}$. 
\end{pf}
\begin{sbpara}\label{Poin} In the rest of this section (Section \ref{ss:Aan}), we assume that $X$ is a log smooth fs log analytic space, and we discuss that the sheaf $A^{\an}_2$ has good properties with respect to the Poincar\'e base of  differential forms in degeneration.

In the theory of degeneration of Hodge structures, the invariant metric (Poincar\'e metric) of the upper half plane $\frak H$ considered at the infinity of $\frak H$ is important.  If we denote the coordinate function of $\frak H$ as $x+iy$ with $x$, $y$ being real, then for the Poincar\'e metric, the tangent vectors $\partial/\partial x$ and $\partial/\partial y$ are orthogonal to each other and have length $1/y$. That is, $y\partial/\partial x$ and $y\partial/\partial y$ form an orthonormal base of the tangent bundle of $\frak H$, and hence $dx/y$ and $dy/y$ form an orthonormal base of the cotangent bundle of $\frak H$.  Via the isomorphism 
$$(2\pi i)^{-1}\log: \Delta\smallsetminus \{0\} \overset\sim\to 
\left(\begin{smallmatrix} 1 & \Z\\ 0& 1\end{smallmatrix}\right)\bs \frak H,$$
the coordinate function $q$ of $\Delta$ is understood as $\exp(2\pi i(x+iy))$, and $0\in \Delta$ is regarded as a point at infinity of 
$\left(\begin{smallmatrix} 1 & \Z\\ 0& 1\end{smallmatrix}\right)\bs \frak H$ at which $y=\infty$.  The Poincar\'e metric pulled back to $\Delta$ (which has a singularity at $0\in \Delta$) is important for the study of degeneration of variations of Hodge structure on $\Delta$ at $0\in \Delta$ (cf.\ \cite{Z} Section 3). 
  Further, iterated derivatives by $y\partial/\partial x$ and $y\partial /\partial y$ on $\Delta$ become important in such studies.

Endow $\Delta$ with the log structure at $0$. 
 By 
$$y\partial(\Re(q))/\partial x= -2\pi \Im(q)y, \quad y\partial(\Re(q))/\partial y= -2\pi \Re(q)y,$$ $$ y\partial(\Im(q))/\partial x= 2\pi \Re(q)y, \quad y\partial(\Im(q))/\partial y= -2\pi \Im(q)y,  $$
we see that $A_2^{\an}$ coincides with the sheaf of functions which are obtained from 
 functions in 
$A_1^{\an}$ by taking the iterated derivatives by $y\partial/\partial x$ and $y\partial/\partial y$.

\end{sbpara}

  \begin{sbpara}\label{dRe} The formulas on $y\partial(\Re(q))/\partial x$ etc.\ in the above are generalized to the following.

 For $f\in M_X$, $$d\Re(f)= \Re(f)d\log|f| - \Im(f)d\arg(f)=\Re(f)d\log|f|+ \Re(if)d\arg(f),$$
 where we denote $\Re(\alpha(f))$ and $\Im(\alpha(f))$ simply by $\Re(f)$ and $\Im(f)$, respectively.

  \end{sbpara}

 \begin{sbpara}\label{Aan1}  Let $j:U=X_{\triv}\to X$ be as in \ref{A_1top}. Let $A^{\an,1}_U$ be the sheaf of real analytic $1$-forms on $U$. 
 
 We  define an $A_2^{\an}$-submodule $A_2^{\an,1}$ of $j_*(A^{\an,1}_U)$ on $X_{[:]}$ to be 
 the $A_2^{\an}$-submodule locally generated by $\frac{d\log|f|}{\log|f|}$ and 
$\frac {d\arg(f)}{\log|f|}$ for $f\in M$ such that $|f|<1$. 
\end{sbpara}

 \begin{sbprop}\label{difbase1} 
  The $A_2^{\an}$-module $A_2^{\an, 1}$ is locally free of rank $2n$, where $n$ is the dimension of $X$. Locally we have the following base. Let $p\in X_{[:]}$, and let $f_{j,k}\in M_{X,p}$ ($1\leq j\leq m$, $1\leq k\leq r(j)$) be as in $\ref{m,f_jk}$. 
  Then there are elements $h_k\in \cO_{X,p}$ ($1\leq k\leq n-r$ with $r=\sum_{j=1}^m r(j)$) such that
  $d\log(f_{j,k})$ ($1\leq j\leq m$, $1\leq k\leq r(j)$) and $dh_k$ ($1\leq k\leq n-r$) form a local base of $\omega^1_X$. 
 We have$:$
 
 \medskip
 
 $\dfrac{d\log|f_{j,k}|}{\log|f_{j,k}|}$, $\dfrac{d\arg(f_{j,k})}{\log|f_{j,k}|}$ ($1\leq j\leq m, 1\leq k\leq r(j)$) and $d\Re(h_k)$, $d\Im(h_k)$ ($1\leq k\leq n-r$)  form an $A_2^{\an}$-base of $A_2^{\an,1}$ at $p$.
\end{sbprop}

\begin{pf} 
We may assume that $X=\Spec(\C[\cS])^{\an} \times \Delta^{n-r}$ for some sharp fs monoid $\cS$ (which defines the log structure) and that the image of $p$ is the origin of $X$. 
  We may assume that $f_{j,k} \in \cS$ and $h_k$ are the coordinate functions on $\Delta^{n-r}$. 

  Then the proposition is reduced to the following two facts (1) and (2). 

(1) For $f,g \in M_X$ such that $|f|<1, |g|<1$, we have 
$\frac{d\log|fg|}{\log|fg|}=P\frac{d\log|f|}{\log|f|} + Q\frac{d\log|g|}{\log|g|}$ and
$\frac{d\arg(fg)}{\log|fg|}=P\frac{d\arg(f)}{\log|f|} + Q\frac{d\arg(g)}{\log|g|}$ 
with $P=\frac{\log|f|}{\log|fg|}$, $Q=\frac{\log|g|}{\log|fg|} \in A_1^{\an}$. 
  This is straightforward. 

(2) For $f \in \cO_X^{\times}$ such that $|f|<1$, $d\log|f|$ and $d\arg(f)$ are $A_2^{\an}$-linear combinations of the base given in the statement. 

  This is reduce to 

$(2)'$  For $h \in \cO_X$, $d\Re(h)$ is an $A_2^{\an}$-linear combination of the base given in the statement. 

  To show this, we may assume $h=\alpha(f_{j,k})$ or $h=\alpha(i f_{j,k})$ for some $j$ and $k$. 
  Since the both cases are similar, we assume $h=\alpha(f_{j,k})$. 
  Then, by the formula in \ref{dRe}, we have 
$d\Re(h) = P\frac{d\log|f_{j,k}|}{\log|f_{j,k}|}+Q\frac{d\arg(f_{j,k})}{\log|f_{j,k}|}$ with $P=\Re(\alpha(f_{j,k}))\log|f_{j,k}|$, $Q=-\Im(\alpha(f_{j,k}))\log|f_{j,k}| \in A_2^{\an}$. 
\end{pf}

\begin{sbpara}  {\it Example.} 
  Let $X=\Delta^n$ and let the notation be as in \ref{Ex2}. 
Then on the standard open set (\ref{Ex2}) of $\Delta_{[:]}^n$, $dx_j/y_j$, $dy_j/y_j$ ($1\leq j\leq r$), $d\Re(q_j), d\Im(q_j)$ ($r+1\leq j\leq n$) form an $A_2^{\an}$-base of $A_2^{\an,1}$. 
\end{sbpara}

\begin{sbprop}
\label{d(A2an)} 

$d(A_2^{\an})\subset A_2^{\an,1}$.

\end{sbprop}

\begin{pf} We first prove  $d(A_1^{\an})\subset A_2^{\an,1}$. 
 If $f,g\in M$ such that $|f|<1$, $|g|<1$ 
 and such that $\log|f|/\log|g|$ does not have value $\infty$, then for $F=\log|f|/\log|g|$, we have 
 $d(F)= F\frac{d\log|f|}{\log|f|}- F\frac{d\log|g|}{\log|g|}$, and 
$d(F^{1/2})= \frac 1{\,2\,} F^{1/2}\frac{d\log|f|}{\log|f|}- F^{1/2}\frac{d\log|g|}{\log|g|}$. 
  Hence we have $d(A_1^{\an})\subset A_2^{\an,1}$.

  The proposition  follows from this and from 
 $$d(\Re(f)(\log|f|)^e)= \Re(f)(\log|f|)^ed\log|f| +\Re(if)(\log|f|)^ed\arg(f)+e\Re(f)(\log|f|)^{e-1}d\log|f|$$
for $f\in M_X$ such that $|f|<1$ (\ref{dRe}).  
\end{pf}

\subsection{Log real analytic geometry of the spaces of ratios}\label{ss:grat}

We describe some log real analytic geometry of the spaces of ratios. We will use it in the proofs of Theorems \ref{A2an} and \ref{sharp=}.

\begin{sbpara}\label{SpI} Let $\cS$ be an fs monoid (we denote the  semi-group law of $\cS$  in the multiplicative way).

Recall that an ideal of $\cS$ is a subset of $\cS$ which is stable under multiplications by elements of $\cS$. Recall that a prime ideal $\frak p$ of $\cS$ is an ideal of $\cS$ such that the complement $\cS\smallsetminus \frak p$ is a submonoid of $\cS$. 
Let $\Spec(\cS)$ be the set of all prime ideals of $\cS$. We have a bijection 
$\frak p\mapsto \cS\smallsetminus \frak p$ 
from $\Spec(\cS)$ to the set of all faces of $\cS$.

For an ideal $I$ of $\cS$, let  $\sqrt{I}= \{f\in \cS\;|\; f^n\in I \;\text{for some}\;n\geq 1\}$. Then $\sqrt{I}$ is an ideal of $\cS$. It coincides with the intersection of all prime ideals of $\cS$ which contain $I$ (\cite{O} Chapter I, Corollary 1.4.3.3). 
\end{sbpara}

\begin{sbpara}\label{cS}

Let $X= \Spec(\C[\cS])^{\an}$ be the complex analytic toric variety associated to $\cS$ with the canonical log structure.

For an ideal $I$ of $\cS$, we have the closed analytic subspace $V(I)$ of $X$ defined by the ideal $\cO_XI$ of $\cO_X$.

 Then $V(I)$ is reduced if and only if $I=\sqrt{I}$.

For a prime ideal $\frak p$ of $\cS$, $V(\frak p)$ is isomorphic to $\Spec(\C[\cS'])^{\an}$, where $\cS'=\cS\smallsetminus \frak p$, as a complex analytic space, and we have an open set $U(\frak p)= \Spec(\C[(\cS')^{\gp}])^{\an}$ of $V(\frak p)$. The log structure of $U(\frak p)$ induced from that of $X$ is constant, that is, $M_{U(\frak p)}/\cO^\times_{U(\frak p)}$ is a constant sheaf, and it is isomorphic to $\cS_{\frak p}:=\cS(\cS')^{\gp}/(\cS')^{\gp}$. 
We have a set-theoretic decomposition $X=\coprod_{\frak p\in \Spec(\cS)}\; U(\frak p)$ (decomposition to log strata). For $x\in X$ and $\frak p\in \Spec(\cS)$, $x$ belongs to $U(\frak p)$  if and only if the kernel of the map $\cS\to (M_X/\cO^\times_X)_x$ 
coincides with  $\cS \smallsetminus \frak p$. 
\end{sbpara}

\begin{sbpara}\label{ILS}

An fs log analytic space is said to be {\it ideally log smooth} if it is locally isomorphic to an open set of $V(I)$ for some $\cS$ and $I$ as above endowed with the log structure induced from that of $\Spec(\C[\cS])^{\an}$.

\end{sbpara}

\begin{sbpara}\label{Cp} In the rest of this Section \ref{ss:grat}, let $X=\Spec(\C[\cS])^{\an}$ for an fs monoid $\cS$. We consider log real analytic natures of the space $X_{[:]}$ of ratios. 

For $\frak p\in \Spec(\cS)$, let $C(\frak p)$ be the closure of the inverse image $U(\frak p)_{[:]}$ of $U(\frak p)$ in $X_{[:]}$. 
We will prove that $C(\frak p)$ is the set of common zeros in $X_{[:]}$ of a locally principal ideal $J(\frak p)$  of $A^{\an}_1= A^{\an}_{1,X}$ defined in \ref{Jp}. This ideal $J(\frak p)$ is $0$ if $\frak p=\emptyset$ and is an invertible ideal of $A^{\an}_1$ if $\frak p\neq \emptyset$. Thus, $C(\frak p)$ for $\frak p\neq \emptyset$ is like a Cartier divisor in the usual analytic  geometry (but the topological shape of $C(\frak p)$ is not like a divisor as the example $C(\frak m)$ in \ref{eCp} shows). This is compared with the fact that  the closure $V(\frak p)$ of $U(\frak p)$ in $X$ is the set of common zeros in $X$ of an ideal $\cO_X\frak p$ which need not be principal. 

\end{sbpara}

\begin{sbpara}\label{eCp}
Example. Let $\cS=\N^2$, and hence $X=\C^2$. For $j=1, 2$, let  $q_j$ be the $j$-th  generator of $\N^2$,  that is, the $j$-th coordinate function of $X$. We have
$$\Spec(\cS)=\{\emptyset, \frak p_1, \frak p_2, \frak m\},$$ where  $\frak p_j$ is the prime ideal of $\cS$ generated by $q_j$ and $\frak m$  is the maximal ideal $\cS\smallsetminus \{1\}$ of $\cS$ generated by $q_1$ and $q_2$. 

We have 
$V(\emptyset)= \C^2\supset U(\emptyset)= (\C\smallsetminus \{0\})\times (\C\smallsetminus \{0\}),$
$V(\frak p_1)= \{0\}\times \C\supset U(\frak p_1)= \{0\}\times (\C\smallsetminus \{0\}),$
$V(\frak p_2)= \C\times \{0\}\supset U(\frak p_2)= (\C\smallsetminus \{0\})\times \{0\},$
$V(\frak m)= U(\frak m)= \{(0,0)\}.$  

We have $C(\emptyset)=X_{[:]}$. 
 $C(\frak m)$ is the inverse image of $\{(0,0)\}\subset X$ in $X_{[:]}$ and is identified with  the interval $ [0, \infty]$ by $\log|q_2|/\log |q_1|$. $C(\frak p_1)$ is the union of $U(\frak p_1)=U(\frak p_1)_{[:]}$ and the point $0$ of $[0, \infty]$, and $C(\frak p_2)$ is the union of $U(\frak p_2)=U(\frak p_2)_{[:]}$ and the point $\infty$ of $[0, \infty]$. 
For $j=1,2$, we have $C(\frak p_j) \overset{\cong}\to V(\frak p_j)$  and  $V(\frak p_j)_{[:]}=C(\frak p_j)\cup C(\frak m)$. 
\end{sbpara}

\begin{sbpara}\label{Jp}  Let $\frak p\in \Spec(\cS)$. We define the ideal $J(\frak p)$ of $A^{\an}_1$ as the ideal generated locally by $(\log|f|/\log|g|)^{1/2}$ for $f, g\in \cS$ such that $f\notin \frak p$, $g\in \frak p$, $|f|<1$, $|g|<1$. 

If $\frak p$ is the empty ideal, $J(\frak p)=0$. In the case where $\frak p$ is not empty, as we show in \ref{jp2} below, $J(\frak p)$ is an invertible ideal of $A^{\an}_1$.

\end{sbpara}

\begin{sbpara}\label{jp2}
Let $p\in X_{[:]}$. We describe the stalk of $J(\frak p)$ at $p$.  Let the filtration $\cS=\cS^{(0)}\supsetneq\cS^{(1)}\supsetneq \cdots \supsetneq \cS^{(m)}$ be the inverse image of the filtration $M_{X,p}= M^{(0)} \supsetneq M^{(1)}  \supsetneq \cdots \supsetneq M^{(m)}= \cO^\times_{X,p}$ in \ref{m,f_jk}. So $\cS^{(j)}$ are faces of $\cS$. Then $J(\frak p)_p\neq  A^{\an}_{1,p}$ if and only if $\cS\smallsetminus \frak p$ appears as one of $\cS^{(j)}$. Assume $\cS\smallsetminus \frak p=\cS^{(j)}$. If $j=0$, then $J(\frak p)_p=0$. If $j\geq 1$, $J(\frak p)_p$ is generated by $(\log|f_{j+1,1}|/\log|f_{j, 1}|)^{1/2}$. Here $f_{m+1,1}$ denotes $e^{-1}$. 

\end{sbpara}

\begin{sbpara} Consider the example in \ref{eCp}. 

At $p\in C(\frak m)= [0, \infty]$, $J(\frak m)$ is generated by $(-1/\log|q_2|)^{1/2}$ if $p\neq \infty$, and is generated by $(-1/\log|q_1|)^{1/2}$ if $p\neq 0$. On the interval $(0, \infty)$, $(-1/\log|q_1|)^{1/2}$ and $(-1/\log|q_2|)^{1/2}$ generate the same ideal $J(\frak m)$ of $A^{\an}_1$. 

At $p\in C(\frak p_1)$, if $p\in U(\frak p_1)$, $J(\frak p_1)$ is generated by $(-1/\log|q_1|)^{1/2}$ at $p$. If $p=0\in [0,\infty]$, $J(\frak p_1)$ is generated by $(\log|q_2|/\log|q_1|)^{1/2}$ at $p$.

\end{sbpara}

\begin{sbprop}\label{Cpprop} $(1)$ We have 
$$C(\frak p)= \{p\in X_{[:]}\;|\; (A^{\an}_1/J(\frak p))_p \neq 0\}.$$

$(2)$  $C(\frak p)\subset V(\frak p)_{[:]}$. 

$(3)$ If $f\in \frak p$, $\Re(f)$ and $\Im(f)$ belong to  $A^{\an}_2J(\frak p)$ in $A^{\an}_2$. 
\end{sbprop}

\begin{pf} 
  (1)  Let $C'(\frak p):= \{p\in X_{[:]}\;|\; (A^{\an}_1/J(\frak p))_p \neq 0\}.$  
We  prove $C(\frak p)= C'(\frak p)$. 

Let $Y=\Spec(\bC[\cS\smallsetminus \frak p])^{\an}$. 
  Then we have 
\begin{equation*}\tag{$*$}
U(\frak p)_{[:]}= Y_{\triv} \times R(\cS_{\frak p}),
\end{equation*}
where  
$R(P):=R(\overline P)$ is the space of ratios of an fs monoid $P$ (\cite{KNU4} 4.1.3). 
  This is because $U(\frak p)$ is $\Spec(\bC[(\cS\smallsetminus \frak p)^{\gp}])^{\an}$ endowed with the constant log structure associated to $\cS_{\pe}$. 
  Further we have 
\begin{equation*}\tag{$**$}
C'(\frak p)= Y_{[:]}\times R(\cS_{\frak p}).
\end{equation*}
  In fact, as is stated in \ref{jp2}, a point of $X_{[:]}$ belongs to $C'(\frak p)$ if and only if $\cS\smallsetminus \frak p$ appears as one of $\cS^{(j)}$ at $p$.
  Hence, for a point $y\in Y$, its fiber of $C'(\frak p) \to Y$ is parametrized by the equivalence classes of the sequence $\cS=\cS^{(0)} \supsetneq \cdots 
\supsetneq \cS^{(m)}=\cS_0$ with an interior $N_j$ of the dual cone of $\cS^{(j-1)}/\cS^{(j)}$ ($1 \le j \le m$) such that some $\cS^{(j)}$ is  $\cS\smallsetminus \frak p$, 
where $\cS_0$ is the pullback of $\cO^{\times}_{Y,y}$ (cf.\ \cite{KNU4} 4.1.6). 
  Then $\cS=\cS^{(0)} \supsetneq \cdots \supsetneq \cS\smallsetminus \pe=\cS^{(j)}$ with $N_1,\ldots, N_j$ gives an element of $R(\cS_{\pe})$ and 
$\cS^{(j)} \supsetneq \cdots \supsetneq \cS_0$ with $N_{j+1},\ldots, N_m$ gives an element of $R((\cS\smallsetminus{\pe})/\cS_0)$. 
  Thus we have a bijection $C'(\pe) \to Y_{[:]} \times R(\cS_{\pe})$ which is a homeomorphism as is easily seen. 
  The equality $(**)$ follows. 
  Since $Y_{\triv}$ is dense in $Y_{[:]}$, $(*)$ and $(**)$ imply that $U(\frak p)_{[:]}$ is dense in $C'(\frak p)$.

Hence $U(\frak p)_{[:]} \subset C'(\frak p)\subset C(\frak p)$. On the other hand, $C'(\frak p)$ is closed in $X_{[:]}$. 
  Therefore, $C(\frak p)= C'(\frak p)$.

  (2) Since $V(\frak p)$ is the closure of $U(\frak p)$, we have that $C(\frak p)$ is contained in the inverse image $V(\frak p)_{[:]}$ of $V(\frak p)$ in $X_{[:]}$. 

To prove (3),  we may assume $|f|<1$. We have $1/\log |f|\in J(\frak p)$. We have also $\Re(f)\log |f|\in A^{\an}_2$ and $\Im(f)\log |f|\in A^{\an}_2$. 
  These prove (3). 
\end{pf}

\begin{sbprop}\label{poles} $$A^{\an}_3= {\textstyle\bigcup}_{e}   {\textstyle\prod}_{\frak p} J(\frak p)^{-e},$$
where $\frak p$ ranges over all non-empty prime ideals of $\cS$ and $e$ ranges over all integers $\geq 0$. 
\end{sbprop}

\begin{pf}Let $p\in X_{[:]}$ and assume that the log structure of $X$  at $p$ is not trivial. Let the notation be as in \ref{m,f_jk}. Then $A^{\an}_{3,p}$ is generated over $ A^{\an}_{1, p}$ by $\log |f_{1,1}|$. For $1\leq j\leq m$, let $\frak p_j$ be the non-empty prime ideal $\cS\smallsetminus \cS^{(j)}$ of $\cS$. Then $(\log |f_{1,1}|)^{1/2}= \prod_{j=1}^m (\log|f_{j,1}|/\log|f_{j+1,1}|)^{1/2}$ ($f_{m+1,1}$ denotes $e^{-1}$) and $(\log |f_{j,1}|/\log |f_{j+1,1}|)^{1/2}$ generates $J(\frak p_j)^{-1}$. 
\end{pf}

\subsection{Properties of log real analytic functions}\label{ss:Aanp}

\begin{sbpara}\label{A*val} Value at a point of an element of the stalk of $A^{\an}_{(2+)}$. 

Let $p\in X_{[:]}$ and let $f\in A^{\an}_{(2+),p}$. Then the value $f(p)\in \R$ is defined because the image of $f$ in the stalk at $p$ of 
$A_{3,\sig(p)}^{\an}= (A_3^{\an}$ of $\sig(p)\in X)$ belongs to  $A^{\an}_1$. 
\end{sbpara}

We prove the following two theorems.

\begin{sbthm}\label{A2an}
  Let $X$ be a log smooth fs log analytic space. 
  Let the subsheaves of rings $A^{\an}_{3,c}$ and  $A^{\an}_{3,b}$ of $A^{\an}_3$ and the subsheaves of rings   $A^{\an}_{(3),c}$ and $A^{\an}_{(3),b}$ of $A_{(3)}^{\an}$ be as follows.  
 
  For an open set $V$ of $X_{[:]}$, $A^{\an}_{3,c}(V)$ (resp.\ $A^{\an}_{3,b}(V)$) is the set of elements of $A^{\an}_3(V)\subset A_U^{\an}(U\cap V)$ which extend to continuous functions on $V$ (resp.\ which are locally bounded on $V$).
  Here $U=X_{\triv}$. Define $A^{\an}_{(3), c}$ and $A^{\an}_{(3), b}$ similarly, by replacing $A_3$ by $A_{(3)}$. 
 
 Then we have
\begin{alignat*}{3}
A^{\an}_2&=A^{\an}_{2+}&&=A^{\an}_{3,c}&&=A^{\an}_{3,b},\ \text{and} \\
A^{\an}_{(2)}&= A^{\an}_{(2+)}&&= A^{\an}_{(3),c}&&= A^{\an}_{(3), b}.
\end{alignat*}
\end{sbthm}

\begin{sbthm}\label{sharp=} Let $X$ be an ideally log smooth fs log analytic space.

$(1)$ 
 Let the sheaf 
$\overline{A}^{\an}_2$ (resp.\ $\overline{A}^{\an}_{2+}$, resp.\ $\overline{A}^{\an}_{(2)}$, resp.\ $\overline{A}^{\an}_{(2+)}$) be the image  of the canonical map from $A^{\an}_2$ (resp.\ $A^{\an}_{2+}$, resp.\ $A^{\an}_{(2)}$, resp.\ $A^{\an}_{(2+)}$) to $A^{\an}_{1,\mathrm{pts}}$ ($\ref{pts}$). Then 
 $\overline{A}^{\an}_2= \overline{A}^{\an}_{2+}$ and  $\overline{A}^{\an}_{(2)}= \overline{A}^{\an}_{(2+)}$. 
  
$(2)$ Assume that $X$ is reduced.  Let $E$ be an $A^{\an}_2$ (resp.\ $A^{\an}_{(2)}$)-module on $X_{[:]}$ locally free of finite rank. Then 
  $$\sig_*(E)\overset{\cong}\to  \sig_*(E\otimes_{A^{\an}_2} A^{\an}_{2+}) \quad (\text{resp.\ } \sig_*(E)\overset{\cong}\to  \sig_*(E\otimes_{A^{\an}_{(2)}} A^{\an}_{(2+)})).$$

\end{sbthm}

\begin{sbrem} (1) The part $A^{\an}_2=A^{\an}_{2+}$ of Theorem \ref{A2an} is contained in (1) of Theorem \ref{sharp=} because in the case where $X$ is log smooth, we have  $\overline{A}^{\an}_2= A^{\an}_2$ and $\overline{A}^{\an}_{2+}=A^{\an}_{2+}$. 

(2) A reformulation of Theorem \ref{sharp=} (2) is the following. Assume that $X$ is an ideally log smooth fs log analytic space and assume that $X$ is reduced. For a sheaf of rings $R$ on $X$, let $\text{lff}(R)$ be the category of $R$-modules on $X_{[:]}$ which are locally free of finite rank. Then for $R=A^{\an}_2$ and $R_+=A^{\an}_{2+}$, and also for  $R=A_{(2)}^{\an}$ and $R_+=A^{\an}_{(2+)}$, 
the functor $R_+\otimes_R : \text{lff}(R) \to\text{lff}(R_+)$ is fully faithful.

(3)  Here is an example of $X$ as in (2) of Theorem \ref{sharp=} for which we have $A^{\an}_2\neq A^{\an}_{2+}$. 
Let $Y$ be the fs log point of log rank $2$ whose log structure is generated by $\bN^2=\langle q_1, q_2 \rangle$.
  Let $X=Y\times \Delta$ with the log structure given by $q_1, q_2, q$, where $q$ is the coordinate function of $\Delta$. So, $X$ is ideally log smooth and reduced.  Let $U$ be the open set of $X_{[:]}$ consisting of all points at which the values of  $\log|q|/\log|q_1|$ and $\log|q_2|/\log|q|$ are not $\infty$.  Then $q\log|q_1|$ belongs to $A^{\an}_{2+,\C}$ (because the value of $q$ at every point of $U$ is $0$), but not to $A^{\an}_{2,\C}$ at the point where the value of $\log|q|/\log|q_1|$ is $0$. 

\end{sbrem}

\begin{sbpara}\label{A2an1}
We give a preparation  for the proof of  Theorem \ref{A2an}. We assume $X=\Spec(\C[\cS])^{\an}\times \C^s$ with $\cS$ being a sharp fs monoid, $s \ge0$. 
  Let $p \in X_{[:]}$ be a point over the origin of $X$. 

Let $\frak p$ be a non-empty prime ideal of $\cS$ and let $p\in C(\frak p) \times \bC^s$ (\ref{Cp}). 
  Let the notation $f_{j,k}$ etc.\ be as in \ref{m,f_jk}, let $y_{j,k}=-(2\pi)^{-1}\log|f_{j,k}|$, let $\cS^{(j)}$ be the pullback of $M^{(j)}$ by $\cS \to M_{X,p}$, and let $\ell$ be the integer such that $\cS\smallsetminus \frak p=\cS^{(\ell)}$.

We define  rings and ring homomorphisms $P\to Q\to Q_+\to R$ as follows.

We define $P, Q, Q_+$, and $R$ as follows.

Let $Q= (A^{\an}_2/A^{\an}_2J(\frak p))_{\C,p}$, $Q_+=(A^{\an}_{2+}/A^{\an}_{2+}J(\frak p))_{\C, p}$.

Let $R$ be the stalk at $p$ of the sheaf on $C(\frak p)\times \bC^s$ of all $\C$-valued functions.

Next we define $P$. 
 First, let 
$$P'= \C\{(\cS\smallsetminus\frak p) \times \overline{\cS\smallsetminus \frak p}, S_j, \overline{S}_j\  (1\leq j\leq s), 
 T_{j,k}\ (1\leq j\leq m, 1\leq k\leq r(j), (j,k)\neq (\ell, 1))\}.$$
Let $P$ be the ring generated over $P'$ algebraically by 
$q(\prod_{j=a}^m T_{j,1})^{-e}$ and $\overline q(\prod_{j=a}^m T_{j,1})^{-e}$
($\ell <a\leq m$, $q\in \cS\smallsetminus \frak p= \cS^{(\ell)}, q\notin \cS^{(a)}$, $e\in \Z_{\geq 1}$). 
By Proposition \ref{anst} (4) and \ref{jp2}, $Q$ is identified with $P$.

We have 
$C(\frak p)\times \bC^s= Y_{[:]}\times \bC^s\times R(\cS_{\frak p})$ 
with $Y=\Spec(\C[\cS\smallsetminus \frak p])^{\an}$ by the proof of Proposition \ref{Cpprop}. 
  We have an open immersion of topological spaces from an open neighborhood of the image of $p$ in 
$R(\cS_{\frak p})$ to $\R_{\geq 0}^{\ell-1}\times \prod_{j=1}^{\ell} \R_{>0}^{r(j)-1}$, 
where the map to $\R_{\geq 0}^{\ell-1}$ is   given by $(y_{j+1,1}/y_{j,1})^{1/2}$  ($1\leq j <\ell$) 
and the map to $\prod_{j=1}^{\ell} \R^{r(j)-1}_{>0}$ is given by $(y_{j,k}/y_{j,1})^{1/2}$ 
($1\leq j \leq \ell$, $2\leq k \leq r(j)$).  The stalk of $A^{\an}_2$ of $Y_{[:]}\times \bC^s$ at the image of $p$ is 
identified with the following ring $P_0$ by Proposition \ref{anst} (4). Let 
$$P'_0=\C\{(\cS\smallsetminus \frak p)\times \overline{\cS\smallsetminus \frak p}, S_j, \overline{S}_j\ (1\leq j\leq s), 
 T_{j,k}\ (\ell<j\leq m, 1\leq k\leq r(j))\}.$$
Let $P_0$ be the ring generated over $P'_0$ by $q(\prod_{j=a}^m T_{j,1})^{-e}$ and $\overline{q}(\prod_{j=a}^m T_{j,1})^{-e}$ ($\ell < a \leq m$, $q\in \cS^{(\ell)}\smallsetminus \cS^{(a)}$, $e\in \Z_{\geq 1}$).  

These show that the map $P\to R$ is injective. Hence 
 $Q\to Q_+$ and $Q\to R$ are injective.

\end{sbpara}

\begin{sbpara}\label{A2an2} We prove Theorem \ref{A2an}. We prove $A^{\an}_{2,p}= A^{\an}_{2+,p}$ for all $p\in X_{[:]}$. We may assume $X=\Spec(\C[\cS])^{\an}\times \C^s$ for a sharp fs monoid $\cS$, $s \ge0$, and $p$ is over the origin of $X$.

For each non-empty prime ideal $\frak p$ of $\cS$, let $t_{\frak p}$ be a generator of the ideal  $J(\frak p)_p$ of $A^{\an}_{2,p}$. 

We prove $A^{\an}_{2,p}= A^{\an}_{2+, p}$. Let $f\in A^{\an}_{2+,p}$. By Proposition \ref{poles}, $g=(\prod_{\frak p} t_{\frak p}^{e(\frak p)})f$ belongs to $A^{\an}_{2, p}$ for some integers $e(\frak p)\geq 0$. 
Assume $e(\frak p_0)\geq 1$ for some $\frak p=\frak p_0$ and denote $t_{\frak p_0}$ by $t_0$.   Because the image of $g$ in $A^{\an}_{2+,p}/t_0A^{\an}_{2+, p}$ is zero, the image of $g$ in $A^{\an}_{2,p}/t_0A^{\an}_{2,p}$ is zero by the injectivity of $A^{\an}_{2,p}/t_0A^{\an}_{2, p}\to A^{\an}_{2+,p}/t_0A^{\an}_{2+, p}$ in \ref{A2an1} for $\frak p=\frak p_0$. Hence $g=t_0h$ for some $h\in A^{\an}_{2, p}$. 
 Hence $(\prod_{\frak p} t_{\frak p}^{e'(\frak p)})f\in A^{\an}_{2, p}$, 
where $e'(\frak p)=e(\frak p)-1$ for $\frak p=\frak p_0$ and $e'(\frak p)=e(\frak p)$ otherwise. By downward induction on $(e(\frak p))_{\frak p}$,  we have $f\in A^{\an}_{2,p}$.  

Next we prove $A^{\an}_2=A^{\an}_{3,c}=A^{\an}_{3,b}$. Because $A^{\an}_2\subset A^{\an}_{3,c} \subset A^{\an}_{3,b}$, it is sufficient to prove $A^{\an}_{3,b}\subset A^{\an}_2$. Let $p\in X_{[:]}$ and
 let $f\in A^{\an}_{3,b,p}$. By Proposition \ref{poles}, $g=(\prod_{\frak p} t_{\frak p}^{e(\frak p)})f$ belongs to $A^{\an}_{2, p}$ for some integers $e(\frak p)\geq 0$. 
Assume $e(\frak p_0)\geq 1$ for some $\frak p=\frak p_0$.  Then since $f$ is locally bounded, the image of $g$ in $R$  in \ref{A2an1} for $\frak p=\frak p_0$ is zero. 
Hence by the injectivity of $Q\to R$, the image of $g$ in $Q$ is zero. By downward  induction on $(e(\frak p))_{\frak p}$, we obtain $f\in A^{\an}_{2, p}$. 

For the proof of Theorem \ref{A2an}, it remains to prove $A^{\an}_{(3), b}\subset A_{(2)}$. 
A local section of $A^{\an}_{(3),b}$ is written locally as $f/g$, where $f\in A_3^{\an}$, $g\in A^{\an}_{2+}=A^{\an}_2$ and the values of $g$ are non-zero. Since $f/g$ is locally bounded, $f$ is locally bounded. Hence $f\in A^{\an}_2$. This proves $f/g\in A^{\an}_{(2)}$. 

\end{sbpara}

\begin{sbpara}\label{sharp=0}
Next we consider Theorem \ref{sharp=}. 

Let $Z=\Spec(\C[\cS])^{\an}\times \C^s$ with $\cS$ being a sharp fs monoid, $s \ge0$.
Let $I$ be an ideal of $\cS$ such that $I\neq \cS$, and let $X=V(I)\times \C^s$. Let $x\in X$ be the origin of $Z$. 
  Let $p$ be a point of $X_{[:]}$ lying over $x$. 
  Let the notation be as in \ref{m,f_jk}. 
  Let $\cS^{(j)}$ be the pullback of $M^{(j)}$ by $\cS \to M_{X,p}$ $(1 \le j \le m)$. 
\end{sbpara}

\begin{sbprop}\label{sharp=1} Let the notation be as above. Let $\frak p$ be a prime ideal of $\cS$ ($\frak p$ can be the empty ideal) and assume $p\in C(\frak p)\times \bC^s$. 
  Let $\ell$ be the integer such that $0\leq \ell \leq m$ and $\frak p=\cS\smallsetminus \cS^{(\ell)}$ , and
let $\mu$ be the smallest integer such that $1\leq \mu \leq m$ and $\cS^{(\mu)} \cap (I\cup \frak p)=\emptyset$ (so we have $\mu \geq \ell$). Then $(\overline{A}^{\an}_{2, X}/\overline{A}^{\an}_{2,X}J(\frak p))_{\C, p}$ is identified with the ring which is algebraically  generated over 
$$\C\{\cS^{(\mu)}  \times \overline{\cS}^{(\mu)}, S_j, \overline{S}_j\ (1\leq j\leq s), T_{j,k}\ (1\leq j\leq m, 1\leq k\leq r(j), (j,k)\neq (\ell, 1))\}$$ 
by $q(\prod_{j=a}^m T_{j,1})^{-e}$ 
and  
$\overline{q}(\prod_{j=a}^m T_{j,1})^{-e}$ ($\mu<a\leq m$, $q\in \cS^{(\mu)}\smallsetminus \cS^{(a)}$, $e\in \Z_{\geq 1}$) 
in the field of fractions.
Here  $T_{j,1}$ corresponds to $(\log|f_{j+1,1}|/\log|f_{j,1}|)^{1/2}$ ($f_{m+1}$ denotes $e^{-1}$) and $T_{j,k}$ for $k\neq 1$ corresponds to $(\log|f_{j,k}|/\log|f_{j,1}|)^{1/2}-a_{j,k}$, where $a_{j,k}\in \R$ denotes the value of $(\log|f_{j,k}|/\log|f_{j,1}|)^{1/2}$ at $p$. 
\end{sbprop}

  This is  deduced from Proposition \ref{anst}.

\begin{sbcor}\label{sharp=1.5}  $(A^{\an}_{2, X}/A^{\an}_{2,X}J(\frak p))_p=(\overline{A}^{\an}_{2, X}/\overline{A}^{\an}_{2,X}J(\frak p))_p$ if there is a prime ideal $\frak q$ of $\cS$ such that $p\in C(\frak q)$ and such that $\frak q$ is minimal in the set of prime ideals of $\cS$ which contain $I\cup \frak p$.

\end{sbcor}

\begin{sbpara}\label{sharp=2} We give a simple example of a reduced ideally log smooth fs log analytic space for which $A^{\an}_2\neq \overline{A}^{\an}_2$. Let $Z=\Delta^2$ with the log structure given by the coordinate functions $q_1, q_2$, let $\cS$ be the multiplicative commutative monoid generated by $q_1$ and $q_2$ (so $Z$ is an open subspace of $\Spec(\C[\cS])^{\an}$), let  $\frak p_1$  be the ideal of $\cS$ generated by $q_1$ and let $X=V(\frak p_1)$. Let $p$ be an element of $X_{[:]}$ which lies over $(0,0)$ of $\Delta^2$ such that the value of $\log|q_2|/\log|q_1|$ at $p$ is not $0$. Then $q_2$ is not zero in $A^{\an}_{2,X, p}$ but is zero in $\overline{A}^{\an}_{2, X,p}$. 

\end{sbpara}

\begin{sbpara}\label{sharp=3}  We prove Theorem \ref{sharp=} (1). Let the notation be as in Proposition \ref{sharp=1}. We define rings and ring homomorphisms $Q\to Q_+\to R$ as follows. 
  Let $Q= (\overline{A}^{\an}_{2,X}/\overline{A}^{\an}_{2,X}J(\frak p))_{\C, p}$, let $Q_+= (\overline{A}^{\an}_{2+,X}/\overline{A}^{\an}_{2+,X}J(\frak p))_{\C, p}$, and let $R=  \varinjlim_V  \prod_{x' \in X} \prod_{V\ni p' \mapsto x'}  (A^{\an}_1/J(\frak p) \; \text{of} \; x')_{\C,p'}$, 
where $V$ ranges over all open neighborhoods of $p$.
  Then there are canonical ring homomorphisms as above. 
  The map $Q\to R$ is injective.  In fact, by Proposition \ref{sharp=1}, an element $f$ of $Q$ can be regarded as a convergent power series of the variables in that proposition times some power of $T_{j,1}^{-1}$ ($1 \le j \le m$). 
  If its image in $R$ is zero, all the specializations of $f$ near the origin are zero, which implies $f=0$. 
 Hence the map $Q\to Q_+$ is injective. This proves  $\overline{A}^{\an}_2=\overline{A}^{\an}_{2+}$  by the arguments in \ref{A2an2}.

Next we prove 
 $\overline{A}^{\an}_{(2)}=\overline{A}^{\an}_{(2+)}$. Any element $a$ of $\overline A^{\an}_{(2+)}$ is the image of 
$g/f$, where $f, g\in A^{\an}_{2+}$ such that all values of $f$ are non-zero.
  Let $\overline f$ and $\overline g$ be the images of $f$ and $g$ in $A^{\an}_{1,\mathrm{pts}}$, respectively.   
  Then $a=\overline g/\overline f$. 
  Since $\overline A^{\an}_{2+}=\overline A^{\an}_{2}$, 
both $\overline g$ and $\overline f$ are in $\overline A^{\an}_{2}$. 
  Then $a$ is in $\overline A^{\an}_{(2)}$. 

\end{sbpara}

\begin{sbprop}\label{sharp=4}  Let the notation be as in $\ref{sharp=0}$. Assume $I=\sqrt{I}$. Let $\frak p$ be a prime ideal of $\cS$ ($\frak p$ can be the empty ideal) and let 
 $\frak p_j$ ($1\leq j\leq h$) be all the minimal elements in the set of prime ideals of $\cS$ which contain $I\cup \frak p$. 
  For $1\leq j\leq h$,  let $p_j$ be a point of $X_{[:]}$ lying over $x$ such that $p_j\in C(\frak p_j)\times \bC^s$ (such a $p_j$ exists). Let $R=A^{\an}_2$ (resp.\ $R=A^{\an}_{(2)}$), let $\overline{R}=\overline{A}^{\an}_{2}$ (resp.\ $\overline{R}=\overline{A}^{\an}_{(2)}$), and let $E$ be an $R$-module on $X_{[:]}$ locally free of finite rank. Then the map
$$(\sig_*(E/J(\frak p)E))_x \to {\textstyle\prod}_{j=1}^h (E/J(\frak p)E)_{p_j}\cong {\textstyle\prod}_{j=1}^h ((E/J(\frak p)E) \otimes_R \overline{R})_{p_j}$$ is injective.  
 (The last $\cong$ is by Corollary $\ref{sharp=1.5}$.)

\end{sbprop}

\begin{pf}
  We may assume $I \supset \pe$. 
  Let $f$ be a section of $E/J(\frak p)E$ over an open neighborhood of $\sigma^{-1}(x)$. 
  Assume that it is zero around each $p_j$. 
  Then, by the description of $A^{\an}_2$ by using convergent power series (Proposition \ref{anst}) and by the usual rigidity of real analytic functions, the support of $f$ does not intersect with $C(\pe_j)\times \bC^s$ for any $j$. 
  In particular, $f$ is zero in a nonempty open set of $V:=\sig^{-1}(x)\smallsetminus(\bigcup_j (C(\pe_j)\times \bC^s))$.
  Again by the rigidity of real analytic functions, $f$ is zero also on $V$. 
  Hence $f=0$. 

\end{pf}

\begin{sbprop}\label{sharp=5} Let the notation be as in $\ref{sharp=0}$, $\ref{sharp=3}$, and Proposition $\ref{sharp=4}$. The map $$\sig^{-1}\sig_*(E/J(\frak p)E) \to (E/J(\frak p)E) \otimes_R \overline{R}$$ is injective.

\end{sbprop}

This follows from Proposition \ref{sharp=4}. 

\begin{sbcor}\label{sharp=5.5} Let $R_+=A^{\an}_{2+}$ (resp.\ $R_+=A^{\an}_{(2+)}$). 
The map $$\sig_*(E/J(\frak p)E) \to \sig_*((E/J(\frak p)E) \otimes_R R_+)$$ is injective.

\end{sbcor}

\begin{pf}
Let $\overline{R}_+=\overline{A}^{\an}_{2+}$ (resp.\ $\overline{R}_+= \overline{A}^{\an}_{(2+)}$). 
We have a commutative diagram
$$\begin{CD}  \sig^{-1}\sig_*(E/J(\frak p)E) @>>> \sig^{-1}\sig_*((E/J(\frak p)E) \otimes_R R_+) \\
@VVV @VVV \\
(E/J(\frak p)E) \otimes_R \overline{R} @>\cong>> 
(E/J(\frak p)E)\otimes_R \overline{R}_+.
\end{CD}$$
The left vertical arrow is injective by Proposition \ref{sharp=5}. The lower horizontal arrow is an isomorphism by Theorem \ref{sharp=} (1).
Hence the upper horizontal arrow is injective. Since $\sig$ is a surjective map, this proves Corollary \ref{sharp=5.5}. 
\end{pf}

\begin{sbpara} Now
Theorem \ref{sharp=} (2) follows from Corollary \ref{sharp=5.5} by the argument in \ref{sharp=3}. 

\end{sbpara}

\begin{sbthm}\label{A*cont}
Let $X$ be an object of $\cB(\log)$. For an open set $U$ of $X_{[:]}$  and for $f\in A^{\an}_{(2+)}(U)$, the map $U \to \R\;;\; p\mapsto f(p)$ is continuous.
\end{sbthm}

\begin{pf} In the case where $X$ is ideally log smooth, this follows from Theorem \ref{sharp=} (1)  because sections of $A^{\an}_{(2)}$ gives continuous functions. 

 The general case is reduced to this case as follows. By the definition of strong topology, we may assume that $X$ is an fs log analytic space. Working locally on $X$, take a strict immersion $X\to Z$  into log smooth. We may assume that $X$ is of codimension $\geq 1$ in $Z$. Let $Y$ be the complement of $Z_{\triv}$ in $Z$.  By using Hironaka resolution,  locally on $Z$, take a finite composition of blowing-ups $Z'\to Z$ such that $Z'$ is non-singular, the inverse image $Y'$ of $X\cup Y$ with the reduced structure  is a divisor with normal crossings, and the inverse image $X'$ of $X$ with the reduced structure is a Cartier divisor. Endow $Z'$ with the log structure associated to $Y'$ and endow $X'$ with the inverse image of the log structure of $Z'$. 
Then we have morphisms $Z'\to Z$ and $X'\to X$ of fs log analytic spaces. Since $X'$ is ideally log smooth and reduced and $X'\to X$ is proper and surjective,   the continuity in problem is reduced to the case where $X$ is ideally log smooth and reduced. 
\end{pf}

\begin{sbpara}  

  Log manifolds (\cite{KU} Definition 3.5.7) are objects of $\cB(\log)$ of special type. Log manifolds 
are important because the toroidal partial compactifications of period domains ($=$ the extended period domains of nilpotent orbits) studied \cite{KU} and \cite{KNU1}--\cite{KNU5}, which are classifying spaces of log Hodge structures, are log manifolds. They play important roles in Section \ref{s:newa}. Log smooth fs log analytic spaces are log manifolds. A log manifold is something like ``a log smooth fs log analytic space with slits''.

For a log manifold, the inclusions $$A^{\an}_2\subset A^{\an}_{2+}\subset A^{\an}_3\cap(\text{continuous functions})$$ can be strict. Here we give an example for which the second inclusion is strict. 
  In \ref{Ex(1)(2)} below, we give an example for which the first inclusion is strict. 

 Let $Z=\Delta^2$ with the log structure given by the coordinate functions $q_1,q_2$, and let $X$ be the log manifold $(\Delta^2\smallsetminus (\{0\}\times \Delta))\cup \{(0,0)\}$. 
Let $p$ be the point of $X$ at which $y_1/y_2$ and $y_2$ have value $\infty$.

Then  $y_1y_2^{-2}$ belongs to the intersection of $A_3^{\an}$ and the sheaf of continuous functions. But it does not belong to $A^{\an}_{2+}$ 
at the point $p$ of $X_{[:]}$.  
 
 We prove that the function $y_1y_2^{-2}$ on $X\smallsetminus \{(0,0)\}$ extends to $X$ as a continuous function with value $0$ at $(0,0)\in X$. 
 Consider the sets $U(\epsilon)=\{(0,0)\}\cup 
\bigcup_n U_n(\epsilon_n)$ of $X$, where  $\epsilon=(\epsilon_n)_n$ is a sequence of positive real numbers and 
$U_n(\epsilon_n)= \{(q_1, q_2)\;|\; |q_2|^n<|q_1|, |q_j|<\epsilon_n\}$, which form a base of open neighborhoods of $(0,0)$ in $X$ for the strong topology (see \cite{KU} 3.1.3). 
  Let $c$ be any positive real number. 
  It is enough to show $y_1/y_2^2 < 1/c$ on $\bigcup_n U_n(\epsilon_n)$ for some $\epsilon$. 
  In the definition of $U_n(\epsilon_n)$, $|q_2|^n<|q_1|$ is equivalent to $ny_2>y_1$, that is, $y_1/y_2<n$. 
  Take $\epsilon_n= \exp(-2\pi cn)$. 
  Then $|q_2|<\epsilon_n$ in the definition of $U_n(\epsilon_n)$ is equivalent to $y_2>cn$. 
  Hence $y_1/y_2^2<1/c$ on $U_n(\epsilon_n)$ for any $n$.

\end{sbpara}

\begin{sbpara}\label{Ex(1)(2)} Examples.  We give an example (1) (resp.\ examples (2)) of a log manifold  (resp.\ fs log analytic spaces) for which  $A^{\an}_{2+}\neq A^{\an}_2$ and more strongly, $\overline{A}^{\an}_{2+}\neq \overline{A}^{\an}_2$ and $\sig_*(A^{\an}_{2+})\neq \sig_*(A^{\an}_2)$. 

Consider the following (1) and (2).

(1)  Let $X$ be the log manifold $(\Delta^2\smallsetminus (\{0\}\times \Delta))\cup  \{(0,0)\}$ with the coordinate functions $q_1, q_2$ and with the log structure given by $q_1$, or by $q_1$ and $q_2$.

(2)  Let $X$ be the fs log analytic space $\Spec(\C[q_1,q_2, q_3]/(q_2^2-q_1q_3))^{\an}$ 
 with the log structure given by $q_1$, or by $q_1$ and $q_2$.

Then for an integer $e\geq 1$, $y_1^{e/2}q_2$ belongs to $\sig_*(A^{\an}_{2+,\C})$.
But it does not belong to $\sig_*(A^{\an}_{2,\C})$ and its image in $\overline{A}^{\an}_{2+,\C}$ does not belong to $\overline{A}^{\an}_{2,\C}$ as is seen by taking the strict immersion  into log smooth $X\to (\C^2 \;\text{or}\; \C^3)$  (with the log structure given by $q_1$ or by $q_1, q_2$). 
\end{sbpara}

\subsection{Log $C^{\infty}$ functions on $X_{[:]}$}\label{ss:A}

  Let $X$ be an object of $\cB(\log)$.
  We define log $C^{\infty}$ functions on $X_{[:]}$. 

  We will have sheaves 
$$
A_1 \subset A_2 \subset A_{2*} \subset A_3
$$
on $X_{[:]}$ of rings over $\R$. 
  The sheaves $A_2, A_{2*}$, $A_3$ and their $(\cdot)^{\log}$-versions on $X_{[:]}^{\log}$ play essential 
roles in our integration theory and the log Poincar\'e lemma in Section \ref{s:integration}. 

 \begin{sbpara} First we 
assume that $X$ is a log smooth fs  log analytic space. 
We define the sheaves of log $C^{\infty}$ functions $A_2$, $A_{2*}$, $A_3$ on $X_{[:]}$ as follows. Let $j:U=X_{\triv}\to X_{[:]}$  be as in \ref{A_1top}.

  Thinking that $A^{\an}_2$ matches well the Poincar\'e base of the differential form (\ref{Poin}--Proposition \ref{d(A2an)}), we define its log $C^{\infty}$ version $A_2$ to be the sheaf of subrings of $j_*(A_U)$ defined as follows. For an open set $V$ of $X_{[:]}$, $A_2(V)$ is the set of all elements $f$ of $(j_*(A_U))(V)=A_U(U\cap V)$ satisfying the following condition: 
 
 \medskip
 All iterated derivatives of $f$ with respect to a local base of $A_2^{\an, 1}$ are extended to continuous functions on $V$. 

 \medskip

Here for a local base $(w_j)_j$ of $A_2^{\an,1}$, the iterated derivatives of $f$ with respect to $(w_j)_j$ mean $f, f_j, f_{j,k}, \dots$ defined by $df=\sum_j\; f_jw_j$, $df_j= \sum_k \; f_{j,k}w_k, \cdots$.

       Let  $A_3$ on $X_{[:]}$  be the sheaf of subrings of $j_*(A_U)$ generated locally over $A_2$ algebraically by $\log|f|$ for $f$ in $M_X$ such that $|f|<1$. That is, $A_3$ is obtained from $A_2$ by inverting $1/\log|f|$ for $f\in M_X$ such that $|f|<1$. 
  We have  $$A_3= A_2\otimes_{A_2^{\an}} A^{\an}_3.$$

  Let $A_{2*}$ on $X_{[:]}$  be the sheaf of subrings of $j_*(A_U)$ generated locally over $A_2$ algebraically by $\log(-\log|f|)$ for $f$ in $M_X$ such that $|f|<1$. 
  We have $A_{2*}\subset A_3$ because $\log(-\log|f|)/\log|f|$ belongs to $A_2$. 
 
\end{sbpara}

\begin{sbprop}\label{122ls}  Assume that $X$ is a log smooth fs log analytic space. 
We have the following inclusions. 

$(1)$ $A^{\an}_{(2)}\subset A_2$.

$(2)$ $A_1\subset A_2$. 
\end{sbprop}

\begin{pf}
We have $A^{\an}_2\subset A_2$    by Proposition \ref{d(A2an)}.  (1) follows from this. 
  
  To prove (2), it is enough to show that for a local section $f$ of $A_1$ and for a local $A_2^{\an}$-base $(w_j)_j$ of $A_2^{\an,1}$, we have $df= \sum_j g_j w_j$ for some  local sections $g_j$ of $A_2$. This is proved in the same way
 as the proof of $d(A_1^{\an})\subset A_2^{\an,1}$ in the proof of Proposition \ref{d(A2an)}.
\end{pf}

\begin{sbpara} Example. 
  Endow $\Delta$ with the log structure by the origin, and let $q$ be the coordinate function. 
  Let $y=-(2\pi)^{-1}\log|q|$.
  On $\Delta$, if $a\in \R$ and $a\leq 0$, then $y^a\in A_2$. This is because $y\partial (y^a)/\partial y=ay^a$ and $y\partial (y^a)/\partial x=0$. Note that $y^a\in A_1$ if and only if $a\in 2^{-1}\Z_{\leq 0}$. 

More generally, for any log smooth fs log analytic space $X$ and for $f\in M(X)$ such that $|f|<1$ and for $a\in \R$ such that $a\leq 0$, we have $(-\log|f|)^a\in A_{2,X}$. 
  
\end{sbpara}

  \begin{sbpara}\label{Ageneral}
  Now let $X$ be an object of $\cB(\log)$.
  Using a strict immersion $X\subset Z$ into log smooth (which always exists locally), we define $A_k$ ($k=0,1,2,2*,3$) on $X_{[:]}$ as follows. 
  
  We define $A_3$ 
as the quotient of $A_3$ of $Z$   by  the ideal generated by $\Re(f)$ for $f\in \Ker(\cO_Z \to \cO_X)$.   

Since $\Im(f)= - \Re(if)$, we have $\Im(f)=0$ in $A_3$ of $X$  for $f\in \Ker(\cO_Z \to \cO_X)$. 

For $k=0,1,2,2*$, we define $A_k$ of $X$ as the image of $A_k$ of $Z$ in  $A_3$ of $X$.

 \end{sbpara}

\begin{sbpara}\label{indep} We prove that $A_k$ for $k=0,1,2, 2*, 3$ is independent of the choice of a strict immersion $X\to Z$ into log smooth, which gives the well-defined sheaf globally. 

  Denote the sheaf $A_k$  given by a strict immersion into log smooth $X\to Z$ by $A_{X,Z}$. 
  By the same arguments as in \ref{anindep}, 
  we are reduced to the following claim. 

\medskip

{\bf Claim 1.}  For  a log smooth fs log analytic space $X$ and for $X \to Z=X \times \C^n;\;\;s\mapsto (s, 0)$, the map $A_{X,Z}\to A_k$ is an isomorphism. 

\medskip

This is reduced to the case $n=1$. Let $t$ be the coordinate function of $\C$ (with the trivial log structure). Then 

\medskip

{\bf Claim 2.} If $k=0,1,2$, a local section $f=f(s, t)$ of $A_{k,Z}$ satisfies locally on $X_{[:]}$, $$f(s, t)= {\sum}_{(a,b)\in I} \frac{(\partial_1^a\partial_2^bf)(s,0)}{a!b!}\Re(t)^a\Im(t)^b  + O(|t|^{n+1}) \quad (s \in U=X_{\triv})$$  for any $n \ge 0$. 
  Here $I$ is the set of pairs of integers $(a,b)$ such that $a\geq 0$, $b\geq 0$, and $a+b\leq n$. $\partial_1$ denotes the derivation by $\Re(t)$ and $\partial_2$ denotes the derivation by $\Im(t)$. 

\medskip

  This Claim 2 is a direct consequence of Taylor's theorem. 
  In fact, for a fixed $s \in U$ and sufficiently small $x,y \in \bR_{>0}$, by Taylor's theorem, $f(s,x+iy)=\sum_{a+b \leq n}\frac{x^a}{a!}\frac{y^b}{b!}(\frac \partial {\partial x})^a(\frac \partial {\partial y})^bf(s,0)+R_{n+1}$ with $R_{n+1}=\sum_{a+b= n+1}\frac{x^a}{a!}\frac{y^b}{b!}(\frac \partial {\partial x})^a(\frac \partial {\partial y})^bf(s, \theta x+i\theta y)$ for some $0<\theta<1$. 
  Since $(\frac \partial {\partial x})^a(\frac \partial {\partial y})^bf$ is continuous on $X_{[:]}$ with respect to $s$ and $t=x+iy$, we have $R_{n+1}=O(|t|^{n+1})$. 

  Applying this Claim 2 also to the iterated derivatives of $f$ with respect to the $s$-direction,  
we see that if $f(s, 0)=0$, then $f=\frac{f(s,\Re(t))}{\Re(t)}\cdot \Re(t)+
\frac{f(s,t)-f(s,\Re(t))}{\Im(t)}\cdot \Im(t)$ belongs to the ideal of $A_{k,Z}$ generated by $\Re(t)$ and $\Im(t)$.

Thus we have proved the independence for $k=0,1,2$. That for $k=3$ follows from it by the construction of  $A_3$ as a localization of $A_2$. The independence for $k=2*$ follows from the construction of $A_{2*}$ as a sheaf of subrings of $A_3$ over $A_2$.
\end{sbpara}

\begin{sbpara}
For an object  of $\cB(\log)$, by Proposition \ref{122ls} (1), we have canonical homomorphisms $$A^{\an}_{(2)}\to A_2, \quad A^{\an}_{(3)}\to A_3.$$

\end{sbpara}

\begin{sbprop}\label{A0toA1}  
  Let $X$ be an object of $\cB(\log)$. Assume that $q_1, \dots, q_n\in \cO(X)$ and assume that  the map $(q_j)_j: X\to \C^n$ induces a surjection from the inverse image of $\cO_{\C^n}$ to $\cO_X$. 
Assume that we are given a chart $\cS\to M_X$ and assume that $\cS$ is generated by $f_1, \dots, f_{m-1}$ such that $|f_j|<1$ on $X$ for all $j$. 
  Let $f_m=e^{-1}$. 
Consider the maps
\begin{equation}
X_{[:]}\overset{\subset}\to \C^n \times [0, \infty]^{m^2}\overset{\subset}\to  \C^n \times \C^{m^2} \tag{$\ast$}
\end{equation}
in which the first map is $$((q_j)_{1\leq j\leq n}, ((\log|f_j|/\log|f_k|)^{1/2})_{1\leq j\leq m, 1\leq k\leq m})$$
and the second map is induced by $$[0, \infty] \to \C\;;\; t\mapsto t/(t+1).$$
    Then the  induced map from the inverse image of $A_0^{\an}$ (resp.\ $A_0$) of $\C^n \times \C^{m^2}$ 
via the composite of these two maps 
to $A_1^{\an}$ (resp.\ $A_1$) of $X_{[:]}$ is surjective.

\end{sbprop}

\begin{pf} We may assume that $X$ is a log smooth fs log analytic space. 
Assuming it, let $p \in X_{[:]}$ and we work around $p$.
  By definition of $A_1^{\an}$ (resp.\ $A_1$), the pullbacks of the sections of $A_0^{\an}$ (resp.\ $A_0$) belong to it. 
  To prove the converse, we apply Proposition \ref{Ex3}.
  When we do so, we choose $f_{j,k}$s in Proposition \ref{Ex3} from $f_j$s.
  Then Proposition \ref{Ex3} says that any section of $A_1^{\an}$ (resp.\ $A_1$) can be written as the real analytic (resp.\ $\bR$-valued $C^{\infty}$) function of $\Re(h)$ for $h \in \cO_X$, some $(\log|f_j|/\log|f_k|)^{-1/2}$, and $(-\log|f_j|)^{-1/2}$. 
  From this, we see that it belongs to the pullback of $A_0^{\an}$ (resp.\ $A_0$).
\end{pf}

\begin{sbprop}\label{barA2} 
Let $X$ be an object of $\cB(\log)$ and let $p\in X_{[:]}$. Then 
the stalk $A_{2,p}$ of $A_2$ at $p$ is a local ring with the residue field $\R$. That is, for $f\in A_{2,p}$, $f$ is invertible in $A_{2,p}$ if and only if the value $f(p)\in \R$ is not zero. 

\end{sbprop}

\begin{pf} In the case where $X$ is a log smooth fs log analytic space, this is an easy consequence of the definition of $A_2$ for the case where $X$ is log smooth. 
For general $X$, for a local strict immersion $X\to Z$ into log smooth, since $A_{2,X, p}$ is a quotient of $A_{2,Z,p}$ by an ideal and since $A_{2,Z,p}$ is a local ring, it is sufficient to prove that $A_{2,X,p}$ is not a zero ring, that is, $1\neq 0$ in $A_{2,X,p}$. We prove this first in the case where $X$ is the standard log point, that is, $X$ is the point $0$ of $\Delta$ whose log is given by the coordinate function $q$ of $\Delta$. In this case, the fact $1\neq 0$ in $A_{2,X,p}$ follows from  the fact that we do not have $1=y^e(\Re(q)f+\Im(q)g)$ for any $e\geq 0$, $f,g\in A_{2, \Delta, 0}$ (in fact the right-hand-side converges to $0$ at $0\in \Delta$). For general $X$, if $1=0$ in $A_{2,X,p}$, there is an open neighborhood $U$ of $p$ in $X_{[:]}$  such that $1=0$ in $A_2(U)$. By Proposition \ref{stpt} below, 
we have a morphism $s\to X$ from the standard log point $s$ such that the image of $s=s_{[:]}$ is contained in $U$. Hence  $1=0$ in $A_2$ of $s$, but this does not happen as we have seen. 
\end{pf}

\begin{sbprop}\label{stpt}
 Let $X$ be an object of $\cB(\log)$ and let $U$ be a non-empty open set of $X_{[:]}$. Then there is a  morphism $s\to X$ from the standard log point $s$ such that the image of $s=s_{[:]}$ in $X_{[:]}$ is contained in $U$.
\end{sbprop}

\begin{pf} 
We may assume that $X$ is an fs log point and is the origin of $\Spec(\C[\cS])^{\an}$ for a sharp fs monoid $\cS$. 
  Let $p\in U$ and let $f_{j,k}\in \cS$ at $p$ ($1\leq j\leq m, 1\leq k\leq r(j)$) be as in \ref{m,f_jk}.  
  Let $a_{j,k}$ be the value of $\log|f_{j,k}|/\log|f_{j,1}|$ at $p$. We prove that if  $c(j,k)$ are integers $\geq 1$ such that 
\medskip

($*$) $c(j,k)/c(j,1)$ are sufficiently near to $a_{j,k}$  and $c(j+1, 1)/c(j,1)$ are sufficiently near to $0$, 

\medskip
\noindent
then there is a morphism $s\to X$
  from the standard log point $s$
 which sends $f_{j,k}$ to $q^{c(j,k)}$.  
  Here $q$ is a generator of the log structure of $s$. 
(Then the image of $s=s_{[:]}$ in $X_{[:]}$ is sufficiently near to $p$ and is contained in $U$.) 
  To prove this, it is sufficient to prove that for such $c(j,k)$,  there is a homomorphism $\cS\to \N$ which sends $f_{j,k}$ to $c(j,k)$.  

 Let $N_{j,k}: \cS^{\gp}\to \Z$ be the homomorphism which sends $f_{j,k}$ to $1$ and $f_{j',k'}$ to $0$ for all $(j',k')\neq (j,k)$. Let $\cS^{(j)}=M_p^{(j)}\cap \cS$. Then for each $j$,  $\sum_k a_{j,k}N_{j,k}$ sends $\cS^{(j-1)}$ to $\R_{\geq 0}$ and the kernel of this homomorphism coincides with $\cS^{(j)}$.  Take a finite subset $S$ of $\cS$ which generates $\cS$. Let $s\in S$, $s\neq 1$, $s\in \cS^{(\ell-1)}$, $s\notin \cS^{(\ell)}$ (here $1\leq \ell\leq m$).
Then $N_{j,k}s=0$ if $j<\ell$, and $\sum_k a_{\ell,k}N_{\ell,k} s>0$. Hence if $c(j,k)$ satisfies the above ($*$), we have
$\sum_{j,k} c(j,k)N_{j,k}s\geq 0$. 
Hence $ \sum_{j,k} c(j,k)N_{j,k}$ gives a homomorphism  $\cS \to \N$ if $c(j,k)$ satisfies ($*$).
 \end{pf}

\begin{sbpara} For an element $f$ of the stalk of $A_2$ at $p\in X_{[:]}$, we define the value $f(p)\in \R$ of $f$ at $p$ as the residue class of $f$ (Proposition \ref{barA2}). This is compatible with the value for $A_{(2+)}^{\an}$ in \ref{A*val}. 
\end{sbpara}

\begin{sbprop}\label{A2cont} Let $X$ be an object of $\cB(\log)$, let $U$ be an open set of $X_{[:]}$, and let $f\in A_2(U)$. Then the map $U\to \R\;;\;p\mapsto f(p)$ is continuous.

\end{sbprop}

\begin{pf}
Using a strict immersion into log smooth, this is reduced to the case where $X$ is a log smooth fs log analytic space.
\end{pf}

\begin{sbprop}\label{loggrow}  Let $X$ be a log smooth fs log analytic space. 
  Working locally on $X$, take   $g\in M_X$ such that $|g|<e^{-1}$, $-\log|g|>1$, and such that 
$U=X_{\triv}$ coincides with the set of all points of $X$ at which $g$ belongs to $\cO_X^\times$, and also take a local base ${\bf e}$ of $A_2^1$ over $A_2$. 
Then for a local section $f$ of $j_*(A_U)$, the following two conditions {\rm (i)} and {\rm (ii)} are equivalent. 

{\rm (i)} $f$ belongs to $A_3$.

{\rm (ii)} Locally on $X_{[:]}$, there is an $n\geq 1$ such that for any iterated derivative $h$ of $f$ with respect to $\bf e$, locally on $X_{[:]}$, there is a constant  $C>0$ such that $|h|\leq C(-\log|g|)^n$. 
\end{sbprop}

\begin{pf} If (i) is satisfied, $f=\sum_{j=1}^m h_j \prod_{k=1}^{a(j)} \log|g_{j,k}|$ with $h_j\in A_2$, $a(j)\geq 1$, $g_{j,k}\in M_X$, $|g_{j,k}|<1$. Since $d\log|g_{j,k}| \in \log|g_{j,k}|A_2^1$, iterated derivatives of $\log|g_{j,k}|$ belong to $A_2\log|g_{j,k}|$. Hence iterated derivatives of $\prod_{k=1}^{a(j)} \log|g_{j,k}|$ belong to $A_2 \prod_{k=1}^{a(j)} \log|g_{j,k}|$. Working locally, we have $-\log|g_{j,k}|\leq -\log|g^c|$ for some $c\geq 1$. Hence we have (ii). 

Conversely, assume that (ii) is satisfied. Let $V$ be an open set of $X_{[:]}$  and assume $f\in j_*(A_U)(V)$. 
Let $f_2= (\log|g|)^{-n-1}f$. Then an iterated derivative $h_2$ of $f_2$ with respect to ${\bf e}$  is a linear combination of $(\log|g|)^{-n-1}h$ over $A_2$, where $h$ is an iterated derivative of $f$ with respect to ${\bf e}$.
Hence locally at each point $p$ of $V$, there is a constant $C>0$ such that $|h_2|\leq C(-\log|g|)^{-1}$. Hence when a point $p'$ of $U \cap V$ converges to $p$, $h_2(p')$ converges to $0$. It implies that 
$h_2$ is continuous on $V$. Therefore $f_2\in A_2$. This implies $f\in A_3$. 
\end{pf}

\begin{sbrem}
\label{r:KMN}
 In \cite{KMN}, for a log smooth fs log analytic space $X$, we used  a sheaf of log $C^{\infty}$ functions on $X$. By Proposition \ref{loggrow}, we see that that sheaf is similar to $A_3$.

 \end{sbrem}

\begin{sbpara} 
$A_{2+}^{\an}$ of a log manifold  is not necessarily contained in  the sheaf of log $C^{\infty}$ functions $A_2$. Let $X=(\Delta^2\smallsetminus (\{0\}\times \Delta))\cup\{(0,0)\}$ with the log structure given by the first coordinate function $q_1$. 
Then $y_1\Re(q_2)$ belongs to $A^{\an}_{2+}$ but it does not belong to $A_2$ because $d(y_1\Re(q_2))= y_1\Re(q_2)dy_1/y_1+y_1d\Re(q_2)$ and the coefficient $y_1$ of $y_1d\Re(q_2)$ is not continuous on $X$. 
\end{sbpara}

\subsection{Sheaves on $X_{[:]}^{\log}$}\label{ss:shlog}

\begin{sbpara}\label{Aklog}
Let $X$ be an object of $\cB(\log)$. 
We define the sheaves of rings 
$$A^{\an,\log}_k, A^{\an,\log*}_k\quad (k=1, 2, 2+, 3, (2), (2+), (3)),$$
$$A_k^{\log}, A_k^{\log*} \quad (k=1,2,2*, 3)$$ 
on $X^{\log}_{[:]}$ as follows.

Let $R=A^{\an}_k$. In the cases $k=1,2,2+,3$, let $R'=A^{\an}_3$. In the cases $k=(2), (2+), (3)$, let $R'=A^{\an}_{(3)}$. Then $R^{\log}$ (resp.\ $R^{\log*}$) is  the sheaf of subrings of $\cO_X^{\log}\otimes_{\cO_X} R'_{\C}$ generated over $R$ locally and algebraically by $\arg(f)/\log|f|$ (resp.\ $\arg(f)$), where $\arg(f):=i^{-1}(\log(f)\otimes 1-1\otimes \log|f|)$, for $f\in M_X$ such that $|f|<1$. 

We define $$A^{\log}_k = A_k \otimes_{A^{\an}_1} A^{\an, \log}_1, \quad A^{\log*}_k=A_k \otimes_{A_1^{\an}} A_1^{\an,\log*}.$$

We have $$\tau_X^{-1}(A^{\an}_k)\subset A^{\an,\log}_k\subset A^{\an,\log*}_k, \quad \tau_X^{-1}(A_k)\subset A^{\log}_k\subset A^{\log*}_k \quad \text{for all $k$},$$ where $\tau_X: X^{\log}_{[:]}\to X_{[:]}$. Because $1/\log|f|$ ($f\in M_X$, $|f|<1$) are invertible in $A_3^{\an}$, we have
$$A^{\an,\log}_k=A^{\an,\log*}_k\;\text{for $k=3, (3)$}, \quad A^{\log}_3=A^{\log*}_3.$$

By the following Proposition \ref{polstalk}, we have 
$$A^{\an,\log}_k = A_k^{\an} \otimes_{A^{\an}_1} A^{\an, \log}_1, \quad A^{\an,\log*}_k=A_k^{\an} \otimes_{A_1^{\an}} A_1^{\an,\log*}.$$

\end{sbpara}

\begin{sbprop}\label{polstalk}  Let $X\in \cB(\log)$. Consider the following cases {\rm (i)} and {\rm (ii)}.

{\rm (i)} $R=A^{\an}_k$ or $R= A_k$, and $R'=R^{\log}$.

{\rm (ii)} $R=A^{\an}_k$ or $R= A_k$, and $R'=R^{\log*}$.

Then the stalk of $R'$  
 is a polynomial ring over the stalk of $R$. 
More precisely, if $p\in X_{[:]}$ and if $\tilde p$ is a point of $X_{[:]}^{\log}$ lying over $p$, then in the case {\rm (i)} (resp.\ {\rm (ii)}), 
 with the notation in $\ref{m,f_jk}$, the stalk of $R'$ at $\tilde p$ is isomorphic to the polynomial ring over the stalk of $R$ at $p$ in $r$ variables $T_{j,k}$ ($1\leq j\leq m, 1\leq k\leq r(j)$), where $T_{j,k}$ corresponds to   $\frac{\arg(f_{j,k})}{\log|f_{j,k}|}$  (resp.\ $\arg(f_{j,k})$). 

\end{sbprop}

\begin{pf} This is by the fact that for $f,g \in M_X$ such that $|f|<1, |g|<1$, we have 
$\frac{\arg(fg)}{\log|fg|}=P\frac{\arg(f)}{\log|f|} + Q\frac{\arg(g)}{\log|g|}$ 
with $P=\frac{\log|f|}{\log|fg|}$, $Q=\frac{\log|g|}{\log|fg|} \in A_1^{\an}$
and the fact that for $f \in \cO_X^{\times}$ such that $|f|<1$, we have $\arg(f), \frac{\arg(f)}{\log|f|} \in A_0^{\an}$. 
\end{pf}

\begin{sbcor}
\label{c:pol}
  Let the notation be as in Proposition $\ref{polstalk}$. For the canonical map $\tau_X: X^{\log}_{[:]}\to X_{[:]}$, we have$:$

$(1)$ $R=\tau_{X*}R'$. More generally, for an open set $U$ of $X_{[:]}$  and for a sheaf $\cF$ of  $R$-modules on $U$, we have 
$\cF= \tau_{X*}(R'\otimes_R \cF)$.

$(2)$ In the case {\rm (ii)},  
   we have the isomorphism
$$R \overset{\cong}\to R\tau_{X*}R'$$ 
in the derived category. More generally, for an open set $U$ of $X_{[:]}$ and for a sheaf $\cF$ of $R$-modules  on $U$, we have the isomorphism $\cF \overset{\cong}\to R\tau_{X*}(R'\otimes_R \cF)$.

\end{sbcor}

\begin{pf} Since the map $\tau_X$ is proper, by the proper base change theorem, the stalk of the higher direct image $R^m\tau_{X*}(R'\otimes_R \cF)$ at $p\in X_{[:]}$ is $H^m(\tau_X^{-1}(p), R'\otimes_R \cF)$. Since $\tau_X^{-1}(p)$ is homeomorphic to $(\R/\Z)^r$,  $H^m(\tau_X^{-1}(p), R'\otimes_R \cF)$ is isomorphic to the group cohomology $ H^m(\Gamma, R_p[T_{j,k}\; (1\leq j\leq m, 1\leq k\leq r(j))]\otimes_{R_p} \cF_p)$, where $\Gamma= \pi_1((\R/\Z)^r)\cong \Z^r$.  The $(j,k)$-th generator of this $\Z^r$ acts 
by $T_{j,k}\mapsto T_{j,k}+ \frac{2\pi}{\log|f_{j,k}|}$ in the case (i), and 
by $T_{j,k}\mapsto T_{j,k}+2\pi$ in the case (ii).
 The results follow from this. (The computation of the group cohomology for (2) is reduced to the fact that for a ring $R$ over $\Q$, we have the splitting exact sequence of $R$-modules $0\to R \to R[T]\overset{h}\to R[T]\to 0$, where $h$ sends $f(T)$ to $ f(T+1)-f(T)$.) 
\end{pf}

\begin{sbpara}\label{cOlog}
For an object $X$ of $\cB(\log)$, we
have a homomorphism $\cO_X^{\log}\to A_{3,\C}^{\an,\log}\;;\; \log(f)\mapsto  \log |f|+ i\arg(f)$ ($f\in M_X$, $|f|<1$)  which extends the homomorphism $\cO_X\to A^{\an}_{1,\C}$ in Corollary \ref{c:cO}. This induces an isomorphism $$ \cO_X^{\log} \otimes_{\cO_X} A^{\an}_{3,\C}\overset{\cong}\to A^{\an,\log}_{3,\C}.$$

\end{sbpara}

\begin{sbpara}\label{OlogAlog} As a sheaf of subrings of $A^{\an,\log}_{3,X,\C}$, $\cO_X^{\log}$ is understood as a sheaf of rings over the inverse image of $\cO_X$ generated locally and algebraically by $\log|f|+  i \arg(f)$ ($f\in M_X$, $|f|<1$).

\end{sbpara}

\begin{sbpara}
  Assume that $X$ is a log smooth fs log analytic space. Let $U=X_{\triv}\subset X$, let the sheaves $A_U^{\an}$ and $A_U$ on $U$ be as in \ref{A_1top}, and let $j^{\log}: U\to X^{\log}_{[:]}$ be the inclusion map. Let the notation be as in Proposition \ref{polstalk}.
  
  Then  in (i) of Proposition \ref{polstalk}, $R'$ is the sheaf of subrings of  $j^{\log}_*(A^{\an}_U)$ or  $j^{\log}_*(A_U)$
generated over $R$ locally and algebraically by  $\frac{\arg(f)}{\log|f|}$ for  $f\in M_X$ such that $|f|<1$.
 In (ii) of Proposition \ref{polstalk}, $R'$ is the sheaf of subrings of  $j^{\log}_*(A^{\an}_U)$ or  $j^{\log}_*(A_U)$
generated over $R$ locally and algebraically by $\arg(f)$ for $f\in M_X^{\gp}$.

\end{sbpara}

 \begin{sbpara}\label{fganlog} For an open set $U$ of $X_{[:]}^{\log}$ and for $f\in M_X^{\gp}(U)$, we have 
 $\arg(f)\in A^{\an,\log*}_1(U)$ if we fix  a branch of the continuous map $\arg(f)$ (\ref{fgcont}, a branch means a continuous lifting to $\R$ of the continuous  map to $\R/2\pi \Z$). This is already given in the case where $X$ is log smooth, and in general, it is obtained by the pullback from the log smooth case by strict immersions into log smooth. 

\end{sbpara}

\begin{sbpara}\label{logval} We have the evaluation map $h\mapsto h(p)$  from the stalk of $A^{\log*}_2$ at $p\in X^{\log}_{[:]}$ to $\R$ extending the evaluation map for  $A_2$. This is because $A_2^{\log*}$ is generated over $A_2$ locally and algebraically by branches of $\text{arg}(f)$ for $f\in M^{\gp}_X$ (\ref{fganlog}). 
It induces evaluation maps from the stalks of 
 $A^{\an,\log}_{(2+)}$, $A^{\an,\log*}_{(2+)}$, $A^{\log}_2$ to $\R$. 
For an open set $U$ of $X^{\log}_{[:]}$ and for $h\in A_2^{\log*}(U)$, the map $U\to \R\;;\;p\mapsto h(p)$ is continuous.

 \end{sbpara}

\begin{sbpara}\label{ringedmor} For a morphism $X\to Y$ in $\cB(\log)$, we have the associated morphisms of ringed spaces
$(X_{[:]}, R_X)\to (Y_{[:]}, R_Y)$ for $R=A_k^{\an}$, $A_k$ ($k=1, 2$, etc.), 
 and 
$(X_{[:]}^{\log}, R_X)\to (Y_{[:]}^{\log}, R_Y)$ for $R=A_k^{\an,\log}$, $A_k^{\an,\log*}$, $A_k^{\log}$, $A_k^{\log*}$. 
These are the evident morphisms in the case where $X$ and $Y$ are log smooth. In general, these are obtained locally from that case by using a commutative diagram 
$$\begin{CD} X @>>> Y \\
@VVV @VVV \\
Z @>>> Z'
\end{CD}$$
in which the vertical arrows are strict immersions into log smooth.
  Such a diagram exists locally on $X$ and on $Y$. 
 It can be shown that the result  is independent of the choice of such a commutative diagram and hence  defined globally.

\end{sbpara}

\subsection{Soft property of log $C^{\infty}$ functions}\label{ss:soft} 
In the usual theory of $C^{\infty}$ functions, it is important that the sheaf of $C^{\infty}$ functions is soft. Here we prove a log version Theorem \ref{softhm}. 

\begin{sbpara}\label{softreview} We review some basic facts about soft sheaves.  

Assume that $T$ is a paracompact topological space. (Recall that we adopt the definition of paracompactness which includes Hausdorffness. Cf.\ Convention.) 

(1) Definition of softness. 
  A sheaf $\cF$ on $T$ is {\it soft} if 
for every closed subset $C$ of $T$ with the inclusion map $i: C\to T$, the map $\cF(T) \to (i^{-1}(\cF))(C)$ is surjective.

(2) Relation with the partition of unity. 

A sheaf $\cF$ of rings on $T$ is soft if and only if for every locally finite open covering $(U_{\la})_{\la}$ of $T$, there is a partition of unity subordinate to it, that is, there is $s_{\la}\in \cF(T)$ for each $\la$ satisfying the following (i) and (ii).

(i) For each $\la$, the support of $s_{\la}$ is contained in $U_{\la}$.

(ii) $1=\sum_{\la} s_{\la}$. 

(3) If $\cF$ is a soft sheaf of rings on $T$, any sheaf of $\cF$-modules on $T$ is soft.

(4) If $\cF$ is a soft sheaf of abelian groups on $T$, $H^m(T, \cF)=0$ for all $m>0$.

(5) Example. If $T$ is a paracompact $C^{\infty}$-manifold, the sheaf of $C^{\infty}$-functions on $T$ is soft.  
(Usually the paracompactness is  assumed in the definition of a $C^{\infty}$-manifold, but we do not impose it.) 

\end{sbpara}

\begin{sbpara}\label{(*)}

Recall that a locally compact space satisfying the second axiom of countability is paracompact.

 Let $X$ be a complex analytic space.  Assume that 
 
 \medskip

$(*)$ $X$ is Hausdorff and satisfies  the second axiom of countability. 

\medskip

\noindent Then $X$ is  paracompact because it is locally compact and satisfies the second axiom of countability. 

\end{sbpara}

\begin{sblem}
Let $X$ be an fs log analytic space and assume that the underlying complex analytic space of $X$ satisfies the above condition $(*)$. Then any open subset of $X_{[:]}$ and any open subset of $X_{[:]}^{\log}$ are paracompact.  
\end{sblem}

\begin{pf}
Since $X$ is Hausdorff and since $X_{[:]}\to X$ and $X_{[:]}^{\log}\to X$ are proper, $X_{[:]}$ and $X_{[:]}^{\log}$ are Hausdorff. 
  Furthermore, locally on $X$, $X_{[:]}$ and $X_{[:]}^{\log}$ are homeomorphic to a closed subspace of $X \times \R^n$ for some $n$  (cf.\ Section \ref{ss:srsa}). 
    Hence they are locally compact and satisfy the second axiom of countability. 
  The lemma follows. 
\end{pf}
  
\begin{sbpara}\label{pos} 
In the classical theory of  partition of unity $1=\sum_{\la} s_{\la}$, it is often important to take $s_{\la}\geq 0$. In our theory of integration in log geometry, to prove the positivity of some integrals, for example, as in Proposition \ref{I1} (2),  we need to consider such a condition $s_{\la}\geq 0$. However, for log $C^{\infty}$ functions in this paper,  how to define the condition $\geq 0$ is rather delicate. Let $X$ be an object of $\cB(\log)$, let $U$ be an open set of $X_{[:]}$, let $k=0,1,2,2*,3$ and let $f\in A_k(U)$.  We think that the following condition (P) is the most natural formulation of $f\geq 0$ in log geometry.

(P) For each $p\in U$, if $X'$ is an open set of $X$ which contains the image of $p$ in $X$ and if $X'\to Y$ is a strict immersion into log smooth, there exist an open neighborhood $V$ of $p$ in $Y_{[:]}$ and $g\in A_k(V)$ satisfying the following (i), (ii), and (iii).

(i) $V\cap X'_{[:]}\subset U$.

(ii) The restriction of $f$ to $V\cap X'_{[:]}$ coincides with the pullback of $g$. 

(iii) $g(p')\geq 0$ for all $p'\in V$ at which the log structure is trivial.  

\medskip

This condition (P) is equivalent to the condition obtained by rewriting the first part of (P) as 
``For each $p\in U$, there are an open set $X'$ of $X$ which contains the image of $p$ in $X$, a strict immersion $X'\to Y$  into log smooth, an open neighborhood $V$ of $p$ (the rest is the same). The equivalence of this condition with (P) is proved by the method of \ref{anindep}. In the case $k=0,1,2$, these conditions are equivalent to the conditions in which we remove  ``at which the log structure is trivial'' in (iii). (In the case $k=0,1,2$, we can take the value $g(p')$ at every $p'\in V$).

For a locally finite open covering $(U_{\la})_{\la}$ of $U$, a partition $1=\sum_{\la} s_{\la}$ of unity in $A_k(U)$ 
subordinate to it is said to be {\it with} (P)  if each $s_{\la}$ satisfies (P).

\end{sbpara}

\begin{sblem}\label{poslem} 
  Let $X$ be an object of $\cB(\log)$, let $U$ be an open set of $X_{[:]}$, let $k=0,1,2,2*, 3$, and let $f\in A_k(U)$.

$(1)$ Assume that $X$ is a log smooth fs log analytic space.  
  Then $f$ satisfies {\rm (P)} if and only if $f(p)\geq 0$ for every $p\in U$ at which the log structure is trivial.

$(2)$ Assume that $X$ is log smooth fs log analytic space and assume $k=0,1,2$. Then $f$ satisfies {\rm (P)} if and only if $f(p)\geq 0$ for all $p\in U$.

$(3)$ If $g\in A_k(U)$ and if $f$ and $g$ satisfy {\rm (P)}, then $fg$ and $af+bg$ for every $a,b\in \R_{\geq 0}$ satisfy {\rm (P)}.

$(4)$ If $Y\to X$ is a morphism in $\cB(\log)$ and if $V$ is an open set of $Y_{[:]}$ whose image in $X_{[:]}$ is contained in $U$ and $f$ satisfies {\rm (P)}, 
 then the pullback of $f$ in $A_k(V)$ also satisfies {\rm (P)}.

\end{sblem}

\begin{pf} (1) and (2) are clear from the definition. (3) is reduced to the log smooth case and to (1) by using a strict immersion into log smooth. (4) is also reduced to the case where $X$ and $Y$ are log smooth and to (1) by using the following commutative diagram which exists locally on $X$ and $Y$ and in which the vertical arrows are strict immersions into log smooth: 
$$\begin{CD}
Y @>>> X \\
@VVV @VVV \\
Y' @>>> X'. 
\end{CD}
$$\end{pf}
\begin{sbpara}\label{pos3}

In the case $k=0,1,2$, we have the following another version of $f\geq 0$.

(P0)  $f(p)\geq 0$ for all $p\in U$.

We have the implication  (P) $\Rightarrow$ (P0).

(P) and (P0) are not equivalent. For example, let $X=U=\{0\}\subset \Delta$ with the log structure given by the coordinate function $q$ of $\Delta$. Then for $k=1,2$ and for $y=-(2\pi)^{-1}\log|q|$, $-y^{-1}\in A_k(U)$ satisfies the condition (P0) because its value at the unique point of $U$ is $0$, but it does not satisfy  (P).

If $X$ is log smooth, we have (P) $\Leftrightarrow$ (P0) for $k=0,1,2$ by Lemma \ref{poslem} (2).

\end{sbpara}

\begin{sbthm}\label{softhm} Let $X$ be an fs log analytic space. Assume that the underlying complex analytic space of $X$ satisfies the condition $(*)$ in $\ref{(*)}$. Let $U$ be an open set of $X_{[:]}$.
 Let $k=1,2,2*,3$. Then $A_k$ on $U$ is soft. As consequences, all sheaves of $A_{k,U}$-modules on $U$ are soft.  
 Furthermore, for every locally finite open covering of $U$, $A_k$ on $U$ has partition of unity subordinate to it with $(\mathrm{P})$ 
($\ref{pos}$). 
\end{sbthm} 

\noindent {\it Remark.} Actually, the case $k=1$ of this theorem implies this theorem for all $k$.

We first prove

\begin{sbprop}\label{A0soft}  Let $X$ be a complex analytic space. Assume that $X$ is paracompact. Then $A_0$ of $X$ is soft. Furthermore, for every locally finite open covering of $X$, 
$A_0$ of $X$ has a partition of unity $1=\sum_{\lam} s_{\lam}$ subordinate to it such that each $s_{\la}$ is 
written as $s_{\la}=\sum_{\mu} s_{\la,\mu}^2$ ($s_{\la,\mu} \in A_0$), where the sum is locally finite in the sense that locally $s_{\la,\mu}=0$ for all but finitely many $\mu$. 
\end{sbprop}

\begin{pf} Let $(U_{\la})_{\la\in \La}$ be a locally finite open covering of $X$. 
  We prove the existence of a partition of unity subordinate to this open covering satisfying the required condition. 

  Since $X$ is paracompact, we can take locally finite open coverings $(V_{\mu})_{\mu\in I}$ and $(W_{\mu})_{\mu\in I}$ of $X$ for an index set $I$ and a map $a:I\to \La$ satisfying the following conditions (i)--(iii).

(i) $V_{\mu}\subset U_{a(\mu)}$ for all $\mu\in I$.

(ii) For each $\mu\in I$, the closure of $W_{\mu}$ in $X$ is contained in $V_{\mu}$.

(iii) $V_{\mu}$ is a closed analytic subspace of $\C^{n(\mu)}$ for some $n(\mu)\geq 0$. 

For each $\mu\in I$, by \cite{PP} Lemma 1.3.2 applied to $\C^{n(\mu)}$, there is an element $g_{\mu}\in A_0(V_{\mu})$ 
which has value $0$ at all points of $V_{\mu}\smallsetminus W_{\mu}$ and has values $\neq 0$ at all points of $W_{\mu}$. 
We extend $g_{\mu}$ to an element of $A_0(X)$ by defining its restriction to $X\smallsetminus \overline W_{\mu}$ to be zero. Then
   $\sum_{\mu} g_{\mu}^2$ 
has values $>0$ at every point of $X$, and hence it is the square of 
an invertible element $h$ of $A_0(X)$. 
  For each $\la$, let $s_{\la}= \sum_{\mu\in a^{-1}(\la)} \;(g_{\mu}h^{-1})^2\in A_0(X)$. Then $(s_{\la})_{\la}$ is a partition of unity subordinate to $(U_{\la})_{\la}$,   $s_{\la}$ is the  sum  of $(g_{\mu}h^{-1})^2$ , and $g_{\mu}h^{-1}$ has support in $U_{\la}$. 
\end{pf}

\begin{sbpara} We prove Theorem \ref{softhm}. Let $U$ be an open set of $X_{[:]}$ and let $(U_{\la})_{\la}$ be a locally finite open covering of $U$. We prove the existence of a partition of unity $1=\sum_{\la} s_{\la}$ subordinate to $(U_{\la})_{\la}$ with (P) (\ref{pos}). 

Take an open covering $(V_{\mu})_{\mu}$ of $X$ such that on each $V_{\mu}$, we have a chart $\cS\to M_{V_{\mu}}$ and generators $f_1^{(\mu)}, \dots, f_{m(\mu)}^{(\mu)}$ of $\cS$ such that $|f_j^{(\mu)}|<1$ for all $j$ and a closed immersion $V_{\mu}
\overset{\subset}\to\C^{n(\mu)}$ of complex analytic spaces. Then $(V_{\mu})_{[:]}$ is identified with a closed subset of $\C^{n(\mu)} \times \C^{m(\mu)^2}$ by Proposition \ref{A0toA1}. 
 (Note that the composition of the maps in $(*)$ in Proposition \ref{A0toA1} is a closed immersion of topological spaces if the map $X\to \C^n$ there is a closed immersion, as is seen from $\ref{chratio}$ $(2)$.) Since $X$ is paracompact, we can take such an open covering  $(V_{\mu})_{\mu}$ of $X$ which is locally finite.

Hence $U\cap (V_{\mu})_{[:]}$ is identified with a closed subset of an open set 
$Z'_{\mu}$ of $Z_{\mu}$. 
When we identify a sheaf on $U\cap (V_{\mu})_{[:]}$ as its direct image on $Z'_{\mu}$,   by Proposition \ref{A0toA1},  we have a  homomorphism from $A_0$ of  $Z'_{\mu}$ to  $A_k$ ($k=1, 2, 2*, 3$) of $U\cap (V_{\mu})_{[:]}$. 

 Take an open set $W_{\mu, \la}$ of $Z'_{\mu}$ such that $W_{\mu,\la}\cap (U \cap (V_{\mu})_{[:]})= U_{\la}\cap (V_{\mu})_{[:]}$. Then we have an open covering of 
$Z'_{\mu}$ consisting of $W_{\mu, \la}$ for all $\la$ and $W_{\mu, e}:= Z'_{\mu}\smallsetminus (U\cap (V_{\mu})_{[:]})$. 
  Take a partition of unity $1= (\sum_{\la} s_{\mu, \la}) + s_{\mu, e}$ in $A_0(Z'_{\mu})$ subordinate to this open covering 
such that each $s_{\mu,\lam}$ is a locally finite sum of the squares of elements in $A_0(Z'_{\mu})$ 
by  Proposition \ref{A0soft}.
  Then the pullback of this to $A_k(U\cap (V_{\mu})_{[:]})$ gives a partition of unity  $1= \sum_{\la} s_{\mu, \la}$ subordinate to the open covering $(U_{\la}\cap (V_{\mu})_{[:]})_{\la}$ of $U\cap(V_{\mu})_{[:]}$ with (P). 
  Note that a section of $A_k$ satisfies (P) if it is locally a finite sum of the squares of sections of $A_k$.

  On the other hand, by Proposition \ref{A0soft}, we have a partition of unity $1=\sum_{\mu} s_{\mu}$ ($s_{\mu}\in A_0(X)$) subordinate to $(V_{\mu})_{\mu}$ with (P).  We have $\sum_{\la} s_{\mu}s_{\mu, \la}=s_{\mu}$ in $A_k(U\cap (V_{\mu})_{[:]})$  and these $s_{\mu}s_{\mu, \la}$ extend by zero to elements in $A_k(U)$. Thus $1=\sum_{\mu, \la} s_{\mu}s_{\mu, \la}$ in $A_k(U)$. For each $\la$, let $s_{\la}= \sum_{\mu}  s_{\mu}s_{\mu,\la}$. Then $1=\sum_{\la} s_{\la}$, and this is a partition of unity subordinate to $(U_{\la})_{\la}$ with (P).
\end{sbpara}

\begin{sbprop}\label{soft3} Let $X$ and $U$ be as in Theorem $\ref{softhm}$. Let $k=1,2,2*,3$  and let $\cF$ be a sheaf of  $A_{k,U}$-modules on $U$.
\medskip

$(1)$ 
$H^m(U, \cF)=0$ for $m>0$.

$(2)$ Let $\tau_X$ be the canonical map 
$X_{[:]}^{\log}\to X_{[:]}$. 
Then $H^m(\tau_X^{-1}(U), A_k^{\log*}\otimes_{A_k} \cF)=0$ for $m>0$.
\end{sbprop}

\begin{pf} (1) follows from Theorem \ref{softhm}. (2) follows from (1) and Corollary \ref{c:pol} (2).
\end{pf}

\subsection{Differential forms}\label{ss:diff} 

\begin{sbpara}\label{diff1}
For $R=A_k$ (resp.\ $ A_k^{\log}$) with $k=2,2*,3$, we define the sheaf $R^1=R^1_X$ of differential forms on $X_{[:]}$ (resp.\ $X_{[:]}^{\log}$) as follows. 

In the case where $X$ is log smooth, define
$$R^1= R \otimes_{A_2^{\an}} A_2^{\an,1},$$
where $A_2^{\an,1}$ is as in \ref{Aan1}.

For an object $X$ of $\cB(\log)$, by using a local strict immersion $X\to Z$ into log smooth, we define $R^1_X$ as the cokernel of
$$I/I^2\overset{d}\to R_Z^1/IR_Z^1,$$
where $I=\text{Ker}(R_Z\to R_X)$. The fact that $R^1_X$ is independent of the choice of a local strict immersion and is defined globally is proved in the same way as in \ref{anindep}, \ref{indep}. 
\end{sbpara}

\begin{sbpara} For a morphism $X\to Y$ in $\cB(\log)$, we have an associated homomorphism  $R_X \otimes_{R_Y} R_Y^1\to R_X^1$. This is the evident map if $X$ and $Y$ are log smooth, and in general, it is defined by using strict immersions to log smooth and by the method in \ref{ringedmor}.

\end{sbpara}

\begin{sbpara}\label{diff2} 
Let $X\to Y$ be a morphism in $\cB(\log)$. 
  For $R$ as in \ref{diff1}, let the sheaf $R^1_{X/Y}$ be the cokernel of 
$$R_X\otimes_{R_Y}  R_Y^1\to R_X^1.$$
Define
$$R^p_{X/Y}= \bigwedge^p_{R_X}\; R^1_{X/Y}.$$
\end{sbpara}

\begin{sblem}\label{d:} Assume that $X$ is a log smooth fs log analytic space. 
  Let $R$ be as in $\ref{diff1}$. Then $d(R)\subset R^1$. 
\end{sblem}

\begin{pf} The case $R=A_2$ is clear.
  In the following, let $f \in M_X$ such that $|f|<1$. 
  By $$d\log(-\log|f|)=\frac{d\log|f|}{\log|f|},$$ the case $R=A_{2*}$ follows. 
  Since $d\log(f)$ belongs to $A_3^{\an,1}$, the case $R=A_3$ follows. 
  \end{pf}
  
  \begin{sbpara}
\label{d:1.5}
 Let $X$ be an object of $\cB(\log)$ and let $R$ be as in $\ref{diff1}$. 
Then we have a canonical map $d: R \to R^1$ which is obtained from the log smooth case Lemma \ref{d:} via  strict immersions into log smooth.

  \end{sbpara}

\begin{sblem}\label{d:2} Let $X$ be an object of $\cB(\log)$. 
 For $R=A_k^{\log}$ with $k=2*,3$, the $R$-module $R^1$ is generated by $d(R)$. 
\end{sblem}
\begin{pf}
The equality  in the proof of Lemma \ref{d:} and 
$$d\left(\frac{\arg(f)}{\log|f|}\right)=
\frac{d\arg(f)}{\log|f|}-
\frac{\arg(f)}{\log|f|}\cdot
\frac{d\log|f|}{\log|f|}$$
show that $d(R)$ includes 
$\frac{d\log|f|}{\log|f|}$ and $\frac{d\arg(f)}{\log|f|}$. 
  This proves Lemma \ref{d:2}. 
\end{pf}

\begin{sbpara} Let $X\to Y$ be  a morphism in  $\cB(\log)$ and let $R$ be as in \ref{diff1}. 
  By \ref{d:1.5}, we have $d: R \to R^1_{X/Y}$. This extends naturally to a de Rham complex
$$R^{\bullet}_{X/Y}= [R \overset{d}\to R^1_{X/Y}\overset{d}\to R^2_{X/Y}\overset{d}\to \cdots].$$

\end{sbpara}

 \begin{sbpara} Let $\omega^1_X$ be the sheaf of holomorphic differential forms with log poles. Then we have a unique $\cO_X$-homomorphism $\omega^1_X\to A_{3,\C}^{\an,1}$ which sends
 $d\log(f)$ to $d\log|f|+ id\arg(f)$ for $f \in M_X$ such that $|f|<1$.
 \end{sbpara}

\subsection{Log points}\label{ss:pt}

In the usual complex geometry, a point is a simple object. However, an fs log point in the complex log geometry has interesting geometry. 

In this Section \ref{ss:pt}, let $X$ be an fs log point. Let $x$ be the unique point of $X$.

\begin{sbpara}\label{ptb} Let $\R_{\geq 0}^{\mult}$ (resp.\ $\C^{\mult}$) be the set $\R_{\geq 0}$ (resp.\ $\C$) regarded as a multiplicative monoid. 

Let $\cT$ (resp.\ $\cT_{\triv}$)  be the set of all homomorphisms $M_{X,x}\to \C^{\mult}$ (resp.\ $M_{X,x}^{\gp}\to \C^\times$) of commutative monoids (resp.\ groups) whose restriction to $\cO_{X,x}^\times= \C^\times$ is the inclusion (resp.\ identity) map. 

Then $\cT$ is regarded as a complex analytic space and $\cT_{\triv}$ is regarded as an open set of $\cT$.
Endow $\cT$ with the log structure consisting of holomorphic functions whose restrictions to $\cT_{\triv}$ are invertible.

If we fix  a chart 
$\cS\to M_X$ of the log structure $M_X$ with $\cS$ being an fs monoid such that $\cS\overset{\cong}\to (M_X/\cO_X^\times)_x$, $\cT$ is identified with the analytic toric variety  $\Hom(\cS, \C^{\mult})=  \Spec(\C[\cS])^{\an}$ and the log structure of $\cT$ coincides with the one given by $\cS\to \cO_{\cT}$ and $\cT_{\triv}$ is identified with the torus $\Hom(\cS^{\gp}, \C^\times)$.

Hence we call $\cT$ the {\it canonical toric variety associated to the fs log point $X$}.

The space $\cT^{\log}$ is identified with the set of homomorphisms $ M_{X,x}\to \R/\Z \times \R^{\mult}_{\geq 0}$ whose restrictions to $\cO_{X,x}^\times =\C^\times$ coincide with $z\mapsto ((2\pi)^{-1}\arg(z), |z|)$.

\end{sbpara}

\begin{sbpara}\label{ptab} We regard  $M_{X, x}$ as a constant submonoid of $M_\cT$, as follows. For $q\in M_{X,x}$, we regard $q$ as the section of  $M_{\cT}$ which  is the holomorphic function on $\cT$ whose value $q(a)\in \C$ at $a\in \cT$ is $a(q)$,  where in $a(q)$, we regard $a$ as a homomorphism $M_{X,x}\to \C^{\mult}$. 

\end{sbpara}

\begin{sbpara}\label{xinT}
We have a canonical morphism $X\to \cT$ of fs log analytic spaces which sends the unique point $x$ of $X$ to the point of $\cT$ corresponding to the structural map $M_{X,x}\to \cO_{X,x}=\C$ and for which  the pullback map $\cO_{\cT, x}\to \cO_{X, x}=\C$ is the evaluation $f\mapsto f(x)$ and the pullback map $M_{\cT,x}= M_{X,x}\times^{\C^\times} \cO_{\cT, x}^\times \to M_{X,x}$ is $(q, f)\mapsto f(x)q$.

\end{sbpara}

\begin{sbpara}\label{XT1}
Let $\cO_{X, \cT}^{\log}$ be the sheaf of  subrings of $\cO_{\cT}^{\log}$ on $\cT^{\log}$ generated over $\C$ algebraically  by logarithms of elements of $M_{X,x}^{\gp}$. 

Then $\cO_{X, \cT}^{\log}$ is a locally constant sheaf and its inverse image on $X^{\log}$ is $\cO_X^{\log}$. 

\end{sbpara}

\begin{sbpara}\label{pta} Let $t\in X^{\log}$. As in \cite{KU} 2.4.6, 
let $\spe(t)$ be the set of all ring homomorphisms $\cO_{X,t}^{\log}\to \C$ whose restriction to $\cO_{X,x}=\C$ is the identity map.

We have a canonical map $\spe(t) \to \cT_{\triv}$ which sends $s\in \spe(t)$ to the homomorphism $M^{\gp}_{X,x}\to \C^\times\;;\;q\mapsto \exp(s(\log(q)))$.

Let $s\in \cT_{\triv}$, $t\in X^{\log}$. If we fix a path between $s$ and $t$ in $\cT^{\log}\supset X^{\log}$,  the stalks of $\cO^{\log}_{X, \cT, s}$ and $\cO^{\log}_{X,\cT,t}=\cO^{\log}_{X,t}$ of the locally constant sheaf $\cO^{\log}_{X,\cT}$ are identified through this path, and the evaluation map $\cO^{\log}_{X,\cT,s}\subset \cO_{\cT, s}\to \C\;;\; f\mapsto f(s)$ is regarded as an element of $\spe(t)$. The image of the last element in $\cT_{\triv}$ under the above canonical map coincides with $s$. When we vary the path, we obtain all elements of $\spe(t)$ with image $s$ in $\cT_{\triv}$.

\end{sbpara}

\begin{sbpara}\label{XT2} 

Let $A^{\an}_{1, X,\cT}$ be the sheaf of subrings of $A^{\an}_{1, \cT}$ consisting of real analytic functions of $(\log|f|/\log|g|)^{1/2}$, where $f,g\in M_{X,x}$ such that $|f|<1$, $|g|<1$ and such that $\log|f|/\log|g|$ does not have value $\infty$.

Let $A^{\an}_{3, X,\cT}$ be the sheaf of subrings of $A^{\an}_{3, \cT}$ generated over  $A^{\an}_{1,X, \cT}$ algebraically by $\log|q|$ for $q\in M_{X,x}$ such that $|q|<1$.

For $k=1,3$, let $A^{\an,\log}_{k, X,\cT}$ (resp.\ $A^{\an,\log*}_{k, X, \cT}$) be the sheaf of subrings of $A^{\an,\log}_{k, \cT}$ (resp.\ $A^{\an,\log*}_{k, \cT}$) generated over the inverse image of $A^{\an}_{k,X,\cT}$  algebraically by $\arg(q)/\log|q|$ for $q\in M_{X,x}$, $|q|<1$ (resp.\ by $\arg(q)$ for $q\in M^{\gp}_{X,x}$). We have $A^{\an,\log}_{3,X,\cT}= A^{\an,\log*}_{3,X,\cT}$.

For $k=1,3$, the inverse image of $A^{\an}_{k, X,\cT}$ on $X_{[:]}$ is $A^{\an}_{k,X}$, and 
the inverse image of $A^{\an,\log}_{k, X,\cT}$ (resp.\ $A^{\an,\log*}_{k, X, \cT}$) on $X^{\log}_{[:]}$ is $A^{\an,\log}_{k,X}$ (resp.\ $A^{\an,\log*}_{k,X}$).  

For $p\in X_{[:]}$ (resp.\ $X^{\log}_{[:]}$, resp.\ $X^{\log}_{[:]}$) 
and $f\in A^{\an}_{k,X,p}$ with $k=1,3$ (resp.\ $A^{\an,\log}_{3,X,p}$, resp.\ $A^{\an,\log*}_{k, X,p}$ with $k=1,3$), $f$ 
 extends uniquely to a section of $A^{\an}_{k,X, \cT}$ (resp.\ $A^{\an,\log}_{3, X,\cT}$, resp.\ $A^{\an,\log*}_{k,X,\cT}$)  on some open neighborhood  of $p$ in $\cT_{[:]}$ (resp.\ $\cT^{\log}_{[:]}$, resp.\ $\cT^{\log}_{[:]}$).  

We have $$\cO_{X,\cT}^{\log}\subset A^{\an,\log}_{3,X,\cT,\C}.$$

\end{sbpara}

\begin{sbprop}\label{pt13} The map $A^{\an}_{3,X} \to A_{3,X}$ is injective. If we regard $A^{\an}_{3,X}$ as a subsheaf of $A_{3,X}$ via this injection, we have  $A^{\an}_{3,X}\cap A_{2,X}=A^{\an}_{1,X}$.

\end{sbprop}

\begin{pf}

We first prove  the following part of Proposition \ref{pt13}: Let $p\in X_{[:]}$. If $f\in A^{\an}_{3,X,p}$ and if the image of $f$ in $A_{3,X,p}$ belongs to $A_{2,X,p}$, then $f\in A^{\an}_{1,X,p}$.

Let $I$ be the ideal of $A_{3,\cT}$ locally generated by $\Re(q)$, where $q\in \cO_{\cT}$ such that the image of $q$ in $\cO_X$ is $0$.
By the assumption on $f$, 
for some  open neighborhood $V$ of $p$ in $\cT_{[:]}$, $f$ extends to 
$\tilde f \in A^{\an}_{3,X,\cT}(V)$  and there are $g_1\in I(V)$ and $g_2\in A_{2,\cT}(V)$  such that $\tilde f= g_1+g_2$ in $A_{3,\cT}(V)$. 

Assume that $f$ does not belong to $A^{\an}_{1,X,p}$. 

We will prove that for some $\epsilon>0$,  there is a continuous map $\varphi:[0,\epsilon)\to \cT_{[:]}$ 
 having the following properties (i)--(iv). (i) $\varphi(0)=p$. (ii) $\varphi$ sends $(0,\epsilon)$ into $V\cap \cT_{\triv}$. (iii) The pullback of $\tilde f$ on $(0,\epsilon)$ under $\varphi$ does not extend to $[0,\epsilon)$ as a continuous function. (iv) The pullback of $g_1$ on $(0,\epsilon)$ extends to $[0,\epsilon)$ as a continuous function. Then since  $g_2$ is continuous on $V$, the pullback of $g_1+g_2$ to $(0,\epsilon)$ extends to a continuous function on $[0,\epsilon)$, and we have a contradiction.

We construct $\varphi$ as follows.

 By using the presentation of $A^{\an}_{1,X,p}$ as convergent power series in Corollary \ref{lgpt}, we have that for some $\ell$ ($1\leq \ell \leq m$), (briefly speaking) $T_{\ell,1}$ appears in the denominator of $f$. That is, precisely speaking,  $f$ is expressed as $(h_1+h_2)/h_3$ with $h_1, h_2, h_3\in A^{\an}_{1,X,p}$ such that $h_3$ is a product of powers of $T_{j,1}$  ($1\leq j\leq m$),  $h_2$ and $h_3$ are divisible by $T_{\ell,1}$, and  $h_1=h_1((T_{\la})_{\la\in \La})$  is a non-zero convergent series in $T_{\la}$ for $\la\in \La=\{(j,k)\;|\; 1\leq j\leq m, 1\leq k\leq r(j), (j,k)\neq (\ell,1)\}$. 
  Take   $c_{\la}\in \R_{>0}$ ($\la\in \La$) such that the convergent power series $h_{1,t}:=h_1((c_{\la}t)_{\la\in \La})$ in one variable  $t$ is non-zero. Take an integer $d>0$ which is strictly bigger than the order of series $h_{1,t}$. Then writing $f=f(T_{\ell,1}, (T_{\la})_{\la\in \La})$, the Laurent series $f_t:=f(t^d, (c_{\la}t)_{\la \in \La})$ in $t$ has a pole. Let $\epsilon>0$ be a sufficiently small number, and let $\varphi:[0,\epsilon) \to \cT_{[:]}$ be the continuous map which sends $t\in \R_{\geq 0}$ to the unique point of $\cT_{[:]}$ at which the value of $(\log|f_{j+1,1}|/\log|f_{j,1}|)^{1/2}$ is $c_{j,k}t$ if $j\neq \ell$ and is $t^d$ if $j=\ell$ ($f_{m+1,1}$ denotes $e^{-1}$) and the value of 
$(\log|f_{j,k}|/\log|f_{j,1}|)^{1/2}-a_{j,k}$ with $k\neq 1$ is $c_{j,k}t$, and the values of $f_{j,k}$ belong to $\R_{\geq 0}$ for all $j,k$. Then the pullback of $\tilde f$ to $[0, \epsilon)$ is $f_t$. This map $\varphi$ has the desired properties.

The remaining part of Proposition \ref{pt13} follows from this part easily. Let $p\in X_{[:]}$ and let $t$ be a non-zero element of the maximal ideal of the local integral domain $A^{\an}_{1,X,p}$. Since $A^{\an}_{1,X,p}$ is a Noetherian local ring, $\bigcap_n \; t^n A^{\an}_{1,X,p}=0$. If $f\in A^{\an}_3$ goes to $0$ in $A_3$,  $t^{-n}f$ for $n\geq 0$ goes to $0$ and hence to $A_2$.  By the above part of Proposition \ref{pt13} which we have proved, $t^{-n}f\in A^{\an}_{1,p}$. Hence $f\in \bigcap_n t^nA^{\an}_{1,p}=0$.

\end{pf}

\section{Local coordinates of log smooth morphisms in log $C^{\infty}$ geometry}
\label{s:ls}

If $X\to Y$ is a smooth morphism of complex analytic spaces, then locally on $X$ and on $Y$, $X$ is isomorphic to an open subspace of $Y \times \Delta^n$  over $Y$ for some $n\geq 0$. 

In the log world, for a log smooth saturated morphism, we give an analogue for this fact in the log $C^{\infty}$ category (using $A_2$ and $A_2^{\log}$) for the morphisms $X_{[:]}\to Y_{[:]}$ and $X^{\log}_{[:]}\to Y^{\log}_{[:]}$. 

  Recall that the underlying morphism of the original log smooth saturated morphism $X\to Y$ can be far from being a smooth morphism, and, though it is well-known that taking the space of arguments improves the morphism topologically (in fact, $X^{\log} \to Y^{\log}$ is topologically smooth (\cite{NO})), it does not improve it differential geometrically. 
  We show here that taking the space of ratios improves the morphism differential geometrically: Theorems \ref{SEisS} and \ref{SEmore}.
  Note that there is a description of absolute value portions of space of ratios in \cite{KNU4} Propositions 4.2.14 and 4.2.19.
  We add to it here the argument portions.

\subsection{Local properties of log smooth saturated morphisms}
\label{ss:locls} 

  In this section, we review necessary facts on the local properties of log smooth and saturated morphisms. 

\begin{sbpara}
  For the definition of log smoothness, see \cite{IKN} Definition (2.2) (3). 
  For its characterization by a toroidal chart, see \cite{IKN} Theorem (2.3), which is adopted as the definition in \cite{KU} 2.1.11.
  The next proposition is a variant of this characterization including the saturatedness.  
\end{sbpara}

\begin{sbprop}
\label{p:lse}
  Let $f: X \to Y$ be a morphism of fs log analytic spaces.
  Then $f$ is log smooth and saturated if and only if, locally on $X$ and on $Y$,  there is a chart by a saturated homomorphism of fs monoids $\cS' \to \cS$ for $f$ such that the induced morphism $X \to Y \times_{\Spec \bC[\cS']}\Spec \bC[\cS]$ is a strict open immersion. 
  Further, if this condition is satisfied, then, for any $x \in X$, the induced homomorphism $\overline M_{Y,f(x)} \to \overline M_{X,x}$ is saturated.
\end{sbprop}

\begin{pf}
  This is by the basic properties of saturated homomorphisms (\cite{T} Propositions I.3.16 and I.3.18).  Cf.\ \cite{IKN} Proposition (A.3.3).
\end{pf}

\begin{sblem}
\label{l:Qint}
  Let $h\colon \cS' \to \cS$ be a local ({\rm \cite{O}} {\rm Chapter I}, $1.4$, {\rm p.15}) and saturated homomorphism. 
  Let $A\supset B$ be faces of $\cS$, and let $A'\supset B'$ be their pullbacks, which are faces of $\cS'$. 
  Then the induced homomorphism $A'{}^{\gp}/B'{}^{\gp} \to A{}^{\gp}/B{}^{\gp}$ is a split injection.
\end{sblem}

\begin{pf} 
  We may assume that $A=\cS$. 
  By \cite{T} Propositions I.3.16 and I.3.18, the induced homomorphism $\cS'/B' \to \cS/B$ is saturated. 
  Hence, we may assume further that $B=B'=\{1\}$. 
  Then, $h$ is injective because it is a local and exact homomorphism from a sharp fs monoid to an fs monoid.
  Since the cokernel of $h^{\gp}$ is torsion-free by \cite{IKN} Lemma (A.4.1), the conclusion follows.
\end{pf}

\begin{sbprop}
\label{p:fiber}
  Let $f: X \to Y$ be a log smooth and saturated morphism of fs log analytic spaces. 
  Then the following holds. 

$(1)$ The set of points where $f$ is strict is open and dense in each fiber.

$(2)$ Any connected component of each fiber is equi-dimensional.

$(3)$ The dimensions of the fibers are locally constant on $Y$.
\end{sbprop}

\begin{pf}
  (1) is by \cite{NO} Proposition 5.8 (1) and \cite{IKN} Lemma (A.4.1). 

  (2) is proved similarly to \cite{N2} Proposition 14.3. 

  (3) By (1), we may assume that $f$ is strict. 
  Then the statement is classical. 
  See also \cite{FN} Proposition 2.8. 
\end{pf}

\subsection{Log $C^{\infty}$ nature of log smooth saturated morphisms}
\label{ss:infls} 

We show that a log smooth saturated morphism $X\to Y$ locally has  good coordinate functions for the log $C^{\infty}$ (more precisely, $A_2$) structures of $X$ and $Y$.

\begin{sbthm}
\label{SEisS}
Let $X\to Y$ be a log smooth and saturated morphism. 
Let $p \in X_{[:]}$. 
Let $R=A_2$ or $A_3$. 
Then there are an open neighborhood $X' \subset X$ of the image of $p$ and a morphism $X'\to Y \times \Delta^n$ over $Y$, where $\Delta^n$ has the log structure defined by $q_1, \dots, q_{n'}$ for some $n' \leq n$ ($q_j$ are the standard coordinate functions), such that the induced $X'_{[:]}\to (Y \times \Delta^n)_{[:]}$ is locally an isomorphism of ringed spaces with the sheaves of rings $R$ around $p$ and such that the induced $X^{\prime\log}_{[:]}\to (Y \times \Delta^n)^{\log}_{[:]}$ is locally an isomorphism of ringed spaces with sheaves of rings  $R^{\log}$ around the inverse image of $p$.
\end{sbthm}

We will call $q_j$ ($1\leq j\leq n$, $q_j\in M_X$ for $1\leq j\leq n'$, $q_j \in \cO_X$ for $n'<j\leq n$)  which give the local morphism from $X'$  to $\Delta^n$ in Theorem \ref{SEisS}, {\it  local coordinate functions} of $X$ over $Y$ at $p$.

\begin{pf}
  Let $x \in X$ and $y \in Y$ be the images of $p$. 
  Then 
by using Proposition \ref{p:lse}, as in the beginning of the proof of the holomorphic log Poincar\'e lemma \cite{IKN} Theorem (5.1) (i.e., the part after \cite{IKN} Lemma (5.3)), we may assume that $X\to Y$ with $x \in X$ is $(u,v): \Spec(\C[\cS])^{\an} \times \C^b \to \Spec(\C[\cS'])^{\an} \times \C^a$  
($0\leq a\leq b$) with the origin (i.e., the pair of the origin 
of the toric variety and the origin of $\bC^b$). 
  Here $u:\Spec(\C[\cS])^{\an} \to \Spec(\C[\cS'])^{\an}$ is the morphism induced by a local saturated homomorphism 
  $\cS'\to \cS$ of sharp fs monoids, and $v$ is the projection $(z_1, \dots, z_b)\mapsto (z_1, \dots, z_a)$. 

In the above reduction step to the case where $Y$ is log smooth,  we take care of the following delicate point which did not appear in \cite{IKN}. If $Y\to \tilde Y$ is a strict immersion into log smooth and if $X\to Y$ comes from a log smooth saturated $\tilde X\to \tilde Y$, $R_X=R_{\tilde X}\otimes_{R_{\tilde Y}} R_Y$ is true for $R=A_3, A_3^{\log}$ but may not be true for $R=A_2, A_2^{\log}$. The cases $R=A_3, A_3^{\log}$ are hence reduced to the case where $Y$ is log smooth, but for $R=A_2, A_2^{\log}$, first we only get that the map from the inverse image of $R_{Y\times \Delta^n}$ to $R_{X'}$  is surjective. But because $A_2\subset A_3$ and $A_2^{\log}\subset A_3^{\log}$, the cases $R=A_3$, $A_3^{\log}$ show that this surjection is in fact an isomorphism.

  Let $M^{(j)}$ $(0\le j \le m)$ be as in \ref{m,f_jk} at $p\in X_{[:]}$. 
  Let $\cS^{(j)}$ be the inverse image of $M^{(j)}$ by $\cS \to M_{X,x}$, and 
let $\cS^{\prime(j)}$ be the inverse image of $\cS^{(j)}$ by $\cS' \to \cS$. 
  Take $f_{j,k}$ 
as in \ref{m,f_jk} satisfying that for each $j$ ($1 \le j \le m$), 
for some $c(j)$ ($0\leq c(j)\leq r(j))$, 
the $f_{j,k}$ for $1\leq k\leq c(j)$ form a 
$\bQ$-base of $\bQ \otimes (\cS^{\prime(j-1)}/\cS^{\prime(j)})^{\gp}$, 
and the $f_{j,k}$ for $c(j)< k\leq r(j)$ give a $\bZ$-base of the cokernel of $(\cS^{\prime (j-1)}/\cS^{\prime (j)})^{\gp} \to (\cS^{(j-1)}/\cS^{(j)})^{\gp}$. 
  This is possible by Lemma \ref{l:Qint}. 
  Let $n'=\dim \cS - \dim \cS'$ and $n=n'+b-a$. 
  Let $\{q_1, \dots, q_{n'}\}$ be the set of elements of $M_{X,x}$ which 
coincides with the set of $f_{j,k}$ ($1\leq j\leq m, c(j)<k \leq r(j)$). 
  Let 
$q_j$ for $n'< j \leq n$ be 
the $(j+b-n)$-th coordinate functions of $\C^b$. 
  Then $(q_j)_{1\leq j\leq n}$ gives a morphism $f:X\to Y\times  \Delta^n$ locally on $X$.

  We prove that this $f$ has the desired property around $p$. 
  For simplicity, in the following, we assume $a=b=0$, so $n=n'$. 
  The general case is similar. 
  We redefine $Y$ by $Y \times \Delta^n$ and $\cS'$ by $\cS'\times \bN^n$ (thus the saturatedness of $\cS' \to \cS$ may be lost), and we prove that the morphisms $f_{[:]}$ and $f^{\log}_{[:]}$ ($f$ is now $X \to Y$) of ringed spaces are local isomorphisms around $p$. 
  We define the new $\cS^{\prime(j)}$ again by the inverse image of $\cS^{(j)}$  by $\cS' \to \cS$, which is isomorphic to the product of the old one and $\bN^{\sum_{j \le j' <m}(r(j')-c(j'))}$. 
  Then $(\cS^{\prime(j-1)})^{\gp}/(\cS^{\prime(j)})^{\gp} \to (\cS^{(j-1)})^{\gp}/(\cS^{(j)})^{\gp}$ are isomorphisms. 
  Hence, $(\cS')^{\gp} \to \cS^{\gp}$ is an isomorphism. 
  Thus we reduce to the following. 
\end{pf}  

\begin{sbprop} 
  Let $h: \cS' \to \cS$ be a homomorphism of sharp fs monoids such that $h^{\gp}: (\cS')^{\gp} \to \cS^{\gp}$ is an isomorphism. 
  Let 
$\cS=\cS^{(0)}\supsetneq \cS^{(1)}\supsetneq \cdots \supsetneq \cS^{(m)}=\{1\}$ be a sequence of faces of $\cS$. 
  Let $\cS^{\prime(j)}$ be the inverse image of $\cS^{(j)}$ by $h$ $(0 \le j \le m)$. 
  Assume that $(\cS^{\prime(j-1)})^{\gp}/(\cS^{\prime(j)})^{\gp} \to (\cS^{(j-1)})^{\gp}/(\cS^{(j)})^{\gp}$ $(1 \le j \le m)$ are isomorphisms. 
  Let $\Phi=\{\cS^{(j)}\,|\,0 \le j \le m\}$ and 
$\Phi'=\{\cS^{\prime(j)}\,|\,0 \le j \le m\}$.

  Let $X$ (resp.\ $Y$) be the inverse image of an open neighborhood $X_0$ (resp.\ $Y_0$) of the origin of 
$\Hom(\cS,\bR_{\ge0}^{\mult})$ (resp.\ $\Hom(\cS',\bR_{\ge0}^{\mult})$) by the projection $\Spec(\C[\cS])^{\an} \to \Hom(\cS,\bR_{\ge0}^{\mult})$ (resp.\ $\Spec(\C[\cS'])^{\an} \to \Hom(\cS',\bR_{\ge0}^{\mult})$).
  Assume that $h$ induces a morphism $X \to Y$, and that for any $f \in \cS \smallsetminus \{1\}$, we have $|f|<1$ on $X$. 
  Then the induced morphism
$X_{[:]}(\Phi)\to Y_{[:]}(\Phi')$ (see {\rm \cite{KNU4} 4.2.11} for the definitions of $X_{[:]}(\Phi)$ and $Y_{[:]}(\Phi')$) with the rings $A_2$ (resp.\ $A_3$)  
is an open immersion. 

  Further, let $X^{\log}_{[:]}(\Phi)$ be the inverse image of $X_{[:]}(\Phi)$ by $X^{\log}_{[:]} \to X_{[:]}$ and let 
$Y^{\log}_{[:]}(\Phi')$ be the inverse image of $Y_{[:]}(\Phi)$ by $Y^{\log}_{[:]} \to Y_{[:]}$.
  Then the induced morphism
$X_{[:]}^{\log}(\Phi)\to Y_{[:]}^{\log}(\Phi')$ with the rings $A_2^{\log}$ (resp.\ $A_3^{\log}$)  
is an open immersion. 
\end{sbprop}

\begin{pf}
  Let $P$ be as in \cite{KNU4} 4.2.12 for $X$ and $\cS$, and $P'$ for $Y$ and $\cS'$.
  Then $P\to P'$ is an open immersion because $P$ and $P'$ are open sets of the same space by \cite{KNU4} Proposition 4.2.19. 
  Let $P_X$ be the inverse image of $X_0$ by $P \to \Hom(\cS,\bR_{\ge0}^{\mult})$ and $P'_Y$ be that of $Y_0$ by $P' \to \Hom(\cS',\bR_{\ge0}^{\mult})$.
  Then \cite{KNU4} Proposition 4.2.14 implies that $X_{[:]}(\Phi)\to P_X$ and $Y_{[:]}(\Phi')\to P'_Y$ are proper because $\C\to \R_{\geq 0}; z \mapsto |z|$ is proper.
  Hence $X_{[:]}(\Phi)\to Y_{[:]}(\Phi') \times_{P'_Y} P_X$ is a continuous map between proper spaces over $P_X$.
  We claim that it is also bijective, and hence a homeomorphism. 
  Then their $A_2$ (resp.\ $A_3$) are the same because, by the assumption, we can choose a common set of $f_{j,k}$s as in \ref{m,f_jk} for $X$ and $Y$.
  Thus we conclude that $X_{[:]}(\Phi)\to Y_{[:]}(\Phi')$ is an open immersion. 

  To prove the above claim, we give an inverse map.
  Let $(a,b) \in Y_{[:]}(\Phi') \times_{P'_Y} P_X$. 
  Let $\cS^{(j)}$ be the inverse image of $\R^{\times}$ by the homomorphism $b':\cS \to \bR_{\ge0}^{\mult}$ induced by $b$.
  Since $(\cS^{\prime(j)})^{\gp} \to (\cS^{(j)})^{\gp}$ is an isomorphism by assumption, there is a unique homomorphism $s:\cS \to \bC^{\mult}$ which coincides on $\cS^{\prime(j)}$ with the homomorphism induced by $a$ and which is zero outside $\cS^{(j)}$. 
  Then we show that the image of $s$ in $\Hom(\cS,\bR_{\ge0}^{\mult})$ coincides with $b'$.
  From this, we have $s \in X$. 
  By \cite{KNU4} Proposition 4.2.14, there is a unique point of $X_{[:]}(\Phi)$ which is sent to $s$ and $b$. 
  It is straightforward to see that it gives an inverse map. 

  As for the latter part of the proposition, we have 
$X^{\log}_{[:]}(\Phi)= \Hom(\cS^{\gp}, \R/\Z)\times P_X$ and $Y^{\log}_{[:]}(\Phi')= \Hom(\cS^{\prime\gp}, \R/\Z)\times P'_Y$, which are shown similarly as in \cite{KNU4} Proposition 4.2.14.
  From this, we have $X^{\log}_{[:]}(\Phi)=Y^{\log}_{[:]}(\Phi') \times_{P'_Y} P_X$ similarly.
\end{pf}

\begin{sbrem}
  The definition of the upper horizontal arrow $S_{[:]}(\Phi) \to P$ in \cite{KNU4} Proposition 4.2.14 is wrong. 
  The corrected one is as follows. 
  It sends $(s,r)\in S_{[:]}(\Phi)$ $(s\in S$, $r\in R((M_S/\cO_S^\times)_s))$ to $((t_j)_{1\le j\le n}, (h_j)_{0 \le j \le n})$ defined as follows. 
  Let $\cS^{(j)} \subset \cS$ be the inverse image of $\cO^{\times}_{S,s}$ by $\cS\to M_{S,s}$.
  Let $t_k=0$ if $k \le j$. 
  Let $t_k= \log|q_{k+1}(s)|/\log|q_k(s)|$ if $j<k<n$.
  If $j<n$, let $t_n= -1/\log|q_n(s)|$.
  Let $h_k(f)=r(f, q_{k+1})$ for $0\leq k<j$, 
where $f$ and $q_{k+1}$ are identified with their images in $(M_S/\cO^{\times}_S)_s$.
  Let $|\alpha(s)|\colon \cS^{(j)} \to \bR_{>0}$ be the homomorphism induced by the chart and the homomorphism $\cO_{S,s}^{\times}\to \bR_{>0}; f \mapsto |f(s)|$. 
  Let $h_k(f)=\log|\alpha(s)|(f)/\log|q_{k+1}(s)|$ for $j \leq k<n$. 
  Let $h_n(f) = -\log|\alpha(s)|(f)$.

  Additionally, the definition of the right vertical arrow in that proposition is miswritten: 
it is not \lq\lq$a\mapsto a(f)$ $(f \in \cS)$'' but $(t,h)\mapsto a$ in \cite{KNU4} Lemma 4.2.13. 
\end{sbrem}

\begin{sbpara}\label{SEEx} Examples of Theorem \ref{SEisS}.  
 
(1) Let $Y$ be the point $\Spec(\C)$ and $X=\Spec(\C[\cS])^{\an}$, where $\cS$ is the fs monoid with the generators $q_1, q_2, q_3, q_4$ and with the relation $q_1q_3=q_2q_4$.
  Then $X$ at the origin and $\Delta^3$ at the origin are not locally homeomorphic. 
But  $X_{[:]}$ is locally isomorphic to $(\Delta^3)_{[:]}$ as a space with $A_2$ as Theorem \ref{SEisS} tells. 
Here $\Delta^3$ has the log structure defined by the coordinate functions. In fact,  consider the open set $U$ of $X_{[:]}$ corresponding to the sequence $\cS=\cS^{(0)}\supsetneq \cS^{(1)}\supsetneq \cS^{(2)}\supsetneq \cS^{(3)}=\{1\}$ of faces of $\cS$, where $\cS^{(1)}$ is generated by $q_2,q_3$, $\cS^{(2)}$ is generated by $q_3$. 
  (All points of $X_{[:]}$ belong to one of the open sets defined similarly to $U$ by changing the variables.) 
  Then $(q_1, q_2, q_3): X\to \Delta^3$ induces an open immersion from $U$ to the standard open set of $(\Delta^3)_{[:]}$.

(2) Let $Y=\Delta$ with the coordinate function $q$, and $X=\Delta^2$ with $q_1, q_2$. 
  Define the log structure on $Y$ and $X$ by $q$ and $q_1, q_2$, respectively. 
  Let $X \to Y$ be the log smooth and saturated morphism defined by $q_1q_2=q$.
  (See \ref{badE2} for a related example.)
  Then, $q_1$ gives a morphism $X \to \Delta \times Y$ which induces an open immersion from the place of $X_{[:]}$ where $y_1/y_2$ does not have value $\infty$ to the standard open set, 
and $q_2$ gives another morphism $X \to \Delta \times Y$ which induces an open immersion from the place of $X_{[:]}$ where $y_2/y_1$ does not have value $\infty$ to the standard open set, where $|q_j|=e^{-2 \pi y_j}$ $(j=1,2)$.  
  The image of the former open immersion consists of the points where $y/y_1 >1$, and that of the latter consists of the points where $y/y_2>1$, where $|q|=e^{-2 \pi y}$.

(3) An example of an explicit description (in the case $n=1$) of a morphism $X\to Y \times \Delta$ which gives a local isomorphism on the space of ratios is given in \ref{E[:]} in the case where $Y$ is the standard log point and $X$ is a degenerate elliptic curve over $Y$. 
\end{sbpara}

About differential forms, we have the following.

\begin{sbcor}\label{SELF}  Let $X\to Y$ be a log smooth and saturated morphism of fs log analytic spaces. 
Then for  $R=A_k$ or $R= A_k^{\log}$ with $k=2,2*,3$, 
the $R$-module $R^1_{X/Y}$ ($\ref{diff2}$) is locally free of rank $2n$, where $n$ is the relative dimension of $X$ over $Y$ (cf.\ Proposition $\ref{p:fiber}$ $(3)$). 
\end{sbcor}

\begin{pf} 
Let $q_j$ ($1\leq j\leq n$) be local coordinate functions of $X$ over $Y$ (Theorem \ref{SEisS}). 
  Then we can prove that $dx_j/y_j$ and $dy_j/y_j$ for $1\leq j\leq n'$ ($q_j=e^{2\pi i(x_j+iy_j)}$ with $x_j$ and $y_j$ being real) together with $d\Re(q_j)$ and $d\Im(q_j)$ for $n'<j\leq n$ form a base of $R^1_{X/Y}$. 
  To see it, we may assume that $Y$ is also log smooth. 
  Then,  by the description of $A_2^{\an,1}$ given in Proposition \ref{difbase1}, 
we can take a base of $R^1_X$ and a base of $R^1_Y$ 
such that the latter is a subset of the former, and such that the complement of the latter in the former coincides with the above $dx_j/y_j$, $dy_j/y_j$, $d\Re(q_{j'})$, $d\Im(q_{j'
})$ $(1 \leq j \leq n'<j' \leq n)$. 
\end{pf}

\begin{sbthm}\label{SEmore} 
  Let $f:X\to Y$ be a log smooth and saturated morphism of fs log analytic spaces. 
  Assume that its relative dimension is constant $n$ in the strong sense that the dimension of any connected component of any fiber is $n$ (cf.\ Proposition $\ref{p:fiber}$ $(2)$, $(3)$). 
  Then, there are an open covering of $X^{\log}_{[:]}$ whose each member  $V$ has the following property. There are an integer $n'$ such that $0\leq n'\leq n$, $q_j\in M_X(V)$ ($1\leq j\leq n'$) with $|q_j|<1$, $s_j\in  A_1^{\an}(V)$ ($1\leq j\leq n'$), $q_j\in \cO_X(V)$ for $n'<j\leq n$, 
and an integer $a$ with $0\leq a \leq n'$ satisfying the following conditions$:$ The values of $s_j$ for $1\leq j\leq a$ are $\geq 0$, 
and we have an open embedding of topological spaces
$$(f^{\log}_{[:]}, (x_j)_{1\leq j\leq n'}, (s_j)_{1\leq j\leq a}, (s_j)_{a<j\leq n'}, (\Re(q_j))_{n'<j\leq n}, (\Im(q_j))_{n'<j\leq n}):\;\:$$ $$V\to Y^{\log}_{[:]}\times (\R/\Z)^{n'} \times \R_{\geq 0}^a\times \R^{n'-a}\times \R^{n-n'} \times \R^{n-n'},$$
where  $x_j=(2\pi)^{-1}\arg(q_j)$ for $1\leq j\leq n'$. 
\end{sbthm}

\begin{sbpara}\label{verrat} 
  Before proving Theorem \ref{SEmore}, we will give some remarks (Remark \ref{r:rounding}). 
  Here is a preparation to state them. 

 Let $f: X\to Y$ be a morphism in $\cB(\log)$. 

Recall that we say that $f$ is vertical if for each $x\in X$ and $g \in  M_{X,x}$, there exists an $h\in M_{Y, f(x)}$ such that $hg^{-1}\in M_{X, x}$ in $M^{\gp}_{X,x}$. 

 For $p\in X_{[:]}$ or $p\in X^{\log}_{[:]}$, we say that $p$ is {\it vertical over $Y$} if for each $g\in M_{X,p}$ such that $|g|<1$, there exists an $h\in M_{Y,p}$ such that $|h|<1$ and such that the value of $\log|g|/\log|h|$ at $p$ is not $\infty$. 

  We list some basic properties of the verticality in this sense. 

\medskip

\ref{verrat}.0. For $p\in X^{\log}_{[:]}$,  $p$ is vertical over $Y$ if and only if the image of $p$ in $X_{[:]}$ is  vertical  over $Y$. 

\medskip

\ref{verrat}.1. If $f$ is vertical, then all points of $X_{[:]}$ are vertical over $Y$. 

\medskip

In fact in this case, if $p\in X_{[:]}$ and $g\in M_{X,p}$, there exists an $h\in M_{Y, p}$ and $g_1\in M_{X,p}$ such that $gg_1=h$. We can take $g_1, h$ satisfying $|g_1|<1$ and $|h|<1$. Then we have $\log|g|/\log|h|+ \log|g_1|/\log|h|=1$ and hence the value of $\log|g|/\log|h|$ at $p$ is $\leq 1$. 

\medskip

\ref{verrat}.2. Let $X\to Y\to Z$ be morphisms in $\cB(\log)$. 
  Then, $p\in X_{[:]}$ is  vertical over $Z$ if and only if 
$p$ is vertical over $Y$ and the image of $p$ in $Y_{[:]}$ is vertical over $Z$.

\medskip

This is clear. 

\medskip

\ref{verrat}.3.  The vertical part in $X_{[:]}$ (or of $X_{[:]}^{\log}$) is open. 

\medskip

  This is seen as follows.  
  It is enough to see the case of $X_{[:]}$. 
  Let $p \in X_{[:]}$, and we prove that if $p$ is vertical over $Y$, then any point of $X_{[:]}$ near $p$ is vertical over $Y$. 
  If $M_{X,p}=\cO^{\times}_{X,p}$, the log structure is trivial around $p$ and the claim follows. 
  Otherwise, take an fs chart $\cS$ around the image of $p$ in $X$. 
  Let $\cS^{(1)}$ be the inverse image of $M^{(1)}$ associated to $p$ in \ref{m,f_jk}. 
  Then $p$ is vertical over $Y$ if and only if the image of $M_Y$ in $(M_X/\cO^{\times}_X)_p$ is not contained in the image of $\cS^{(1)}$. 
  If $p$ satisfies this condition, then any point $p'$ near $p$ also satisfies the same condition because $\cS^{(1)}$ for $p'$ is smaller than $\cS^{(1)}$ for $p$ unless the log structure is trivial at $p'$. 
\end{sbpara}

\begin{sbrem}
\label{r:rounding}
 Let $f:X\to Y$ be as in Theorem \ref{SEmore}.  We have 

(1) $f^{\log}_{[:]}$ is a topological submersion; and the proof of Theorem \ref{SEmore} below shows 

(2) If $U$ denotes the open set of $X^{\log}_{[:]}$ consisting of all points which are vertical over $Y$ (\ref{verrat}), we can take $a=0$ on $U$ so that any fiber of $U\to Y^{\log}_{[:]}$ is a topological manifold without boundaries. 

That is, this proposition is an analogue of the relative rounding result in \cite{NO} (which is about $X^{\log}\to Y^{\log}$) for $X^{\log}_{[:]}\to Y^{\log}_{[:]}$. 
  Note that the proof of Theorem \ref{SEmore} below is easier than the proof of the result in \cite{NO} 
in the sense that, though in \cite{NO}, we have to show that $\Hom(h,\R_{\ge0}^{\mult})$ is submersive for any $\bQ$-integral homomorphism $h$ of fs monoids, the proof in this time need in essence only the case where $h$ is the diagonal map $\bN \to \bN^n$ for some $n$, that is, the submersivity of the product map $\bR_{\ge0}^n \to \bR_{\ge0}; (y_1,\ldots, y_n) \mapsto y_1\cdots y_n$.

\end{sbrem}

\begin{sbpara}  
  We prove Theorem \ref{SEmore}.
  Let $n'$, $f_{j,k}$, and the local coordinate functions $(q_{\ell})_{1\le \ell \le n}$ be as in the proof of Theorem \ref{SEisS}. 
  We may assume that if $f_{j,k}$ mod $\cO_X^{\times}$ belongs to the image of $M_Y$, then $f_{j,k}$ is the pullback of an element of $M_Y$, which is denoted by the same symbol by abuse of notation. 
  Let $B$ be the set of such $f_{j,k}$s. 

  We define $s_1, \dots, s_{n'}\in A_1^{\an}$  as follows. 
  Renumbering $q_{\ell}$, we may assume that if $q_{\ell}=f_{j,k}$, $q_{\ell'}=f_{j',k'}$, and if either $k=1<k'$ or $k=k'=1$ and $j<j'$, then $\ell<\ell'$.
    Let $1\leq \ell\leq n'$.  

  If $q_{\ell}$ is $f_{j,k}$ with $1\leq j\leq m$ and $k\neq 1$, then 
$s_{\ell}=\log|f_{j,k}|/\log|f_{j,1}|$.

  If $q_{\ell}$ is $f_{j,1}$ and if there is no $j'<j$ such that $f_{j',1} \in B$, then $s_{\ell}= \log|f_{j+1, 1}|/\log |f_{j,1}|$.  
  (Here and hereafter, $f_{m+1,1}$ means $e^{-1}$.) 
  Let $a$ be the number of such an $\ell$.

  If $q_{\ell}$ is $f_{j,1}$ and if there is a $j'<j$ such that $f_{j', 1} \in B$ and if $j'$ is the largest among such $j'$, then 
$s_{\ell}= \log|f_{j+1, 1}|/\log|f_{j,1}| - \log|f_{j'+1,1}|/\log|f_{j',1}|$.
  Let $b$ be the number of such an $\ell$. 

  We prove that the above $(s_j)_j$ give the desired open immersion. 
  First, let $t_j=\log|f_{j+1, 1}|/\log |f_{j,1}|$ ($1 \le j \le m$).  
  Let $j_1<j_2<\cdots<j_{m-a-b}$ be the set of $j$ such that $f_{j,1} \in B$. 
  Then, by the map $(f^{\log}_{[:]}, (x_j)_{1\leq j\leq n'}, (t_j)_{1\leq j\leq m}, 
(s_j)_{a+b<j\leq n'}, (\Re(q_j))_{n'<j\leq n}, (\Im(q_j))_{n'<j\leq n})$, we see that $X^{\log}_{[:]}$ is locally homeomorphic to the subspace of 
$Y^{\log}_{[:]}\times (\R/\Z)^{n'} \times \R_{\geq 0}^m\times \R^{n'-a-b}_{>0}\times \R^{n-n'} \times \R^{n-n'}$ consisting of the points $(p',(x_j)_j, (t_j)_j, (s_j)_j, (r_j)_j, (i_j)_j)$ satisfying that for any $j_{\alpha}$ $(1 \le \alpha \le m-a-b)$, $t_{j_{\alpha}}\cdots t_{j_{\alpha+1}-1}$ coincides with the value of $\log|f_{j_{\alpha+1},1}|/\log|f_{j_{\alpha},1}|$ at $p'$, where $j_{m-a-b+1}$ means $m+1$. 

  If we replace $(t_j)_{j_{\alpha}\leq j <j_{\alpha+1}}$ by $(t_j-t_{j_{\alpha}})_{j_{\alpha}< j <j_{\alpha+1}}$ for each $\alpha$, we can see that this subspace is homeomorphic to an open subspace of 
$Y^{\log}_{[:]}\times (\R/\Z)^{n'} \times \R_{\geq 0}^a \times \bR^{m-a-(m-a-b)}\times \R^{n'-a-b}_{>0}\times \R^{n-n'} \times \R^{n-n'}$ by the fact that for any $c \ge 1$, the product map $\bR_{\ge0}^c \to \bR_{\ge0}; (t_1,\ldots, t_c) \mapsto t_1\cdots t_c$ is trivialized in the sense that it is homeomorphic to the projection $\bR^{c-1} \times \bR_{\ge0} \to \bR_{\ge0}$ by $(t_1,\ldots,t_c) \leftrightarrow (t_2-t_1,...,t_c-t_1, t_1\cdots t_c)$.
  Since $s_j=t_j$ if $j \le a$ and the set of the functions $t_j-t_{j_{\alpha}}$ ($1 \le \alpha \le m-a-b$, $j_{\alpha}< j <j_{\alpha+1}$) coincides with $\{s_{a+1},\ldots,s_{a+b}\}$, it gives the desired open immersion. 
  \end{sbpara}

  In the rest of this section (Section \ref{ss:infls}), we explain the necessity of the exactness condition, which is contained in the saturatedness, in the above results.

\begin{sbpara}
\label{Blbad1}

  We need the exactness in the assumption of Corollary \ref{SELF}. 
  For example, let $Y=\Delta^2$ endowed with the log structure defined by the coordinate functions $q_1, q_2$, and let $X\to Y$ be the log blowing up at the origin $(q_1,q_2)$. 
  Then $f: X\to Y$ is log \'etale (that is, log smooth \lq\lq of relative dimension $0$''). But $R^1_{X/Y}$ is a non-zero torsion module over $R$. 
 
  In fact, if $C\subset X_{[:]}$ or $X^{\log}_{[:]}$ denotes the inverse image of the point $p$ of $Y_{[:]}$ at which $q_1=q_2=0$ and  $y_2/y_1=1$, the support of $R^1_{X/Y}$ is contained in $C$ because $f_{[:]}$ is an isomorphism outside $p$. 
  On the other hand, 
consider the open set of $X_{[:]}$ or $X^{\log}_{[:]}$ on which $q_2|q_1$ and $q_1|q_2$, and let $q_3= q_1/q_2$, which is invertible there. 
  We show that $R^1_{X/Y}$ is a non-zero torsion module over $R$ there. 
  We observe that $dx_2/y_2$, $dy_2/y_2$, $dx_3$, $dy_3$  form a base of $R^1_X$ there, 
and that $dx_1/y_1$, $dy_1/y_1$, $dx_2/y_2$, $dy_2/y_2$ form a base of $R^1_Y$. 
  Since 
$dy_1/y_1 -dy_2/y_2 = (dy_2+dy_3)/(y_2+y_3) -dy_2/y_2
= y_2^{-1}(Pdy_2/y_2+Qdy_3)$ 
and $dx_1/y_1 -dx_2/y_2 = y_2^{-1}(Pdx_2/y_2+Qdx_3)$, where $P=-\frac{y_3}{1+y_2^{-1}y_3} \in R$, $Q= \frac1{1+y_2^{-1}y_3}%
\in R^{\times}$, 
$R^1_{X/Y}$ is isomorphic to $(R/y_2^{-1}R)^{\oplus2}$ there. 
  A similar calculation shows that 
$R^1_{X/Y}$ is isomorphic to $(R/(y_3/y_2)R)^{\oplus2}$ near a point of $C$ where $y_2 \gg y_3$ and that 
$R^1_{X/Y}$ is isomorphic to $(R/(-y_3/y_1)R)^{\oplus2}$ near a point of $C$ where $y_1 \gg -y_3$.
\end{sbpara}

\begin{sbpara}
\label{Blbad2}

  In \ref{Blbad1}, we do not have $f_{[:]*}(A_{2,X})=A_{2,Y}$.  
  In fact, take a non-zero $C^{\infty}$ function on $X$ whose support does not intersect with the strict transformations of $\{q_1=0\}$ and $\{q_2=0\}$ but intersect with the exceptional fiber.
  We regard this function as a global section of $A_{2,X}$. 
  If $p$ denotes the point of $Y_{[:]}$ at which $q_1=q_2=0$ and  $y_2/y_1=1$, this function in $A_{2,X}(X_{[:]})$ is not constant in the fiber of $p$ in $X_{[:]}$ and hence this does not belong to $A_{2,Y}(Y_{[:]})$.
\end{sbpara}

\begin{sbpara}
  In \ref{Blbad1}, we also do not have $f_{[:]*}(A_{3,X})=A_{3,Y}$.
  We prove it. 
  Let $V$ be the open set of $X_{[:]}$ where $y_3>1$. 
  Let $q$ be the point of $X$ where $q_2=q_3=0$. 
  Let $W \subset V$ be an open neighborhood of the fiber of $q$ such that the closure of $W$ is contained in $V$. 
  Let $g$ be a usual $C^{\infty}$ function on the real line such that $g^{(n)}$ is bounded near $+\infty$ for any $n \ge 0$. 
  Then, by Proposition \ref{loggrow}, the function $g(\log(y_3))$ belongs to $A_3(V)$.
  By the softness of $A_3$, we have an element $h \in A_3(X_{[:]})$ which is $g(\log(y_3))$ on $W$ and zero outside $V$. 
  We prove that, for some choice of $g$, the section $h$ is not necessarily come from $A_{3,Y}$. 
  (We guess that $g=\sin(x)$ will do, though we could not prove it.)
  To see it, we work around the point $p$ (\ref{Blbad2}).
  Let $V'$ be the image of $V$ in $Y_{[:]}$. 
  Then the closure of $V'$ contains $p$. 
  If $h$ belongs to $A_{3,Y}$, by Proposition \ref{loggrow}, there is an $n$ such that all iterated derivations $l$ of $h$ satisfy $|l|<C(y_1+y_2)^n$ near $p$. 
  Consider $(y_1\frac \partial{\partial y_1})^{n+1}g(\log(y_1-y_2))$. 
  As is easily seen, it is calculated as 
$$(y_1-y_2)^{-n-1}\bigl(y_1^{n+1}g^{(n+1)}(\log(y_1-y_2))+{\textstyle\sum}_{j=1}^n p_{j,n+1}(y_1,y_2)g^{(j)}(\log(y_1-y_2))\bigr),$$ where $p_{j,n+1}$ is a homogeneous polynomial of $y_1$ and $y_2$ of degree $n+1$. 
  Let $a_{j,n+1}$ be the sum of the absolute values of the coefficients of $p_{j,n+1}$. 
  Then, under the condition $y_1>y_2>0$, 
$|\sum_{j=1}^n p_{j,n+1}(y_1,y_2)g^{(j)}(\log(y_1-y_2))|\le
y_1^{n+1}\sum_{j=1}^n a_{j,n+1}|g^{(j)}(\log(y_1-y_2))|$.
  Let $b_k=\overline \lim_{x \to \infty} |g^{(k)}(x)|$ for each $k \ge0$. 
  (Recall that we assumed that $g^{(k)}$ is bounded near the infinity.) 

\smallskip

{\bf Claim 1.} If $b_{n+1} > \sum_{j=1}^n a_{j,n+1}b_j$, then it yields a contradiction. 

\smallskip

  We prove Claim 1. 
  Let $\varepsilon=(b_{n+1}-\sum_{j=1}^n a_{j,n+1}b_j)/3$.
  Near $p$ on $V'$, we have 
$$(\tfrac{y_1}{y_1-y_2})^{n+1}\bigl(|g^{(n+1)}(\log(y_1-y_2))|-\varepsilon-{\textstyle\sum}_{j=1}^n a_{j,n+1}b_j\bigr)\le C(y_1+y_2)^n.$$
  Considering the locus where $y_2=y_1-\log{y_1}$, it implies that 
$$(\tfrac{y_1}{\log{y_1}})^{n+1}\bigl(|g^{(n+1)}(\log\log{y_1})|-\varepsilon-{\textstyle\sum}_{j=1}^n a_{j,n+1}b_j\bigr)\le 2^nCy_1^n.$$ 
  Then we have 
$(\frac{y_1}{\log y_1})^{n+1}\varepsilon\le 2^nCy_1^n$ for infinitely many $y_1$ such that $y_1 \to \infty$, which is a contradiction. 

  By Claim 1, it is enough to show the existence of $g$ which satisfies 
$b_{n+1} > \sum_{j=1}^n a_{j,n+1}b_j$ for any $n \ge0$. 
  Let $a_{n+1}=1+n\underset {1 \le j \le n} \max a_{j,n+1}$ for any $n\ge0$. 

\smallskip

{\bf Claim 2.} It is enough to show the existence of $g$ which satisfies 
$b_{n+1}>a_{n+1}b_n$ for any $n \ge0$. 

\smallskip

  In fact, for such a $g$, $b_j$ is strictly increasing so that 
$\sum_{j=1}^n a_{j,n+1}b_j\le n(\max a_{j,n+1})b_n \le a_{n+1}b_n$.

  Thus the rest is to see the following. 

\smallskip

{\bf Claim 3.} Let $a_n$ ($n \ge 1$) be a real number which is greater than 1.
  Then there is a smooth function $g$ on $\bR$ such that $b_n=
\overline \lim_{x \to \infty} |g^{(n)}(x)|$ ($n \ge0$) are finite and satisfy $b_{n+1}>a_{n+1}b_n$ for any $n \ge0$. 

\smallskip

  Consider the Fourier series $$g(x)={\textstyle\sum}_{k \ge0} \sin(c_k x)/(3^k c_k^{k-1}),$$ where $c_k$ is an increasing sequence with $c_0=1$. 
  Then, it uniformly converges and we have 
$$g^{(n)}(x)={\textstyle\sum}_{k \ge0} c_k^{n-k+1}\sin(c_k x+n\pi/2)/3^k.$$
From this, we have 
$$|g^{(n)}(x)| \le {\textstyle\sum}_{k \le n} {c_k^{n-k+1}}/{3^k}+(2\cdot 3^n)^{-1}$$
and 
$$\overline \lim |g^{(n+1)}(x)| \ge c_{n+1}/3^{n+1}-{\textstyle\sum}_{k \le n} {c_k^{n-k+2}}/{3^k}-(2\cdot 3^{n+1})^{-1}.$$

  Now we can choose $c_k$ iteratively satisfying 
$$\bigl({\textstyle\sum}_{k \le n} {c_k^{n-k+1}}/{3^k}+(2\cdot 3^n)^{-1}\bigr)a_{n+1} < 
c_{n+1}/3^{n+1}-{\textstyle\sum}_{k \le n} {c_k^{n-k+2}}/{3^k}-(2\cdot 3^{n+1})^{-1}$$
for any $n$.
  Then, we have $b_na_{n+1} < b_{n+1}$, as desired. 
\end{sbpara}

\section{Log integration}
\label{s:integration}

In this Section \ref{s:integration}, we develop a theory of integration in log geometry. We also prove a  Poincar\'e lemma (Theorem \ref{LPL}) for log $C^{\infty}$ differential  forms.

\subsection{Our integration}\label{ss:ourint}

  Here is an overview of Section \ref{s:integration}. 

\begin{sbpara} In this Section \ref{s:integration}, we develop an  integration theory for a log smooth and saturated morphism $f: X\to Y$ of fs log analytic spaces, by using log $C^{\infty}$ differential forms. 
We also prove a Poincare lemma for log $C^{\infty}$ differential forms for such an $f:X\to Y$.
\end{sbpara}

\begin{sbpara}\label{2n}  We will discuss  the following theory of  integration of $2n$-forms, where $n$ is the relative dimension of $X$ over $Y$. See Section \ref{ss:2nint}. 
Let $U$ be an open set of $X_{[:]}$ (resp.\ $X_{[:]}^{\log}$) and let $V$ be an open set of $Y_{[:]}$ (resp.\ $Y^{\log}_{[:]}$). Assume that $f_{[:]}(U)\subset V$ (resp.\ $f^{\log}_{[:]}(U)\subset V$) and that every point of $U$ is vertical over $Y$ in the sense of \ref{verrat}. Let $k=2*$ or $3$,  let $w$ be an element of $A^{2n}_{k,X/Y}(U)$ (resp.\ $A^{\log,2n}_{k,X/Y}(U)$),  and assume that the support of $w$ is proper over $V$.  In Section \ref{ss:2nint}, we define  $$\int_{U/V} w \in A_{k,Y}(V).$$ The following  (i)--(iii) are the main properties of this integration. 

(i) This integration is $A_{k,Y}(V)$ (resp.\ $A^{\log}_{k,Y}(V)$)-linear in $w$.

(ii) In the case where the log structure of $Y$ is trivial and $Y$ is smooth (then $U$ and $V$ are smooth complex analytic spaces and the morphism $U\to V$ is smooth), 
this is the classical integration along the fibers. That is, 
$\int_{U/V} w$ is  the function $v\mapsto \int_{U(v)} w(v)$ on $V$. 
  Here, $U(v)$ is the fiber of $v$ in $U$ and $w(v)$ is the pullback of $w$ to $U(v)$ which is a differential form of degree $2n$ with compact support on the smooth complex analytic space $U(v)$.

\medskip

 (iii) Assume that we have a cartesian diagram of fs log analytic spaces
$$\begin{CD}
X' @>>> X \\
@Vf'VV @VVfV \\
Y' @>>> Y.
\end{CD}
$$
Let $V'$ be  an open set of  $Y'_{[:]}$ (resp.\ $(Y')^{\log}_{[:]}$) whose image in  $Y_{[:]}$ (resp.\ $Y^{\log}_{[:]}$) 
is contained in $V$,  let $U'$ be the intersection of the inverse image of $U$ and the inverse image of $V'$ in   $X'_{[:]}$ (resp.\ $(X')^{\log}_{[:]}$), and let 
$w'$ be the pullback of $w$ in $A^{2n}_{k,X'/Y'}(U')$ (resp.\ $A^{\log,2n}_{k, X'/Y'}(U')$).  Then  $\int_{U'/V'} w'$ coincides with the
  pullback of $\int_{U/V} w$.

\medskip

In particular, if we apply (iii) to the case where $Y$ is log smooth and $Y'=Y_{\triv}$, we see that $\int_{U/V} w$ is 
 the unique extension  to $V$ of the map given on $V \cap Y_{\triv}$ by (ii). 
 
In general, locally on $X$ and $Y$, we have a strict immersion $Y\to Y'$ into a log smooth $Y'$ and $X\to Y$ comes from $X' \to Y'$ which is log smooth and saturated as $X=X'\times_{Y'} Y$. 
  Then the integration for $X/Y$ is compatible with the integration for $X'/Y'$ by (iii) (applied  exchanging the notation $X$ and $X'$,  and $Y$ and  $Y'$) and hence is related to the classical integration in (ii).

\end{sbpara}

\begin{sbpara}  In this Section \ref{s:integration}, we also treat integrals of $1$-forms. See  Sections \ref{ss:inty}, \ref{ss:intx}, and \ref{ss:intcoh2}.

\end{sbpara}

\begin{sbpara}

We plan that in a forthcoming paper \lq\lq 
Integrations in log geometry'' (\cite{KNUi}), we discuss the following things.

(i) We discuss the integration of $m$ forms for arbitrary $m$.

(ii) We will unify the theories of integrations of $1$-forms in the above three sections  \ref{ss:inty}, \ref{ss:intx}, and \ref{ss:intcoh2}
which are not well-related to each other in the present paper. 

(iii) In the present paper, we assume the verticality of $U/V$  in the integration of $2n$ forms. In the paper \cite{KNUi},  we will discuss it without this assumption (cf.\ \ref{nonvert}). 

\end{sbpara}
\begin{sbpara}

In this Section \ref{s:integration}, we prove a log Poincar\'e lemma (Theorem \ref{LPL}) in Section \ref{ss:LPL} based on our path integrals in Sections \ref{ss:inty}, \ref{ss:intx}. The integrations of $2n$ forms introduced in \ref{2n} is constructed by using this log Poincar\'e lemma.

\end{sbpara}

\subsection{The logarithmic  path $|\Delta|$}
\label{ss:abDel}

\begin{sbpara}
\label{path}
  In the following, we need the log $C^{\infty}$ functions on the interval $|\Delta|:=[0,1)\subset \Delta$, which is like the path in log geometry. 
  For example, the log $C^{\infty}$ functions on $|\Delta|$ in the sense of $A_2$ are $C^{\infty}$ functions on $(0, 1)$ whose iterated derivatives by $y\partial/\partial y$ extend to the boundary (i.e., to the origin) as continuous functions. 
  Here $y=-(2 \pi)^{-1}\log q$, where $q$ is the coordinate function of $|\Delta|$. 
  More generally, we consider $X\times |\Delta|$ for any fs log analytic space $X$, introduce the spaces  $(X\times |\Delta|)_{[:]}$ and $(X\times|\Delta|)_{[:]}^{\log}$, and 
   introduce log $C^{\infty}$ functions on these two spaces.  These will be used in Section \ref{ss:inty} and in the proof of log Poincar\'e lemma in Section \ref{ss:LPL}. 
\end{sbpara}

\begin{sbpara}
\label{Xxpath}
  Let $X$ be an fs log analytic space. Let the topological space  $(X \times |\Delta|)_{[:]}$ be inverse image of $|\Delta|$ under $(X \times \Delta)_{[:]}\to \Delta$.  
  
Let $(X \times |\Delta|)_{[:]}^{\log}$ be the topological space $(X \times |\Delta|)_{[:]}\times_{X}X^{\log}$.

We have a canonical homeomorphism 
$$(X \times \Delta)^{\log}_{[:]}\cong (X \times |\Delta|)^{\log}_{[:]}\times \R/\Z$$
induced by $(2\pi)^{-1}\arg(q): (X \times \Delta)^{\log}_{[:]}\to \R/\Z$  with $q$ being the coordinate function of $\Delta$. 
\end{sbpara}

\begin{sbpara}\label{Xxpath2} Let $X$ be as above, let $p\in X_{[:]}$, and let $f_{j,1}$ ($1\leq j\leq m$) be as in \ref{m,f_jk}. Then replacing $X_{[:]}$ by an open neighborhood $U$ of $p$, we have a  commutative cartesian diagram
$$\begin{matrix}  (X \times |\Delta|)_{[:]} & \to & X_{[:]} \times |\Delta| \times [0, \infty]^m \\
\cap && \cap\\
(X \times \Delta)_{[:]} & \to & X_{[:]} \times \Delta \times [0, \infty]^m
\end{matrix}$$
in which the maps to $[0, \infty]^m$ are given by $(\log|f_{j,1}|/\log|q|)_{1\leq j\leq m}$, the horizontal arrows are injective, and  the topology of the left-hand-side of each horizontal arrow is that as a subspace of the right-hand-side. 

Each horizontal arrow induces a homeomorphism from the fiber of $(p,0)\in X_{[:]}$ in the space of the left-hand-side to the subspace of $[0, \infty]^m$ consisting of all elements $(a_j)_{1\leq j\leq m}$ satisfying the following conditions (i) and (ii). (i) $a_1 \geq a_2 \geq \cdots \geq a_m$.
(ii) There is at most one $j$ such that $a_j\neq 0, \infty$. 

\end{sbpara}

\begin{sbpara}\label{Xxpath3} We have a canonical continuous map $(X \times \Delta)_{[:]}\to (X \times |\Delta|)_{[:]}$ such that  we have the commutative cartesian diagram
$$\begin{matrix}  (X \times \Delta)_{[:]} & \to & X_{[:]} \times \Delta \times [0, \infty]^m \\
\downarrow  && \downarrow \\
(X \times |\Delta|)_{[:]} & \to & X_{[:]} \times |\Delta| \times [0, \infty]^m
\end{matrix}$$
with the  horizontal arrows  as in \ref{Xxpath2} in which the right vertical arrow is given by the identity maps on $X_{[:]}$ and on $[0, \infty]^m$ and the map $\Delta\to |\Delta|\;;\;q\mapsto |q|$. 

The vertical arrows in this diagram are proper, and hence the topology of $(X \times |\Delta|)_{[:]}$ is the image of the topology of $(X \times \Delta)_{[:]}$.

\end{sbpara}

\begin{sbpara}\label{thetac1} Let $c\in \R_{>0}$. Then we have a canonical homeomorphism $\nu_c: (X \times |\Delta|)_{[:]}\to (X\times |\Delta|)_{[:]}$ which is compatible with $X\times |\Delta|\to X \times |\Delta|; 
(u, q)\mapsto (u, q^c)$. It is the map which is,  via the upper horizontal arrow in \ref{Xxpath2},  compatible with
the map $X_{[:]}\times |\Delta|\times [0, \infty]^m\to X_{[:]}\times |\Delta|\times [0, \infty]^m$ induced by the identity map of $X_{[:]}$, the map $|\Delta|\to |\Delta|\;;\; q \mapsto q^c$, and the map $[0, \infty]^m \to [0, \infty]^m\;;\; (a_j)_j\mapsto (c^{-1}a_j)_j$.

\end{sbpara}

\begin{sbprop}\label{3.1fiber} 
  Let $X$ be an fs log analytic space. 
  Let the notation be as in $\ref{Xxpath}$.
  If $p\in X^{\log}_{[:]}$ (resp.\ $X_{[:]}^{\log}$, resp.\ $X_{[:]}$), there is an open neighborhood $V$  of $p$ such that if we denote by $U$ the inverse image of $V$ in $(X\times \Delta)^{\log}_{[:]}$ (resp.\ $(X\times |\Delta|)^{\log}_{[:]}$, resp.\ $(X \times|\Delta|)_{[:]}$), we have a homeomorphism $$U \cong V\times \R/\Z\times [0,1)\quad (\text{resp.\ } U\cong V \times [0,1), \quad \text{resp.\ } U \cong V \times [0,1))$$ over $V$ given by $x:U \to \Delta^{\log} \to \bR/\bZ$ and (resp.\ by, resp.\ by) 
$$\ell:=
\tfrac1{m+1}{\textstyle\sum}_{j=1}^{m+1} \log|f_{j,1}|/\log|qf_{j,1}|
= \tfrac1{m+1}{\textstyle\sum}_{j=1}^{m+1} (y/y_{j,1}+1)^{-1}; 
U\to [0,1).$$ 
  Here $m$ and $f_{j,1}$ ($1\leq j\leq m$) are for $M_{X,p}$ as in $\ref{m,f_jk}$, $f_{m+1,1}$ denotes $e^{-1}$, 
$q$ is the coordinate function of $|\Delta|$, $y=-(2\pi)^{-1}\log q$, and 
$y_{j,1}=-(2\pi)^{-1}\log |f_{j,1}|$. 

In the case of $(X\times \Delta)^{\log}_{[:]}$, a point $p'$ of $U$ is vertical over $X$ ($\ref{verrat}$) if and only if $\ell(p')\neq 0$.
\end{sbprop}

\begin{sbrem} 
 The map $\ell: U \to [0, 1)$ in Proposition \ref{3.1fiber} has the following nice properties. Let $p\in X^{\log}_{[:]}$ (resp.\ 
$p\in X_{[:]}^{\log}$, resp.\ $p\in X_{[:]}$),  let $\mu\in (X\times \Delta)^{\log}_{[:]}$ (resp.\ $(X \times |\Delta|)_{[:]}^{\log}$, resp.\ $(X \times |\Delta|)_{[:]}$) be a point in the fiber of $p$ and let $j$ be an integer such that
  $1\leq j\leq m+1$. Then:
  
  (i)  $\ell(\mu)\geq j/(m+1)$ if and only if $y/y_{j, 1}=0$ at $\mu$.
  
  (ii)  $\ell(\mu) \leq (j-1)/(m+1)$ if and only if $y/y_{j,1}=\infty$ at $\mu$.
  
  (iii) $(j-1)/(m+1) < \ell(\mu) <j/(m+1)$ if and only if $0<y/y_{j,1}<\infty$ at $\mu$.
  
  (iv) Assume $1\leq j\leq m$. Then $\ell(\mu)=j/(m+1)$ if and only if $y/y_{j,1}=0$ and $y/y_{j+1, 1}=\infty$ at $\mu$.

\end{sbrem}

\begin{pf} 
  The proof for $(X \times \Delta)_{[:]}^{\log}$ is 
 reduced to that for   $(X \times |\Delta|)_{[:]}^{\log}$ by taking care of $x$-component, and the proof for $(X \times |\Delta|)^{\log}_{[:]}$ is reduced to that for $(X \times |\Delta|)_{[:]}$ taking $X^{\log}_{[:]} \times_{X_{[:]}}-$. Consider  $(X \times |\Delta|)_{[:]}$.

  To see that the map in the statement is bijective, it is enough to show it in any fiber of $U \to V$. 
  First we see from the proof of Theorem \ref{SEmore} that each fiber of $U \to V$ is homeomorphic to $[0,1)$.
  Second we observe that if $\mu\in U$ moves along a fiber from the point where $y/y_{1,1}=\infty$ to 
the other end where $y \to 0$, 
the function $\mu \mapsto (y/y_{j,1}+1)^{-1}$  is increasing for any $1 \le j \le m+1$ whose set of values is $[0,1]$ for $j \le m$ and is $[0,1)$ for $j=m+1$. 
  Further, since each of these functions $(y/y_{j,1}+1)^{-1}$ is strictly increasing whenever $y/y_{j,1}$ is non-zero finite, 
one of $(y/y_{j,1}+1)^{-1}$ is strictly increasing in any part of the fiber. 
  Since $\ell$ is the average of the above functions, $\ell$ restricted to each fiber gives a bijection from each fiber to $[0,1)$. 
  Hence the map in the statement is bijective. 

  We claim that it is also an open map. 
  In fact, since the proof of 
Theorem \ref{SEmore} implies that $U \to V$ is locally a product map (cf.\ also Remark \ref{r:rounding}), the claim follows from 
the fact that for any open $V'$ of $V$, and for any pair of continuous functions $a$ and $b$ from $V'$ to $[0,1)$ such that $a(v')<b(v')$ ($v' \in V'$), the subset $\{(v',t)\,|\,v' \in V', t$ is between $a(v')$ and $b(v')\}$ is open in $V \times [0,1)$. 
  Thus the map in the statement is open and hence a homeomorphism. 
\end{pf}

\begin{sbpara}\label{A|Del|}

  Let the notation be as in \ref{Xxpath}.
  For $k=2,2*,$ or $3$,
  we define a sheaf of rings  $A_{k,X\times |\Delta|}$ on $(X \times |\Delta|)_{[:]}$ (resp.
$A_{k, X\times |\Delta|}^{\log}$ on $(X \times |\Delta|)_{[:]}^{\log}$).

  First, we assume that $X$ is log smooth. Let $k=2,2*,3$, let $U$ be an open set of $(X \times |\Delta|)_{[:]}$ (resp.\ $(X \times |\Delta|)_{[:]}^{\log}$), and let $U'$ be the inverse image of $U$ in $(X\times \Delta)_{[:]}$ (resp.\ $(X\times \Delta)_{[:]}^{\log}$) under the canonical map $(X\times \Delta)_{[:]}\to (X \times |\Delta|)_{[:]}$ (\ref{Xxpath3}) (resp. \; $(X\times \Delta)^{\log}_{[:]}\to (X \times |\Delta|)^{\log}_{[:]}$ (\ref{Xxpath})). 
Define $A_{k,X\times |\Delta|}(U)$ (resp.\ $A_{k, X\times |\Delta|}^{\log}(U)$) to be the set of all $\R$-valued functions on $U\cap (X_{\triv} \times (0,1))$ 
whose pullback to $U' \cap (X_{\triv}\times (\Delta\smallsetminus \{0\}))$ 
 extends to an element of $A_k(U')$ (resp.\ $A_k^{\log}(U')$).

  For a general $X$, we define them as $R_{X\times|\Delta|}=R_X\otimes_{R_Z} R_{Z\times|\Delta|}$ ($R=A_k, A_k^{\log}$) by using a strict immersion  $X\to Z$ into log smooth (\ref{strim}). 
  
  In the case where $X$ is log smooth, for an open set $U$ of $(X\times |\Delta|)_{[:]}$, $A_2(U)$ is equal to the set of all $C^{\infty}$ functions on $U\cap (X_{\triv}\times (0,1))$ whose all iterated derivatives with respect to the local base of $A^{\an,1}_{1, X}$  and $dy/y$ ($y=-(2\pi)^{-1}\log(q)$ with $q$ being the coordinate function of $|\Delta|$) extend to continuous functions on $U$. For general $X$, $A_{2*, X\times |\Delta|}$ (resp.\ $A_{3, X \times |\Delta|}$) of $(X\times |\Delta|)_{[:]}$ is generated over  $A_{2*, X}\otimes_{A_{2,X}} A_{2, X\times|\Delta|} $ algebraically by $\log(y)$ (resp.\ $y$), and $A_{k, X}^{\log}\otimes_{A_{k,X}} A_{k, X\times|\Delta|} \overset{\cong}\to A^{\log}_{k,X\times |\Delta|}$ for $k=2,2*, 3$. 
  
  For $R=A_k, A_k^{\log}$, the sheaf of  differential forms (1) $R^1_{X \times |\Delta|}$,  
  (2) $R^1_{X \times |\Delta|/X}$ and 
(3) $R^1_{Y/ X \times |\Delta|}$  for an fs log analytic space $Y$ over $X \times \Delta$, are defined as follows.

(1) is the direct sum of $R_{X \times |\Delta|} \otimes_{R_X} R^1_X$ and the free $R_{X \times |\Delta|}$-module with base $dy/y$. 

(2) is the cokernel of $R_{X \times  |\Delta|}\otimes_{R_X}  R^1_X\to R_{X \times |\Delta|}^1$.  
It is a free $R_{X\times |\Delta|}$-module of rank $1$ with base $dy/y$. 

(3) is the cokernel of $R_Y\otimes_{R_{X \times |\Delta|}}  R^1_{X \times |\Delta|}\to R_Y^1$.  In the case $Y=X\times \Delta$, it is a free $R_Y$-module of rank $1$ with base $dx/y$.

\end{sbpara}

\begin{sbpara}\label{thetac2} Let $c\in \R_{>0}$. Then the map $\nu_c: (X\times |\Delta|)_{[:]}\to (X \times |\Delta|)_{[:]}$ (\ref{thetac1}) and the induced map $(X\times |\Delta|)_{[:]}^{\log}\to (X \times |\Delta|)_{[:]}^{\log}$ underlie isomorphisms $((X \times |\Delta|)_{[:]}, A_k)\overset{\cong}\to ((X \times |\Delta|)_{[:]}, A_k)$ and $((X \times |\Delta|)_{[:]}^{\log}, A_k^{\log})\overset{\cong}\to ((X \times |\Delta|)_{[:]}^{\log}, A_k^{\log})$
($k=2,2*,3$) defined as follows, respectively. We will denote these morphisms also as $\nu_c$.

In the case where $X$ is log smooth, we have these morphisms because the pullback of $dy/y$ by $\nu_c$ is $cdy/y$. In the general case, these morphisms are obtained from the log smooth case by local strict immersions into log smooth.

\end{sbpara}

\subsection{Path integrals in the direction of $y$}\label{ss:inty}
  We discuss path integrals in the direction of $y$ (Proposition \ref{inty0}).

  The following Theorem \ref{intyp} and Theorem \ref{intxp} in the next section are used for the proof of log Poincar\'e lemma (Theorem \ref{LPL}). 
 
\begin{sbthm}\label{intyp}
  Let $X$ be an fs log analytic space.
  For $k=2*,3$, we have exact sequences 
\begin{align*}
0&\to A_{k,X}\to A_{k, X\times |\Delta|} \overset{d}\to A^1_{k, (X \times |\Delta|)/X}\to 0\quad \text{on}\;\;  X_{[:]}, {\text and} \\
0&\to A^{\log}_{k,X}\to A^{\log}_{k, X\times |\Delta|} \overset{d}\to A^{\log,1}_{k, (X \times |\Delta|)/X}\to 0 \quad \text{on}\;\;  X_{[:]}^{\log}.
\end{align*}
\end{sbthm}

\begin{sbpara}\label{intyrem} The second sequence is obtained as 
$A^{\log}_{k,X}\otimes_{A_{k,X}}-$ of the first, and the exactness of the second sequence  follows from that of the first by the flatness of $A_{k,X}\to A^{\log}_{k,X}$ (Proposition \ref{polstalk}) (or by the locally splitting property of the first (the case (1) of Proposition \ref{inty0} given below)).

\end{sbpara}

\begin{sbpara}\label{sec1} 
Let  $V$ be an open set of $X_{[:]}$ (resp.\ $X^{\log}_{[:]}$). 

By a {\it $V$-point} of $|\Delta|$, we mean a section of the inverse image of the sheaf $\Q\otimes M^{\gp}_X/\{z\in \C^\times\;|\; |z|=1\}$ on $V$ which is locally of the form $c\otimes f$  
 with $c\in \Q_{>0}$ and $f\in M_X$ such that $|f|<1$.  A simple example of a $V$-point (called a {\it constant $V$-point}) of $|\Delta|$ is the class of $1\otimes z$ for $z\in \C^\times$ such that $|z|<1$.

Let $\tilde V$ be the inverse image of $V$ in $(X\times |\Delta|)_{[:]}$ (resp.\ $(X \times |\Delta|)^{\log}_{[:]}$). A $V$-point $b$ of $|\Delta|$ gives a morphism of ringed spaces
$$(V, A_{k, X}|_V) \to (\tilde V, A_{k, X \times |\Delta|})$$
$$(\text{resp}.\;  (V, A_{k, X}^{\log}|_V) \to (\tilde V, A^{\log}_{k, X \times |\Delta|}))$$
($k=2, 2*, 3$)
defined as follows. Locally on $V$, $b$ is the class of $c \otimes f$, where $c\in \Q_{>0}$, $f\in M_X(V')$ for some open set $V'$ of $X$ which contains the image of $V$ such that $|f|<1$. Then the morphism in problem is given by the composition
$$(V'_{[:]}, A_k) \to ((X\times \Delta)_{[:]}, A_k) \to ((X\times |\Delta|)_{[:]}, A_k) \to ((X\times |\Delta|)_{[:]}, A_k)$$ 
$$(\text{resp}.\; ((V')^{\log}_{[:]},  A^{\log}_k) \to ((X\times \Delta)_{[:]}^{\log}, A_k^{\log}) \to ((X\times |\Delta|)_{[:]}^{\log}, A_k^{\log}) \to ((X\times |\Delta|)_{[:]}^{\log}, A_k^{\log})),$$ where the first morphism is induced by the morphism $V'\to X \times \Delta$ defined by $f$, the second morphism is  the canonical projection, and the third morphism is  $\nu_c$ (\ref{thetac2}). These local constructions are canonical and give a global morphism whose image is in $\tilde V$.

We often identify a $V$-point of $|\Delta|$ with these morphisms of ringed spaces. 

We will use a $V$-point of $|\Delta|$ as the base point of our indefinite integral in this Section \ref{ss:inty}. 
\end{sbpara}

\begin{sbpara}\label{Ub}  We consider a triple $(V, U, b)$ in the following two situations (1) and (2).  In the situation (1) (resp.\  (2)), $V$ is an open set of $X_{[:]}$ (resp.\ $X^{\log}_{[:]}$), $U$ is an open set of the inverse image $\tilde V$ of $V$ in  $(X\times |\Delta|)_{[:]}$ (resp.\  $(X\times |\Delta|)^{\log}_{[:]}$) such that all the fibers of the projection  $p: U\to V$ are connected, and $b$ is a $V$-point of $|\Delta|$ such that the image of the morphism $V\to \tilde V$ 
in \ref{sec1} is contained in $U$. 

Then by Proposition \ref{3.1fiber} and the existence of the section $b:V\to U$, the map $p$ 
is a homotopical isomorphism over $X^{\log}_{[:]}$. (This is a generalization of the fact that every non-empty connected open set of the interval $[0, 1)$ is contractible.) 
  Hence, by \cite{KS} Corollary 2.7.7 (i), we have in the situation (1) (resp.\  (2))

\ref{Ub}.1.\;   $R^1p_*p^{-1}(A_{k,X}|_V)=0$ and $p_*p^{-1}(A_{k,X}|_V)=A_{k,X}|_V$. (resp.\ 
$R^1p_*p^{-1}(A^{\log}_{k,X}|_V)=0$ and $p_*p^{-1}(A^{\log}_{k,X}|_V)=A^{\log}_{k,X}|_V$.) 
\end{sbpara}

\begin{sbpara}\label{interval}
We have the following basic examples (i)--(iii) of $(V, U, b)$. 

Let $V$ be an  open set of $X_{[:]}$ (resp.\  $X^{\log}_{[:]}$).

For  $V$-points $b_1, b_2$ of $|\Delta|$,  by $b_1<b_2$, we  mean that  locally on $V$, $b_j=c(j) \otimes f_j$ $(j=1,2)$ and the value of  $\log|f_1|/\log|f_2|$ (\ref{fgan}) at every point  is  $> c(2)/c(1)$. For example, if $b_j$ are constant $V$-points given by $z_j\in \C^\times$ such that $|z_j|<1$ ($j=1,2$), $b_1<b_2$ is equivalent to $|z_1|<|z_2|$.

For a $V$-point $b$ of $|\Delta|$ and for $p\in (X\times |\Delta|)_{[:]}$ (resp.\ $(X\times |\Delta|)^{\log}_{[:]}$), by $p<b$, we mean that locally on $V$, $b$ is the class of $c\otimes f$ ($c \in \bQ_{>0}$, $|f|<1$) and every value of $\log |q|/\log |f|$ is $>c$. Here $q$ is the coordinate function of $\Delta$. By $b<p$, we mean that locally on $V$, $b$ is the class of $c\otimes f$ ($c\in \bQ_{>0}$, $|f|<1$) and every value of $\log |q|/\log |f|$ is $<c$.

 (i) $U$ is $\tilde V$ 
and $b$ is a $V$-point of $|\Delta|$.  

 (ii) We have $V$-points $b, b_2$ of $|\Delta|$ such that  $b<b_2$  and $U$ is the set of all points $p$ of $\tilde V$ 
such that $p<b_2$.

(iii) We have $V$-points $b_1, b, b_2$ of $|\Delta|$ such that $b_1<b<b_2$ and $U$ is the set of all points $p$ of $\tilde V$ 
such that $b_1<p<b_2$.

In (ii), if $b_2$ is given first,  locally on $V$, we have a $V$-point
$b<b_2$ of $|\Delta|$ defined as $b_2= c(2)\otimes f_2$, $b= a \otimes f_2$ with $a>c(2)$. In (iii), if $b_1$ and $b_2$ are given first, locally on $V$, we have $b$ such that $b_1<b<b_2$ defined as 
$b_j=c(j) \otimes f_j$,  $b= \frac1{a(1)+a(2)}\otimes f_1^{a(1)c(1)}f_2^{a(2)c(2)}$ with $a(1), a(2) >0$ such that $a(j)c(j)\in \Z$.

\end{sbpara}

\begin{sbpara}\label{sec2}

Let   $p\in (X \times |\Delta|)_{[:]}$ (resp.\ $(X \times |\Delta|)^{\log}_{[:]}$) and let $\bar p$ be the image of $p$ in  $X_{[:]}$ (resp.\ $X^{\log}_{[:]}$). Assume that the image of $p$ in  $(X \times \Delta)_{[:]}$  (resp.\ $(X \times \Delta)^{\log}_{[:]}$) 
 is vertical over $X$ in the sense of \ref{verrat}.

Then if  $E$ is an open neighborhood of $p$, there are an open neighborhood $V$ of $\bar p$ and $U$ as in (ii) in \ref{interval}   such that $U\subset E$.  This follows from Proposition \ref{3.1fiber}. 
\end{sbpara}

\begin{sbprop}\label{inty0}
  Let $k$ be either $2*$ or $3$. 
Let $(V,U, b)$ be as in $\ref{Ub}$ $(1)$ (resp.\  $(2)$) and let $w$ be an element  of $A^1_{k,(X\times |\Delta|)/X}(U)$ (resp.\  $ A^{\log, 1}_{k, (X\times |\Delta|)/X}(U)$). Then there is a unique element $F$ of $A_{k,X\times |\Delta|}(U)$ (resp.\   $ A^{\log}_{k, X\times |\Delta|}(U)$) such that $dF=w$ and $F(b)=0$. Here $F(b)$ denotes the pullback of $F$ in $A_{k, X}(V)$ by $b$. 

Denote this $F$ by $\int_b w$. Then for every element $f$ of $A_{k, X\times |\Delta|}(U)$ (resp.\  $A^{\log}_{k, X\times |\Delta|}(U)$), we have 
$$f= \int_b df+f(b).$$
\end{sbprop}

This proposition follows from Theorem \ref{intyp} by \ref{Ub}.1.
 In actual, we prove this proposition together with Theorem \ref{intyp}.

\begin{sbpara}\label{inty2}  

Assuming that $X$ is log smooth, we prove Proposition \ref{inty0} for $(X\times |\Delta|)_{[:]}$ (Case (1)) and  for $k=2*$.

Let $q$ be the coordinate function of $\Delta$ and let $y= -(2\pi )^{-1}\log|q|$. 

   Let $f\in A_{2,X \times |\Delta|}(U)$ and let $n\geq 0$ be an integer. It is enough to prove $\int_b f (\log y)^n dy/y \in A_{2*, X \times |\Delta|}$. 
   
   Let $y_0\in A_{3,X}(V)$ be the pullback of $y$ under $b$.  
Twisting by the isomorphism $\nu_c$ (\ref{thetac1}, \ref{thetac2}) if necessary, we may assume $y>1$ on $U$ and $y_0>1$ on $V$. We prove 
$\int_b f (\log y)^n dy/y \in ((\log y)^{n+1}+(\log y_0)^{n+1})A_{2, X \times |\Delta|}$.

 Below we use $y$ as a new coordinate function of $|\Delta|$, so $|\Delta|$ is regarded as $(0, \infty]$.

 Let $\mu \in U$.
We prove that on some open neighborhood of $\mu$, $U_{\triv}\ni (p, y)\mapsto h(p,y)= ((\log y)^{n+1}+(\log y_0)^{n+1})^{-1}\int^y_b f(p,u)du/u$ is in $A_{2,X \times |\Delta|}$.

  We first show that $h$ extends to a continuous function on $U$ (not only on $U_{\triv}$).

  It is sufficient to prove that when $U_{\triv}\ni (p,y)\to \mu$ and $U_{\triv}\ni (p',y') \to \mu$, we have the following (1)--(3).

(1)  $A-B_1\to 0$, (2) $A-B_2\to 0$, (3) $A-B_3\to 0$, where $$A=((\log y)^{n+1}+(\log y_0)^{n+1})^{-1}\int_b^y f(p, u)(\log u)^ndu/u,$$ $$B_1=((\log y)^{n+1}+(\log y_0)^{n+1})^{-1}\int_b^y f(p', u)(\log u)^ndu/u,$$.
 $$B_2=((\log y)^{n+1}+(\log y_0)^{n+1})^{-1}\int_b^{y'} f(p, u)(\log u)^ndu/u, $$ $$B_3=((\log y')^{n+1}+(\log y_0)^{n+1})^{-1}\int_b^y f(p, u)(\log u)^ndu/u.$$

We prove (1). Let $\nu\in X_{[:]}$ be the image of $\mu$. Let $\epsilon>0$. By Proposition \ref{3.1fiber},  if $p, p'\to \nu$, we have $|f(p,u)-f(p',u)|\leq \epsilon$ on the route of integration. 
We have $|\int_b^y f(p, u)(\log u)^ndu/u- \int_b^y f(p', u)(\log u)^ndu/u|\leq \epsilon \int_b^y (\log u)^ndu/u \leq  \epsilon (n+1)^{-1}((\log y)^{n+1}+(\log y_0)^{n+1})$. This proves (1).

We prove (2). Take $C$, which is independent of $p$,  such that $|f(p,u)|\leq C$. We have $|\int_b^y f(p, u)(\log u)^ndu/u- \int_b^{y'} f(p, u)(\log u)^ndu/u| \leq 
C|\int_{y'}^y (\log u)^ndu/u|
\leq (n+1)^{-1}C((\log y)^{n+1}+(\log y')^{n+1})$. This proves (2). 

We prove (3). Let $C$ be as above. Then $|\int_b^y f(p, u)(\log u)^ndu/u|
\leq 
(n+1)^{-1}C((\log y)^{n+1}+(\log y_0)^{n+1})$ and hence $|A-B_3|\leq |((\log y)^{n+1}+(\log y_0)^{n+1})^{-1}-((\log y')^{n+1}+(\log y_0)^{n+1})^{-1}|(n+1)^{-1}C((\log y)^{n+1} +(\log y_0)^{n+1})\to 0$.

Let $D$ be an iterated derivative on the direction along 
$X$. 
  Since $D \int=\int D$ and we can do the same thing as above for $D(f)$, 
$D(h)$ is continuous on $U$. To prove that $(y\partial/\partial y)^mD(h)$ is continuous on $U$, it is sufficient to prove that $(y\partial/\partial y)^mh$ is continuous on $U$.
  
  By $(y\partial/\partial y)(h ((\log y)^{n+1}+(\log y_0)^{n+1}))=f(\log y)^n$, we have for $m\geq 1$ 
  \begin{align*}
(y\textstyle\frac{\partial}{\partial y})^{m-1}(f (\log y)^n) =&
(y\textstyle\frac{\partial}{\partial y})^m(h)((\log y)^{n+1}+(\log y_0)^{n+1})\\
&+\textstyle\sum_{j=1}^m
\textstyle\binom{m}{j}
(y\textstyle\frac{\partial}{\partial y})^{m-j}(h)
(y\textstyle\frac{\partial}{\partial y})^j((\log y)^{n+1}).
\end{align*}
By this and by the fact that $((\log y)^{n+1}+(\log y_0)^{n+1})^{-1}$ times the left-hand-side extends to a continuous function on $U$, 
$(y\partial/\partial y)^m(h)$ also extends by induction on $m$.

\end{sbpara}

\begin{sbpara}\label{cint1}  Assume still that $X$ is log smooth.

 Let $c\in \R$, $c\not=0$ and let $n\in \Z_{\geq 0}$.  
    We prove that if $f\in A_{2, X \times |\Delta|}$, then $\int_b f y^c(\log y)^n dy/y \in y^c(\log(y))^n  A_{2,X\times|\Delta|}+  y_0^c((\log y_0))^n  A_{2,X\times|\Delta|}$.
    
    We may assume $y>1$, $y_0>1$.

  It is enough to prove 
that $U_{\triv} \ni (p,y)\mapsto h(p,y)=
(y^c(\log y)^n+y_0^c(\log y_0)^n)^{-1}\int_b^y fu^c(\log u)^ndu/u$  satisfies the Cauchy condition  when $(p,y)$ converges to $\mu$, 
for, then, $h$ extends over $U$, and, 
by 
\begin{align*}
(y\textstyle\frac{\partial}{\partial y})^{m-1}(f y^c (\log y)^n)=&
(y\textstyle\frac{\partial}{\partial y})^m(h)(y^c (\log y)^n+y_0^c(\log y_0)^n)\\
&+\textstyle\sum_{j=1}^m
\textstyle\binom{m}{j}
(y\textstyle\frac{\partial}{\partial y})^{m-j}(h)
(y\textstyle\frac{\partial}{\partial y})^j(y^c (\log y)^n+y_0^c(\log y_0)^n)
\end{align*}
together with the fact that $(y^c(\log y)^n+y_0^c(\log y_0^n))^{-1}$ times the left-hand-side extends, 
$(y\partial/\partial y)^m(h)$ also extends by induction on $m$.

  It is sufficient to prove that when $U_{\triv}\ni (p,y)\to \mu$ and $U_{\triv}\ni (p',y') \to \mu$, we have the following (1)--(3).

  (1)  $A-B_1\to 0$, (2) $A-B_2\to 0$, (3) $A-B_3\to 0$, where 
  $$A= (y^c(\log y)^n+y_0^c(\log y_0)^n)^{-1}\int_b^y f(p, u)u^c(\log u)^ndu/u,$$ $$B_1=(y^c(\log y)^n+y_0^c(\log y_0)^n)^{-1}\int_b^y f(p', u)u^c(\log u)^ndu/u,$$
  $$B_2=(y^c(\log y)^n+y_0^c(\log y_0)^n)^{-1}\int_b^{y'} f(p, u)u^c(\log u)^ndu/u,$$ $$B_3= ((y')^c(\log y')^n+y_0^c(\log y_0)^n)^{-1}\int_b^y f(p, u)u^c(\log u)^ndu/u.$$
 
 To prove these, we use the elementary fact that $$\int^y u^c(\log u)^ndu/u= y^c{\textstyle\sum}_{j=0}^n a_j(\log y)^{n-j}+a_{n+1}$$ for some $a_j\in \R$ with $a_0= c^{-1}$.

We prove (1). Let $\nu\in X_{[:]}$ be the image of $\mu$. Let $\epsilon>0$. By Proposition \ref{3.1fiber},  if $p, p'\to \nu$, we have $|f(p,u)-f(p',u)|\leq \epsilon$ for $1\leq u\leq y$.
We have $|\int_b^y f(p, u)u^c(\log u)^ndu/u- \int_b^y f(p', u)u^c(\log u)^ndu/u|\leq \epsilon \int_b^y u^c(\log u)^ndu/u \leq  \epsilon O(y^c(\log y)^n+y_0^c(\log y_0)^n)$. This proves (1).

We prove (2). Take $C$, which is independent of $p$,  such that $|f(p,u)|\leq C$. We have $|\int_b^y f(p, u)u^c(\log u)^ndu/u- \int_b^{y'} f(p, u)u^c(\log u)^ndu/u|\leq 
 C |\int_{y'}^y u^c(\log u)^ndu/u|  
\leq 
C|y^c\sum_{j=0}^n a_j (\log y)^{n-j}- (y')^c\sum_{j=0}^n a_j (\log y')^{n-j}|=o(y^c(\log y)^n)$.  
This proves (2). 

We prove (3). We have $|\int_b^y f(p, u)u^c(\log u)^ndu/u|=O(y^c (\log y)^n+y_0^c(\log y_0)^n)$ and hence

$|A-B_3|= |(y^c(\log y)^n+y_0^c(\log y_0)^n)^{-1} -((y')^{-c}(\log y')^n+ y_0^c(\log y_0)^n)^{-1} |O(y^c(\log y)^n+y_0^c(\log y_0)^n)  \to 0$.

  Together with \ref{inty2}, 
for a log smooth $X$, Proposition \ref{inty0} for $(X\times |\Delta|)_{[:]}$ (Case (1))  is proved.

\end{sbpara}

\begin{sbpara}
For a log smooth $X$, Proposition \ref{inty0} for $(X\times |\Delta|)_{[:]}^{\log}$ (Case (2)) follows from this by taking $A_{2,X}^{\log} \otimes_{A_{2,X}}-$. 
\end{sbpara}

By this, for a log smooth $X$, Proposition \ref{inty0}  is proved and hence Theorem \ref{intyp} is proved.

\begin{sbpara}\label{pfyp}
We complete the proof of Theorem \ref{intyp}.  
  We give the proof for $(X\times |\Delta|)_{[:]}^{\log}$. The proof for  $(X\times |\Delta|)_{[:]}$ is similar.

The log smooth case is already shown. We consider the general case by using a strict immersion into log smooth. Let $k=2*$ or $3$. 

1. The surjectivity of  $d: A^{\log}_{k, X\times |\Delta|} \to A^{\log, 1}_{k, (X \times |\Delta|)/X}$. Locally an element of $A^{\log, 1}_{k, (X \times |\Delta|)/X}$ comes from the log smooth. Hence this is reduced to the log smooth case.

2. The injectivity of $A^{\log}_{k,X}\to A^{\log}_{k, X\times |\Delta|}$. For each non-empty open neighborhood $U$ of $p'\in (X \times |\Delta|)^{\log}_{[:]}$ with image $p\in X_{[:]}^{\log}$,  we have a $V$-point whose image is contained in $U$ for some open neighborhood $V$ of $p$ by \ref{sec2}. This proves the injectivity.

3. The exactness at the middle. Let $p'\in (X\times |\Delta|)^{\log}_{[:]}$. Let $f\in A^{\log}_{k, X\times |\Delta|}$ around $p'$ and assume $df=0$. We use a strict immersion $X\to Z$ into a log smooth $Z$ around $p'$. Take an open neighborhood $\tilde U$ of $p'$ in $(Z\times |\Delta|)_{[:]}^{\log}$ and let $U$ be the inverse image of $\tilde U$ in $(X\times |\Delta|)_{[:]}^{\log}$. We may assume that $f$ comes from $\tilde f$ on $\tilde U$. Since $df=0$, we have $d\tilde f =\sum_j a_jw_j$, where $a_j$ is an element of the kernel of $A^{\log}_{k, Z}\to A^{\log}_{k,X}$, $w_j\in A^{\log, 1}_{3, (Z\times |\Delta|)/Z}$, and we shrink $\tilde U$ if necessary. 
By using \ref{sec2}, take a $V$-point $b$ of $|\Delta|$  into $\tilde U$ for some open subset $V$ of $X^{\log}_{[:]}$. 
  By Proposition \ref{3.1fiber},  we may assume that the fibers of $\tilde U \to V$ are connected as in \ref{Ub}. 
By the log smooth case, $\tilde f = \int_b d\tilde f+\tilde f(b)$. Since $\int_b$ is $A^{\log}_{3, Z}$-linear, we have $\tilde f = (\sum_j a_j \int_b w_j) + \tilde f(b)$. Hence $f= f(b)$.

\end{sbpara}

\begin{sbpara}\label{definty1} We have considered indefinite integrals. We now consider definite integrals. Let $(V, U, b)$ be as in \ref{Ub} (1) (resp.\ (2)), and let $a$ be a $V$-point of $|\Delta|$ such that the image of $a: V \to (X\times |\Delta|)_{[:]}$ (resp.\ $V \to (X\times |\Delta|)_{[:]}^{\log}$) is contained in $U$.
  For $w\in A^1_{3, (X\times|\Delta|)/X}(U)$ (resp.\ $A^{\log,1}_{3, (X\times |\Delta|)/X}(U)$), we define
$$\int_b^a w = F(a) \in A_{3,X}(V) \;(\text{resp}.\; A^{\log}_{3,X}(V)),\quad \text{where}\;\; F=\int_b w.$$
\end{sbpara}

\begin{sbprop}\label{definty2}  
Let the situation be as in $\ref{definty1}$. 
Assume  $a<b$ ($\ref{interval}$) and let $y_0\in A_{3,V}$ be the pullback of $y$ under $a$. Assume $y_0>1$. Let $c\in \R$, $d\in \Z_{\geq 0}$. Let $R=A_2$ (resp. $R=A_2^{\log}$).

$(1)$ If $c>0$, we have $$\int_b^a \;(y^c\log(y)^d R^1_{(X\times |\Delta|)/X})(U)\subset (y_0^c\log(y_0)^dR_X)(V).$$

$(2)$ If $c<0$ or if $c=d=0$, we have $$\int_b^a\;(y^c\log(y)^dR^1_{(X\times |\Delta|)/X})(U) \subset R_X(V).$$ 
\end{sbprop}

\begin{pf} 

By using a strict immersion into log smooth, this is reduced to the case where  $X$  is log smooth.  
  Then the statements follow from the calculation in \ref{cint1}.

\end{pf}

\begin{sbprop}\label{definty3} 
Let the notation be as in the situation $(1)$ in $\ref{definty1}$.   
Assume  $a<b$ ($\ref{interval}$). 
  Let $f\in M_Y(V)$ such that $|f|<1$  and assume that $\log|q|/\log|f|$ does not have value $\infty$ on $U$.
  Let $g\in A_2(U)$ and assume that $g$ satisfies {\rm (P)} ($\ref{pos}$). 
  Let $c\in \R_{>0}$. Let $\varphi:= \int_b^a  (\log|q|/\log|f|)^cg dy/y\in A_2(V)$ (Proposition $\ref{definty2}$).

$(1)$  $\varphi$ satisfies {\rm (P)}.

$(2)$ 
Let $p\in V$ and  assume that for some $p'\in U$ lying over $p$, the value of $\log|q|/\log|f|$ at 
$p'$ is  not zero and $g(p')>0$. Then $\varphi(p)>0$.
\end{sbprop}

\begin{pf} By using a strict immersion into log smooth, these are reduced to the case where $Y$ is log smooth and then reduced to the classical theory on $Y_{\triv}$.
\end{pf}

\subsection{Path integrals in the direction of $x$}\label{ss:intx}

\begin{sbthm}\label{intxp} For $k=2, 2*, 3$, 
we have an exact sequence $$0\to A^{\log}_{k,X\times |\Delta|}\to A^{\log}_{k, X\times\Delta} \overset{d}\to A^{\log, 1}_{k, (X \times \Delta)/(X\times |\Delta|)}\to 0.$$
\end{sbthm}

\begin{sbpara}\label{Ub2} Like in \ref{Ub}, we consider a triple $(V, U, b)$, where $V$ is an open set of $(X \times|\Delta|)_{[:]}^{\log}$, $U$ is an open set of the inverse image of $V$ in $(X \times \Delta)^{\log}_{[:]}= (X \times |\Delta|)_{[:]} \times \R/\Z$  such that $U=V \times \{x\bmod \Z\;|\; x\in\R, \theta_1<x<\theta_2\}$ for some $\theta_1, \theta_2\in \R$ such that $\theta_1<\theta_2\leq \theta_1+1$, and $b$ is an element of $\R$ such that $\theta_1<b<\theta_2$.  

For $k=2, 2*, 3$,  $b$ induces a morphism of ringed spaces $(V, A^{\log}_{k, X\times |\Delta|}|_V)\to (U, A^{\log}_{k, X \times \Delta}|_U)$. For the projection $p:U\to V$, we have 
$$R^1p_*p^{-1}(A^{\log}_{k, X\times |\Delta|}|_V)=0, \quad p_*p^{-1}(A^{\log}_{k,X\times |\Delta|}|_V)=A^{\log}_{k, X\times |\Delta|}|_V.$$

\end{sbpara}

\begin{sbprop}\label{intx0} Let $k$ be either $2$, $2*$ or $3$. 
Let $(V, U, b)$ be as in $\ref{Ub2}$ and let $w\in A^{\log, 1}_{k, (X\times \Delta)/(X\times |\Delta|)}(U)$. Then there is a unique element $F 
\in A^{\log}_{k, X\times \Delta}$ such that $dF=w$ and $F(b)=0$. 

  Denote this $F$ by $ \int_b w$.
Then every $f\in A^{\log}_{k, X\times \Delta}$ can be written as $\int_b df+f(b)$.

\end{sbprop}

This proposition follows from Theorem \ref{intxp}.
 In actual, we prove this proposition together with Theorem \ref{intxp}.

These results for $A_k^{\log}$ with $k=2*, 3$ are obtained from the results for $A_2^{\log} $ by taking $A^{\log}_k \otimes_{A^{\log}_2}-$. 

\begin{sbpara}\label{pfintx} In the case where $X$ is with trivial log structure and is smooth, Theorem \ref{intxp} is classical. 

We prove Proposition \ref{intx0} in the case where $X$ is log smooth. (The log smooth case of Theorem \ref{intxp} follows from this.)

We may assume $b=0$. 

  Let $Z=X\times \Delta$.
  Let $q$ be the coordinate function of $\Delta$ and let $q=e^{2 \pi i (x+iy)}$ with $x$ and $y$ being real as usual. 
  Let $p\in Z^{\log}_{[:]}$ be a point of $U$. 
  Since $A^{\log}_{k,Z}= A^{\log}_{k,X\times |\Delta|} \otimes_{A_{2, X\times |\Delta|}} A_{2,Z}[x/y]$, 
  it is sufficient to  consider  $F=\int_0 f\cdot (x/y)^n dx/y$ for $f\in A_{2,Z}$ and $n\geq 0$ 
on $Z_{\triv}$ near $p$ and show that $F$ belongs to $A^{\log}_{2, Z}$ at $p$. 

  First we consider the case where the image of $p$ in $\Delta$ is $0\in \Delta$. 

  Then the continuity of $F$ at $p$ follows from the continuity of $f$. 

  We work on the fiber product of $Z^{\log}_{[:]} \to\R/\Z \leftarrow \R$ and regard $x$ as the coordinate function of this $\R$. 
Since the image of $p$ is $0$ and $f\in A_{2,Z}$, 
$$g_n(x) := \int_0^1 f(x+\theta)\theta^nd\theta$$ is defined. This $g_n(x)$ depends only on $x\bmod \Z$ and 
belongs to $A_{2,Z}$ as is seen by taking iterated derivatives. 

  Further, we have
\medskip

(1) $F(x+1)-F(x)=y^{-n-1} \int_0^1 f(x+\theta)(x+\theta)^n d\theta=\sum_{r=0}^n \Bigl(\begin{matrix} n\\r\end{matrix}\Bigr) g_r(x) x^{n-r}/y^{n+1}$. 

\medskip
\noindent
Let $P_n(x)$ be the polynomial of degree 
$n$ such that $P_n(x+1)-P_n(x)= x^{n-1}$ and $P_n(0)=0$. Let 

\medskip

 (2) $h(x):= F-G \quad \text{with}\;\; G=\sum_{r=0}^n \Bigl(\begin{matrix} n\\r\end{matrix}\Bigr) g_r(x) P_{n+1-r}(x)/y^{n+1}.$

\medskip
\noindent
By (1), we have $h(x+1)=h(x)$. To prove $F\in A_{2,Z}^{\log}$, it is sufficient to prove $h(x)\in A_{2,Z}$. 
By the continuity of $F$, we have that $h$ is continuous. Applying $y\partial/\partial x$ to (2), we obtain

\medskip

(3) $y\partial h(x)/\partial x= f -y\partial G/\partial x$.

\medskip
\noindent
Hence the left-hand-side is continuous. 
  When we  apply $y\partial/\partial x$ more times, the left-hand-side is still continuous. Derivatives in other directions commute with our integration. Hence $h$ belongs to $A_2$. 
  Thus the case where the image of $p$ in $\Delta$ is $0$ is proved.

   Finally, if the image of $p$ in $\Delta$ is not $0\in \Delta$, 
considering $(X\times(\Delta\smallsetminus \{0\}))_{[:]}^{\log}=X^{\log}_{[:]}\times (\Delta\smallsetminus \{0\})=P \times \R/\Z$, where $P:=(X \times (|\Delta|\smallsetminus \{0\}))^{\log}_{[:]}= X^{\log}_{[:]}\times (|\Delta|\smallsetminus \{0\})$, we see that 
the integral is classical and $F$ clearly belongs to $A_2^{\log}$. 

\end{sbpara}
\begin{sbpara}\label{pfxp} Using the above log smooth case, the proof of Theorem \ref{intxp} is similar to the proof of Theorem \ref{intyp} in \ref{pfyp}. 
\end{sbpara}

\begin{sbpara}\label{defintx1} We have discussed indefinite integrals. Now we discuss definite integrals. Let $k=2, 2*$ or $3$.

(1) Let $(V, U, b)$ and $w$ be as in Proposition \ref{intx0}, and let $a\in \R$ satisfying $\theta_1< a<\theta_2$. We define $$\int_b^a w = F(a)\in A^{\log}_{k,X}(V), \quad \text{where}\;\; F=\int_b w.$$

(2) Let $V$ be an open set of $(X\times |\Delta|)_{[:]}^{\log}$ and let $U$ be the inverse image of $V$ in $(X\times \Delta)^{\log}_{[:]}$. Let $w\in A^{\log, 1}_{k,(X\times \Delta)/(X\times |\Delta|)}(U)$. For $b, a\in \R$ such that $b<a$, we define $\int_b^a w\in A^{\log}_{k,X}(V)$ as the sum $\sum_{j=1}^n \int_{b(j-1)}^{b(j)} w$, where $b(0)=b$, $b(n)= a$, $b(j-1)<b(j)< b(j-1)+1$ for $1\leq j\leq n$. Then $\int_b^a w$ is independent of the choices of $b(j)$. Furthermore $\int_b^{b+1} w$ is independent of the choice of $b\in \R$. We denote this $\int_b^{b+1} w$ by $\oint w$. 

\end{sbpara}

\begin{sbprop}\label{defintx2}  Both in the cases $(1)$ and $(2)$ in $\ref{defintx1}$, for $c\in \R$ and $d \in \Z_{\geq 0}$,  $\int_b^a$ sends $y^c \log(y)^d A^{\log,1}_{2, (X\times \Delta)/(X \times |\Delta|)}(U)$ to $y^{c-1}\log(y)^d A^{\log}_{2, X\times |\Delta|}(V)$. 
\end{sbprop} 

\begin{pf}   This is reduced to the case $c=1$ and $d=0$. By using a strict immersion to log smooth, this is reduced to the case where $Y$ is log smooth. 
   Let $p\in Z_{[:]}^{\log}$, where $Z=X \times \Delta$.
It is sufficient to consider  $h:=\int_b^a f\cdot (x/y)^n dx$ for $f\in A_{2,Z}$ and $n\geq 0$ 
on $Z_{\triv}$ near $p$ and show that $h$ belongs to $A_{2, X\times |\Delta|}$ at $p$. (Note that $dx/y$ in the definition of $F$ in \ref{pfintx} is replaced by $dx$ here.) The continuity of $h$ at $p$ follows from the continuity of $f$.  Since the iterated derivatives with respect to a local base of $A^1_{2,X\times |\Delta|}$ commute with $\int_b^a$, this proves $h\in A_{2,X\times |\Delta|}$ at $p$. 
\end{pf}

\begin{sbprop}\label{defintx3}
 Let $V$ be an open set of $(X\times|\Delta|)_{[:]}$ and let $U$ be the inverse image of $V$ in $(X\times\Delta)_{[:]}$. Let $g\in A_2(U)$ and assume that $g$ satisfies {\rm (P)}. Let 
 $\varphi:= \oint  g dx\in A_2(V)$ (Proposition $\ref{defintx2}$).

$(1)$  $\varphi$  satisfies {\rm (P)}.  

$(2)$ Let $p\in V$ and  assume that $g(p')>0$ for some $p'\in U$ lying over $p$.
Then $\varphi(p)>0$.
\end{sbprop}

\begin{pf}
By using a strict immersion into log smooth, these are reduced to the case where $X$ is log smooth and then to the classical theory on $X_{\triv}$. 
\end{pf}

\subsection{Log Poincar\'e lemma}\label{ss:LPL}

  We prove the log Poincar\'e lemma. 

\begin{sbthm}\label{LPL} 
  For a log smooth saturated morphism $X\to Y$ of fs log analytic spaces,  $A_{2*,Y}^{\log}\to A_{2*,X/Y}^{\log, \bullet}$ and $A_{3,Y}^{\log}\to A_{3,X/Y}^{\log,\bullet}$ are quasi-isomorphisms. 
\end{sbthm} 

\begin{pf}
  This can be reduced to the statements of relative $\R$-dimension one treated in Theorem \ref{intyp} and Theorem \ref{intxp}, as follows.

By Theorem \ref{SEisS}, we are reduced to the case $X=Y\times \Delta^n$, where $\Delta^n$ is endowed with the log structure given by the coordinate functions. 

Let $k=2*$ or $3$. 
  For $0\leq j\leq 2n$, let  $P_j= Y\times \Delta^{n-m}$ if $j$ is an even integer $2m$, and $P_j= Y\times \Delta^{n-m} \times |\Delta|$ if $j$ is an odd integer  $2m-1$. Hence  $P_0=X$ and $P_{2n}=Y$.
Let $S(j)$ be the statement that the map
 $$A^{\log}_{k,P_j} \to A^{\log,\bullet}_{k, X/P_j}$$
 is a quasi-isomorphism,  
 where we use the notation concerning differential modules in \ref{A|Del|}.
Then $S(2n)$ is the log Poincare lemma and  $S(0)$ is evident. 
We prove $S(j)$ by induction on $j$. 

Assume $1\leq j\leq 2n$ and assume that we already know that $S(j-1)$ is valid. Consider the morphisms $X=P_0\to P_{j-1}\to P_j$. 
 We have an exact sequence

\begin{equation}\tag{$*$}
0\to A^{\log, \bullet-1}_{k,X/P_{j-1}} \otimes_{A^{\log}_{k, P_{j-1}}} A^{\log, 1}_{k,P_{j-1}/P_j}\to A^{\log, \bullet}_{k,X/P_j} \to A^{\log, \bullet}_{k,X/P_{j-1}} \to 0.
\end{equation}

From this exact sequence $(*)$ and $S(j-1)$, we have an exact sequence
$0\to \cH^0(A^{\log, \bullet}_{k,X/P_j}) \to A^{\log}_{k,P_{j-1}} \to A^{\log,1}_{k,P_{j-1}/P_j} \to \cH^1(A^{\log, \bullet}_{k,X/P_j})\to 0$ and the other $\cH^p(A_{k,X/P_j}^{\log, \bullet})$ are $0$. Hence 
in the case where $j$ is even (resp.\ odd), we can go from $S(j-1)$ to $S(j)$ by the result Theorem \ref{intyp} (resp.\ Theorem \ref{intxp}) for $P_{j-1}/P_j$. 
\end{pf}

The following is a log version of the classical theorem of de Rham. 
\begin{sbthm}\label{LPLcoh}  Let $f:X\to Y$ be a proper log smooth saturated morphism of fs log analytic spaces. 

$(1)$ Denote the map $X^{\log}_{[:]}\to Y_{[:]}$ by $g$. Then we have a canonical isomorphism of $A_{3,Y}$-modules
$$Rg_*(A_{3,Y}^{\log})\cong f_{[:]*}(A^{\bullet}_{3,X/Y}).$$

$(2)$ We have 
a canonical isomorphism of $A^{\log}_{3,Y}$-modules
$$A^{\log}_{3,Y}\otimes_{\R} Rf^{\log}_{[:]*}(\R) \cong A^{\log}_{3,Y}\otimes_{A_{3,Y}}f_{[:]*}(A^{\bullet}_{3,X/Y}).$$
In particular, for $p\in Y_{[:]}$ and for $\tilde p\in Y_{[:]}^{\log}$ lying over $p$, and for each $m$, we have an isomorphism of $A_{3,Y,p}$-modules between  $A_{3,Y,p}\otimes_{\R} R^mf^{\log}_{[:]*}(\R)_{\tilde p}$ and the stalk at $p$  of the $m$-th cohomology sheaf of 
the complex $f_{[:]*}(A_{3,X}) \overset{d}\to f_{[:]*}(A_{3,X/Y}^1) \overset{d}\to f_{[:]*}(A_{3,X/Y}^2)\to \cdots$.
\end{sbthm}

\begin{pf} (1) We have an isomorphism
$$Rg_*(A_{3,Y}^{\log})\overset{\cong}\to Rg_*(A^{\log,\bullet}_{3,X/Y})$$
by the log Poincar\'e lemma (Theorem \ref{LPL}). On the other hand, we have isomorphisms
$$f_{[:]*}(A^{\bullet}_{3,X/Y})\overset{\cong}\to Rf_{[:]*}(A^{\bullet}_{3,X/Y})\overset{\cong}\to Rg_*(A^{\log,\bullet}_{3,X/Y}).$$
Here the first homomorphism 
is an isomorphism by Proposition \ref{soft3} (1), and the second one 
is an isomorphism by Corollary \ref{c:pol} (2). Here to prove that these homomorphisms are isomorphisms, we may work locally on $Y$ and hence we can assume that $Y$ (and hence $X$) satisfies the condition $(*)$ in \ref{(*)}. 
  From these isomorphisms, we obtain the desired isomorphism. 

(2)
We have  $Rg_*(A^{\log}_{3,Y})= R\tau_{Y*}Rf^{\log}_{[:]*}(A_{3,Y}^{\log})= R\tau_{Y*}(A_{3,Y}^{\log}\otimes_{\R} Rf^{\log}_{[:]*}(\R))$.
By  Claim 2 below, the  first half of (2) is obtained from (1) by applying 
$A^{\log}_{3,Y}\otimes_{A_{3,Y}}-$ to the isomorphism of (1).

{\bf Claim 1.} Let $H$ be a local system of finite-dimensional $\R$-vector spaces on $Y^{\log}$ such that for each $s\in Y$, the action of $\pi_1(s^{\log})$ on a stalk of $H$ at a point of $s^{\log}$ is unipotent. Then the 
canonical homomorphism 
$A_{3,Y}^{\log}\otimes_{A_{3,Y}}  R\tau_{Y*}(A_{3,Y}^{\log}\otimes_{\R} H) \to A^{\log}_{3,Y}\otimes_{\R} H$ is an isomorphism. 

In fact, this is reduced to the case where $H$ is the constant sheaf $\R$ and is reduced to $A_{3,Y}\overset{\cong}\to R\tau_{Y*}(A_{3,Y}^{\log})$ (Corollary \ref{c:pol} (2)). 

{\bf Claim 2.} The canonical homomorphism $A_{3,Y}^{\log}\otimes_{A_{3,Y}}  R\tau_{Y*}(A_{3,Y}^{\log}\otimes_{\R} Rf^{\log}_{[:]*}(\R))\to A^{\log}_{3,Y}\otimes_{\R} Rf^{\log}_{[:]*}(\R)$
is an isomorphism.

This follows from Claim 1 because for any $m$, $R^mf^{\log}_*\R$  satisfies the condition on $H$ in Claim 1 (\cite{IKN} Theorem (6.3)(1)).

The latter half of (2) is obtained from the first half by taking a ring homomorphism $A^{\log}_{3, Y, \tilde p}\to A_{3,Y,p}$ over $A_{3,Y,p}$. 
\end{pf}

\begin{sbrem} As above,  we have the log Poincar\'e lemma by using log $C^{\infty}$ functions. 
For our log real analytic functions, we do not have the log Poincar\'e lemma even for $X=\Delta$ with the log structure given by the coordinate function $q$. In fact, for $\omega:= q \bar q dy/y\in A^{\an,1}_{1,X}$, we have  $d\omega=0$ but locally at $0\in \Delta$, there is no $f\in A^{\an}_{(3), X}$ such that $ df=\omega$.

\end{sbrem}

\subsection{Integration along fibers}\label{ss:2nint}

\begin{sbpara}\label{2nint1}
Let the notation be as in \ref{2n}.
  In particular, $U$ is an open set of $X_{[:]}$ (resp.\ $X_{[:]}^{\log}$), $V$ is an open set of $Y_{[:]}$ (resp.\ $Y^{\log}_{[:]}$), and every point of $U$ is vertical over $Y$ in the sense of \ref{verrat}. 
  Let $u: U \to V$ be the induced map. Let $k=2*$ or $3$. We define  a homomorphism
$$\int_{U/V} : u_! (A^{2n}_{k,X/Y}|_U)\to A_{k,Y}|_V \quad (\text{resp}.  \;\;  u_!(A^{\log,2n}_{k,X/Y}|_U) \to A^{\log}_{k,Y}|_V)$$
discussed in \ref{2n} as below. Here $(-)_!$ is the direct image with proper support.

Call the case $U\subset X_{[:]}$ and $V\subset Y_{[:]}$ Situation 1, and call the case $U\subset X_{[:]}^{\log}$ and $V\subset Y^{\log}_{[:]}$ Situation 2. 

The map $\int_{U/V}$ in Situation 1 is obtained from the map $\int_{U/V}$ in Situation 2 as follows. Assume that we are in Situation 1,  let $U^{\log}\subset X^{\log}_{[:]}$ and let $V^{\log}\subset Y^{\log}_{[:]}$ be the inverse images of $U$ and $V$, respectively, and denote the map $U^{\log}\to V^{\log}$ by $u^{\log}$. Let $\tau_U: U^{\log}\to U$ and $\tau_V: V^{\log}\to V$ be the canonical maps.
By Corollary \ref{c:pol} (1), the map $\int_{U/V}$ in Situation 1 is obtained from that in Situation 2  as 
$$u_!(A_{k,X/Y}^{2n}|_U)
= u_!\tau_{U*}(A_{k, X/Y}^{\log, 2n}|_{U^{\log}})
= (u \tau_U)_!(A_{k, X/Y}^{\log, 2n}|_{U^{\log}}) 
= (\tau_V u^{\log})_!(A_{k, X/Y}^{\log, 2n}|_{U^{\log}})$$ $$= \tau_{V*}u^{\log}_!(A_{k,X/Y}^{\log,2n}|_{U^{\log}}) \to \tau_{V*}(A_{k, Y}^{\log}|_{V^{\log}})=A_{k, Y}|_V.$$

\end{sbpara}

\begin{sbpara}\label{tr0}  In \ref{tr0}--\ref{tr2}, assume that we are in Situation 2.

We define the cohomology class map $u_!(A^{2n}_{k,X/Y}|_U) \to A^{\log}_{k,Y}|_V\otimes_\R R^{2n}u_!\R$ as follows: 
$$u_!(A^{\log, 2n}_{k, X/Y}|_U) \to  R^{2n}h_!(A^{\log, \bullet}_{k, X/Y}|_U)= R^{2n}u_!(A^{\log}_{k,Y}|_U)=A^{\log}_{k,Y}|_V \otimes_\R R^{2n}u_!\R.$$ 
Here, the first equality is by the log  Poincar\'e lemma (Theorem \ref{LPL}).
\end{sbpara}

\begin{sbpara}\label{trmap}  
  By Remark \ref{r:rounding} (2), 
the map $h: U\to V$ is topologically submersive whose all fibers are topological manifolds without boundaries. 
  Hence we have a canonical isomorphism 
$u^!\R\cong  \R[2n]$ 
by the Verdier duality (cf.\ \cite{Ve}). 
  This induces the trace map $$R^{2n}u_!\R\to \R.$$

\end{sbpara}

\begin{sbpara}\label{tr2}

The composite map $u_!(A^{\log, 2n}_{k,X/Y}|_U)\to A^{\log}_{k,Y}|_V\otimes_\R R^{2n}u_!\R\to A^{\log}_{k,Y}|_V$ is called the {\it integration along the fibers} and denoted by $\int_{U/V}$.

 \end{sbpara}

\begin{sbpara}\label{char0}
     By the construction of the integration $\int_{U/V}$, it has the properties (i), (ii), and (iii) stated in \ref{2n}.

In Situation 1, as is proved in \ref{intchar} below,  our integration is characterized as a map by 
 the properties  (ii), (iii) in \ref{2n} and  the following weaker version (i)$'$ of (i) in \ref{2n},  and the property (iv) below. 

\medskip

(i)$'$   $\int_{U/V}$ is additive in $w$. 

(iv) If $U'$ is an open set of $U$ and the support of $w$ is contained in $U'$, and if $w'$ denotes the restriction of $w$ to $U'$, we have $\int_{U/V} w=\int_{U'/V} w'$.

\end{sbpara}

\begin{sbpara}\label{XYZ}
  Let $f:X \to Y$ and $g:Y \to Z$
be log smooth saturated morphisms of fs log analytic spaces of relative dimension $n$ and $m$, respectively. Let $U, V, W$ be open sets of $X_{[:]}$, $Y_{[:]}$, $Z_{[:]}$ (resp.\ $X^{\log}_{[:]}$,  $Y^{\log}_{[:]}$,  $Z^{\log}_{[:]}$), respectively, such that these morphisms induce maps $U\to V$ and $V\to W$ which we denote by $u$ and $v$, respectively. 
Assume that $U$ is  vertical over $Y$ and $V$ is vertical over $Z$ (and hence $U$ is vertical over $Z$ by \ref{verrat}.2). Then we have
$$\int_{U/W}=\int_{V/W} \circ \int_{U/V}$$
as homomorphisms $(vu)_!(A^{2(m+n)}_{3, X/Z}|_U) \to A_{3, Z}|_W$ (resp.\ $(vu)_!(A^{\log, 2(m+n)}_{3, X/Z}|_U) \to A^{\log}_{3, Z}|_W$). Here $\int_{U/V}$ means 
\begin{align*}
\int_{U/V}\otimes 1 : (vu)_!(A^{2(m+n)}_{3, X/Z}|_U) &= (vu)_!(A^{2n}_{3, X/Y}|_U \otimes_{A_{3,Y}|_V} A^{2m}_{3,Y/Z}|_V) \\
&= v_!(u_!(A^{2n}_{3, X/Y}|_U) \otimes_{A_{3,Y}|_V} A^{2m}_{3,Y/Z}|_V)\to v_!(A^{2m}_{3,Y/Z}|_V),
\end{align*}
where the first equality follows from
 Corollary  \ref{SELF}.

This follows from the fact that the trace map for $U/W$ is the composition of the trace maps for $U/V$ and $V/W$. 
\end{sbpara}

\begin{sblem}\label{lift} Assume that we are  in Situation $1$ with $X,Y, U, V$ as in $\ref{2nint1}$. 
  Assume $X=X' \times_{Y'} Y$ for a strict closed immersion $Y\to Y'$ and for a log smooth saturated morphism $X'\to Y'$ of fs log analytic spaces. 
  Assume that we are given an open set $V'$ of $Y'_{[:]}$ such that $V = V'\cap Y_{[:]}$. Let $\cF$ be a sheaf of $A_{1,X'}$-modules on $X'_{[:]}$, let $\cG$ be a sheaf of abelian groups on  $X_{[:]}$,  and assume that we are given a surjective homomorphism from the inverse image of $\cF$ on $X_{[:]}$ to $\cG$. Let $w\in \cG(U)$,  and assume  that the support of $w$ is proper over $V$.  Then locally on $V'$, there are an open set $U'$ of $(X')_{[:]}$ such that $U=U' \cap X_{[:]}$,  such that the image of $U'$ in $(Y')_{[:]}$ is contained in $V'$, and such that $U'$ is vertical over $Y'$, and an element $w'$ of $\cF(U')$ whose image in  $\cG(U)$ is $w$ and whose support is proper over $V'$.

\end{sblem}

\begin{pf} Let $C$ be the support of $w$. 
  Working locally on $Y'$, we may assume that $Y'$ satisfies the condition $(*)$ in \ref{(*)}. Since $C$ is proper over $V$ and hence over $V'$, we may assume that $C$ is covered by the inverse images of finitely many open sets $D_{\nu}$ of $X'$ which are isomorphic to a closed analytic subspace of $\C^{n(\nu)}$ for some $n(\nu)$. Replacing $X'$ by the union $\bigcup_{\nu}\; D_{\nu}$,  we may assume that $X$ satisfies the condition $(*)$ in \ref{(*)}. 

Take an open set $U'$ of $X'_{[:]}$ such that $U=U'\cap X_{[:]}$ and such that the image of $U'$ in $(Y')_{[:]}$ is contained in $V'$. 

Since $C$ is proper over $V'$, locally on $V'$, we have a finite family of open sets  $E_{\la}$ of $U'$
whose closure is proper over $Y'$ such that $\bigcup_{\la} E_{\la}$ contains $C$  and such that for each $\lambda$, 
we have a section  $w'_{\la} \in \cF(E_{\la})$ 
whose pullback to $\cG(X_{[:]}\cap E_{\la})$ coincides with the pullback $w_{\la}$ of $w$.  Let  $\La'=\La\;{\scriptstyle \coprod}\;\{e\}$ 
and $E_e= U'\smallsetminus C$. 
By Theorem \ref{softhm}, we have  a partition of unity $(s_{\la})_{\la\in \La'}$ in $A_1(U')$ such that the support of $s_{\la}$ is contained in $E_{\la}$ for each $\la\in \La'$. Let $w'\in \cF(U')$  be  the sum of the $0$-extensions of $s_{\la}w'_{\la}$ to $U'$ for all $\la\in \La$ (not for all $\la\in \La'$). Since the pullback of $s_{\la} w'_{\la}$ to $U$ is $s_{\la}w_{\la}$, the image of $w'$ in $\cG(U)$ is $w$. Furthermore, the support of $w'$ is proper over $V$. 
\end{pf}

\begin{sbpara}\label{intchar}We give the proof of the characterization of the  integration in Situation 1 stated in \ref{char0}.

By the argument in the first part of the proof of Lemma \ref{lift}, we see that we may assume that $X$ satisfies the condition $(*)$ in \ref{(*)}. 
  Locally on $X$ and $Y$, we have a strict immersion $Y\to Y'$ into log smooth and a log smooth saturated morphism $X'\to Y'$ such that $X=X' \times_{Y'}Y$. By the partition of unity (Theorem \ref{softhm}), we are reduced to the case where we have such $Y\to Y'$ and $X'\to Y'$. 

Then by Lemma \ref{lift} and by using the comments after (iii) in \ref{2n}, we are reduced to the classical case (ii) in \ref{2n}.

\end{sbpara}

We have the following relation between  the integration in Situation 1 in this Section \ref{ss:2nint} and  the integrations in Sections \ref{ss:inty} and \ref{ss:intx}.

\begin{sbprop}\label{6=34} In Situation $1$ in $\ref{2nint1}$, assume 
$X=Y\times \Delta$. Let $w\in A^2_{3,X/Y}(U)$  and assume that the support of $w$ is proper over $V$. 

$(1)$ Locally on $Y$, there are $V$-points $a,b$ of $|\Delta|$ such that $a<b$ ($\ref{interval}$) and such that the inverse image of $\{p\in (Y \times |\Delta|)_{[:]}\;|\; a<p<b\}$ ($\ref{interval}$) in $(Y\times \Delta)_{[:]}$ contains the support of $w$.  

$(2)$ Let $a, b$ be $V$-points of $|\Delta|$ satisfying the condition in $(1)$. Then 
$$\int_{U/V} w= \int_b^a   \oint w.$$
Here 
$\int_b^a$ is the integration in Section $\ref{ss:inty}$ defined as in $\ref{definty1}$, and 
$\oint$ is  the integration in Section $\ref{ss:intx}$ defined as in $\ref{defintx1}$. 
  To define $\oint w$, we identify $w$ with its $0$-extension from $U$ to the inverse image of $V$ in $(Y \times \Delta)_{[:]}$.
\end{sbprop}

\begin{pf} (1) follows from Proposition \ref{3.1fiber}. 

We prove (2). 
By the arguments in the first part of the proof of Lemma \ref{lift}, we see that we can assume that $X$ and $Y$ satisfy the condition $(*)$ in \ref{(*)}. Locally on $Y$, we have a strict immersion from $Y$ into log smooth. 
  Hence, by the partition of unity (Theorem \ref{softhm}), we can assume that we are given a strict immersion $Y\to Y'$ into log smooth. By applying Lemma \ref{lift} replacing $X$ and $X'$ there by $Y \times \Delta$ and $Y'\times \Delta$ here, respectively, and by the comment after the condition (iii) in \ref{2n}, we are reduced to the case where $Y$ is log smooth and then to the case where the log structure of $Y$ is trivial and $Y$ is smooth. In the last case, the statement is classical. 
\end{pf}

\begin{sbprop}\label{2*to2} 
Let $\cE$ be the $A_{2,X}$-submodule of $A^{2n}_{3,X/Y}$ defined as follows. For an open set $W$ of $X_{[:]}$, $\cE(W)$ is the set of all $w\in A^{2n}_{3,X/Y}(W)$ satisfying the following condition. For every $p\in W$, if $q_1, \dots, q_n$ denote local coordinate functions of $X$ over $Y$ (Theorem $\ref{SEisS}$), then for  $y_j=-(2\pi)^{-1}\log|q_j|$, 
the image of $w$  in the stalk $A^{2n}_{3,X/Y,p}$ belongs to $(\prod_{j=1}^n y_j^{a(j)})A^{2n}_{2,X/Y,p}$ for some  $a(j)<1$ ($1\leq j\leq n$). 

Then $\int_{U/V}$ sends $u_!(\cE|_U)$ to $A_{2,Y}|_V$. 
\end{sbprop}

 For the proof, see Proposition \ref{tildeJ} below.

\begin{sbpara}\label{2*to2b}
In the case $n>0$, $A^{2n}_{2*, X/Y}\subset \cE$ and hence Proposition \ref{2*to2} shows that $\int_{U/V}$ sends $A^{2n}_{2*,X/Y}$ to $A_{2,Y}$.

\end{sbpara}

The following Proposition \ref{tildeJ} is  slightly stronger than Proposition \ref{2*to2}. (If we take only  $h_j=1$ and $c(j)<0$ for all $j$ in the definition of $\tilde \cE$ below, then we get $\cE$.)  

\begin{sbprop}\label{tildeJ}
Let $\tilde \cE\supset \cE$ be the $A_{2,X}$-submodule of $A^{2n}_{3,X/Y}$ defined as follows. For an open set $W$ of $X_{[:]}$, $\tilde \cE(W)$ is the set of all $w\in A^{2n}_{3,X/Y}(W)$  satisfying the following condition. For every $p\in W$, if $q_1, \dots, q_n\in M_{X,p}$ denote local coordinate functions of $X$ over $Y$ (Theorem $\ref{SEisS}$), then for $y_j=-(2\pi)^{-1}\log|q_j|$, 
the image of $w$  in the stalk $A^{2n}_{3,X/Y,p}$ belongs to $(\prod_{j=1}^n h_jy_j^{c(j)})(\prod_{j=1}^n y_j)A^{2n}_{2,X/Y,p}$ for some  $h_j\in A_{3,Y, p}$ and $c(j)\in \R\smallsetminus \{0\}$ ($1\leq j\leq n$) such that $ h_jy_j^{c(j)}\in A_{2,X,p}$ for all $j$.  

 Then $\int_{U/V}$ sends $u_!(\tilde \cE|_U)$ to $A_{2,Y}|_V$. 
\end{sbprop}

\begin{pf}
By the argument in the first part of the proof of Lemma \ref{lift}, we see that we 
 may assume that $X$ and $Y$ satisfy the condition $(*)$ in \ref{(*)}. 

By Theorem \ref{SEisS} and by the partition of unity, we are reduced to the case $X=Y \times \Delta^n$ and we may assume that we can take $h_j\in A_3(U)$ and $c(j)\in \R\smallsetminus\{0\}$ which satisfy the condition of $h_j$ and $c(j)$ at every point of $U$. 
For $1 \leq j \leq n-1$, regard $\Delta^j$ as the quotient of $\Delta^n$ forgetting the $k$-th components for $j<k\leq n$. 
  Let the open set $U_j$ be the image of $U$ in $(Y \times \Delta^j)_{[:]}$, 
every point of which is vertical over $Y$ by \ref{verrat}.2. 
  Then by \ref{XYZ}, we have
$\int_{U/V}= \int_{U_1/V} \circ \dots \circ \int_{U/U_{n-1}}$. 
  For $u\in A^2_{2, (Y\times \Delta^n)/(Y\times \Delta^{n-1})}(U)$, we have  
$$\int_{U/U_{n-1}} ({\textstyle\prod}_{j=1}^n h_j y_j^{c(j)})({\textstyle\prod}_{j=1}^n y_j) u = ({\textstyle\prod}_{j=2}^n h_jy_j^{c(j)})({\textstyle\prod}_{j=2}^n y_j) \int_{U/U_{n-1}} h_1y_1^{c(1)}y_1u.$$
By this and by induction on $n$, 
we are reduced to the case $X=Y \times \Delta$. 
Working locally on $Y$, take $V$-points $a,b$ of $|\Delta|$ as in Proposition \ref{6=34} (1). 
  Again by the partition of unity, we may assume that the image of $a: V \to (Y \times|\Delta|)_{[:]}$ is contained in the image of $U$. 
  Then by 
 Proposition \ref{6=34} (2), writing $h_1$ and $c(1)$ by $h$ and $c$, respectively, we have  
$$\int_{U/V} hy^cyA^2_{2,X/Y}|_U= h \int_b^a y^c \oint A_{2,X} dx\wedge dy/y \subset h\int_b^a y^cA_{2, Y\times |\Delta|}dy/y \subset A_{2,Y},$$
where the first $\subset$ follows from 
Proposition \ref{defintx2} and the second $\subset$ follows from 
 Proposition \ref{definty2}. 
  Note here that, when $c>0$, we have $hy_0^c \in A_{2,Y}$ by the above assumption on $a$.
\end{pf}

\begin{sbprop}\label{tildeJ2} 
 Let $w\in A^{2n}_{3,X/Y}(U)$ and assume that the support of $w$ is proper over $V$. Fix $p\in V$. Assume that for each $p'\in U$ lying over $p$, there are 
local coordinate functions $q_1, \dots, q_n\in M_{X,p'}$ of $X$ over $Y$, $f_1, \dots, f_n\in M_{Y,p}$ such that $|f_j|<1$ and such that  $y_j/y_{f_j}$ do not have value $\infty$, $c(j)>0$, and $g \in A_{2,X,p'}$ satisfying {\rm (P)} ($\ref{pos}$),  for which 
$$w= ({\textstyle \prod}_{j=1}^n (y_j/y_{f_j})^{c(j)})({\textstyle\prod}_{j=1}^n y_j) gdx_1/y_1\wedge dy_1/y_1\wedge \dots \wedge dx_n/y_n \wedge dy_n/y_n.$$  
  Here $q_j=e^{2 \pi i(x_j+iy_j)}$ with $x_j$ and $y_j$ being real and 
$y_{f_j}:= -(2\pi)^{-1}\log |f_j|$ ($1 \le j \le n$). 
   Let $\varphi:=\int_{U/V} w\in A_{2,Y,p}$ (Proposition $\ref{tildeJ}$).  

$(1)$ 
$\varphi$ satisfies {\rm (P)}.

$(2)$ %
  Assume that for some $p'\in U$ lying over $p$, the values of $y_j/y_{f_j}$  at 
$p'$ are  not zero for $1\leq j\leq n$, and $g(p')>0$. Then $\varphi(p)>0$.
\end{sbprop}

\begin{pf} The proof is similar to the proof of Proposition \ref{tildeJ}. 
  By using Theorem \ref{SEisS} and partition of unity with (P)  (Theorem \ref{softhm}),  we are reduced to the case $X=Y\times \Delta^n$. By using the induction on $n$ as in the proof of Proposition \ref{tildeJ}, we are reduced to  the case $X=Y\times \Delta$. Then we are reduced to Propositions \ref{definty3} and \ref{defintx3}.    
\end{pf}

\begin{sbrem}\label{nonvert}
In a forthcoming paper \cite{KNUi}, without assuming the verticality in \ref{2nint1}, we define $\int_{U/V}: u_!(\cE) \to A_{2,Y}$, where $\cE$ is as in Proposition \ref{2*to2}. 
 The method is to start from the integrations of   Sections \ref{ss:inty} and \ref{ss:intx}.  In this method, the formula in Proposition \ref{6=34} (2) becomes the definition of the left-hand-side by using the right-hand-side.
\end{sbrem}

\subsection{Integration of closed $1$-forms}\label{ss:intcoh2}

\begin{sbpara}\label{d=0}
 Let $T$ be a $C^{\infty}$ manifold and let $w$ be a  $C^{\infty}$  $1$-form on $T$ such that $dw=0$. Then:

(1) If $T$ is contractible, then for $b \in T$,  we have a  $C^{\infty}$ function $F= \int_b w$ on $T$ characterized by $dF=w$ and $F(b)=0$. 

(2.1) We have a homomorphism $H_1(T,\Z)\to \R\;;\; \gamma \mapsto \int_{\gamma} w$ called a period integral of $w$. This implies:

(2.2) We have the integral  $\int_{\gamma} w$ for a loop $\gamma: \R/\Z\to T$  which depends only on the homotopy class of $\gamma$.  

We discuss  log relative versions of these  by using our log Poincar\'e lemma (Theorem \ref{LPL}). 

In this Section \ref{ss:intcoh2}, let $f: X\to Y$ be a log smooth saturated morphism, and let $k$ be $2*$ or $3$.  Let

\medskip

$A=A^{\log}_{k, Y}$ and $A^{\log,1}_{k, X/Y,d=0}= \text{Ker}(d: A^{\log,1}_{k, X/Y}\to A^{\log,2}_{k, X/Y})$.

\medskip

By the log Poincar\'e lemma (Theorem \ref{LPL}), we have 
an exact sequence

\medskip

\noindent \ref{d=0}.1.   $0\to A\to A^{\log}_{k,X}\to A^{\log, 1}_{k,X/Y,d=0}\to 0$. 
\end{sbpara}

\begin{sbpara}\label{H1.5} We discuss a log relative version of \ref{d=0} (1). Let $U$ be an open set of $X_{[:]}^{\log}$ and let $V$ be the image of $U$ in $Y^{\log}_{[:]}$. (Then $V$ is an open set of $Y^{\log}_{[:]}$ by Remark \ref{r:rounding} (1).) Assume that the projection  $p:U\to V$ is a homotopical isomorphism over $V$. Assume that we are given a section $b : Y\to X$ of $f$ and assume that $b$ induces $V\to U$.

 Then for $w\in A^{\log, 1}_{k,X/Y, d=0}(U)$, by the exact sequence \ref{d=0}.1, we have a unique element $F$ of $A^{\log}_{k, X}(U)$ such that $dF=w$ and $F(b)=0$. 

We denote this $F$ by $\int_b w$. 

If we have another section $a:Y\to X$ which induces $V\to U$, we define $\int_b^a w:= F(a)\in A^{\log}_{k, Y}(V)$.

\end{sbpara}

\begin{sbpara}\label{H1.1} We discuss a log relative version of \ref{d=0} (2.1). Assume that $f:X\to Y$ is proper. Then 
for $f^{\log}: X^{\log}_{[:]}\to Y^{\log}_{[:]}$, $R^1f^{\log}_*\Z$ is a locally constant sheaf of finitely generated $\Z$-modules by \cite{KN}. Let  $V$ be an open set of $Y^{\log}_{[:]}$ and let $U$ be the inverse image of $V$ in $X_{[:]}^{\log}$ and let  $\gamma$ be a homomorphism $R^1f^{\log}_*\Z|_V\to \Z$ on $V$. We define 
$$\int_{\gamma}  : A^{\log,1}_{k, X/Y, d=0}(U)\to A^{\log}_{k,Y}(V)$$
by the composition
$$f^{\log}_*A^{\log,1}_{k, X/Y, d=0}|_V\to R^1f_*^{\log}A _V= A \otimes_{\Z} R^1f^{\log}\Z|_V \to A|_V=A^{\log}_{k,Y}|_V,$$
in which the first arrow comes from the exact sequence \ref{d=0}.1.

 \end{sbpara}

\begin{sbpara}\label{H1.0} We discuss a log relative version of \ref{d=0} (2.2). 

Let $U$ be an open set of $X^{\log}_{[:]}$, let $V$ be an open set of $Y^{\log}_{[:]}$, and let $\gamma: V \times \R/\Z\to U$ be a continuous map over $Y^{\log}_{[:]}$. Let $w\in A^{\log,1}_{k, X/Y,d=0}(U)$. Then we define $$\int_{\gamma} w\in A^{\log}_{k,Y}(V)$$ as follows. 
By the exact sequence \ref{d=0}.1,  $w$ gives an element of $H^1(U, A)$ and the map  $\gamma$ gives $H^1(U, A)\to H^1(V \times \R/\Z, A) \to H^0(V, A)$. 
 
We have $\int_{\gamma} w=\int_{\gamma'} w$ if $\gamma$ and $\gamma'$ are homotope over $Y^{\log}_{[:]}$. 

In the case where $f:X\to Y$ is proper and $U$ is the inverse image of $V$ in $X^{\log}_{[:]}$, $\gamma$ induces a homomorphism $h: R^1f^{\log}_*\Z\to \Z$, and $\int_{\gamma} w$ coincides with $\int_h w$ in \ref{H1.1}. 
 
 \end{sbpara}

 \begin{sbpara}\label{H1.4}
We consider the relation of \ref{H1.5} to  \ref{H1.1} and \ref{H1.0}. Let $V$ be an open set of $Y_{[:]}^{\log}$. 
 Let  $U_j$ ($1\leq j \leq m$) be  open sets of $X^{\log}_{[:]}$ with the image $V$ in $Y_{[:]}^{\log}$. Assume that the map $U_j\to V$ is a homotopy equivalence over $V$ for every $j$.  
   Let $U=\bigcup_{j=1}^m U_j$ and let $w\in A^{\log,1}_{k, X/Y, d=0}(U)$.

  Assume that we have sections $s_j:Y \to X$ which induce  $V\to U_j\cap U_{j+1}$ for $0\leq j\leq m-1$, where $U_0$ denotes $U_m$. Let $s_m=s_0$. Let $\int_{U_j, s_{j-1}}^{s_j} w$ be the $\int_{s_{j-1}}^{s_j} $ of the restriction of $w$ to $U_j$ (\ref{H1.5}). 
  
  Let $A=A^{\log}_{k, Y}$. By the vanishing of $\cH^1(U_j,A)$ ($\cH^m$ denotes the  higher direct image to $V$), $\cH^1(U,A)$ is equal to the \v{C}ech cohomology $\check{\cH}^1((U_j)_j, A)$. Consider the homomorphism $h: H^1(U, A)\to A(V)$ given by the homomorphism $\check{\cH}^1((U_j)_j, A)\to A$
 induced by $\bigoplus_{j<k}  \cH^0(U_j\cap U_k,A)\to A $ whose $(j,k)$-component sends $g\in  H^0(U_j\cap U_k, A)$  to $s_j^*(g)\in A$ if $m\neq 2$, $1\leq j\leq m-1$  and $k=j+1$, to $-s_0^*(g)\in A$ if $m\neq 2$, $j=1$ and $k=m$, to $s_1^*(g)-s_0^*(g)$ if $m=2$, $j=1$ and $k=2$, and to $0$ otherwise.

  We consider two cases Case 1, Case 2, to treat \ref{H1.1} and \ref{H1.0}, respectively.

    In Case 1, assume that $f:X\to Y$ is proper and assume that $U$ is the inverse image of $V$ in $X^{\log}_{[:]}$, and let $\gamma: R^1f^{\log}_*\Z|_V\to \Z$  be the above homomorphism $h$. 
    
   In Case 2, assume that we are given a continuous map $\gamma: V \times \R/\Z\to U$ over $V$. 
  Assume that  for $1\leq j\leq m$, the map $\gamma(\cdot, j/m): V \to U$  coincides with the map given by the above $s_j$ and that the image of $V\times \{v/m\;|\;v\in \R,  j-1\leq v \leq j\}$ in $U$ under $\gamma$ is contained in $U_j$.
 \end{sbpara}

\begin{sbprop}\label{H1.2} Under the above assumptions in $\ref{H1.4}$, in both Case $1$ and Case $2$, we have $$\int_{\gamma}w= {\textstyle\sum}_{j=1}^m \int_{U_j, s_{j-1}}^{s_j} w.$$

\end{sbprop}

\begin{pf} 
 Let $F_j\in A^{\log}_{k, X}(U_j)$ be $\int_{s_{j-1}} w$.   The image of $w$ in  $\cH^1(U, A)$ coincides with the class of the  $1$-cocycle 
in $\bigoplus_{j<k}  \cH^0(U_j\cap U_k,A) $ whose $(j,k)$ component coincides with $F_j-F_k$. In Case 2, the homomorphism $\cH^1(U, A)\to A$ induced by $\gamma$ coincides with $h$  in \ref{H1.4}. In both Case 1 and Case 2, the homomorphism $h$ sends the $1$-cocycle $(F_j-F_k)_{1\leq j<k\leq m}$ to $\sum_{j=1}^{m-1} (s_j^*(F_j)-s_j^*(F_{j+1})) + s_0^*(F_m)- s_0^*(F_1)= \sum_{j=1}^m \int_{U_j, s_{j-1}}^{s_j} w$. 
\end{pf}

\subsection{Example: Degenerate elliptic curve}\label{ss:intex}
\label{s:newex}

\begin{sbpara}\label{badE1}  Let $S$ be an fs log analytic space and let $q$ be a section of $M_S$ such that $|q|<1$. Then we have an fs log analytic space $E_{S,q}$ over $S$  regarded as an elliptic curve with degeneration over $S$ with $q$-invariant $q$, which represents the following functor. For an fs log analytic space $T$ over $S$, $\text{Mor}_S(T, E_{S,q})$ is the set of all sections of the sheaf 
$$\{t\in M_T^{\gp}\; |\; \text{locally}\; q^n|t|q^{n+1}\;\text{for some $n\in \Z$}\}/q^{\Z}.$$

In the case where $q$ is a section of $\cO_S^\times$, $E_{S, q}$ represents ${\bf G}_m/q^{\Z}$ and it  is an elliptic curve over $S$ whose fiber on $s\in S$ is  the elliptic curve $\C^\times/q(s)^{\Z}\cong \C/(\Z+ \Z\tau)$ over $\C$, where $\tau\in \C$ is such that $q(s)=e^{2\pi i\tau}$.  

\end{sbpara}

\begin{sbpara} The morphism $f: E_{S, q}\to S$ is log smooth and saturated, and for the induced map $f^{\log}_{[:]}: E_{S, q, [:]}^{\log}\to S^{\log}_{[:]}$, $R^1f^{\log}_{[:]*}\Z$ is locally isomorphic to $\Z^2$. 

\end{sbpara}

\begin{sbpara}\label{badE2}
In this section (Section \ref{ss:intex}), 
we  consider the case where $S$ is an fs log point of log rank $1$ 
and $q$ is a generator of the log structure of $S$. Thus $q$ goes to $0$ in $\cO_S=\C$. We will denote $E_{S,q}$ simply by $E$. In this case, as a complex analytic space, $E$ is identified with the quotient ${\bf P}^1(\C)/(0\sim \infty)$ obtained by identifying $0$ and $\infty$ of ${\bf P}^1(\C)$. Let $X$ be the fs log analytic space over $S$ which represents the functor $T\mapsto \{(t_1, t_2)\; |\; t_1,t_2\in \Gamma(T, M_T),  t_1t_2=q\}$.
Then as a complex analytic space, $X$ is identified with $\{(t_1, t_2)\in \C^2\;|\; t_1t_2=0\}$. We have the morphism $X\to E$ over  $S$ corresponding to $\{(t_1, t_2)\;|\; t_1, t_2\in \Gamma(T, M_T), t_1t_2=q\}\to
\{t\in M_T^{\gp}\; |\; \text{locally}\; q^n|t|q^{n+1}\;\text{for some $n\in \Z$}\}/q^{\Z};\; (t_1, t_2)\mapsto t=t_1$.  
By this morphism, as a complex analytic space, $E$ is identified with the quotient of $X$ obtained by identifying $(c, 0)$ with $(0, c^{-1})$ for $c\in \C^\times$. 
\end{sbpara}

\begin{sbpara}\label{E[:]} Let $S$ and $E$ be as in \ref{badE2}. 

We describe $E_{[:]}$ and describe what Theorem \ref{SEisS} says for $E\to S$. 

We can identify $E_{[:]}$ with the union of ${\bf P}^1(\C)$ and the interval $[0, \infty]$ in which $0\in {\bf P}^1(\C)$ is identified with $0\in [0, \infty]$ and $\infty\in {\bf P}^1(\C)$ is 
identified with $\infty\in [0, \infty]$. In this identification, the inverse image in $E_{[:]}$ of the unique singular point of $E$ is identified with $[0,\infty]$ by the value $\log|t_1|/\log|t_2|$. That is, $E_{[:]}$ is homeomorphic to 
the union of the surface of the earth and a path whose one end point is identified with the south pole and the other end point is identified with the north pole. On the other hand, $(S\times \Delta)_{[:]}$ 
is
 identified with the union of $\Delta$ and the closed interval $[0, \infty]$ in which $0\in \Delta$ is identified with $0\in [0, \infty]$. In this identification, a point of $(S \times \Delta)_{[:]}$ lying over $0\in S\times \Delta= \Delta$ 
is identified with the value $\in [0,\infty]$ of $\log|q_*|/\log|q|$, where $q_*$ is the coordinate function of $\Delta$ (recall that $q$ is the generator of the log structure of $S$). 

We have an isomorphism from an open neighborhood of the class of $0$ in $E_{[:]}$ to an open neighborhood of $0$ in $(S \times \Delta)_{[:]}$ induced by $t_1$, and an isomorphism from an open neighborhood of the class of $\infty$ in $E_{[:]}$ to an open neighborhood of $0$ in $(S \times \Delta)_{[:]}$ induced by $t_2$. Thus $t_1$ (resp.\ $t_2$) is  a local coordinate function of $E$ over $S$ at the class of $0$ (resp.\ $\infty$) in $E_{[:]}$ in the sense of Theorem \ref{SEisS}. 

We describe the fiber of $E^{\log}_{[:]}\to S^{\log}_{[:]}$ ($S^{\log}_{[:]}$ is homeomorphic to a circle). It is obtained from $E_{[:]}$ by replacing the south pole by a circle $S^1$ and the north pole by another circle $S^1$, and the path between the two poles by $S^1 \times $ (the path).

\end{sbpara}

\begin{sbpara}\label{badE3}
Let $S$, $E$ and $X$ be as in \ref{badE2}. In Proposition  \ref{per} below, we compute the period integral on $E$ mentioned in \ref{Intro2}.

Here we define the loops $\alpha$ and $\beta$ along which we perform the period integrals. 

Write $$\log(q)=2\pi i(x_0+iy_0), \quad \log(t_1)=2\pi i(x_1+iy_1), \quad \log(t_2)=2\pi i(x_2+iy_2),$$
$x_1+x_2=x_0$, $y_1+y_2=y_0$ ($x_j$ and $y_j$ being real).

Let $V$ be a contractible open set of $S^{\log}$. Take  a continuous map $\tilde x_0: V\to \R$ such that the composite map $V \overset{\subset}\to S^{\log}\underset{x_0}{\overset{\cong}\to} \R/\Z$ factors through $\tilde x_0$.  

We have the following homeomorphism $$\nu: V\times (\R/\Z)^2 \overset{\cong}\to (f^{\log}_{[:]})^{-1}(V)\subset E^{\log}_{[:]}$$
 over $S^{\log}$.

If $0\leq b\leq 2/3$, let $\nu(p,a, b)$ be the image of the unique point of $X^{\log}_{[:]}$ lying over $p$ at which $x_1=a$, $y_1= b/(2/3-b)$, $y_2=\infty$. 

If $2/3\leq b\leq 5/6$, let $\nu(p,a, b)$ be the image of the unique point of $X^{\log}_{[:]}$ lying over $p$ at which $x_1=a+ 6\tilde x_0(p)(b-2/3)$, $y_1=\infty$, $y_2=\infty$, $y_1/y_2= (b-2/3)/(5/6-b)$. 

If $5/6\leq b\leq 1$, let $\nu(p,a, b)$ be the image of the unique point of $X^{\log}_{[:]}$ lying over $p$ at which $x_1=x_0(p)$,  $y_1=\infty$, $y_2= b/(5/6-b)$. 

Consider the continuous maps  $$\alpha, \beta: V \times \R/\Z\to E^{\log}_{[:]}, \quad \alpha(p,b)=\nu(p,0,b), \;\; \beta(p, a)= \nu(p,a,2/3).$$  Then 
the image of the loop $\alpha$ on ${\bf P}^1(\C)/(0\sim \infty)$ is  $[0,\infty]/(0\sim \infty)$ (with the orientation to go from $\infty$ to $1$ and then from $1$ to $0$). The loop $\beta$ is a vanishing cycle which crushes to the point $0=\infty$ of ${\bf P}^1(\C)/(0\sim \infty)$.
In the usual geometry, the two integrals in Proposition \ref{per} below are like $(2\pi i)^{-1}\int_0^{\infty} d\log(t)$ and $(2\pi i)^{-1}\int_0^0d\log(t)$, respectively,  which are not defined, but we can compute them correctly by our log geometry. 
\end{sbpara}

\begin{sbprop}\label{per} Denote by $d\log(t)$  the global log differential form on $E$ over $S$ whose pullback to $\C^\times$ is $d\log(t)$ and whose pullback to $X$ in $\ref{badE2}$ is $d\log(t_1)$. 
  Then we have $$(2\pi i)^{-1}\int_{\alpha} d\log(t)= \tilde x_0+ iy_0, \quad (2\pi i)^{-1}\int_{\beta} d\log(t)=1.$$
  
\end{sbprop}

\begin{pf} We write the proof of the first formula. The proof of the second formula is similar and easier. 

We use the following open sets $U_1$ and $U_2$ of $E^{\log}_{[:]}$ and sections $s_0, s_1: S\to E$ and use Proposition \ref{H1.2}. 

Let $U_1$ be the set of $\nu(p, a, b)$ such that $a\neq 1/2$ and $b$ is in the image of $(-1/6, 2/3)\to \R/\Z$. Let  $U_2$ be the set of $\nu(p,a,b)$ such that $a\neq 1/2$ and $b$ is in the image of $(1/3, 7/6)\to \R/\Z$. 
Then $U_1$ and $U_2$ are contractible, and the image of the map $\alpha$ is contained in $U_1\cup U_2$.

Let $s_0: S\to E^{\log}$ be the morphism over $S$ given by $1\in {\bf G}_m\subset {\bf G}^{\log}_m/q^{\Z}$, and let $s_1:S\to U_1\cap U_2$ be the morphism over $S$ given by $e^{-6\pi}\in {\bf G}_m\subset {\bf G}^{\log}_m/q^{\Z}$. Then $s_0$ sends $p\in V$ to $\alpha(p, 0)$ and $s_1$ sends $p\in V$ to $\alpha(p,1/2)$.

By Proposition \ref{H1.2}, 
$\int_{\alpha} d\log(t)=\int_{U_1, s_0}^{s_1}d\log(t) +\int_{U_2, s_1}^{s_2} d\log(t)$, where $s_2=s_0$. 

We have $(2\pi i)^{-1}d\log(t) = dx_1+idy_1$. 

On $U_1$, we have $(2\pi i)^{-1}d\log(t)= d(\tilde x_1+ iy_1)$, where $\tilde x_1$ is the branch of $x_1$ which has value $0$ at $s_0$.  We have $s_0^*(\tilde x_1)=0$,  $s_0^*(y_1)=0$, $s_1^*(\tilde x_1)=0$, $s_1^*(y_1)=3$.  Hence $(2\pi i)^{-1}\int_{U_1, s_0}^{s_1}d\log(t)= 3i$. 

On $U_2$, we have $(2\pi i)^{-1}d\log(t)= d(\tilde x_1+ iy_1)$, where $\tilde x_1$ is the branch of $x_1$ which has value $0$ at $s_1$. We have
$s_1^*(\tilde x_1)=0$,  $s_1^*(y_1)=3$, $s_2^*(\tilde x_1)=\tilde x_0$, $s_2^*(y_1)=y_0$.  Hence $(2\pi i)^{-1}\int_{U_2, s_1}^{s_2} d\log(t)=\tilde x_0+ iy_0-3i$. 

Hence $(2\pi i)^{-1}\int_{\alpha} d\log(t)= \tilde x_0+iy_0$. 
\end{pf}

\begin{sbpara}\label{another}    The first and the second formulas in Proposition \ref{per} are the degenerate versions of the period integrals
$$\int_0^{\tau} dz=\tau, \quad \int_0^1 dz=1$$
of the elliptic curve $\C/(\Z+\Z\tau)$, respectively. 

Proposition \ref{per}  can be proved also as follows, by using these classical period integrals.   

The space $E$ over $S$ is the fiber at $S=0\in \Delta$ in the  degenerating family of elliptic curves $\cE:=E_{\Delta, q}$ over $\Delta$.
Let $\tilde V$ be a contractible open set of $\Delta^{\log}$ such that $V\subset \tilde V$. Then the lifting $\tilde x_0: V\to \R$ of $x_0$ extends to a continuous lifting $\tilde x_0: \tilde V\to \R$ of $x_0: \tilde V\to \R/\Z$. We have a homomorphism $\tilde \alpha$ (resp.\  $\tilde \beta): R^1g^{\log}_*\Z|_{\tilde V} \to \Z$ on $\tilde V$ whose restriction to the fiber $\C^\times/q^{\Z}\cong \C/(\Z+ \Z\tau)$ at $q\in \tilde V\cap(\Delta\smallsetminus \{0\})$ is induced by the loop which is the image of the path from $0$ to $\tau$ (resp.\  $0$ to $1$) on $\C$. Here in the definition of  $\tilde \alpha$, $\tau$ is the branch of $(2\pi i)^{-1}\log(q)$ whose real part is $\tilde x_0(q)$.
Then the restriction of $\tilde \alpha$ (resp.\ $\tilde \beta$) to $V$ coincides with the homomorphism induced by $\alpha$ (resp.\ $\beta$). 

Let $d\log(t)$ be the global log differential form on $\cE$ over $\Delta$ whose pullback to $\C^\times/q^{\Z}$ is $d\log(t)$, where $t$ is the coordinate function of $\C^\times$. The pullback of $(2\pi i)^{-1}d\log(t)$ to $\C/(\Z+\Z\tau)$ is $dz$, where $z$ is the coordinate function of $\C$. Since we have the classical period integrals

\medskip

(1) $(2\pi i)^{-1}\int_{\tilde \alpha} d\log(t)= \tau=\tilde x_0+iy_0$ on $\tilde V\cap (\Delta\smallsetminus \{0\})$, and 

(2)  $(2\pi i)^{-1}\int_{\tilde \beta} d\log(t)= 1$ on $\tilde V\cap (\Delta\smallsetminus \{0\})$, 

\noindent
we have (1) (resp.\ (2)) on $\tilde V$, and hence have (1) (resp.\  (2)) on $V$, that is, the first (resp.\ second) formula in Proposition \ref{per}.

A good thing of our theory of integration is that the integration can be done not using the family $\cE$ over $\Delta$, but just by using the log structure of $E={\bf P}^1(\C)/(0\sim \infty)$. Thus the log structure of $E$ can play the role of the family  $\cE$ and has the power to give the period integral in degeneration (which can be given by $\cE$) without the help of $\cE$.

\end{sbpara}

\begin{sbpara}

Integrations in Sections \ref{ss:inty} and \ref{ss:intx} are not used here. They will be related to the above integral
  in another paper \cite{KNUi}.

\end{sbpara}

\section{Log real analytic functions and $\SL(2)$-orbits}\label{s:newa}

In degeneration of Hodge structure, 

\medskip

(1) nilpotent orbits appear in the holomorphic world (\cite{Sc} Section 4) , and 

\medskip

(2) $\SL(2)$-orbits appear in the real analytic world (\cite{Sc} Section 5,  \cite{CKS1}, \cite{KNU08}).

\medskip
 In \cite{KU} and \cite{KNU3}, we understand (1) geometrically as follows. In degeneration, a variation of Hodge structure extends to the boundary as a log Hodge structure, and the nilpotent orbit which appears  is understood as the fiber of the log Hodge structure at a point of the boundary. In Section \ref{ss:HA}, we understand (2) geometrically as follows. The pullback of a log Hodge structure on $X$ to the space of ratios $X_{[:]}$  is treated nicely by log real analytic functions, and the SL(2)-orbit which appears is the log real analytic fiber of
  the log Hodge structure at a point of $X_{[:]}^{\log}$. 

The Hodge decomposition of a variation of  Hodge structure $H$ and the Hodge metric of a variation of polarized Hodge structure $H$ are  real analytic. But when a degeneration happens, these become not real analytic and diverge. In Section \ref{ss:HA}, we show that we have a nice $A^{\an}_{(2+)}$-structure ${}^{\prime}H^{\natural}_A$ of $H$ on the space $X_{[:]}$ for which the coefficients of  the Hodge decomposition and the Hodge metric  belong to $A_{(2+)}^{\an}$ and hence converge. 
In  Section \ref{ss:HA}, we also show that the associated $\SL(2)$-orbit in degeneration can be nicely understood by using $A^{\an}_{(2+)}$ on $X_{[:]}$.
The proofs of the theorems in Section \ref{ss:HA}  are given in Section \ref{ss:pfSL}. 
They 
are proved by using the known SL(2)-orbit theorems, but they will be  more useful than the original SL(2)-orbit theorems in our forthcoming studies of higher direct images.

In Section \ref{ss:CKS}, as an application to extended period domains, we prove the log real analyticity of the CKS map $D^{\sharp}_{\Sig,[:]}\to D_{\SL(2)}$, an important map in the theory of extended period domains which connects the world of nilpotent orbits and the world of $\SL(2)$-orbits. This gives a new proof of the continuity of this CKS map which we proved in \cite{KU},\cite{KNU3}, and \cite{KNU4}.

\subsection{Log Hodge theory}\label{ss:LMH}

We give a short introduction to  log Hodge theory.

\begin{sbpara} 
There are two log versions for the notion mixed Hodge structure.

\medskip

(1) Log mixed Hodge structure (LMH, in short),

(1)$^*$ Logarithmic variation of mixed Hodge structure (LVMH, in short).

\medskip

There are two log versions for the notion polarized Hodge structure.

\medskip

(2) Polarized log Hodge structure (PLH, in short),

(2)$^*$ Logarithmic variation of polarized Hodge structure (LVPH, in short).

\medskip

We have

\medskip

(1)$^*$ (resp.\ (2)$^*$)=   ((1) (resp.\ (2)) satisfying the big Griffiths transversality).

\medskip
\noindent
These are explained in \cite{KU} 2.4, 2.6, \cite{KNU3} 1.3.1, and also explained  in  \ref{LMH1}--\ref{LMH3} below. 

(1) and (2) are useful for the study of moduli (see \ref{PLH1}). We think that (1)$^*$ and (2)$^*$ are useful for the study of higher direct images (see Section \ref{ss:hdi}). 

\end{sbpara}
\begin{sbpara}\label{LMH1} Let $X$ be an object of $\cB(\log)$. Let $R$ be a subring of $\R$ (important cases are $R=\Z, \Q, \R$).

A {\it pre-$R$-LMH}  on $X$ is a triple $$H=(H_R, H_{\cO}, \iota),$$ where $H_R$ is a locally constant sheaf on $X^{\log}$ of finitely generated $R$-modules endowed with a locally constant increasing filtration $W$ on $H_R \otimes_{\Z}\Q$ such that $W_w=H_R \otimes_\Z \Q$ if $w\gg 0$ and $W_w=0$ if $w\ll 0$, $H_{\cO}$ is a sheaf of  $\cO_X$-modules on $X$ which is locally free of finite rank endowed with a decreasing filtration $F$ by $\cO_X$-submodules which splits locally such that $F^p=0$ if $p\gg 0$ and $F^p=H_\cO$ if $p\ll 0$, and $\iota$ is an isomorphism of $\cO_X^{\log}$-modules 

\begin{equation}\tag{$*$} 
\cO_X^{\log} \otimes_R H_R \cong \cO_X^{\log} \otimes_{\cO_X} H_{\cO}\quad\text{ on }X^{\log}.
\end{equation}

This $(*)$ implies that we have $H_{\cO}= (\tau_X)_*(\cO^{\log}_X\otimes_R H_R)$ (\cite{KU} Proposition 2.2.10). 

Let $w\in \Z$. A {\it pre-$R$-PLH of weight $w$}  is  a pre-$R$-LMH of pure weight $w$ (this means $W_w=H_R \otimes_{\Z} \Q$ and $W_{w-1}=0$) endowed with an $R$-bilinear form  $\langle\cdot,\cdot\rangle: H_R\times H_R\to R\otimes_{\Z} \Q$ satisfying the following conditions (i)--(iii).

(i) $\langle\cdot,\cdot\rangle$  is symmetric if $w$ is even and anti-symmetric if $w$ is odd.

(ii) $\langle\cdot, \cdot\rangle$  induces an isomorphism between  $H_R \otimes_\Z \Q$ and the dual $R\otimes_\Z \Q$-module of it. 

(iii) For the induced pairing $H_{\cO} \times H_{\cO}\to \cO_X$ of $\cO_X$-modules,  we have $\langle F^p, F^q\rangle=0$ if $p+q>w$.

\end{sbpara}

\begin{sbpara}\label{Gri} We consider the Griffiths transversality of a pre-$R$-LMH $H$. 

 The map $d: \cO_X^{\log}\to \omega^{1,\log}_X$ (here $\omega_X^{1,\log}= \cO_X^{\log} \otimes_{\cO_X} \omega^1_X$) induces  $d: \cO_X^{\log}\otimes_R H_R \to \omega^{1,\log}_X\otimes_R H_R\;;\;f\otimes h\mapsto df\otimes h$ ($f\in \cO_X^{\log}$, $h\in H_R$). By taking $\tau_*$ of this, we obtain a  connection $\nabla: H_{\cO}\to \omega^1_X \otimes_{\cO_X} H_\cO$. We say that $H$ satisfies the {\it big Griffiths transversality} (or simply the {\it Griffiths transversality}) if $\nabla(F^p)\subset \omega^1_X \otimes_{\cO_X} F^{p-1}$ for all $p$. 

The reason why we say big here will be explained in Remark \ref{LMHrem} (1)  below. 

\end{sbpara}

\begin{sbpara}\label{LMH2} Assume that $X=\{x\}$ is an fs log point. Let $\cT\supset X$ be the canonical toric variety associated to $X$ (Section \ref{ss:pt}). 

Let $H=(H_R, H_{\cO}, \iota)$ be a pre-$R$-LMH on $X$. Then $H$ extends canonically 
to a pre-$R$-LMH $H_{\cT}$ on $\cT$ as follows. Since $X^{\log}\to \cT^{\log}$ is a homotopy equivalence, 
$H_R$ with the weight filtration $W$ on $X^{\log}$ extends uniquely to a locally constant $R$-module $H_{\cT, R}$ 
with the weight filtration $W$ on $\cT^{\log}$. Also, $H_{\cO}$ on $X$ extends to the $\cO_{\cT}$-module $H_{\cT, \cO}= \cO_{\cT} \otimes_{\C} H_{\cO}$  
with the filtration $F_{\cT}:=(\cO_{\cT}\otimes_{\C} F^p)_p$ on $\cT$, and the isomorphism $\iota: \cO_X^{\log}\otimes_R H_R \cong \cO_X^{\log}\otimes_{\C} H_{\cO}$  on $X^{\log}$
extends to an isomorphism 
$\iota: \cO_{\cT}^{\log}\otimes_R H_{\cT, R} \cong \cO_{\cT}^{\log} \otimes_{\cO_{\cT}} H_{\cT, \cO}$ via $\cO_{X,\cT}^{\log}\overset{\subset}\to  \cO_{\cT}^{\log}$ (\ref{XT1}) on $\cT^{\log}$. 

For $s\in \cT_{\triv}$, let $H_{\cT}(s)$ be  the stalk $H_{\cT, R, s}$ with the filtration $W(s)$ on $H_{\cT,R,s}\otimes_\Z \Q$ induced by $W$ and with the filtration $F(s)$ on $H_{\cT,R,s}\otimes_R \C$ induced by the filtration $F_{\cT}$ via the evaluation map $\cO_{\cT,s}^{\log}=\cO_{\cT,s}\to \C$ at $s$.

We call $H$ an {\it $R$-LMH} on $X$ if it satisfies the following conditions (i)--(iii).

\medskip

(i) If $s\in \cT_{\triv}$ is sufficiently near to $x$ in $\cT$, $H_{\cT}(s)$ ($=H_{\cT, R,s}$ with $W(s)$ and $F(s)$) 
is a mixed Hodge structure in the usual sense. 

(ii) The Griffiths transversality \ref{Gri} is satisfied.

(iii) The local monodromy of $H_\R=H_R\otimes_R \R$ is admissible. This is explained below in \ref{adm}. 

\medskip

Let $H$ with $\langle\cdot, \cdot\rangle$ be a pre-$R$-PLH on $X$ of weight $w$. Then $H$ extends canonically to a pre-$R$-PLH $H_{\cT}$ of weight $w$ on $\cT$. The pairing $\langle\cdot,\cdot\rangle$ on $H_{\cT, R}$ is the unique pairing which extends $\langle\cdot,\cdot\rangle$ on $H_R$. 

  We call $H$ an {\it $R$-PLH of weight $w$} on $X$ if it is an $R$-LMH and of pure weight $w$ satisfying the following condition (iv).

\medskip

(iv)  If $s\in \cT_{\triv}$ is sufficiently near to $x$ in $\cT$, $H_{\cT}(s)$ (with the induced $\langle\cdot,\cdot\rangle$)  is a polarized Hodge structure of weight $w$ in the usual sense. 

\end{sbpara}

\begin{sbpara}\label{LMH3} Let $X$ be an object of $\cB(\log)$.

Let $H$ be a pre-$R$-LMH on $X$. We call $H$ an {\it $R$-LMH} if its pullbacks to the fs log point $x$ for all $x\in X$ are $R$-LMH in the sense of \ref{LMH2}.

Let $H$ be a pre-$R$-PLH of weight $w$ on $X$. We call $H$ an {\it $R$-PLH} if its pullbacks to the fs log point $x$ for all $x\in X$ are $R$-PLH in the sense of \ref{LMH2}.  

An  $R$-LMH (resp.\ $R$-PLH) $H$ satisfying the  big Griffiths transversality (\ref{Gri})  is called an {\it $R$-LVMH} (resp.\ {\it $R$-LVPH}). 
\end{sbpara}

\begin{sbrem}\label{LMHrem}
(1) By the definition of $R$-LMH, an $R$-LMH (and hence an $R$-PLH) $H$ satisfies the following weaker version of Griffiths transversality which we call the {\it small Griffiths transversality}:  For every $x\in X$, the pullback of $H$ to the fs log point $x$ satisfies the Griffiths transversality.  Note that the log differential module $\omega^1_x$ of an fs log point $x$ need not be trivial. If the log structure of $X$ is trivial, $\omega^1_x=0$ for every $x\in X$ and hence the small Griffiths transversality is the empty condition.

(2) In the case $X=\Delta^n$ with the log structure given by the coordinate functions $q_j$ ($1\leq j\leq n$) and $x$ is the origin $(0,\dots, 0)$ of $X$, a pre-$R$-PLH on $X$ is written classically as $\exp(\sum_{j=1}^n z_jN_j)F(q)$ ($q=(q_j)_j$, $q_j=\exp(2\pi iz_j)$, $N_j$ is the $j$-th local monodromy logarithm at $x$, $F$ is the Hodge filtration), the canonical toric variety $\cT$  associated to the fs log point $x$  is identified with $\C^n$, and then $H_{\cT}$ can be written as $\exp(\sum_{j=1}^n z_jN_j)F(0)$ in $\C^n$, where $F(0)$ denotes the fiber of $F$ at $x$.

(3) In \cite{KU} 2.4, 2.6 and \cite{KNU3}, the conditions (i) and (iv)  in \ref{LMH2} are written by using the set $\spe(t)$ in \ref{pta}, not by using the canonical toric variety $\cT$. We can write these conditions in both ways by the  relation of $\spe(t)$ and $\cT$  explained in \ref{pta}.

(4) If $H$ is an $R$-PLH of weight $w$, then for every $p\in \Z$, the annihilator of $F^p$ for $\langle\cdot,\cdot\rangle$ is $F^{w+1-p}$.

(5)  For a pre-$R$-PLH on an fs log point, the admissibility of  monodromy (the condition (iii) of \ref{LMH2}) is a consequence of the other conditions (i), (ii) and (iv) of \ref{LMH2}  (\cite{CK},  \cite{KU} 2.6.6).
\end{sbrem}

\begin{sbpara}\label{PLH1}
We have nice moduli spaces of $\Z$-PLH as in \cite{KU}, and we have nice moduli spaces of $\Z$-LMH with polarized $\gr^W$ as in \cite{KNU3}. These are partial toroidal compactifications of classical period domains. 

If $X$ is an fs log point, we have $\{$R-LMH on $X\}=\{$R-LVMH on $X\}=\{$the classical notion ``nilpotent orbit''$\}$ (\cite{KU} Section 2.5, \cite{KNU3} 2.2.2). 
  The above moduli spaces are regarded as spaces of nilpotent orbits.
\end{sbpara}

\subsection{Local monodromy and relative monodromy filtration}\label{ss:mono}

\begin{sbpara} Let $X$ be an object of $\cB(\log)$. For $x\in X$, we have a canonical isomorphism $\pi_1(x^{\log})\cong\Hom((M^{\gp}_X/\cO_X^\times)_x, \Z)$ (\cite{KU} 2.2.9).

For a pre-$R$-LMH $H$ on $X$, for $t\in X^{\log}$ lying over $x$, we have the action of $\pi_1(x^{\log})$ on the stalk $H_{\R, t}$, which we call the {\it local monodromy  of $H$ at $t$}.  This action is unipotent
 by \cite{KU} Proposition 2.3.3 (ii).

\end{sbpara}

\begin{sbpara}\label{adm} We say that the local monodromy of $H$ is {\it admissible} if the following conditions (i) and (ii) are satisfied for each $x\in X$ and each $t\in X^{\log}$ lying over $x$. 

For $h\in \Hom((M^{\gp}_X/\cO_X^\times)_x, \R)$, let $\gamma_h$ be the corresponding element of $\pi_1(x^{\log})\otimes_\Z \R$, and let $N_h=\log(\gamma_h)$ acting on $H_{\R, t}$. Let $\sig(x)$ be the cone $\Hom((M_X/\cO_X^\times)_x, \R^{\add}_{\geq 0})\subset  \Hom((M^{\gp}_X/\cO_X^\times)_x, \R)$.

(i) For each $h\in \sig(x)$, the relative monodromy filtration of $N_h$ with respect to $W$ exists. Furthermore, this relative monodromy filtration depends only on the smallest face of $\sig(x)$ which contains $h$.

(ii) For each face $\tau$ of $\sig(x)$, let $W(\tau)$ be the relative monodromy filtration of $N_h$ of an element $h$ of the interior of $\tau$ (that is, an element $h$ of $\tau$ such that $\tau$ is the smallest face of $\sig(x)$ which contains $h$), which is independent of such an $h$ by (i). Then for $h\in \sig(x)$ and for faces $\tau$, $\tau'$ of $\sig(x)$ such that $\tau'$ is the smallest face of $\sig(x)$ which contains $\tau$ and $h$, $W(\tau')$ is the relative monodromy filtration of $N_h$ with respect to $W(\tau)$. 

\end{sbpara}

\begin{sbpara}\label{admrem} Assume that (i) is satisfied at $t\in X^{\log}$ lying over $x\in X$. 

(1) For every $h\in \sig(x)$ and every face $\tau$ of $\sig(x)$, we have $N_hW(\tau)_w\subset W(\tau)_w$ for every $w\in \Z$. This follows from the fact  $N_hN_{h'}= N_{h'}N_h$  for an element $h'$ of the interior of $\tau$. 

(2) For every face $\tau$ of $\sig(x)$ and for every $h\in \tau$, we have $N_hW(\tau)_w\subset W(\tau)_{w-2}$ for every $w\in \Z$. This is because $h=h_1-h_2$ for some elements $h_1, h_2$ of the interior of $\tau$. 

(3)  By (1),   for each face $\tau$ of $\sig(x)$, $W(\tau)$ on $H_{\R, t}$ is stable under the action of $\pi_1(x^{\log})$. Hence $W(\tau)$ is a locally constant filtration on the restriction of $H_\R$ to $x^{\log}$. 

\end{sbpara} 

\begin{sbpara}\label{nearx} Let $x\in X$. The local monodromy at points near to $x$  are related  to that of $x$ as follows.

Take a chart $\cS\to M_U$ of the log structure on an open neighborhood $U$ of $x$ with $\cS$ being an fs monoid such that $\cS\overset{\cong}\to (M_X/\cO_X^\times)_x$. 
Let $\cT=\Spec(\C[\cS])^{\an}$ with the canonical log structure. We have a morphism $U\to \cT$ and we can regard $x\in \cT$ and $x^{\log}\subset \cT^{\log}$.  Since $\pi_1(x^{\log}) \overset{\cong}\to \Gamma:=\pi_1(\cT^{\log})$, there is a unique local system $\tilde H_\R$ on $\cT^{\log}$ whose pullback to $x^{\log}$ is the restriction of $H_\R$ to $x^{\log}$. By the properness of $X^{\log}\to X$, by shrinking $U$ if necessary, the pullback of $\tilde H_\R$ on $U$ and the restriction of $H_\R$ on $U$ are identified. For $t' \in U^{\log}\subset X^{\log}$, the stalk $H_{\R, t'}$ is identified with the stalk $\tilde H_\R$ at the image of $t'$, and hence we have an isomorphism $H_{\R,t'}=\tilde H_{\R,t'}\cong\tilde H_{\R,t}=H_{\R,t}$ which is canonical modulo the action of $\pi_1(x^{\log})$. 

For $x'\in U$, via $\cS\to (M_X/\cO_X^\times)_{x'}$, we have an isomorphism from $\pi_1((x')^{\log})$ to a subgroup of $\Hom(\cS^{\gp}, \Z)\cong \pi_1(x^{\log})$, and we have an isomorphism from $\sig(x')$ to a face of $\Hom(\cS, \R^{\add}_{\geq 0})\cong \sig(x)$. For $t'$ lying over $x'$, via  this map $\pi_1((x')^{\log}) \to \pi_1(x^{\log})$ 
 and the  isomorphism between $H_{\R,t'}$ and $H_{\R,t}$,   the action of $\pi_1((x')^{\log})$ on $H_{\R,t'}$ and the action of $\pi_1(x^{\log})$ on $H_{\R,t}$ are compatible. 

\end{sbpara}

By \ref{nearx}, we have

\begin{sbprop} If the local monodromy of the restriction of $H$ to $x$ is admissible, the local monodromy of $H$ is admissible on some open neighborhood of $x$.

\end{sbprop}

\begin{sbpara} If the local monodromy of $H$ is admissible, in the situation of \ref{nearx}, the relative monodromy filtrations $W(\tau)$ on $H_{\R,t}$  for faces $\tau$ of $\sig(x)$ extend to  local systems on  $\cT^{\log}$, and hence to local systems on $U^{\log}$. For each $x' \in U$ and $t'$ lying over $x'$, the relative monodromy filtrations $W(\tau')$  on $H_{\R,t'}$ for faces $\tau'$ of $\sig(x')$  form a subset of the set of the stalks of these filtrations $W(\tau)$ at $t'$. 
\end{sbpara}

\begin{sbpara}\label{nearp}

Assume that the local monodromy of $H$ is admissible.
 Let $p\in X_{[:]}$ with the image $x$ in $X$. Then the set $\Phi(p)$ of relative monodromy filtrations at $p$ is given as follows (cf.\ \cite{KNU5} 3.3.2). 
  Let $M_{X,x}=M^{(0)}\supsetneq M^{(1)}\supsetneq \dots \supsetneq M^{(m)}=\cO^\times_{X,x}$ be as in \ref{m,f_jk} at $p$, and let $0=\sig_0\subsetneq \sig_1\subsetneq \dots \subsetneq \sig_m=\sig(x)$ be the dual sequence of cones, where $\sig_j $ is the set of all homomorphisms $M_{X,x} \to \R^{\add}_{\geq 0}$ which kills $M^{(j)}$. These $\sig_j$ are faces of $\sig(x)$. Write $W(\sig_j)$ by $W^{(j)}$.  Define $$\Phi(p)=\{W^{(j)}\;|\; 1\leq j \leq m\}.$$ 

Consider the situation in \ref{nearx}.
If $p' \in U_{[:]}$ with the image $x'$ in $U$, we have similarly faces $0=\sig'_0\subsetneq \sig'_1 \subsetneq \dots \subsetneq \sig'_{m'}$ of $\sig(x')$ associated to $p'$. Since $\sig(x')$ is a face of $\sig(x)$, these $\sig'_j$ are faces of $\sig(x)$. If $\Psi$ is a subset of the set of all faces of $\sig(x)$ and if $\Psi$ is  totally ordered for the inclusion, we have the open subset $U_{[:]}(\Psi)$ of $U_{[:]}$ consisting of all $p'\in U_{[:]}$ such that the set of faces of $\sig(x)$ associated to $p'$ is contained in $\Psi$. We have an open covering $U_{[:]}=\bigcup_{\Psi} U_{[:]}(\Psi)$. If $\Psi$ is the set of faces of $\sig(x)$ associated to $p$, then for any point $p'$ of the open neighborhood $U_{[:]}(\Psi)$ of $p$, we can regard $\Phi(p')$ as a subset of $\Phi(p)$.

\end{sbpara}

\begin{sbpara}\label{xi0} Consider the situation in \ref{nearx}.
Since $\pi_1(\cT^{\log})$ is commutative, for each $\gamma\in \pi_1(\cT^{\log})$, we have an automorphism of $\tilde H_\R$ (denoted by $\gamma$) which induces the action of $\gamma$ on each stalk of $\tilde H$.

Let $(q_{\la})_{\la\in \La}$ be a base of the free abelian group  $\cS^{\gp}$ of finite rank, and let $(\gamma_{\la})_{\la}$ be the dual base of 
$\pi_1(\cT^{\log})=\Hom(\cS^{\gp}, \Z)$. Let $N_{\la}:\tilde H_\R\to \tilde H_\R$ be the logarithm  of the above action of $\gamma_{\la}$ on $\tilde H_\R$. Let 
$$R= A^{\an}_k \quad (k=1,2,2+,3,(2), (2+), (3)),\quad\text{or}\quad R=A_k\quad(k=1,2,2*,3).$$
The $\cO_X^{\log}$-linear automorphism   $\xi:=\exp(\sum_{\la} z_{\la}N_{\la})$ of $\cO_{\cT}^{\log}\otimes_\R \tilde H_\R$ and 
the  $R_{\cT}^{\log*}$-linear automorphism  $\xi^{\R}:= \exp(\sum_{\la\in \La} x_{\la}N_{\la})$ of $R^{\log*}_{\cT}\otimes_\R \tilde H_\R$, where $z_{\la}=(2\pi i)^{-1}\log(q_{\la})$ and $z_{\la}=x_{\la}+iy_{\la}$ with $x_{\la}$, $y_{\la}$ being real, are defined locally modulo the linear action of  $\pi_1(S^{\log})$ and is independent of the choice of $(q_{\la})_{\la}$. The subsheaf $\xi \tilde H_\R\subset \cO^{\log}_{\cT}\otimes_\R \tilde H_\R$ and $\xi^{\R} \tilde H_\R\subset R^{\log*}_{\cT}\otimes_\R \tilde H_\R$ are well-defined globally (without modulo anything). Furthermore, we have:

\end{sbpara}

\begin{sblem}\label{cnst}  These sheaves $\xi \tilde H_\R$ and $\xi^{\R} \tilde H_\R$  are constant sheaves. 

\end{sblem}
\begin{pf} It is sufficient to prove that for  $\mu \in \La$, the monodromy action of $\gamma_{\mu}$ on each stalk of $\xi \tilde H_\R$ and that for  $\xi^{\R}\tilde H_\R$ are trivial. Since  $\gamma_{\mu}(z_{\la})= (\gamma_{\mu}^*)^{-1}(z_{\la})$ is equal to $z_{\la}-1$ if $\la=\mu$ and is $z_{\la}$ if $\la\neq \mu$,  
for an element $v$ of the stalk of $H_\R$, we have $\gamma(\xi v)=(\xi \gamma_{\mu}^{-1})(\gamma_{\mu}v)=\xi v$. 
\end{pf}

\begin{sbpara}\label{xi1} Let the situation be as in \ref{nearx}, \ref{nearp}, \ref{xi0}, and Lemma \ref{cnst}. By pulling back 
$\xi\tilde H_{\bR}$ (resp.\ $\xi^{\bR}\tilde H_{\bR}$) 
 from $\cT$ to $U\subset X$ (resp.\ $U_{[:]}(\Psi)\subset X_{[:]}$, where $\Psi=\Phi(p)$), we have a constant sheaf $\xi H_\R\subset \cO_X^{\log}\otimes H_\R$ on $U^{\log}$ and hence on $U$ which we denote by the same notation $\xi H_\R$ (resp.\ $\xi^\R H_\R\subset R^{\log*}_X\otimes_\R H_\R$ on the inverse image $U_{[:]}(\Psi)^{\log}$ of $U_{[:]}(\Psi)$ and hence on $U_{[:]}(\Psi)$ which we denote by the same notation $\xi^{\R} H_\R$).

We have
$$H_{\cO} = \cO_X \otimes_\R \xi H_\R\quad \text{on $U$}.$$
$$(\text{resp.\ }\tau_*(R^{\log*}_X\otimes_\R H_\R)= R_X\otimes_\R \xi^{\R} H_\R \quad \text{on $U_{[:]}(\Psi)$}.)$$

\end{sbpara}

\subsection{Log real analytic structures $H_A$ and $H^{\natural}_A$}\label{ss:twoHA}

Let $X$ be an object of $\cB(\log)$ and let $H$ be an $\R$-LMH  
on $X$ such that the graded quotients $\gr^W_wH$ ($w\in \Z$) of the weight filtration $W$ of $H$ are polarizable. We consider two log real analytic structures  
$$H_A\quad\text{and}\quad  H_A^{\natural}$$
 of $H$. These are sheaves of $A^{\an}_{(2+)}$-modules on $X_{[:]}$ and 
 are locally free of finite rank. 
In the case where $H$ is pure, we also consider 
$$\text{a modified version ${}^{\prime}H^{\natural}_A$ of $H^{\natural}_A$}. $$

We have chosen the notation $\natural$ in $H^{\natural}_A$ and ${}^{\prime}H^{\natural}_A$ because as we will see  in Theorem \ref{HA1} (3),  ${}^{\prime}H_A^{\natural}$ plays the role to make
natural the pitch of the Hodge metric which is raised or lowered  by degeneration.  

\begin{sbpara}\label{H0} We first define $H_A$.

For $R=A^{\an}_k$ ($k=1,2, 2+, 3, (2), (2+), (3)$) or  $R=A_k$ ($k=1,2,2*,3$), define
$$H_R := \tau_{X*}(R^{\log*}\otimes_{\R} H_{\R}),$$
where $\tau_X:X^{\log}_{[:]}\to X_{[:]}$.

This is an analogue of $H_{\cO}= \tau_{X*}(\cO_X^{\log}\otimes_{\C} H_\C)$. 

We denote $H_R$ for  $R=A^{\an}_{(2+)}$ by $H_A$.

\end{sbpara}

Next we define $H^{\natural}_A$. The following Proposition \ref{HA} is a preparation for it. 

\begin{sbprop}\label{HA} Let $R$ be as in $\ref{H0}$. In the case $R=A^{\an}_k$ with $k=1,2,2+,3$ (resp.\ $A^{\an}_k$ with $k=(2), (2+), (3)$, resp.\ $A_k$ with $k=1,2,2*,3$), let $R'=A^{\an}_3$ (resp.\ $A^{\an}_{(3)}$, resp.\ $A_3$).  
Then 
there is a unique  $R^{\log}$-submodule $H_R^{\natural\log}$  of $(R')^{\log} \otimes_\R  H_\R$ such that for every $p\in X_{[:]}$ and every $\tilde p\in X_{[:]}^{\log}$ lying over $p$, we have 
  $$H^{\natural\log}_{R,\tilde p}={\textstyle\sum}_{\ell\in \Z^m} \; R^{\log}_{\tilde p}{\textstyle\prod}_{j=1}^m (y_j/y_{j+1})^{-\ell(j)/2} {\textstyle\bigcap}_{j=1}^m W_{\ell(j)}^{(j)}H_{\R,\tilde p}\subset (R')_{\tilde p}^{\log} \otimes_\R  H_{\R,\tilde p}$$ 
 with $y_j=-(2\pi)^{-1}\log |f_j|$ ($y_{m+1}$ denotes $1$)  for any choices of $f_j\in M_{X,p}$ ($1\leq j\leq m$) such that $f_j\in M^{(j-1)}$ and $f_j\notin M^{(j)}$. 
Here $(W^{(j)})_{1\leq j\leq m}$ is as in $\ref{nearp}$. 
\end{sbprop}

 \begin{pf} Fix $p\in X_{[:]}$ and let $x$ be the image of $p$ in $X$. Let the notation be as in \ref{nearx} and \ref{nearp}.
Take $f_j$ in $\cS$ through $\cS \overset{\cong}\to (M_X/\cO_X^\times)_x$. Let $\Psi$ be the set of faces of $\sig(x)$ associated to $p$, and let $V$ be the open neighborhood of $p$ in $X_{[:]}$ consisting of all points of $U_{[:]}(\Psi)$ at which we have $|f_j|<1$ for all $j$. Let $\tilde V$ be the open set of $\cT_{[:]}$ consisting of all points at which $|f_j|<1$ for all $j$. We have a morphism $V\to \tilde V$. On $\tilde V$, we define 
$$\tilde H_R^{\natural\log}:= {\textstyle\sum}_{\ell\in \Z^m}\; R^{\log}_{\cT} {\textstyle\prod}_{j=1}^m (y_j/y_{j+1})^{-\ell(j)/2} {\textstyle\bigcap}_{j=1}^m W^{(j)}_{\ell(j)} \tilde H_\R,$$  where $\tilde H_\R$ is as in \ref{nearx}. We show that the pullback of $\tilde H_R^{\natural\log}$ to  $V$ has the desired property of $H_R^{\natural\log}$ at every point $p'$ of $V$.  We have $\Phi(p')=\{W^{(a(j))}\;|\; 1\leq j\leq m'\}$ with $1\leq a(1)< \dots<a(m')\leq m$. At $p'$,  instead of $f_j$ ($1\leq j\leq m$), we can use $f_{a(j)}$ ($1\leq j \leq m'$). It is sufficient to prove
that the stalk at $p'$ of the pullback of $\tilde H_R^{\natural\log}$ is the stalk of $R^{\log}_{X} \prod_{j=1}^{m'} (y_{a(j)}/y_{a(j+1)})^{-\ell(a(j))/2} \bigcap_{j=1}^{m'} W^{(a(j))}_{\ell(a(j))}H_\R.$ Here $y_{a(m'+1)}$  means $1$.  But this follows from the fact that 
$y_j/y_{j+1}$ is invertible in $R_{X, p}$ if $j$ is not in the image of $a$. 
\end{pf}

\begin{sbpara}
 Let the notation be as in Proposition \ref{HA}. Define
$$H^{\natural}_R=\tau_{X*}(H_R^{\natural\log}),$$
where $\tau_X: X^{\log}_{[:]}\to X_{[:]}$.

In the case where $H$ is pure of weight $w$, we have a slightly modified version ${}^{\prime}H_R^{\natural}\subset H_{R'}^{\natural}$ 
 of $H_R^{\natural}$ characterized by
$${}^{\prime}H^{\natural}_{R,p}=y_1^{w/2} H^{\natural}_{R,p}\quad  (p\in X_{[:]}),$$ which is useful.

 We denote $H^{\natural\log}_R$,  $H^{\natural}_R$, and ${}^{\prime}H^{\natural}_R$ for $R=A^{\an}_{(2+)}$  by $H_A^{\natural\log}$, $H^{\natural}_A$ and ${}^{\prime}H^{\natural}_A$, respectively. 

\end{sbpara}
\begin{sbpara}\label{HAfiber} Fibers. 
Let $p\in X_{[:]}$. Except the cases $R=A^{\an}_3$, $A^{\an}_{(3)}$, $A_{2*}$, $A_3$ (for which we do not have the evaluation map used below), we define
$$H_R(p):=\R \otimes_{R_p} H_{R,p},\quad  H^{\natural}_R(p):=\R \otimes_{R_p} H^{\natural}_{R,p},\quad {}^{\prime}H^{\natural}_R(p):=\R \otimes_{R_p} {}^{\prime}H^{\natural}_{R,p},$$
 where $R_p\to \R$ is the evaluation $f\mapsto f(p)$ at $p$. We call it the fiber at $p$.  
\end{sbpara}

\begin{sbpara}\label{csp}

Let $V$ be a finite dimensional vector space and let $\Phi$ be a finite set of increasing filtrations on $V$. In the case where we are given a splitting of $W'$  for each $W'\in \Phi$, we say that these splittings are {\it compatible} if  there is a  direct sum decomposition $V=\bigoplus_{w\in \Z^\Phi} V[w]$ such that, for each $k\in\Z$, 
the part of weight $k$ of the given splitting of $W'$ for each $W'\in \Phi$ is $ \bigoplus_{w\in \Z^{\Phi}, w(W')= k} V[w]$. 

This is equivalent to the condition that  
there is a 
homomorphism $\tau: {\bf G}_m^{\Phi}\to \Aut(V)$ of algebraic groups such that for each $W'\in \Phi$, its $W'$-component ${\bf G}_m \to \Aut(V)$ provides the given splitting of $W'$.

If such a compatible family of splittings exists, we say that 
$\Phi$ {\it has a compatible splitting} and we call a homomorphism  $\tau: {\bf G}_m^{\Phi}\to \Aut(V)$ with the above property a {\it compatible splitting of $\Phi$}.
\end{sbpara}

\begin{sbprop}\label{cspthm}  For $p\in X_{[:]}$ and for $\tilde p\in X^{\log}_{[:]}$ lying over $p$, the set $\{W\}\cup \Phi(p)$ of increasing filtrations on $H_{\R,\tilde p}$  has a compatible splitting in the sense of $\ref{csp}$.

\end{sbprop}

This is shown in 10.2 of \cite{KNU08}. See Proposition \ref{twi0}.

\begin{sbprop}\label{HA10} Let $R$ and $R'$ be as in Proposition $\ref{HA}$. 

$(1)$  The $R$-module  $H_R$ is locally free of finite rank, and we have
$$R^{\log*}\otimes_R H_R \overset{\cong}\to R^{\log*} \otimes_{\R} H_\R.$$ We have
$$R\otimes_{A^{\an}_1} H_{A^{\an}_1}\overset{\cong}\to H_R.$$

For $p\in X_{[:]}$ and for $\tilde p\in X_{[:]}^{\log}$ lying over $p$, we have 
$$H_A(p) \overset{\cong}\to H_{\R,\tilde p}.$$

The weight filtration $W$ of $R^{\log*}\otimes_{\R} H_\R$ comes from a locally splitting filtration on $H_R$ over $R$. 
For $p$ and $\tilde p$ as above,  the relative monodromy filtration $W^{(j)}$ ($1\leq j\leq m$) of $R^{\log*}_{\tilde p} \otimes_{\R} H_{\R,\tilde p}$ comes from a splitting filtration on $H_{R,p}$ over $R_p$.

$(2)$ The $R$-module $H^{\natural}_R$ and the $R^{\log}$-module $H_R^{\natural\log}$ are locally free of finite rank, and 
$$R^{\log}\otimes_R H_R^{\natural} \overset{\cong}\to H_R^{\natural\log}.$$ We have
$$R\otimes_{A^{\an}_1} H^{\natural}_{A^{\an}_1}\overset{\cong}\to H^{\natural}_R, \quad R'\otimes_R H_R= R'\otimes_R H^{\natural}_R.$$
The weight filtration $W$ of $(R')^{\log}\otimes_{\R} H_\R= (R')^{\log}\otimes_R H^{\natural}_R$ comes from a locally splitting filtration on $H^{\natural}_R$ over $R$. 
For $p\in X_{[:]}$ and for $\tilde p\in X^{\log}_{[:]}$ lying over $p$, the relative monodromy filtration $W^{(j)}$ ($1\leq j\leq m$) of $(R')^{\log}_{\tilde p} \otimes_{\R} H_{\R,\tilde p}$ comes from a splitting filtration on $H^{\natural}_{R,p}$ over $R_p$.

$(3)$ We have an isomorphism
$$A^{\an}_{3,\C}\otimes_{\cO_X} H_{\cO}\cong H_{A^{\an}_3,\C}=H^{\natural}_{A^{\an}_3,\C}$$
induced from 
the homomorphism $\cO_X^{\log}\to A^{\an,\log}_{3,\C}$ ($\ref{cOlog}$) and the isomorphism  $\cO_X^{\log}\otimes_{\cO_X} H_{\cO}\cong \cO_X^{\log}\otimes_{\R} H_{\R}$. 

\end{sbprop}

\begin{sbpara}\label{HA1(1)} 
  We prove Proposition \ref{HA10}.

Consider the situation in \ref{xi1} and consider the 
constant sheaf $\xi^{\R}H_\R$ on $U_{[:]}(\Psi)$. It is isomorphic to the stalk of $\xi^{\R} H_\R$  at $\tilde p$.  

 (1) follows from \ref{xi1}.

 We prove (2). We fix polarizations of $\gr^W_w H$ ($w\in \Z$).

We use the following notation:  $p\in X_{[:]}$, $\tilde p\in X_{[:]}^{\log}$ lies over $p$, $m$ and $f_{j,k}$ are as in \ref{m,f_jk}, and $(2\pi i)^{-1}\log(f_{j,k})=x_{j,k}+iy_{j,k}$ with $x_{j,k}$ and $y_{j,k}$ being real. 
  We denote $y_{j,1}$ also by $y_j$. 
Let $$\tau: {\bf G}_{m,\R}^m\to \Aut_\R(H_{\R,\tilde p})$$ be a homomorphism which gives a compatible splitting of $(W^{(j)})_{1\leq j\leq m}$ of $H_{\R,\tilde p}$ in Proposition \ref{cspthm}. 
Consider the splitting 
of $H_{\R,\tilde p}=\bigoplus_{\ell\in \Z^m}\;  P_{\ell}$, where $P_{\ell}$ is the part of weight $\ell$ for $\tau$. We have $\xi^{\R} H_{\R, \tilde p}= \bigoplus_{\ell\in \Z^m}\;Q_{\ell}$, where $Q_{\ell}= \xi^{\R} P_{\ell}$. We have this decomposition of 
the constant sheaf $\xi^{\R} H_\R= \bigoplus_{\ell\in \Z^m} \;Q_{\ell}$. 

Since $N_{j,k}$ preserves $W^{(j')}$ 
for  $1\leq j'\leq m$ and $N_{j,k}W^{(j')}_w\subset W^{(j')}_{w-2}$ if $j\leq j'\leq m$, we have 
$${\textstyle\prod}_{j=1}^m (y_j/y_{j+1})^{-\ell(j)/2} Q_{\ell}\subset 
{\textstyle\sum}_{\ell'\leq \ell} R^{\log}{\textstyle\prod}_{j=1}^m (y_j/y_{j+1})^{-\ell'(j)/2} P_{\ell'}.$$
Hence we have 
 $$H^{\natural}_R ={\textstyle \bigoplus}_{\ell\in \Z^m} R{\textstyle\prod}_{j=1}^m (y_j/y_{j+1})^{-\ell(j)/2} Q_{\ell}.$$
 The statements in (2) except the ones about $W^{(j)}$ and $W$ follow from this.

The filtration $W^{(j)}$ over $(R')^{\log}$ comes from the filtration $W^{(j)}H^{\natural}_R$ on $H^{\natural}_R$  defined by $W^{(j)}_a H^{\natural}_R= \bigoplus_{\ell \in \Z^m, \ell(j)\leq a} \; R\prod_{j=1}^m (y_j/y_{j+1})^{-\ell(j)/2}\otimes_\R Q_{\ell}$.

We consider the statement about $W$. By Proposition \ref{cspthm}, there are direct sum decompositions $P_{\ell}=\bigoplus_{a\in \Z} P_{\ell,a}$ for all $\ell\in \Z^m$ such that $W_aH_{\R, \tilde p}=\bigoplus_{\ell\in \Z^m,b\leq a} \; P_{\ell,b}$ for all $a\in \Z$. The filtration $W$ over $(R')^{\log}$ comes from the filtration $W_{\bullet}H^{\natural}_R$ defined by $W_aH^{\natural}_R=  \bigoplus_{\ell\in \Z^m, b\leq a}\;  R\prod_{j=1}^m (y_j/y_{j+1})^{-\ell(j)/2}\otimes_\R Q_{\ell,a}$, where $Q_{\ell,a}$ is the constant sheaf $\xi^{\R}P_{\ell,a}$. This proves the statement about $W$.

(3) is clear. 

\end{sbpara}

By the above proof in \ref{HA1(1)}, we have the following \ref{tortwi} and Lemma \ref{twigr}. 
\begin{sbpara}\label{tortwi}
In  \ref{HA1(1)}, consider the automorphism $$t(y):= \tau\Bigl(\Bigl(\sqrt{y_{j+1}/y_j}\Bigr){}_{1\leq j\leq m}\Bigr)$$ of $(R')^{\log}\otimes_\R H_{\R, \tilde p}$ over $(R')^{\log}$ with $y_j$ as in \ref{HA1(1)}  ($y_{m+1}$ denotes $1$).

(1) We have an isomorphism $$\xi^{\R}t(y)(\xi^{\R})^{-1}: H_{R,p} \overset{\cong}\to H^{\natural}_{R,p}$$ in $(R')^{\log}_p \otimes_{\R} H_{R,\tilde p}$.

(2) We have $\xi^{\R}t(y): R_p\otimes _\R H_{\R,\tilde p} \overset{\cong}\to H^{\natural}_{R,p}$. 

(3) We have $t(y) : R^{\log}_{\tilde p}\otimes_\R H_{\R,\tilde p} \overset{\cong}\to H_{R,\tilde p}^{\natural\log}$.

\end{sbpara}

\begin{sblem}\label{twigr} For $\ell\in \Z^m$, let $\gr_{[\ell]}= (\bigcap_{j=1}^m W^{(j)}_{\ell(j)})/(\sum_{\ell'<\ell} \bigcap_{j=1}^m W^{(j)}_{\ell'(j)})$. Then the maps $\gr_{[\ell]}H_{R,p}\to \gr_{[\ell]}H_R(p)$ are surjective, and the maps $\prod_{j=1}^m (y_j/y_{j+1})^{-\ell(j)/2}: \bigcap_{j=1}^m W^{(j)}_{\ell(j)}H_{R,p} \to H^{\natural}_R$ induce homomorphisms 
$\gr_{\ell} H_{R,p} \to H^{\natural}_R(p)$ and an isomorphism 
$${\textstyle\bigoplus}_{\ell\in \Z^m}\; \gr_{[\ell]}H_R(p) \overset{\cong}\to H^{\natural}_R(p).$$

\end{sblem}

\subsection{Geometric interpretation of the theory of $\SL(2)$-orbits}\label{ss:HA}

In our geometric interpretation,  the theory of $\SL(2)$-orbit is summarized in Theorems \ref{HA1}, \ref{SL2I}, \ref{SL2II}, \ref{ourSL2}, \ref{oldSL20}   by using the sheaves $H_A$, $H^{\natural}_A$, ${}^{\prime}H^{\natural}_A$ in Section \ref{ss:twoHA} on the space of ratios.

Let $X$ be an object of $\cB(\log)$ and let $H$ be an $\R$-LMH  
on $X$ such that the graded quotients  of the weight filtration $W$ of $H$ are polarizable. 

Theorem   \ref{HA1} shows that for $H_A^{\natural}$ (resp.\ for $H^{\natural}_A$ and ${}^{\prime}H_A^{\natural}$ in the case where $H$ is pure, resp.\ for ${}^{\prime}H_A^{\natural}$ in the case where $H$ is a PLH),  the Hodge filtration (resp.\  the Hodge decomposition, resp.\ the Hodge metric) 
does  not degenerate even when $H$ degenerates.

Theorem \ref{SL2I}  shows that in the case where $H$ is a PLH, the relative monodromy filtrations of $H_A$ and those of ${}^{\prime}H_A^{\natural}$  have orthogonal splittings over $A^{\an}_{(2+)}$ for the Hodge metric.

Theorem \ref{SL2II} shows that the weight filtration of $H_A$ and that of $H_A^{\natural}$ have a canonical splitting over $A^{\an}_{(2+)}$.

Theorems \ref{ourSL2} and \ref{oldSL20} give
  our interpretation  of the associated  SL(2)-orbit (Schmid \cite{Sc} and Cattani--Kaplan--Schmid \cite{CKS1}) as the log analytic fiber of $H$ at a point of $X_{[:]}$.

The proofs of the statements in  this Section \ref{ss:HA} are given in Section \ref{ss:pfSL}  after a preparation Section \ref{ss:HA0}. 

\begin{sbthm}\label{HA1} 

$(1)$ (Hodge filtration and $H^{\natural}_A$.) The Hodge filtration of $A^{\an}_{3,\C}\otimes_{\cO_X} H_{\cO}= H^{\natural}_{A^{\an}_3,\C}$ (Proposition $\ref{HA10}$ $(3)$) induced by the Hodge filtration of $H_{\cO}$ comes from a locally splitting filtration on $H^{\natural}_{A^{\an}_{2+},\C}$ over $A^{\an}_{2+,\C}$.  Consequently, the Hodge filtration of $A^{\an}_{(3),\C}\otimes_{\cO_X} H_{\cO}=H^{\natural}_{A^{\an}_{(3)},\C}$ induced by the Hodge filtration of $H_{\cO}$ comes from a locally splitting filtration on $H^{\natural}_{A,\C}$ over $A^{\an}_{(2+),\C}$.  Every fiber $H^{\natural}_A(p)$ ($p\in X_{[:]}$) with the induced Hodge filtration is a mixed Hodge structure.

$(2)$ (Hodge decomposition and $H^{\natural}_A$.) Assume that $H$ is pure of weight $w$. Then for the Hodge filtration $F$ of $A_{(3)}^{\an,\log}\otimes_\R H_\C$, we have the Hodge decomposition $A_{(3)}^{\an,\log}\otimes_\R H_\C= \bigoplus_{p+q=w} F^p\cap \overline{F^q}$. Furthermore, this direct sum decomposition comes from a direct sum decomposition of $H^{\natural}_{A,\C}$ over $A^{\an}_{(2+)}$.

$(3)$ (Hodge metric and ${}^{\prime}H^{\natural}_A$.) Assume that $H$ is a PLH. Then the  duality $\langle \cdot,\cdot\rangle : H_\R \times H_\R \to \R$ induces the perfect self-duality ${}^{\prime}H^{\natural}_A\times {}^{\prime}H^{\natural}_A\to A^{\an}_{(2+)}$ of the $A^{\an}_{(2+)}$-module  ${}^{\prime}H^{\natural}_A$. 
Define the Hodge metric $$(\cdot,\cdot): (A^{\an,\log}_{(3),\C}\otimes_\R H_\R)\times (A^{\an,\log}_{(3),\C}\otimes_\R H_\R)\to A^{\an,\log}_{(3),\C}$$ by $(u,v)= i^{p-q}\langle u, \bar v\rangle $ for $u\in $ 
the $(p,q)$-Hodge component. Note that the $(p,q)$-Hodge component  and $(p',q')$-Hodge component  are orthogonal for $(\cdot,\cdot)$  unless $(p,q)= (p',q')$.
Then this Hodge metric gives a perfect duality of $A^{\an}_{(2+)}$-modules 
$${}^{\prime}H^{\natural}_{A,\C}\times {}^{\prime}H^{\natural}_{A, \C}\to A^{\an}_{(2+),\C}.$$
 Every fiber ${}^{\prime}H^{\natural}_A(p)$ ($p\in X_{[:]}$) with $\langle\cdot,\cdot\rangle$ and with the induced Hodge filtration is a polarized Hodge structure.

\end{sbthm}

\begin{sbthm}\label{SL2I} Assume that $H$ is a PLH. Let $p\in X_{[:]}$ and let $W^{(j)}$ ($1\leq j\leq m$) be the relative monodromy filtrations at $p$. 

$(1)$ The filtration $W^{(j)}$ on ${}^{\prime}H^{\natural}_{A,p}$, which is induced by that on $H^{\natural}_{A,p}$ in  Proposition $\ref{HA10}$ $(2)$, 
has an orthogonal splitting over $A^{\an}_{(2+),p}$ for the Hodge metric. 

$(2)$ The induced orthogonal splitting of $H_{A^{\an}_{(3)},p}^{\natural}=H_{A^{\an}_{(3)},p}$ comes from a splitting of $W^{(j)}$ of $H_{A,p}$ (Proposition 
$\ref{HA10}$ $(1)$) over $A^{\an}_{(2+),p}$. 

\end{sbthm}

This splitting of $W^{(j)}$ in Theorem \ref{SL2I} is denoted by $\spl_{W^{(j)}}^{\BS}$.

\begin{sbpara}\label{canspl} Recall that every $\R$-mixed Hodge structure $H$ has a canonical splitting of the weight filtration defined by the $\R$-split $\R$-mixed Hodge structure $\tilde F_0$ in \cite{CKS1} (3.30), (3.31) (not the $\tilde F$ there) whose $\R$-structure and  weight filtration are the same as those of $H$.

\end{sbpara}

The following result is proved in \cite{logmot}.

\begin{sbprop}\label{mixpure} Let $X$ be an object of $\cB(\log)$,  let $R$ be a subfield of $\R$, and let $H$ be an $R$-LMH (resp.\ $R$-LVMH) on $X$. Let $S$ be the standard log point. 
Then locally on $X$, there is an $R$-PLH (resp.\ $R$-LVPH) $H'$ on $X\times S$ and an injective homomorphism $H_R \to H'_R$ of the local systems of $R$-modules on $(X\times S)^{\log}$ satisfying the following conditions {\rm (i)} and {\rm (ii)}.

{\rm (i)} The Hodge filtration of $H$ on $H_{\cO}= (\tau_X)_*(\cO^{\log}_X\otimes_R H_R)= (\tau_{X\times S})_*(\cO^{\log}_{X\times S} \otimes_R H_R)$ coincides with the restriction of the Hodge filtration of $H'$ on $H'_{\cO}=(\tau_{X\times S})_*(\cO^{\log}_{X\times S} \otimes_R H'_R)$.

{\rm (ii)} For every $t\in (X\times S)^{\log}$,  the  weight filtration of $H$ on $H_{R, t} \otimes_{\Z} \Q$  is the restriction of the relative monodromy filtration on $H'_{R,t}\otimes_{\Z} \Q$ of the logarithm $N:H'_{R, t}\otimes_{\Z} \Q \to H'_{R,t}\otimes_\Z \Q$ of the action of the standard generator of $\pi_1(S^{\log})$.

\end{sbprop}

\begin{sbthm}\label{SL2II}  
$(1)$ There is a unique way to give a splitting $\spl_W$ of $W$ on $H_A$ or on $H^{\natural}_A$ over $A^{\an}_{(2+)}$ for every $\R$-LMH $H$ with polarizable $\gr^W$ 
  on an object of $\cB(\log)$ (called the {\rm canonical splitting}), satisfying the following {\rm (i)} and  {\rm (ii)}.

{\rm (i)} It is compatible with pullbacks of LMH. 

{\rm (ii)} In the case where $X$  is $\Spec(\C)$ with the trivial log structure and $H$ is a usual mixed Hodge structure, $\spl_W$ is the canonical splitting ($\ref{canspl}$).

$(2)$ This $\spl_W$ is also characterized by the following properties {\rm (i)} and  {\rm (ii)}.

{\rm (i)} It is compatible with pullbacks of LMH.

{\rm (ii)}  For  $H'$ as in Proposition $\ref{mixpure}$, $\spl_W$ is the orthogonal decomposition of $W$ by the Hodge metric of $H'$.

\end{sbthm}

Our understanding of the associated $\SL(2)$-orbit $(\rho,\varphi)$ is given in the following Theorem \ref{ourSL2} as a fiber of a log real analytic structure ${}^{\prime}H^{\natural}_A$ of PLH. The relation with the usual associated $\SL(2)$-orbit in \cite{Sc} and \cite{CKS1} is explained in Theorem \ref{oldSL20}.

\begin{sbthm}\label{ourSL2} 

 Let $p\in X_{[:]}$.
 Assume that $H$ is a PLH. 
 
 $(1)$ There is a unique homomorphism 
$$\rho: \SL(2)^m_\R\to \Aut({}^{\prime}H^{\natural}_A(p), \langle\cdot,\cdot\rangle)$$ of algebraic groups over $\R$  characterized by the following properties {\rm (i)} and {\rm (ii)}.
 
{\rm (i)}   For $a=(a_1, \dots, a_m)\in ({\bf G}_m)^m$, let $\Delta(a)$ be the element of $\SL(2)^m$ whose $j$-th component is $\begin{pmatrix} a_j^{-1}& 0 \\0 &a_j\end{pmatrix}$. Then for $\ell\in \Z^m$, on the image of $y_1^{w/2}\gr_{[\ell]}H_A(p) \to {}^{\prime}H^{\natural}_A(p)$, 
$\rho(\Delta(a))$ acts as the multiplication by $\prod_{j=1}^m (a_j/a_{j+1})^{\ell(j)-w}=a_1^{-w}\prod_{j=1}^m (a_j/a_{j+1})^{\ell(j)}$, where $a_{m+1}:=1$.

{\rm (ii)} For $1\leq j\leq m$, the  homomorphism  $\Lie(\rho):{\frak {sl}}_2(\R)^m \to \mathrm{End}_\R({}^{\prime}H^{\natural}_A(p))$ induced by $\rho$ sends $\begin{pmatrix} 0 & 1 \\ 0& 0\end{pmatrix}$ in the $j$-th ${\frak {sl}}_2(\R)$ in ${\frak {sl}}_2(\R)^m$ 
to the following $\hat N_j$.

Let $(f_{j,k})_{j,k}$ be as in $\ref{m,f_jk}$ and 
let $y_{j,k}= -(2\pi)^{-1}\log|f_{j,k}|$. Let $N_{j,k}$ be the monodromy logarithm corresponding to the homomorphism $(M^{\gp}_X/\cO_X^\times)_p\to \R$ 
sending $f_{j,k}$ to $1$ and the other $f_{j',k'}$ to $0$, which acts on the inverse image of $H_\R$ on the inverse image of $p$ in $X^{\log}_{[:]}$. Let  $1\leq j\leq m$. 
Then the action of $\sum_{k=1}^{r(j)} y_{j,k}N_{j,k}$ on  $A^{\an,\log*}_{(3),\tilde p}\otimes_\R H_{\R,\tilde p}$ preserves  ${}^{\prime}H^{\natural}_{A, p}$, and the induced map 
${}^{\prime}H^{\natural}_A(p)\to {}^{\prime}H^{\natural}_A(p)$ is independent of the choice of $(f_{j,k})_{j,k}$.
  Let $\hat N_j$ be this induced map. 

$(2)$  
Let $D$ be the set of all decreasing filtrations $F'$ on  ${}^{\prime}{H^{\natural}_A(p)}_{\C}$ 
such that $({}^{\prime}H^{\natural}_A(p), \langle \cdot,\cdot\rangle, F')$ is a polarized Hodge structure of weight $w$. Then we have a 
 holomorphic map
$$\varphi: \frak H^m \to D,$$
where $\frak H$ denotes the upper half plane, characterized by the following properties {\rm (i)} and {\rm (ii)}.

{\rm (i)}  $\varphi(g\alpha)= \rho(g)\varphi(\alpha)$ for all $g\in \SL(2,\R)^m$ and $\alpha\in \frak H^m$. 

{\rm (ii)}  $\varphi(i, \dots, i)$ is the induced Hodge filtration (Theorem $\ref{HA1}$ $(3)$) on ${{}^{\prime}H^{\natural}_A(p)}_\C$. 

\end{sbthm}

\begin{sbthm}\label{oldSL20} 
 Let $p\in X_{[:]}$.

$(1)$ Assume that $\gr^W_wH$ are polarized. Then 
there is a unique compatible splittings of $W$ and $W^{(j)}$ for $1\leq j\leq m$ of  $H_A(p)$ 
such that the splitting of $W$ comes from the canonical splitting $\spl_W$ of $H_A$ (Theorem $\ref{SL2II}$) and the induced splitting of $W^{(j)}$ of $\gr^W_wH_A(p)$ comes from  the orthogonal splittings $\spl_{W^{(j)}}^{\BS}$ of $W^{(j)}$ of $\gr^W_wH_A$ for the Hodge metric (Theorem $\ref{SL2I}$). 

$(2)$  Let $\tilde p\in X_{[:]}^{\log}$ lying over $p$. 
 Then we have the following triple $(\rho, \varphi, \br)$, where 
$\rho$ is a homomorphism $\SL(2)^m_\R\to \prod_w\; \Aut(\gr^W_wH_{\R,\tilde p}, \langle\cdot,\cdot\rangle_w)$ of algebraic groups over $\R$, 
 $\varphi$ is a holomorphic map $\frak H^m\to \prod_w\;D(\gr^W_w)$, where  $D(\gr^W_w)$ is the space of decreasing filtrations $F'$ on $\gr^W_wH_{\C,\tilde p}$ such that $(\gr^W_w(H_{\R, \tilde p}), \langle\cdot, \cdot\rangle_w, F')$ is a polarized Hodge structure, such that  
$$\varphi(gz)= \rho(g)\varphi(z)\quad \text{for all $g\in \SL(2)^m_\R$ and $z\in \frak H^m$},$$ and $\br$ is a decreasing filtration on $H_{\C,\tilde p}$ such that $\varphi(i, \dots, i)=(\br(\gr^W_w))_w$.

First $\br$ is obtained from the Hodge filtration of $H^{\natural}_A(p)$ (Theorem $\ref{HA1}$ $(1)$) via the isomorphism $$H_{\R,\tilde p} \cong H_A(p)\cong{\textstyle\bigoplus}_{\ell\in \Z^m} \gr_{[\ell]} H_A(p)  \overset{\cong}\to H^{\natural}_A(p),$$
where the middle isomorphism is by the compatible splitting of $(W^{(j)})_{1\leq j\leq m}$ in $(1)$ and 
the last isomorphism is that of Lemma $\ref{twigr}$. 
  Next the $\gr^W_w$-component of $(\rho, \varphi)$ is obtained from $(\rho, \varphi)$ 
 for ${}^{\prime}(\gr^W_wH)_A^{\natural}$
 in Theorem $\ref{ourSL2}$ $(1)$ and $(2)$ via the isomorphism $\gr^W_wH_A(p) \overset{\cong}\to {}^{\prime}\gr^W_wH_A^{\natural}(p)$ which is obtained from the above isomorphism by multiplying it by $y_1^{w/2}$.

This $(\rho, \varphi, \br)$ depends on the choices of $f_{j,k}$. If we change the choices of $f_{j,k}$, $(\rho, \varphi, \br)$ is changed to 
$ (\tau(a)\rho(\cdot)\tau(a)^{-1}, \tau(a)\varphi, \tau(a)\br)$ for some $a\in \R_{>0}^m$, where $\tau$ is the homomorphism ${\bf G}^m_{m,\R} \to \Aut(H_{\bR,\tilde p})$ which gives the compatible splitting of $(W^{(j)})_{1\leq j\leq m}$ of $H_{\R, \tilde p}=H_A(p)$ in $(1)$.
\end{sbthm}

\begin{sbrem}

(1) If $H$ is a PLH, the pair $(\rho, \varphi)$ on $H_{\R,\tilde p}$ in Theorem \ref{oldSL20} (2)  is our interpretation of 
 the   $\SL(2)$-orbit $(\rho, \varphi)$
in \cite{Sc} and \cite{CKS1}. 

(2) We have the $\SL(2)$-orbit $(\rho, \varphi)$ in Theorem  \ref{ourSL2} and that in Theorem \ref{oldSL20}. The former is canonical, though the latter is determined only modulo the action of $\tau(a)$ ($a\in \R^m_{>0})$. In the former,   $\rho$ depends only on the local system $H_\R$, and does not depend on the Hodge filtration nor polarization.

(3) In Theorem \ref{ourSL2}, in the case where $H$ is a PLH, the splitting of ${}^{\prime}H^{\natural}_A(p)$ into the direct sum of the images of $\gr_{[\ell]}H_A(p)$ is the compatible orthogonal splitting of $W^{(j)}$ ($1\leq j\leq m$) for the polarized Hodge structure ${}^{\prime}H_A^{\natural}(p)$, and this splitting of $W^{(j)}$ coincides with the one induced by the orthogonal splitting of $W^{(j)}$ of ${}^{\prime}H^{\natural}_{A,p}$ 
in Theorem \ref{SL2I}. 

\end{sbrem}

\begin{sbpara}\label{6.1Ex} Example.

Let the degenerate elliptic curve  $\cE= E_{\Delta,q}$ over $\Delta$ be as in \ref{another}. Consider the following $\Z$-LVPH $H$ on 
  $X=\Delta$  of weight $1$ which extends the variation of polarized Hodge structure $\{H^1(\cE_q)\}_{q \in \Delta\smallsetminus \{0\}}$   (\cite{KU} 0.2.18). This is the first higher direct image $H^1$ of $\cE/\Delta$ (\ref{PLH2}).
  
  (1) Description of $H$.

The local system $H_\Z$ on $\Delta^{\log}$ 
is locally isomorphic to $\Z^2$. Let  $q\in \Delta\smallsetminus \{0\}$ and let $\tau=(2\pi i)^{-1}\log(q)=x+iy$ ($x, y$ being real). Then 
the fiber $\cE_q$ is identified with $\C/(\Z+\Z\tau)$, $H_1(\cE_q, \Z)$ is identified with $\Z+\Z\tau$, and $H_{\Z,q}=\Hom(H_1(\cE_q, \Z), \Z)=\Z e_1+\Z e_2$, where $e_1$ corresponds to the homomorphism 
which sends $1$ to $0$ and $\tau$ to $1$, and $e_2$ corresponds to the homomorphism which sends $1$ to $1$ and $\tau$  to $0$. 
These $\tau$, $x$,  and $e_2$ are not determined by $q$ (they have branches), but $\tau e_1+e_2$ is canonical and extends to a global section of $H_{\cO}$ on $\Delta$. We have
$$H_{\cO}= \cO_X e_1 + \cO_X (\tau e_1+e_2).$$ The Hodge filtration $F$ on $H_{\cO}$ is given as $$F^0=H_{\cO},\quad F^1= \cO (\tau e_1+e_2), \quad F^2=0.$$

(2) $H_A$, $H^{\natural}_A$, ${}^{\prime}H^{\natural}_A$. 

In what follows, let $R=A^{\an}_{(2+)}$. Then $xe_1+e_2$ at $q\in \Delta\smallsetminus \{0\}$ is canonical and extends to a global section of $H_A$ on $\Delta$. We have 
$$H_A= R e_1 + R(xe_1+e_2), \quad H^{\natural}_A= R e_1+ R y^{-1}(xe_1+e_2), \quad 
 {}^{\prime}H^{\natural}_A= R y^{1/2}e_1+ R y^{-1/2}(xe_1+e_2).$$

(3) The Hodge filtration of $H^{\natural}_A$ (Theorem \ref{HA1} (1)). 

It is given by $$F^0= H^{\natural}_{A,\C},\quad  F^1= R_{\C}\cdot y^{-1}(\tau e_1+e_2), \quad F^2=0.$$

(4) The Hodge decomposition.  

The Hodge decomposition of ${}^{\prime}H^{\natural}_{A,\C}$, which is induced by that of $H^{\natural}_{A,\C} $ in Theorem \ref{HA1} (2), 
is  $${}^{\prime}H^{\natural}_{A,\C}= {}^{\prime}H_{A,\C}^{\natural,1,0}\oplus {}^{\prime}H_{A,\C}^{\natural, 0,1}, \quad{\text {where }}\quad  {}^{\prime}H_{A,\C}^{\natural, 1,0}= R_{\C}u, \quad {}^{\prime}H_{A,\C}^{\natural, 0,1} =R_{\C} \bar u$$ $$\text{with} \;\;u= y^{-1/2}((x+iy)e_1+e_2),  \quad \bar u= y^{-1/2}((x-iy)e_1+e_2).$$

(5) The Hodge metric (Theorem \ref{HA1} (3)).

The bilinear form $\langle\cdot,\cdot\rangle$ which gives the canonical polarization of $H$  is the skew-symmetric form satisfying $\langle e_2,e_1\rangle=1$.

The Hodge metric $(\cdot,\cdot)$ is the Hermitian form on  ${}^{\prime}H^{\natural}_{A,\C}$ given by $(u,u)=(\bar u, \bar u)=2$, $(u, \bar u)= (\bar u, u)=0$, where $u$ is as in (4). 

(6) The relative monodromy filtration and its orthogonal splitting (Theorem \ref{SL2I}). 

Let $p$ be the origin of $\Delta$. Then the 
 local monodromy logarithm $N$ at a point $\tilde p\in \Delta^{\log}$ lying over $p$ is:
$$Ne_2=e_1, \quad Ne_1=0.$$
We have  $\Phi(p)=\{W'\}$, where $$0=W'_{-1}\subset W'_0=\R e_1=W'_1\subset W'_2=H_{\R, \tilde p}.$$ The splitting of $W'$ of ${}^{\prime}H^{\natural}_A$ by the Hodge metric is given by 
$y^{-1}(xe_1+e_2)$. 

(7) The associated $\SL(2)$-orbit (Theorem \ref{ourSL2}, Theorem \ref{oldSL20}). 
 
 The associated $\SL(2)$-orbit $(\rho, \varphi)$ at $p$ in Theorem \ref{oldSL20} (resp.\ Theorem \ref{ourSL2}) is as follows. $\rho$ is the identity map $\SL(2)\to \SL(2)$ with respect to the base $e'_1$, $e'_2$ of $H_{\R,\tilde p}$ (resp.\ ${}^{\prime}H^{\natural}_A(p)$), where $e'_1$, $e'_2$ are the classes of $e_1$, $xe_1+e_2$ (resp.\ $y^{1/2}e_1$, $y^{-1/2}(xe_1+e_2)$), respectively. 
  $\varphi(\alpha)$ for $\alpha\in \frak H$ is the Hodge filtration such that $\varphi(\alpha)^0= H_{\bC,\tilde p}$  (resp.\ ${}^{\prime}H^{\natural}_A(p)_{\C}$), $\varphi(\alpha)^1$ is 
generated by $\alpha e_1'+e_2'$, and $\varphi(\alpha)^2=0$. 
\end{sbpara}

\subsection{Review on fs log points}\label{ss:HA0}
In this section (Section \ref{ss:HA0}), as  a preparation for the proofs of the statements in Section \ref{ss:HA} which are given in the next section (Section \ref{ss:pfSL}), we 
deduce results on the  case where $X$ is an fs log point from the known $\SL(2)$-orbit theorems.

In this section (Section \ref{ss:HA0}), let $X=\{x\}$ be an fs log point. Let $H$ be an $\R$-LMH on $X$ with polarizable $\gr^W$.

Let $\cT$ be the canonical toric variety associated to $X$ (Section \ref{ss:pt}).  Let the notation be as in Section \ref{ss:pt} and Section \ref{ss:LMH}.

\begin{sbpara}\label{ss0} For $s\in \cT_{\triv}$ and  $h\in \Hom((M^{\gp}_X/\cO_X^\times)_x,\C)$, define  $s^{(h)}\in \cT_{\triv}$ by 
$s^{(h)}(q)=\exp(2\pi ih(q))s(q)$ for $q\in M_{X,x}$.

Fix $s_0\in \cT_{\triv}$. Then every point $s$ of $\cT_{\triv}$ is written in the form $s=s_0^{(h)}$ for some  $h\in \Hom(
(M^{\gp}_X/\cO_X^\times)_x,\C)$ and this $h$ is unique modulo $\Hom((M^{\gp}_X/\cO_X^\times)_x, \Z)$.

Let $\tilde p\in X_{[:]}^{\log}$.  Fix  a path from $\tilde p$ to $s_0$ in $\cT_{[:]}^{\log}$. For an $h$ as above and for $s=s_0^{(h)}$, we identify $H_{\cT, \R, s}$ with $H_{\R, \tilde p}$ through the composition of this path and the  path from $s_0$ to $s$ in $\cT_{\triv}$ given by $s_0^{(th)}$ ($0\leq t\leq 1$).

Let $F$ be the Hodge filtration of $H$ which is extended canonically to the Hodge filtration of $H_{\cT}$, and regard the Hodge filtration $F(s)$ on $H_{\cT,\C, s}$ for $s=s_0^{(h)}$ as a filtration on $H_{\C,\tilde p}$ via the above identification.

For an $h$ as above, let the nilpotent operator $N_h:H_{\C,\tilde p}\to H_{\C, \tilde p}$ be 
the image of $h$ under $\Hom((M^{\gp}_X/\cO^\times_X)_x, \C)\cong \pi_1(x^{\log})\otimes_{\Z}\C \to \End_{\C}(H_{\C, \tilde p})$, where the last arrow is the logarithm of the action of $\pi_1(x^{\log})=\pi_1(x^{\log}_{[:]})$ on $H_{\C,\tilde p}$.
\end{sbpara}

\begin{sblem}\label{ss1} In the situation of $\ref{ss0}$, we have 
$$F(s)= \exp(N_h)F(s_0).$$

\end{sblem}

\begin{pf} Let the notation be as in \ref{xi0}.
  Write 
$F=\exp(\sum_{\la\in \La} z_{\la}N_{\la})F_0$, where $F_0$ is a filtration on $H_{\C,\tilde p}$. Here  $z_{\la}$ is a branch of $(2\pi i)^{-1}\log(q_{\la})$ at $\tilde p$. Let $z_{\la}(s)\in \C$ be the branches of $(2\pi i)^{-1}\log(q_{\la}(s))$ defined by the fixed path from $\tilde p$ to $s_0$ and the path from $s_0$ to $s$ 
given by $s_0^{(th)}$ ($0\leq t\leq 1$). 
  We have $F(s)= \exp(\sum_{\la\in \La} z_{\la}(s)N_{\la})F_0$ and $z_{\la}(s)=z_{\la}(s_0)+ h(q_{\la})$. 
Hence $F(s)= \exp(\sum_{\la\in \La} h(q_{\la})N_{\la})F(s_0)$, and $\sum_{\la \in \La} h(q_{\la})N_{\la}=N_h$. 
 \end{pf}

\begin{sbpara}\label{handN} Let ($h_j)_{j\in J}$ be a finite family of generators of the cone $\Hom((M_X/\cO_X^\times)_x, \R_{\geq 0})$. Fix $s_0\in \cT_{\triv}$. Then  
the sets $$\{s_0^{(\sum_{j \in J} \alpha_jh_j)}\;|\; \alpha_j\in \C, \Im(\alpha_j)\geq c\}$$ given for $c\in \R$ form a base of the filter in $\cT$ consisting of all subsets of the form $U\cap \cT_{\triv}$, where $U$ ranges over all neighborhood of $x$.

Hence the study of $H_{\cT}(s)$ when $s\in \cT_{\triv}$ converges to $x$ is the study of $\exp(\sum_{j\in J} \alpha_jN_j)F(s_0)$ when $\Im(\alpha_j)$ tend to $\infty$ for all $j\in J$. 

Let $p\in X_{[:]}$ and let $(f_{j,k})_{j,k}$ be as in \ref{m,f_jk} and let $y_{j,k}=-(2\pi)^{-1}\log|f_{j,k}|$.  Let $b_{j,k}\in \R_{>0}$ be the value of $y_{j,k}/y_{j,1}$ at $p$. We will denote $y_{j,1}$ also by $y_j$. Let $h_{j,k}: (M_X^{\gp}/\cO_X^\times)_x\to \R$ be the homomorphism which sends $f_{j,k}$ to $1$ and sends all the other $f_{j',k'}$ to $0$. 

Then the sets 
$$\{s_0^{(\sum_{j,k} \alpha_{j,k}h_{j,k})}\;|\; \alpha_{j,k}\in \C, \Im(\alpha_{j,k})>0, \Im(\alpha_{j,1})/\Im(\alpha_{j+1,1})\geq c, |\Im(\alpha_{j,k})/\Im(\alpha_{j,1})-b_{j,k}|\leq c^{-1}\}$$ ($\alpha_{m+1,1}$ denotes $i$) given for $c\in \R_{>0}$ form a base of the filter in $\cT$ consisting of all subsets of the form $U\cap \cT_{\triv}$, where $U$ ranges over all neighborhood of $p$.

Hence the study of $H_{\cT}(s)$ when $s\in \cT_{\triv}$ converges to $p$ is the study of $\exp(\sum_{j,k} \alpha_{j,k}N_{j,k})F(s_0)$, where $N_{j,k}$ is the monodromy logarithm corresponding to $h_{j,k}$,  when 
 $\Im(\alpha_{j,k})>0, \Im(\alpha_{j,1})/\Im(\alpha_{j+1,1})\to \infty$, $ \Im(\alpha_{j,k})/\Im(\alpha_{j,1})\to b_{j,k}$.
\end{sbpara}

\begin{sbpara} Hence the $\SL(2)$-orbit theorems in \cite{Sc}, \cite{CKS1}, and \cite{KNU08} concerning the SL(2)-orbits associated to
$(N_1^{(c)}, \dots, N^{(c)}_m, F(s_0))$, 
where $c=(c_{j,k})$  is near to $b=(b_{j,k})$ ($c_{j,j}=1$) and $N_j^{(c)}=\sum_{k=1}^{r(j)} c_{j,k}N_{j,k}$, imply the following Propositions \ref{sW0}, \ref{sBS0}, \ref{twi0}, and \ref{twi1}.  In particular, Proposition \ref{sW0} follows from \cite{KNU08} Theorem 0.5 (1) and (2) and Proposition \ref{sBS0} follows from \cite{KNU08} Proposition 8.5.

 Note that
$\sum_{j,k} y_{j,k}N_{j,k}= \sum_{j=1}^m y_{j,1}(\sum_{k=1}^{r(j)} (y_{j,k}/y_{j,1})N_{j,k}).$ 
  Note also that $W^{(j)}$ is the relative monodromy filtration of $N^{(c)}_1+\dots +N^{(c)}_j$ for $1\leq j\leq m$ if $c$ is near to $b$ and we have the real analytic dependence of the invariants of $\SL(2)$-orbits on the parameter $c$ stated in \cite{KNU08} 10.8 (1) (in the pure Hodge  case, it is written in   the last five lines in (4.65) (ii) of \cite{CKS1}).
\end{sbpara}

\begin{sbprop}\label{sW0}   Let $p\in X_{[:]}$. Then for some open neighborhood $U$ of $p$ in $\cT_{[:]}$, 
 there is a unique  splitting $\spl_W$ of $W$ of $A^{\an,\log*}_{1, X, \cT} \otimes_\R H_{\cT,\R}$ over $A^{\an,\log*}_{1,X,\cT}$ ($\ref{XT2}$) on the inverse image $U^{\log}$ of $U$ in $X_{[:]}^{\log}$ satisfying the following condition$:$

For each $s\in U\cap \cT_{\triv}$, via the evaluation $A^{\an,\log*}_{1,X,\cT,s}\to A^{\an, \log*}_{1, \cT,s}= A^{\an}_{0,\cT,s}\to \R$ 
 (the last arrow is the evaluation $f\mapsto f(s)$), it gives the canonical splitting  ($\ref{canspl}$)  of $W$ of the mixed Hodge structure $H_{\cT}(s)$. 
\end{sbprop}

\begin{sbprop}\label{sBS0} Let  $p\in X_{[:]}$, and assume that $H$ is a PLH. Let $W^{(j)}$ ($1\leq j\leq m$) be as in $\ref{nearp}$. Then for some open neighborhood $U$ of $p$  in $\cT_{[:]}$, 
 there is a unique splitting of $W^{(j)}$ of $A^{\an,\log*}_{1, X, \cT} \otimes_\R H_{\cT,\R}$ over $A^{\an,\log*}_{1,X,\cT}$ on the inverse image $U^{\log}$ of $U$ in $X_{[:]}^{\log}$  satisfying the following condition$:$

For each $s\in U\cap \cT_{\triv}$, via the evaluation $A^{\an,\log*}_{1,X,\cT,s}\to A^{\an, \log*}_{1, \cT,s}= A^{\an}_{0,\cT,s}\to \R$, it gives the orthogonal  splitting of $W^{(j)}$ for the Hodge metric of $H_{\cT}(s)$.

\end{sbprop}

\begin{sbprop}\label{twi0} Let $p\in X_{[:]}$  and let $\tilde p$ be a point of $X_{[:]}^{\log}$ lying over $p$. Assume that $\gr^WH$ is polarized.

Let $s^{(0)}$ be the splitting of $W$  on $H_{\R, \tilde p}$ induced  from  the splitting of $W$ on $A^{\an,\log*}_{1, \tilde p} \otimes_\R H_\R$  in Proposition $\ref{sW0}$ by the evaluation  map $A^{\an,\log*}_{1,\tilde  p}\to \R $ ($\ref{logval}$). Then there is a unique $m$-ple $(s^{(1)}, \dots, s^{(m)})$ of a splitting $s^{(j)}$ of $W^{(j)}$ of $H_{\R,\tilde p}$ such that the $m+1$-ple $(s^{(0)}, s^{(1)}, \dots, s^{(m)})$ is compatible in the sense of $\ref{csp}$ 
 and such that for $1\leq j\leq m$ and $w \in \bZ$, the splitting induced by $s^{(j)}$ on $\gr^W_wH_{\R,\tilde p}$ is obtained by the evaluation map
$A^{\an,\log*}_{1,\tilde  p}\to \R$ ($\ref{logval}$) from the splitting of $W^{(j)}$ on $A^{\an,\log*}_{1, \tilde p} \otimes_\R \gr^W_w(H_\R)$ in Proposition $\ref{sBS0}$.

\end{sbprop}

\begin{sbpara}\label{GgK0}

Assume that $\gr^WH$ is polarized. Let $G$ be the algebraic subgroup of $\Aut_\R(H_{\R,\tilde p})$ consisting of all elements which preserve the weight filtration $W$ and $\langle\cdot,\cdot \rangle$ of $\gr^W$. (In the case where $H$ is  a PLH, $G=\Aut(H_{\R,\tilde p}, \langle\cdot,\cdot\rangle)$). 
\end{sbpara}

\begin{sbprop}\label{twi1}  Let $p$ and $\tilde p$  be as in Proposition $\ref{twi0}$. 
  For $1\leq j\leq m$, let  $\tau_j$ be the action of ${\bf G}_{m,\R}$ on $H_{\R,\tilde p}$ corresponding to  the splitting $s^{(j)}$ of the filtration $W^{(j)}$ on $H_{\R,\tilde p}$ in Proposition $\ref{twi0}$. 
Let $t(y)= \prod_{j=1}^m \tau_j(\sqrt{y_{j+1}/y_j})$. 
   Write $F=\exp(\sum z_{j,k}N_{j,k})F_0$, where $F_0$ is a filtration on $H_{\bC,\tilde p}$, $N_{j,k}$ is as in $\ref{handN}$, and $z_{j,k}$ is a branch of $(2\pi i)^{-1}\log(f_{j,k})$ with $f_{j,k}$ as in $\ref{m,f_jk}$.
Then 
  there is a unique $\br$ such that $(H_{\R,\tilde p}, W, \br)$ is a mixed Hodge structure with  polarizable $\gr^W$ and such that 
$$\exp({\textstyle\sum_{j,k}} iy_{j,k}N_{j,k})F_0= t(y) g_y \br$$
with $g_y\in G(A^{\an}_1)G_u(A^{\an}_{1,\C})\subset G(A^{\an}_{1,\C})$ satisfying $g_y(p)=1$. Here $G_u$ is the unipotent radical of $G$, that is, the group of all elements of $G$ which induce the identity map of $\gr^W$. 
If $H$ is a PLH of weight $w$, $(H_{\R,\tilde p},\langle\cdot,\cdot\rangle, \br)$ is a polarized Hodge structure of weight $w$. 

\end{sbprop}

\begin{sbpara}\label{GgK} Assume that $H$ is a PLH. 
  Let $K_{\br}$ be the algebraic subgroup of $G$ consisting of elements which preserve the Hodge metric of $\br$. 

Let 
$\frak g=\Lie(G)$, which we regard as the set of  linear maps $h:H_{\R,\tilde p}\to H_{\R,\tilde p}$  such that $\langle hu,v\rangle+\langle u, hv\rangle=0$. 

Define
$\frak g_{(A)}:= \sum_{\alpha\in \Z^m} A^{\an}_1\prod_{j=1}^m (y_j/y_{j+1})^{-|\alpha(j)|/2}\frak g_{\alpha}\subset A^{\an}_1\otimes_\R \frak g$.   
Here $\frak g_{\alpha}$ is the part of weight $\alpha\in \Z^m$ for the adjoint action of $(\tau_j)_j$.

\end{sbpara}

\begin{sbprop}\label{twi2}
 Let $p\in X_{[:]}$.  
Assume that $H$ is a PLH. Then we can take $g_y\in G(A^{\an}_1)$ in Proposition $\ref{twi1}$ such that $g_y=\exp(h_y)k_y$ with $k_y\in K_{\br}(A^{\an}_1)$, $h_y\in \frak g_{(A)}$ such that $k_y(p)=1$, $h_y(p)=0$. 

\end{sbprop}
  This is shown in  \cite{KU} Proposition 6.2.2 as a consequence of the SL(2)-orbit theorems in \cite{Sc} and \cite{CKS1}.

\begin{sbprop}\label{oldSL2} Let $p$ and $\tilde p$  be as in Proposition $\ref{twi0}$. Assume that $H$ is a PLH 
on $X$ of weight $w$. 

Let $D$ be the set of all decreasing filtrations $F'$ on $H_{\C,\tilde p}$ 
such that $(H_{\R,\tilde p}, \langle \cdot,\cdot\rangle, F')$ is a polarized Hodge structure of weight $w$. Then we have a pair $(\rho, \varphi)$ of a homomorphism
$$\rho: \SL(2)^m_\R \to \Aut_\R(H_{\R, \tilde p}, \langle\cdot,\cdot \rangle)$$ of algebraic groups over $\R$ characterized by the following properties {\rm (i)}, {\rm (ii)} and a holomorphic map
$$\varphi: \frak H^m \to D$$
characterized by the following properties {\rm (iii)} and {\rm (iv)}. 
Let the notation $\Delta(a)$  be as in Theorem $\ref{ourSL2}$. 

{\rm (i)}  $\rho(\Delta(a_1,\dots, a_m))= a_1^{-w/2}\prod_{j=1}^m \tau_j(a_j/a_{j+1})$, where $\tau_j$ is as in Proposition $\ref{twi1}$ and  $a_{m+1}$ denotes $1$. 

{\rm (ii)}
 $\Lie(\rho):{\frak {sl}}_2(\R)^m \to \End_\R(H_{\R,\tilde p})$ sends $\begin{pmatrix} 0&1\\0&0\end{pmatrix}$ in the $j$-th ${\frak {sl}}_2(\R)$ in ${\frak {sl}}_2(\R)^m$ to $\hat N_j$, where
$\hat N_j$ is 
the component of weight $0$ of $\sum_{k=1}^{r(j)} b_{j,k}N_{j,k}$ for the adjoint action of $(\tau_k)_{1\leq k\leq j-1}$.  Here $b_{j,k}$ and $N_{j,k}$ are  as in $\ref{handN}$.

{\rm (iii)}  $\varphi(g\alpha)= \rho(g)\varphi(\alpha)$ for all $g\in \SL(2,\R)^m$ and $\alpha\in \frak H^m$. 

{\rm (iv)}  $\varphi(i, \dots, i)$ is $\br$ in Proposition $\ref{twi1}$. 

\end{sbprop}
This is by the SL(2)-orbit theorem of \cite{Sc}, \cite{CKS1}. 

\subsection{Proofs of the results in Section \ref{ss:HA}}\label{ss:pfSL}

We prove the results in Section \ref{ss:HA}. Let $X$ be an object of $\cB(\log)$ and let $H$ be an LMH on $X$ with polarizable $\gr^W$. 

\begin{sbpara}\label{tortwist} Let $p\in X_{[:]}$ and let $\tilde p$ be a point of $X^{\log}_{[:]}$ lying over $p$. 
  Write the Hodge filtration $F$ of $\cO^{\log}_{X, \tilde p}\otimes_\R H_{\R,\tilde p}$ as $F=\exp(\sum_{j,k} z_{j,k}N_{j,k})F_0$, where 
$F_0$ is a splitting decreasing filtration of $\cO_{X,p}\otimes_\R H_{\R,\tilde p}$ (\ref{xi1}), $N_{j,k}$ is as in $\ref{handN}$, and $z_{j,k}$ is a branch of $(2\pi i)^{-1}\log(f_{j,k})$ with $f_{j,k}$ as in $\ref{m,f_jk}$ at $p$.
  Let $z_{j,k}=x_{j,k}+iy_{j,k}$ with $x_{j,k}$ and $y_{j,k}$ being real.
  Let $R=A^{\an}_{2+}$. By \ref{tortwi} (2),  as a pair of 
an $R_p$-module $M$ and a decreasing  filtration on $A^{\an}_{(3),p, \C}\otimes_{R_p}  M$, we have
\begin{align*}
(H^{\natural}_{R,p}, F) & \cong (R_p\otimes_\R H_{\R,\tilde p},  t(y)^{-1}\exp(-{\textstyle\sum}_{j,k} x_{j,k}N_{j,k})F)\\
&= (R_p\otimes_\R H_{\R,\tilde p},  t(y)^{-1}\exp({\textstyle\sum}_{j,k} iy_{j,k}N_{j,k})F_0).
\end{align*}

\end{sbpara}

\begin{sbpara}\label{pfHA1(2)}

We prove Theorem \ref{HA1} (1) except the last part which says that $H^{\natural}_A(p)$ with the Hodge filtration is a mixed Hodge structure (this part will be  proved in \ref{HAp}). 

We have the following (i) and (ii). 

(i) By $A^{\an}_{3,\C} \otimes_{\cO_X} H_{\cO} \cong A^{\an}_{3, \C}\otimes_R H^{\natural}_R$ (Proposition \ref{HA10} (3)), 
the Hodge filtration of $H$ gives a locally splitting filtration on $A^{\an}_{3, \C}\otimes_R H^{\natural}_R$. 

(ii) If $X$ is an fs log point, it comes from a locally splitting filtration on $H^{\natural}_R$. 
In fact, by Proposition \ref{twi1}, near $p$, we have  $$(A^{\an}_1\otimes_\R H_\R, t(y)^{-1}\exp({\textstyle\sum}_{j,k} iy_{j,k}N_{j,k})F_0)=(A^{\an}_1\otimes_\R H_\R, g_y\br)$$
with $\br$ being as in Proposition \ref{twi1} and with $g_y\in G(A^{\an}_1)G_u(A_{1,\bC}^{\an})$.

Let $p\in X_{[:]}$ and let $\tilde p\in X_{[:]}^{\log}$ lying over $p$. For the associated SL(2)-orbit at $p$ (Proposition \ref{oldSL2}), we have $H_{\C,\tilde p}= \bigoplus_{a\in \N^m} E^{(a)}$, where $E^{(a)}$ is the part on which the representation of $\SL(2)^m$ is through $\bigotimes_{j=1}^m \text{Sym}^{a(j)}(\rho_j)$, where $\rho_j$ is the $j$-th projection $\SL(2)^m\to \SL(2)$. Let $E^{(a,0)}$ be the part of $E^{(a)}$ on which all $\hat N_j$ are zero. Then $E^{(a)}=\bigoplus_b E^{(a,b)}$, where $b$ ranges over all elements of $\N^m$ such that
$b(j)\leq a(j)$ for all $j$ and $E^{(a,b)}=({\hat N}^+_1)^{b(1)}\cdots ({\hat N}^+_m)^{b(m)}E^{(a,0)}$ with 
${\hat N}^+_j$ being the Lie action of the  $j$-th $\begin{pmatrix} 0& 0\\1 & 0\end{pmatrix}$ in ${\frak {sl}}_2(\R)^m$. Then $E^{(a,b)}$ has a structure of  a pure Hodge structure of weight $w-\sum_j a(j)+2\sum_j b(j)$.
Let $E^{(a,b), (p,q)}$ ($p+q=w-\sum_j a(j)+2\sum_j b(j))$ be the Hodge $(p,q)$-component of $E^{(a,b)}_\C$. For $\alpha\in \Z$, let $E_{\alpha}$ be the direct sum of all $E^{(a,b),(p,q)}$ such that $p<\alpha$.

Let $V$ be the affine scheme over $\C$ such that for every scheme $T$ over $\C$, $V(T)$ is the set of all decreasing filtrations $(\cF^{\alpha})_{\alpha\in \Z}$ on $\cO_T \otimes_{\C} H_{\C, \tilde p}$ such that $\cF^{\alpha}\overset{\cong}\to \cO_T \otimes_{\C} H_{\C, \tilde p}/E_{\alpha}$ for all $\alpha\in \Z$.

By (i) and (ii), $t(y)^{-1}\exp(\sum_{j,k} iy_{j,k}N_{j,k})F_0$ gives a point of the fiber product
$$V(A^{\an}_{3,\C,p})\to V(A^{\an}_{3,\pts,\C,p})\leftarrow V(A^{\an}_{1, \pts,\C,p}),$$ 
where $A^{\an}_{k,\pts}$ $(k=1,3)$ are as in \ref{pts}.
Since $A^{\an}_{2+,\C,p}$ is the fiber product of $A^{\an}_{3,\C,p}\to A^{\an}_{3, \pts,\C,p}\leftarrow A^{\an}_{1, \pts,\C,p}$ and since $V$ is affine, the map from $V(A^{\an}_{2+,\C,p})$ to the above fiber product 
is bijective.  Hence 
$t(y)^{-1}\exp(\sum_{j,k} iy_{j,k}N_{j,k})F_0$ belongs to $V(A^{\an}_{2+,\C,p})$ and hence
comes from a splitting filtration on
 $A^{\an}_{2+,\C,p}\otimes_{\C} H_{\C, \tilde p}$ over  $A^{\an}_{2+,\C,p}$.

 \end{sbpara}
 
\begin{sbpara}\label{HAp}
The fiber $H^{\natural}_A(p)$ of $(H^{\natural}_A, F)$ at $p$ is isomorphic to $(H_{\R, \tilde p}, \br)$ in Proposition \ref{twi1} by \ref{tortwist} and \ref{pfHA1(2)}. Hence by Proposition \ref{twi1}, $H^{\natural}_A(p)$ is a mixed Hodge structure with polarizable $\gr^W$, and if $H$ is a PLH of weight $w$, 
${}^{\prime}H_A^{\natural}(p)$ is a polarized Hodge structure of weight $w$. 

\end{sbpara}

\begin{sbpara}\label{HA1(3)} We prove Theorem \ref{HA1} (2). We may assume that $H$ is a PLH of weight $w$. Let $F$ be the Hodge filtration on $H_{A,\C}$. It is sufficient to prove that $F^r\oplus \overline{F^{w+1-r}}\to H_{A,\C}$ is an isomorphism for every $r$. 
  Let $p \in X_{[:]}$. 
  Because $A^{\an}_{(2+),p}$ is a local ring, 
the bijectivity of  the homomorphism $F^r_p\oplus \overline{F^{w+1-r}_p}\to H_{A,\C,p}$ of free $A^{\an}_{(2+),p}$-modules of finite rank is reduced to the fact that the fiber $H^{\natural}_A(p)$ at $p$ is a Hodge structure of weight $w$ (\ref{HAp}). 
\end{sbpara}

\begin{sbpara} Assume that $H$ is a PLH of weight $w$. We prove the  part of Theorem \ref{HA1} (3) that we have a perfect pairing $\langle \cdot, \cdot\rangle: {^{\prime}}H^{\natural}_A\times {^{\prime}}H^{\natural}_A\to A^{\an}_{(2+)}$. By \ref{tortwi},  this is reduced to the perfectness of the pairing $\langle\cdot, \cdot \rangle: H_{\R, \tilde p} \times H_{\R, \tilde p} \to \R$.   

\end{sbpara}

\begin{sbpara} Assume that $H$ is a PLH. We prove the following part of Theorem \ref{HA1} (3): the perfectness of the  Hodge metric of  ${}'H_{A,\C}^{\natural}$. 
  Because $A^{\an}_{(2+),p}$ is a local ring, this is reduced to the fact that ${}^{\prime}H^{\natural}_A(p)$ is a polarized Hodge structure.

\end{sbpara}

\begin{sblem}\label{(A)metric}  

Let $S$ be  an $A^{\an}_{(2+)}$-submodule of  ${}^{\prime}H_A^{\natural}$ which is locally a direct summand. Then the map $S\to S^*\;;\; a \mapsto (b \mapsto (a,b))$ is an isomorphism. Here 
 $S^*$ is the sheaf of  $A^{\an}_{(2+)}$-homomorphisms $h: S\to A^{\an}_{(2+)}$ 
 and $(\cdot,\cdot)$ is the Hodge metric.
\end{sblem}

\begin{pf} By \ref{HAp}, 
the fiber of the  map in problem at $p$ is an isomorphism by the positivity of the Hodge metric of ${}^{\prime}H^{\natural}_A(p)$. 
 Because $A^{\an}_{(2+), p}$ is a local ring and the stalks $S_p$ and $(S^*)_p$ are free $A^{\an}_{(2+),p}$-modules of finite rank, this map is an isomorphism. 
\end{pf}

 To prove Theorem \ref{SL2I}, we give the following preparations Proposition \ref{prop1} and Proposition \ref{prop2}. 

\begin{sbprop}\label{prop1} Let $V$ 
be a finite dimensional $\R$-vector space. Assume that we are  given  metrics $\mu_1$, $\mu_2$  on $V$ (a metric here is a metric induced by a positive definite symmetric form). Let $\Phi$ be a finite set of increasing filtrations $\fil_{\bullet}$ (such that $\fil_j=0$ if $j\ll 0$ and $\fil_j=V$ if $j\gg 0$) on $V$ and let $G=\Aut(V)$, $G_{\Phi}=\{g\in G\;|\; gW=W\;\text{for all} \; W\in \Phi\}$.  Assume that for each $j=1,2$, $\mu_j$ gives compatible splittings 
of $\Phi$. Then there is 
$g\in G_{\Phi}$ 
such that $\mu_2$ coincides with  $g\mu_1:=\mu_1\circ g^{-1}$.

\end{sbprop}  
 \begin{pf} For $j=1,2$, let $t_j:{\bold G}_m^{\Phi}\to G_\Phi$ be the compatible splitting given by $\mu_j$. Then  the map ${\bold G}_m^{\Phi}\to G_\Phi/G_{\Phi, u}$ induced by $t_1$  and that by $t_2$ coincide. 
Take  a parabolic subgroup $P$ such that $G^{\circ}_{\Phi}\subset P$, $G_{\Phi,u} \subset P_u$ and such that the image of the composite map  ${\bold G}_m^{\Phi}\to G_{\Phi}^{\circ}/G_{\Phi,u}\to P/P_u$ is central. (To see the existence of $P$, write $\Phi=\{\phi_1, \dots, \phi_n\}$ and let $P$ be the group of all automorphisms of $V$ which preserve $V_{\leq a}$ for all $a\in \Z^n$, where $V_{\leq a}$ denotes the sum of the $b$-component of $V$ for all $b\in \Z^n$
 such that $b\leq a$ for the lexicographic order of $\Z^n$ with respect to the splitting of $\Phi$ by $\mu_1$.) Let $K_{\mu_j}= \{g\in G\;|\; g\mu_j=\mu_j\}$ $(j=1,2)$. 
By  $G=PK_{\mu_1}$, we have $K_{\mu_2}= \text{Int}(g)K_{\mu_1}$ and $\mu_2=g\mu_1$ for some $g\in P$.
Since the splitting $t_j:{\bf G}_m^{\Phi}\to P$ is the Borel--Serre lifting  of ${\bold G}_m^{\Phi}\to P/P_u$ at $K_{\mu_j}$, we have $t_2= \text{Int}(g)(t_1)$. On the other hand,  since $t_1$ and $t_2$ are compatible splittings of $\Phi$, $t_2=\text{Int}(h)(t_1)$
for some $h\in G_{\Phi, u}$. 
Hence $\text{Int}(h^{-1}g)(t_1)=t_1$.
Therefore $h^{-1}g\in G_{\Phi}$ and so $g\in G_{\Phi}$. 
\end{pf}

\begin{sbprop}\label{prop2}
Let $V$ be a finite dimensional $\R$-vector space. Assume that we are given a metric $\mu$ on $V$, a finite set $\Phi'$ of increasing filtrations on $V$ which has a compatible splitting, and a subset $\Phi$ of $\Phi'$. Assume that $\mu$ gives compatible splittings of $\Phi$.  Then there is $g\in \Aut(V)$ such that $g$ keeps each element of $\Phi$ and such that $\mu$ gives compatible splittings of $g\Phi'$. (Note that $g\Phi=\Phi$.)
\end{sbprop} 

\begin{pf} In Proposition \ref{prop1}, take this $\mu$ as $\mu_2$. Take as $\mu_1$ a metric which gives a  compatible splitting of $\Phi'$. (To see the existence of $\mu_1$, use the direct sum decomposition of $V$ with index set $\Z^{\Phi'}$. Then take the direct sum of the metric of each component.) By Proposition \ref{prop1}, there is $g\in G_{\Phi}$ such that $\mu_2= g\mu_1$. Since $\mu_1$ gives a  compatible  splitting 
of $\Phi'$,  $\mu_2$ gives a compatible  splitting 
 $g\Phi'$. \end{pf}

\begin{sbpara} We reduce Theorem \ref{SL2I} to the following Proposition \ref{SL2Ib}.
By Lemma \ref{(A)metric}, we have an orthogonal splitting of $W^{(j)}$ on $H^{\natural}_{A,p}$. Consider the situation in \ref{nearx} and \ref{nearp}. There is an open neighborhood $V$ of $p$ in $X_{[:]}$ such that
 all elements of $\Phi(p)$  extend as local systems on the inverse image $V^{\log}$  of $V$  in $X_{[:]}^{\log}$ and such that at each point of $V^{\log}$  lying over each $p'\in V$, $\Phi(p)$ on $H_{\R,p'}$ has a compatible splitting and  $\Phi(p')\subset \Phi(p)$.  The orthogonal splitting of $W^{(j)}$ of ${}^{\prime}H^{\natural}_A$ extends to an open neighborhood of $p$. 
Hence Theorem \ref{SL2I} is reduced to 
Proposition \ref{SL2Ib} by using $W^{(j)}$ as $W'$ in Proposition \ref{SL2Ib}. 

\end{sbpara}

\begin{sbprop}\label{SL2Ib} Assume that $H$ is a PLH. 
Let $U$ be an open set of $X_{[:]}$ and let $W'$ be a locally constant increasing filtration on the restriction of $H_\R$ to the inverse image $U^{\log}$ of $U$ in $X^{\log}_{[:]}$. Assume that for every $p\in U$ and every $\tilde p$ lying over $p$, 
 the set of increasing filtrations on $H_{\R,\tilde p}$ consisting of  $W'$  and all elements of $\Phi(p)$  has a
 compatible splitting in the  sense of $\ref{csp}$. Assume furthermore that on $U$, $W'$ on $A_{(3)}^{\an, \log*} \otimes_{\R} H_{\R}$  comes from a locally splitting filtration $W'$ on ${}^{\prime}H^{\natural}_A$ and the last filtration has an orthogonal splitting for the Hodge metric. 
Then 
 there is a unique orthogonal  direct sum decomposition of  $H_A$ over $A^{\an}_{(2+)}$ on $U$ for the Hodge metric which splits $W'$ on $H_A$.

\end{sbprop}

 \begin{pf}
 We have an orthogonal splitting $\spl^{\BS}_{W'}(F)$ of $W'$ on $A^{\an,\log}_{(3)}\otimes_{\R} H_\R$ for the Hodge filtration $F$ of $H$. It is sufficient to  prove that this splitting comes from an orthogonal splitting of $W'$ of $A^{\an,\log*}_{(2+)}\otimes_\R H_\R$ because then the desired orthogonal splitting of $H_A$ is obtained from it by taking the direct image under $U^{\log}\to U$.  
  For this, by Lemma \ref{()pt}, it  is sufficient to prove that  in the case where $X$ is an fs log point, 
 it comes from an orthogonal splitting of $W'$ of 
$A_1^{\an,\log}\otimes_\R H_\R$. 

This is shown as follows.  
  Assume that $X$ is an fs log point. 
  Write $F=\exp({\textstyle\sum}_{j,k} (x_{j,k}+iy_{j,k})N_{j,k})F_0$ 
as in \ref{tortwist}. 
  It is sufficient to prove that the twist 
  $\spl_{W'}^{\BS}(\exp({\textstyle\sum}_{j,k} iy_{j,k}N_{j,k})F_0)$ of $\spl_{W'}(F)$ 
  by $\exp(-{\textstyle\sum}_{j,k} x_{j,k}N_{j,k})$ comes from an orthogonal splitting of $W'$ of 
$A_1^{\an,\log}\otimes_\R H_\R$.  
  By Proposition \ref{prop2}, 
  for some $g\in \Aut_\R(H_{\R, \tilde p})$ such that $gW^{(j)}=W^{(j)}$ for all $1\leq j\leq m$, the Hodge metric of $\br$ gives a compatible splitting of $\{gW'\}\cup\{W^{(j)}\;|\;1\leq j\leq m\}$.
  Let the notation be as in Proposition \ref{twi1}.
We have 
$$\spl^{\BS}_{W'}(\exp({\textstyle\sum}_{j,k} iy_{j,k}N_{j,k})F_0)=\spl^{\BS}_{W'}(t(y)a_yk_y\br)$$ 
with $a_y\in \exp(\frak g_{(A)})$, $k_y\in K_{\br}(A^{\an}_1)$ as in Proposition \ref{twi2}. 
 Here $\frak g_{(A)}$ is as in \ref{GgK}.
  Let 
$$\frak p_{(A)}:= {\textstyle\sum}_{\alpha\in \Z^m} A^{\an}_1{\textstyle\prod}_{j=1}^m (y_j/y_{j+1})^{-|\alpha(j)|/2}\frak p_{\alpha}\subset A^{\an}_1\otimes_\R \frak p,$$
where $\frak p=\Lie(G_{gW'})$.
Here $\frak p_{\alpha}$ is the part of weight $\alpha\in \Z^m$  of $\frak p$ for the adjoint action of $(\tau_j)_j$.  
  Then we have 
\begin{align*}
\spl^{\BS}_{W'}(t(y)a_yk_y\br)&=g^{-1}\spl^{\BS}_{gW'}(gt(y)a_yk_y\br)= g^{-1}\spl^{\BS}_{gW'}(t(y)g_{1,y}a_yk_y\br)\\
&= g^{-1}\spl^{\BS}_{gW'}(t(y) g_{2,y} k'_y\br) 
 = g^{-1}\spl^{\BS}_{gW'}(g_{3,y}t(y)k'_y\br)\\
&= g^{-1}g_{3,y}\spl^{\BS}_{gW'}(t(y)k'_y\br)= g^{-1}g_{3,y}\spl^{\BS}_{gW'}(\br),
\end{align*}
where  $g_{1,y}\in \exp(\frak g_{(A)})$, 
 $g_{2,y}\in\exp(\frak p_{(A)})$, 
 $k'_y\in K_{\br}(A^{\an}_1)$, and 
 $g_{3,y}\in \exp(A^{\an}_1\otimes_\R \frak p)$.
  This implies the desired statement. 
 
In the above, we use the following three inclusions: 

$\exp(\frak g_{(A)})\subset \exp(\frak p_{(A)})K_{\br}(A^{\an}_1)$ by that $\frak p$ is parabolic;

$\exp(\frak g_{(A)})\exp(\frak g_{(A)})\subset \exp(\frak g_{(A)})$;

$t(y) \frak p_{(A)} t(y)^{-1}\subset A^{\an}_1\otimes \frak p$. 

  In fact, we use the first two inclusions for the third equality, 
and 
the last one for the fourth equality.
 \end{pf}

The next proposition follows from Proposition 4.2 of \cite{KNU08}.

\begin{sbprop}\label{splWBS} Let $X$ be an fs log point and let $H$ be an $\R$-LMH on $X$ with polarizable $\gr^W$.
 Let $S$ be the standard log point. 

$(1)$ There is a unique splitting of $W$ of $A^{\an}_{1,X\times S}\otimes_\R H_\R$ such that for any $H'$ and $H_\R\to H'_\R$ as in Proposition $\ref{mixpure}$, it is an orthogonal splitting for the Hodge metric on $A^{\an,\log}_{(3), X\times S} \otimes_\R H'_\R$. 

$(2)$ Assume that the log structure of $X$ is trivial. Then  
the canonical splitting ($\ref{canspl}$) of $W$ is obtained from the splitting in the above $(1)$ by the evaluation  $A^{\an}_{1, S,p}\to \R\;;\;f\mapsto f(p)$ at the unique point $p$ of $(X\times S)_{[:]}=S_{[:]}=S$.
\end{sbprop}

The following Proposition \ref{fromLM} is \cite{logmot} Proposition 2.2.

\begin{sbprop}\label{fromLM}  Let $X$ be an object of $\cB(\log)$, let $R$ be a subfield of $\R$, and let $H$ be an $R$-LMH  on $X$.
Then locally on $X$, there are a log manifold $Z$ and a morphism $X\to Z$ of $\cB(\log)$ 
 such that $H$ is the pullback of an $R$-LMH on $Z$ with polarizable $\gr^W$.

\end{sbprop}

\begin{sbpara}\label{pfII1} We start the proof of Theorem \ref{SL2II}. 

  We use Proposition \ref{mixpure}. Let $X$ be an object of $\cB(\log)$ and let $H$ be an $\R$-LMH on $X$ with polarizable $\gr^W$. Let $S$ be the standard log point and let $q$ be the standard generator of the log structure of $S$. 
  Then a pure $H'\supset H$ appears on $X\times S$ locally on $X$. The Hodge metric of ${}^{\prime}(H'){}^{\natural}_A$ belongs to $A^{\an}_{(2+), X\times S}$ by Theorem \ref{HA1} (3) and it induces an orthogonal splitting of $W$ on ${}^{\prime}(H'){}^{\natural}_A$ over $A^{\an}_{(2+), X \times S}$ by Theorem \ref{SL2I} (1). Let $J$ be the ideal of $A^{\an}_{(2+), X \times S}$ generated by $(\log|f|/\log|q|)^{1/2}$ for  $f\in M_X$, $|f|<1$. Then $A^{\an}_{(2+), X\times S}/J\overset{\cong}\leftarrow A^{\an}_{(2+), X}$. Hence we get a splitting of $W$ on  $H^{\natural}_A$. Hence we have a splitting of $W$ on $A^{\an,\log}_{(3), X}\otimes_\R H_\R$. 

 By taking the pullbacks to points of $X$ and reducing to the case where $X$ is an fs log point, we see by Proposition \ref{splWBS} (1) that this splitting $\spl_W$ of $W$ on $A^{\an,\log}_{(3)}\otimes_\R H_\R$ comes from a splitting of $W$ on $A^{\an,\log*}_{(2+)}\otimes_\R H_\R$. A problem is to prove that this $\spl_W$ is independent of the choice of $H'\supset H$. Once this problem is solved (as in \ref{pfII2}--\ref{pfII3} below), we have that this $\spl_W$ satisfies the required conditions (use Proposition \ref{splWBS} (2)).

\end{sbpara}

\begin{sbpara}\label{pfII2}
 In the case where $X$ is a log manifold, by restricting to $X_{\triv}$ on which the splitting is the classical splitting for a usual mixed Hodge structure by Proposition \ref{splWBS} (2), we have that the  $\spl_W$ obtained by \ref{pfII1} is independent of  the choice. 

\end{sbpara}

\begin{sbpara}\label{pfII3}
 Assume that we are given two pure $H'_j\supset H$ $(j=1,2)$ as in Proposition \ref{mixpure}. Locally on $X$, there is a morphism $X\to Z$ as in Proposition \ref{fromLM} via which $H$ comes from an $\R$-LMH $\tilde H$ with polarizable $\gr^W$ on $Z$ and via which times $S$, 
$H'_1 \supset H \subset H'_2$ comes from 
$\tilde H'_1 \supset \tilde H \subset \tilde H'_2$ on $Z \times S$.
  This is proved similarly to Proposition \ref{fromLM}. 
  (In the proof of Proposition \ref{fromLM} in \cite{logmot}, we consider the period map $X \to E$, where $E$ is a moduli space of Hodge filtration with Hodge numbers of $H$. 
  Here we consider the period maps $X\to E_j$ $(j=1,2)$, where $E_j$ is the corresponding one for $H'_j$, and the fiber product $E_1 \times_EE_2$.) 
  Then the splitting of $W$ of $H$ induced by $H'_j$ is the pullback of the splitting of $W$ of $\tilde H$ induced by $\tilde H'_j$, which is independent of $j$ 
by \ref{pfII2}. 
  Hence the splitting of $W$ of $H$ is independent of $j$.

\end{sbpara}

\begin{sbpara} By \ref{pfII3}, the part of Theorem \ref{SL2II} on $W$ of $H_A$ was proved. We prove the part on $W$ of $H_A^{\natural}$. By the part on $H_A$, we have a canonical splitting of $W$ of 
 $A^{\an}_{(3)} \otimes_{A^{\an}_{(2+)}} H_A= H^{\natural}_{A^{\an}_{(3)}}$. It is sufficient to prove that this splitting comes from a splitting of $W$ of $H^{\natural}_A$. This is reduced to the case where $X$ is an fs log point. By \ref{tortwist} and \ref{pfHA1(2)} (ii), we are reduced to the fact that the canonical splitting of $W$ of the usual mixed Hodge structure $g_y\br$ is a real analytic function in $g_y$ and hence a real analytic function.

\end{sbpara}

\begin{sbpara}\label{phi(i)} We prove Theorems \ref{ourSL2} and \ref{oldSL20}.

Theorem \ref{oldSL20} (1) follows from Proposition \ref{twi0}. 

The isomorphism $H_{\R,\tilde p} \cong H^{\natural}_A(p)$  in Theorem \ref{oldSL20} (2) is defined as is written there by using 
Theorem \ref{oldSL20} (1) and Lemma \ref{twigr}. Through this isomorphism, Theorem \ref{ourSL2} is reduced to Proposition \ref{oldSL2}. Then Theorem \ref{oldSL20} (2)  follows from Theorem \ref{ourSL2}. 

\end{sbpara}

\subsection{Complements to Section \ref{ss:HA}}\label{ss:comple}

We discuss whether our formulations of log real analytic functions are the best ones for the results in Section \ref{ss:HA}. In \ref{HA(2)}, we give a relation of $H^{\natural}_A$ and log $C^{\infty}$ functions in the case where $X$ is  ideally log smooth and reduced. 

\begin{sbpara} We can formulate the smaller sheaves $B^{\an}_k\subset A^{\an}_k$%
, $B^{\an,\log}_k\subset A^{\an,\log}_k$
,  and $B^{\an,\log*}_k\subset A^{\an,\log*}_k$ ($k=1,2,2+,3, (2), (2+), (3)$)  by replacing $(\log|f|/\log|g|)^{1/2}$ in \ref{A_1} by $\log|f|/\log|g|$ and modifying the subsequent arguments in the natural way. For example, $B_1^{\an}$ of the standard log point is the ring of convergent series in $y^{-1}$, not in $y^{-1/2}$. 

It is probable that we can replace $H_A$ (not ${}^{\prime}H^{\natural}_A$ or $H^{\natural}_A$)
in  Theorem \ref{SL2I} and in Theorem \ref{SL2II} 
 by its $B^{\an}_{(2+)}$-version. By \cite{KNU08}, this is true if $X$ is an fs log point. 

However, we do not pursue this subject in this paper.

\end{sbpara}

\begin{sbpara}\label{ExA(2+)} In the following Example 1 (resp.\ 2), we show that we do not have 
 (3) (resp.\ (1)) 
of Theorem \ref{HA1}  if we replace ${}^{\prime}H^{\natural}_A$ 
 (resp.\ $H^{\natural}_A$)  
by  ${}^{\prime}H^{\natural}_{A^{\an}_{2+}}$ (resp.\ $H^{\natural}_{A^{\an}_{(2)}}$) even for $X=\Delta$ with the log structure given by the coordinate function $q$ (resp.\ for a log manifold $X$ or for an fs log analytic space $X$).

Example 1. $X=\Delta$ with the log structure given by the coordinate function $q$. 

Example 2. Let $X$ be as in (1) or (2) of \ref{Ex(1)(2)}. 

 We give an example of PLH of rank $5$ and of Hodge type $(4,0)+(3,1)+(2,2)+(1,3)+(0,4)$ on $X$.

In both examples, $$F= \exp(zN)F_0,$$
where $F_0$ is as follows, $z= (2\pi i)^{-1}\log(q)$ (resp.\  $z=z_1=(2\pi i)^{-1}\log(q_1)$) in Example 1 (resp.\ 2), and $N$ is the monodromy logarithm coming from $q$ (resp.\ $q_1$) in Example 1 (resp.\ 2). We have an $\R$-base $e_j$ ($1\leq j\leq 5$) with the symmetric $\R$-bilinear form $\langle\cdot,\cdot \rangle $ defined by $\langle e_j, e_k\rangle=0$ unless $j+k= 6$, $\langle e_1, e_5\rangle = \langle e_3, e_3\rangle = 1$, $\langle e_2,e_4\rangle=-1$.
  The action of $N$ is given by 
$$N e_j=e_{j-3}\quad \text{for $1\leq j\leq 5$},$$ where $e_j$ for $j=0,-1,-2$ denote $0$. Let $h=q$ in Example 1, and let $h=q_2$ in Example 2. 
Let $\alpha= e_5-ie_4$, $\beta= e_2-ie_1+he_3$, $\gamma= e_3+he_4$, $\delta= e_5+ie_4$.
  The Hodge filtration $F_0$ is given by
$$0=F_0^5\subsetneq  F_0^4=\C \cdot \alpha \subsetneq  F_0^3= F_0^4+ \C \cdot \beta $$ 
$$\subsetneq F_0^2= F_0^3+ \C \cdot \gamma \subsetneq F_0^1= F_0^2+\C\cdot \delta\subsetneq F_0^0= {\textstyle\bigoplus}_{j=1}^5 \C e_j.$$
Let $R=A^{\an}_{2+}$ in Example 1 and $R=A^{\an}_{(2)}$ in Example 2. Then  $e_{j,R}:=y^{1/2}e_j$ (resp.\ $e_{j,R}:=y^{1/2}_1e_j$) for $j=1,2$ and $e_{3,R}:=e_3$ form a part of an $R$-base of ${}^{\prime}H_R^{\natural}$.

In Example 1,  the Hodge decomposition  comes from  ${}^{\prime}H_R^{\natural}$. The $(3,1)$-Hodge component of ${}^{\prime}H^{\natural}_{R,\C}$ is generated by $\beta':=y^{1/2}\beta+ 2iy^{-1/2}(\alpha+ ze_2)$, and $(\beta', \beta')= 1-yq\overline{q}$ which is invertible in $A^{\an}_{(2+)}$ but not invertible in $R$. 

In Example 2, the Hodge filtration does not come from a locally splitting filtration on  ${}^{\prime}H^{\natural}_{R,\C}$ because $y^{1/2}\beta=e_{2,R}-ie_{1,R}+ y_1^{1/2}q_2e_{3,R}$ and $y_1^{1/2}q_2$ does not belong to $R$.

\end{sbpara}

\begin{sbrem} In the above examples in \ref{ExA(2+)}, the PLH does not satisfy the big Griffiths transversality.
For LVPH (not for PLH), the authors are not sure 

(1) whether (3) of Theorem \ref{HA1} does not hold when we use ${}^{\prime}H^{\natural}_{A^{\an}_{2+}}$ instead of ${}^{\prime}H^{\natural}_A$ even for $X=\Delta$. 

(2) whether (1) of Theorem \ref{HA1} does not hold when we use $H^{\natural}_{A^{\an}_{(2)}}$ instead of $H^{\natural}_A$   even for a log manifold $X$ or for an fs log analytic space $X$.

 \end{sbrem}

\begin{sbpara}\label{HA(2)}  
In the case where  $X$ is ideally log smooth and reduced and $H$ is a PLH, the statements of Theorem \ref{HA1} (1) and (2)  for $H^{\natural}_A=H^{\natural}_{A^{\an}_{(2+)}}$ remains true if we replace $H^{\natural}_A$ by $H^{\natural}_{A^{\an}_{(2)}}$ and Theorem \ref{HA1} (3) remains  true if we replace
${}^{\prime}H^{\natural}_A$ by ${}^{\prime}H^{\natural}_{A^{\an}_{(2)}}$.  

  In fact, by taking $E= \cE nd_{A^{\an}_{(2)}}(H^{\natural}_{A^{\an}_{(2)}})_{\C}$ in Theorem \ref{sharp=} (2), we see that the projection to the $(p,q)$-Hodge component in $E\otimes_{A^{\an}_{(2)}} A^{\an}_{(2+)}$ belongs to $E$. This proves Theorem \ref{HA1} (2) for $H^{\natural}_{A^{\an}_{(2)}}$. Theorem \ref{HA1} (1) for $H^{\natural}_{A^{\an}_{(2)}}$ and Theorem \ref{HA1} (3) for ${}^{\prime}H^{\natural}_{A^{\an}_{(2)}}$ follow from it.

Since we have a homomorphism $A^{\an}_{(2)}\to A_2$, 
for an ideally log smooth and reduced fs log analytic space $X$, Theorem \ref{HA1}  for a PLH $H$ remains true if we replace 
 $H^{\natural}_A$ (resp.\ 
${}^{\prime}H^{\natural}_A$) by the log $C^{\infty}$ version 
$H^{\natural}_{A_2}$ (resp.\ 
${}^{\prime}H^{\natural}_{A_2}$) of it. 
\end{sbpara}

\subsection{Relations to the rings in  \cite{CKS1} Section 5}\label{ss:CKS5}

\begin{sbpara}
An aim of this paper is to give a good ring of functions which is an enlargement of the ring of real analytic functions and with which we can treat degeneration of Hodge structure nicely. 

This subject is  already discussed in Section 5 of the paper \cite{CKS1} by Cattani, Kaplan and Schmid for smooth complex analytic space with a normal crossing divisor.  Their ring ${\cal R}^b_K$ is similar to our $A^{\an}_{(2)}$ and their 
ring $(\cO \otimes \overline{\cO}\otimes \cL)^b$ is similar to our $A^{\an}_2$. More precisely:
 
\end{sbpara}
\begin{sbpara}\label{CKSrel1} Let $X=\Delta^{n+m}$ endowed with the log structure generated by the first $n$ coordinate functions $q_1, \dots, q_n$. 
Fix a subset $I$ of $\{1,\dots, n\}$, let $K>1$,  and let ${\cal R}^b_K$ be the ring defined in (5.17) in \cite{CKS1}. Write $I=\{i_{\alpha}\}$ with $1\leq i_1< \dots < i_r=n$. Then we have
$${\cal R}^b_K\subset \varinjlim_{\varepsilon} A^{\an}_{(2), X,\C}(U_{\varepsilon}),$$
where for $\varepsilon>0$, $U_{\varepsilon}$ is the open set of $X_{[:]}$ consisting of all points $p$ satisfying the following conditions (i)--(iii). 
 (i) $|q_j|<\varepsilon$ at $p$ for all $1\leq j\leq n+m$.
(ii) $y_{i_{\alpha}}/y_{i_{\alpha+1}}>K$ for $1\leq \alpha\leq r-1$, $y_n>K$. (iii) $K^{-1}< y_j/y_{i_{\alpha}}<K$ for $1\leq \alpha\leq r$ and $i_{\alpha-1}<j <i_{\alpha}$ ($i_0$ means $0$). Here $y_j=- (2\pi)^{-1}\log |q_j|$. 

The spaces of ratios are not used in \cite{CKS1}, but some atmosphere of the space of ratios appears in \cite{CKS1}, for example, in the definition of the space $(U^n)^I_K$ in (5.6) in \cite{CKS1} which is used to define ${\cal R}^b_K$.

Proposition (5.19) of \cite{CKS1} shows that the degeneration of the Hodge metric of a variation of polarized Hodge structure $H'$ on $X_{\triv}=(\Delta\smallsetminus \{0\})^n \times \Delta^m$ at the boundary is treated nicely by the ring ${\cal R}^b_K$. By the above inclusion of the rings and by \ref{tortwist}, this proves Theorem \ref{HA1} (3) for the PLH $H$ on $X$ which is the unique extension of $H'$ to $X$. 
Similarly,  Proposition (5.10) of \cite{CKS1}
implies our Theorem \ref{HA1} (1) for this $X$ and for this $H$. 

We have
$$(\cO \otimes \overline{\cO}\otimes \cL)^b \subset \varinjlim_{\varepsilon} A^{\an}_{2,X, \C}(V_{\varepsilon}),$$
where the left-hand-side is the ring defined in (5.17) in \cite{CKS1} and for $\varepsilon>0$, $V_{\varepsilon}$ is the open set of  $X_{[:]}$ consisting of all points  $p$ satisfying the following conditions (i)--(iii). (i) $|q_j|<\varepsilon$ at $p$ for all $1\leq j\leq n+m$. (ii) The value of $y_{j+1}/y_j$ at $p$ is not $\infty$ for $1\leq j\leq n$, where $y_{n+1}$ means $1$. (iii) If $1\leq j\leq n$ and $j\notin I$, the  value of $y_{j+1}/y_j$ at $p$ is not $0$.

\end{sbpara}

\begin{sbpara}\label{differ} Assume $X=\Delta$ with the log structure given by the coordinate function $q$.  Then
any real analytic functions of $y^{-1/2}$ on $\Delta$ belongs to $A_2^{\an}$, but it belongs to $(\cO \otimes \overline{\cO}\otimes \cL)^b$ (resp.\ ${\cal R}^b_K$) if and only if it is a polynomial in $y^{-1/2}$ (resp.\ rational function in $y^{-1/2}$ which has a finite value at $0\in \Delta$).
 This is the essential point of the difference of their rings and our rings. 
 \end{sbpara}

\begin{sbpara} There are  several good points of our rings of log real analytic functions  compared with the rings in \cite{CKS1}  as follows. 

1. The rings $(\cO \otimes\overline{\cO} \otimes \cL)^b$  and  ${\cal R}_K^b$ are not canonical in the sense that they depend on the choices of local coordinate functions. For example, let $X$ be smooth  complex analytic space  of dimension one with the log structure at one point and identify $X$ locally with $\Delta$ by using a local coordinate function $q$ at that point. Then  these rings depend on the choice of $q$ because which function is a polynomial of $y^{-1/2}$, or a rational function of $y^{-1/2}$, depends on the choice of $q$ (see \ref{differ}). 
For this reason,  the rings $(\cO\otimes\overline{\cO}\otimes  \cL)^b$ and ${\cal R}^b$ cannot be globalized as sheaves on smooth complex analytic spaces endowed with a normal crossing divisor. 

2. It is nice to work on $X_{[:]}$ as we do in this paper. For example, the $\SL(2)$-orbit theorem is understood  geometrically there.

3. Our rings of log real analytic functions  work well for any objects of $\cB(\log)$, and so for spaces with singularities.  
\end{sbpara}

\subsection{Log real analyticity of the CKS map}\label{ss:CKS}

Here we obtain the log real analyticity of the CKS map, the map which connects  the space of nilpotent orbits and the space of $\SL(2)$-orbits (Theorem \ref{partR2}).

\begin{sbpara}\label{La0} We review the space $D_{\SL(2)}$ of $\SL(2)$-orbits.

Let $\La=(H_0, W, (\langle\cdot,\cdot\rangle_w)_w, (h^{p,q})_{p,q})$ be as in \cite{KNU3} 2.1.1. (Here $H_0$ is a finitely generated free $\Z$-module, $W$ is an increasing filtration on $H_0\otimes \Q$ called the weight filtration, $\langle\cdot, \cdot\rangle_w$ is a non-degenerate $(-1)^w$-symmetric bilinear form on  $\gr^W_w$, and $(h^{p,q})_{p,q}$ are non-negative integers given for $p, q\in \Z$.)

Let $D$ be the period domain for $\La$. 
Let $D_{\SL(2)}\supset D$ be the space of $\SL(2)$-orbits for $\La$ defined in \cite{KNU2}. Here $D$ is the set of decreasing filtrations $F$ on $H_{0,\C}$ such that $H=(H_0, W, F)$ is a mixed Hodge structure satisfying the following conditions (i) and (ii). (i) The dimension of the Hodge $(p,q)$-component of $\gr^W_{p+q}H$ is $h^{p,q}$ for every $p,q$. (ii) $\gr^W_wH$ is polarized by $\langle\cdot,\cdot\rangle_w$ for each $w$. On the other hand, $D_{\SL(2)}$ is the set of equivalent classes (for a certain equivalent relation) of triples $(\rho, \varphi, \br)$, where $\rho$ is a homomorphism of algebraic groups $\SL(2)^n_\R \to \prod_w \Aut(\gr^W_wH_{0,\R})$ over $\R$, $\varphi$ is a holomorphic map from $\frak H^n$ to $\prod_w D(\gr^W_w)$ for some $n\geq 0$ ($D(\gr^W_w)$ is the period domain for $\La_w= (\gr^W_wH_0, \langle\cdot, \cdot\rangle_w, (h^{p,q})_{p+q=w})$), and $\br$ is a decreasing filtration on $H_{0,\C}$ such that $(F(\gr^W_w))_w=\varphi(i, \dots, i)$, satisfying certain conditions.

As is explained in \cite{KNU2} 3.2, $D_{\SL(2)}$ has two 
structures of locally ringed spaces over $\R$, $D^I_{\SL(2)}$ and $D^{II}_{\SL(2)}$, and we have a morphism $D_{\SL(2)}^I\to D_{\SL(2)}^{II}$ underlain by the identity map of $D_{\SL(2)}$. 
 The structure sheaves of rings are called the sheaves of real analytic functions on $D_{\SL(2)}^I$ and on $D_{\SL(2)}^{II}$, respectively. In this section (Section \ref{ss:CKS}), we will consider $D_{\SL(2)}^I$ and often denote the locally ringed space $D_{\SL(2)}^I$ simply by $D_{\SL(2)}$.
\end{sbpara}

\begin{sbpara}\label{L*} Let $X$ be a log manifold. We define the sheaf of rings ${\hat A}_{(2+)}^{\an, \log*}$  on $X^{\log}_{[:]}$.
Let ${\hat A}_{(2+)}^{\an, \log*}$  be the sheaf of real analytic functions of finite numbers of local sections of $A^{\an,\log*}_{(2+)}$.

\end{sbpara}

\begin{sbpara}\label{La2}
  Let $\Lambda$ be as in \ref{La0}. 
Let $X$ be a log manifold. Let $H$ be 
a $\bZ$-LMH on $X$ with polarized $\gr^W$ such that the dimension of the $(p,q)$-Hodge component of $\gr^W_{p+q}H$ is equal to $h^{p,q}$ for every $p$ and $q$.

\end{sbpara}

\begin{sbpara}\label{La22}  Let $\Gamma$ be a neat subgroup of $\Aut_\Z(H_0, W, (\langle\cdot, \cdot\rangle_w)_w)$. Assume that $H$ is endowed with a $\Gamma$-level structure. 
  Here a {\it $\Gamma$-level structure} on $H$ means a global section of the sheaf $\Gamma\bs I$, where $I$ is the sheaf of isomorphisms from
$(H_\Z, W, (\langle \cdot, \cdot\rangle_w)_w)$ to $(H_0, W,  (\langle \cdot, \cdot\rangle_w)_w)$.

\end{sbpara}
\begin{sbthm}\label{L*thm} Let the assumptions be as in $\ref{La2}$ and $\ref{La22}$. Then the period map $X_{\triv}\to \Gamma\bs D$ associated to the restriction of $H$ to $X_{\triv}$  extends uniquely to a continuous map $\psi: X^{\log}_{[:]}\to \Gamma \bs D^I_{\SL(2)}$ having the following property$:$ For every open set $U$ of $\Gamma\bs D^I_{\SL(2)}$  and for every real analytic function $f$ on $U$, $f\circ \psi$ belongs to  ${\hat A}_{(2+),X}^{\an, \log*}(\psi^{-1}(U))$.
\end{sbthm}

This improves  \cite{KNU4} Theorem 6.3.1 in which $\psi$ is constructed only as a unique continuous extension to $X^{\log}_{[:]}$ of the period map on $X_{\triv}$, and gives a new proof of that theorem. The proof of Theorem \ref{L*thm} is given in 
\ref{La1}--\ref{La13}.

\begin{sbpara}\label{La1} Let $\La$ be as in \ref{La0} and let $X$ and $H$ be as in \ref{La2}. 
  Let $V$ be an open set of $X^{\log}_{[:]}$ on which we are given an isomorphism
$(H_\Z, W, (\langle \cdot, \cdot\rangle_w)_w)\cong (H_0, W,  (\langle \cdot, \cdot\rangle_w)_w)$. 

We define a map $\psi_V: V\to D_{\SL(2)}$ which will be compatible with the map  $\psi$ in Theorem \ref{L*thm}. 
Let $\tilde p\in V$ and let $p$ be the image of $\tilde p$  in $X_{[:]}$. We define  $\psi_V(\tilde p)\in D_{\SL(2)}$ to be   the class of  $(\rho, \varphi, \br)$ on $H_{\R,\tilde p}= H_0$ in Theorem \ref{oldSL20}. 
\end{sbpara}

\begin{sbprop}\label{La12}

For every open set $U$ of $D^I_{\SL(2)}$  and for every real analytic function $f$ on $U$, $f\circ \psi_V$ belongs to  ${\hat A}_{(2+),X}^{\an, \log*}(\psi_V^{-1}(U))$.
\end{sbprop}

\begin{sbpara} This is a preparation for the proof of Proposition \ref{La12}.

Let $\tilde p\in V$ and let $p$ be the image of $\tilde p$  in $X_{[:]}$. We work locally at $\tilde p$. Let $\Phi=\Phi(p)=\{W^{(j)}\;|\;1\leq j\leq m\}$ and let $\tau:{\bf G}_{m,\R}^{\Phi} \to \Aut_\R(H_{\R,\tilde p})= \Aut_\R(H_0)$ be as in Theorem \ref{oldSL20} (2). Let $D_{\SL(2)}(\Phi)$ be the open set of $D_{\SL(2)}$ consisting of points  whose associated set of weight filtrations are contained in $\Phi$. Let $D_{\SL(2)}(\Phi)'=D_{\SL(2)}(\Phi)$ 
(resp.\ $D_{\SL(2)}(\Phi)'=D_{\SL(2)}(\Phi)_{\text{nspl}}$)
 if $W\notin \Phi$ (resp.\ $W\in \Phi$). %
Here $D_{\SL(2)}(\Phi)_{\text{nspl}}$ is the open set of $D_{\SL(2)}$ in \cite{KNU2} 3.2. Replacing $V$ by some open neighborhood of $\tilde p$, we may assume that the image of $V$ in $D_{\SL(2)}$ is containd in $D_{\SL(2)}(\Phi)'$. 
Let $D'=D\cap D_{\SL(2)}(\Phi)'$.

  We briefly review the definition of the sheaf of real analytic functions on  $D_{\SL(2)}(\Phi)'$. 

We have an injection $$(*)\quad  D_{\SL(2)}(\Phi)' \to \R^{\Phi}\times D \times\spl(W) \times {\textstyle\prod}_{j=1}^m \spl(W^{(j)}(\gr^WH_0)) $$
 (\cite{KNU2} Proposition 3.2.7 (i)).  
Here $\spl(W)$ is the real analytic manifold of all splittings of $W$ of 
$H_0$ and $\spl(W^{(j)}(\gr^WH_0))$  is the real analytic manifold of 
all splittings of $W^{(j)}$ on $\gr^WH_0$. To define the maps from $D_{\SL(2)}(\Phi)'$ to  $\R^{\Phi} \times D$, we fix a real analytic map $\beta:D'\to \R_{>0}^{\Phi}$ 
such that $\beta(\tau(a)F')= a\beta(F')$ for all $a\in \R_{>0}^{\Phi}$ and $F'\in D'$. 
Then 
the map 
  $\beta$ 
  extends canonically to a map $\beta: D_{\SL(2)}(\Phi)'\to \R_{\geq 0}^{\Phi}\subset \R^{\Phi}$  and the map $\nu: D'\to D\;;\; F' \mapsto \tau(\beta(F'))^{-1}F'$ 
extends canonically to a map $\nu: D_{\SL(2)}(\Phi)'\to D$ (\cite{KNU2} Proposition 3.2.7 (i)). 
By definition, an $\R$-valued function on an open set of  $D_{\SL(2)}(\Phi)'$ is a real analytic function if and only if locally, it is the restriction of a local section of the sheaf of real analytic functions on the right-hand-side of the above $(*)$.
\end{sbpara}

\begin{sbpara}  We prove Proposition \ref{La12}. 

It is sufficient to prove that the maps (i) $V\to \R^{\Phi}$, (ii) $V\to D$, 
(iii) $V\to \spl(W)$, and (iv) $V\to \spl(W^{(j)}(\gr^WH_0))$ obtained by the composition of $\psi_V$ and the above map $(*)$  have the property that 

\medskip

(A) the pullbacks of  real analytic functions by these maps belong to $\hat A^{\an,\log*}_{(2+)}$. 

\medskip

The maps (iii) and (iv) have the  property (A) by Theorems \ref{SL2II} and \ref{SL2I}, respectively, because these maps are nothing but $\spl_W$ and 
$\spl_{W^{(j)}(\gr^WH_0)}^{\BS}$, respectively. We consider the maps (i) and (ii). Let $F$ be the Hodge filtration of $H$. Define $t(y)$ as in \ref{tortwi} by using $\tau$. 
Let $t(y)^{-1}: V \to D'$ be the map $p'\mapsto (t(y)^{-1}F)(p')$. This map has the  property (A) by Theorem \ref{HA1} (1) and \ref{tortwi}.

Next, the map (i) has the  property (A) because it is the product of the composition  $V \overset{t(y)^{-1}}\to D' \overset{\beta}\to \R_{> 0}^{\Phi}$ and $(\sqrt{y_{j+1}/y_j})_{1\leq j\leq m}$.

Lastly, the map (ii) has the  property (A) by the following commutative diagram:
$$\begin{CD} V @> {\psi_V}>> D_{\SL(2)}(\Phi)'\\
@V t(y)^{-1} VV @VV \nu V \\
D'@> {\nu} >> D.
\end{CD}$$
\end{sbpara}

\begin{sbpara}\label{La13} We complete the proof of Theorem \ref{L*thm}. The space $X$ is covered by open sets $V$ on which the $\Gamma$-level structure comes from an isomorphism $(H_\Z, W, (\langle \cdot, \cdot\rangle_w)_w)\cong (H_0, W,  (\langle \cdot, \cdot\rangle_w)_w)$. We have $\psi_V:V\to D_{\SL(2)}$, and the compositions $V\to D_{\SL(2)}\to \Gamma \bs D_{\SL(2)}$ glue to a map $X^{\log}_{[:]}\to \Gamma \bs D_{\SL(2)}$ which has the property stated in Theorem \ref{L*thm}.
\end{sbpara}

\begin{sbpara}

For a weak fan $\Sig$  (\cite{KNU3} 2.2.3)  for $\La$ which  is strongly compatible with $\Gamma$ (\cite{KNU3} 2.2.6), let
  $\Gamma\bs D_{\Sig}\supset \Gamma \bs D$
 be the space of nilpotent orbits (\cite{KNU3} Section 2). This $\Gamma \bs D_{\Sig}$ is a moduli space of LMH with a $\Gamma$-level structure.

\end{sbpara}

\begin{sbthm}\label{partR2} Let the notation be as above. Then the identity map of $\Gamma\bs D$ extends uniquely to a continuous map $$\psi: (\Gamma\bs D_{\Sig})_{[:]}^{\log}=\Gamma \bs D_{\Sig,[:]}^{\sharp} \to \Gamma \bs D^I_{\SL(2)}$$ 
(the CKS map) having the following property$:$ For every open set $U$ of $\Gamma\bs D^I_{\SL(2)}$ and for every real analytic function $f$ on $U$, $f\circ \psi$ belongs to  ${\hat A}_{(2+),X}^{\an, \log*}(\psi^{-1}(U))$.

\end{sbthm}
This follows from Theorem \ref{L*thm} applied to the case where $X=\Gamma\bs D_{\Sig}$ and $H$ is the universal object on it.

\begin{sbpara} This gives a new proof of the continuity of the CKS map $D^{\sharp}_{\Sig,[:]}\to D^I_{\SL(2)}$  proved in \cite{KU}, \cite{KNU3}, \cite{KNU4}.

\end{sbpara}

\section{Higher direct images (plan and a special case)}\label{s:HDI}

In Section \ref{ss:hdi}, we explain basic things about higher direct images and our plan for a future study. In Section \ref{ss:curve}, we consider the case of relative dimension one.

We consider $\Z$-PLH and $\Z$-LMH, and call them just PLH and LMH, respectively.

\subsection{Higher direct images}\label{ss:hdi}

\begin{sbpara}\label{12}

Let $f:X\to Y$ be a projective log smooth saturated morphism of fs log analytic spaces. Let $H'$  be an LVMH on $X$ with polarizable $\gr^W$. We expect that for each $m$, $H_{\Z}:=R^mf^{\log}_*H'_\Z$, which is a locally constant sheaf of finitely generated $\Z$-modules by \cite{KN},  underlies an LVMH $H$ on $Y$ with polarizable $\gr^W$. If  $f$ is vertical and $H'$ is an LVPH of weight $w$, 
we expect that this $H$ is an LVPH on $Y$ of weight $w+m$.

\end{sbpara}

\begin{sbpara}\label{123}

  In a forthcoming paper \cite{KNU7}, we  plan to study the subject in \ref{12}.
Section \ref{s:newa} and Section \ref{ss:curve} will serve as the first steps for our study of higher direct images.

Let $f: X\to Y$ be a projective log smooth saturated vertical morphism, let $H'$ be an LVPH on $X$ of weight $w$,  and let $H$ be the $m$-th higher direct image of $H'$ which we wish to show to be an LVPH of weight $w+m$. Our method is, very roughly speaking,  as follows.

(1) We first consider the $\SL(2)$-orbits  associated to $H'$ which appear on $X_{[:]}$. This is discussed in Section \ref{ss:HA}. 

(2) Before we prove that $H$ is an LVPH, we show that  the $\SL(2)$-orbits associated to $H$ appear on $Y_{[:]}$, by using (1). 

(3) We deduce that $H$ is an LVPH from (2). 
\end{sbpara}

\begin{sbpara}\label{PLH2} In \ref{12}, 
consider the case where $X\to Y$ is vertical and  $H'$ is the Hodge structure $\Z$. 

Then for each $m$, we have a pre-PLH $H=H^m(X/Y)$ of weight $m$
 whose $H_\Z$ and $H_{\cO}$ are given by $H_\Z=R^mf^{\log}_*\Z$  and  $H_{\cO}=H^m_{\dR}(X/Y):=R^mf_*(\omega_{X/Y}^{\bullet})$
 and whose $\langle\cdot,\cdot\rangle$ is defined by using a polarization of $X$.  

By \cite{IKN} Theorem  (6.3) (3), we have an isomorphism 
$$\iota: \cO_Y^{\log} \otimes_{\Z} H_\Z \cong \cO_Y^{\log} \otimes_{\C} H_{\cO}.$$  
Thus we have  $H^m(X/Y):=(H_\Z, H_{\cO}, F, \iota,\langle\cdot,\cdot \rangle)$. 

  By \cite{FN} Theorem 2.5 (1) and (2), the Hodge to de Rham spectral sequence 
$$E_1^{pq}= R^qf_*(\omega^p_{X/Y}) \Rightarrow E^m_{\infty}=H^m_{\dR}(X/Y)$$ degenerates at $E_1$ and $R^qf_*(\omega^p_{X/Y})$ are locally free and commute with base changes. In particular,  the map $F^p:= R^mf_*(\omega_{X/Y}^{\bullet \geq p})\to H_{\cO}$ is injective and this filtration $(F^p)_p$ splits locally. 

The property $\langle F^p, F^q \rangle=0$ for $p+q>m$ is reduced to $F^{d+1}H^{2d}(X/Y)=0$, where $d$ is the relative dimension of $X$ over $Y$, and  hence to $\omega_{X/Y}^{\bullet,\geq d+1}=0$. 

Thus $H$ is a pre-PLH.

Furthermore, $H$ satisfies  the big Griffiths transversality by \ref{PLH3} below. 

The annihilator of $F^p$ for $\langle\cdot, \cdot\rangle$ is $F^{w+1-p}$ for every $p$. This is proved as follows. By using $\langle F^p, F^{w+1-p}\rangle =0$ and by using the fact that $R^qf_*(\omega^p_{X/Y})$ commute with base changes, we are reduced to the case where $Y$ is the standard log point. In this case, this follows from the Serre duality \cite{T0} Theorem 2.21.

It is known  in certain cases that $H$ is an LVPH  of weight $m$ on $Y$. For example, this is known in the following cases (i) and (ii). 

(i) The case where $Y$ is log smooth (\cite{KMN}).

(ii) The case where $Y$ is of log rank $\leq 1$ (\cite{FN}). 
\end{sbpara} 

\begin{sbpara}\label{PLH3} We prove that the above $H^m(X/Y)$ satisfies the big Griffiths transversality. Since $X\to Y$ is log smooth, we have a locally splitting exact sequence 

(1) $0\to \cO_X \otimes_{\cO_Y}\omega^1_Y \to \omega^1_X \to \omega^1_{X/Y}\to 0.$

From (1), we obtain an exact sequence

(2) $0\to \omega^{\bullet}_{X/Y}[-1]\otimes_{\cO_Y} \omega^1_Y \to \omega_X^{\bullet}/\text{fil}^2 \to \omega^{\bullet}_{X/Y}\to 0,$

\noindent
where $\text{fil}^r$ is the image of $\omega_Y^r \otimes_{\cO_Y}\omega_X^{\bullet}[-r]\to \omega_X^{\bullet}$. 
By the quasi-isomorphism $\cO_Y^{\log} \to \omega_{X/Y}^{\bullet,\log}$ (\cite{O} Chapter V, Theorem 3.3.4; the assumption there is satisfied by \cite{IKN} Lemma (A.4.1)) whose right-hand-side appears in the exact sequence ($\cO_X^{\log} \otimes_{\cO_X}$ of the exact sequence  (2)), we see that the connection $\nabla: H_{\cO}\to \omega^1_Y \otimes_{\cO_Y} H_{\cO}$ coincides with the connecting map of $Rf_*$  of the exact sequence (2). 
The exact sequence (2) gives an exact sequence

$0\to  \omega^{\geq p-1}_{X/Y}[-1]\otimes_{\cO_Y} \omega^1_Y \to \omega_X^{\geq p}/\text{fil}^2 \to \omega^{\geq p}_{X/Y}\to 0,$

\noindent
and hence we have the big Griffiths transversality. 
\end{sbpara}

\begin{sbpara}\label{sp(t)} By \ref{PLH2} and \ref{PLH3} and by the definition of LVPH in  \ref{LMH3}, $H^m(X/Y)$ is an LVPH of weight $m$ on $Y$ if and only if for every $y\in Y$, $H=H^m(X_y/y)$ satisfies the following condition (1). Here $X_y$ is the fiber of $y$ in $X$.  Let $\cT$ be the canonical toric variety associated to the fs log point $y$. 

\medskip

(1) If $s\in \cT_{\triv}$   is sufficiently near to $y$ in $ \cT$,  $H_{\cT}(s)$  is a polarized Hodge structure of weight $m$.

\medskip

By \cite{FN}, this condition  for $y\in Y$  are satisfied if the rank of $(M^{\gp}_Y/\cO_Y^\times)_y$ is $\leq 1$.

\medskip
By the properness of $\cT_{[:]}\to \cT$,  the condition (1) is equivalent to the following condition (2).

\medskip

(2) For each $p\in Y_{[:]}$, if $s\in \cT_{\triv}$ is sufficiently near to $p$ in $\cT_{[:]}$,  $H_{\cT}(s)$ is a polarized Hodge structure of weight $m$.

Thus we can check whether $H_{\cT}(s)$ is a PLH by considering it locally on the space $Y_{[:]}$ of ratios. This is the basic spirit of the step (3) of \ref{123}. 
This condition  (2)  will be useful in Section \ref{ss:curve}. 

\end{sbpara}

\begin{sbpara}\label{1235}

In Section \ref{ss:curve},  in the case where $X\to Y$ is a relative curve, we prove that $H^1(X/Y)$ is an LVPH. The step (1) in \ref{123}  is trivial in this situation. Proposition \ref{I1} and the proof of Proposition \ref{p:basic} given in \ref{S5HA2} are based on  the  spirit of  Steps (2) and (3) in \ref{123}, respectively.
\end{sbpara}

\subsection{Case of relative curves}\label{ss:curve}

In this section (Section \ref{ss:curve}), let $f:X\to Y$ be a 
  projective, vertical, saturated and log smooth morphism of fs log analytic spaces whose every fiber is
one-dimensional. Our aim is to prove the following Proposition \ref{p:basic}
by using 
our theory of integration and by using the geometry of the space of ratios.

\begin{sbprop}
\label{p:basic}
  
  The pre-PLH  $H^1(X/Y)$ of weight $1$ on $Y$ ($\ref{PLH2}$) endowed with the Hodge filtration and the bilinear form defined by the cup product 
   is an LVPH. 
\end{sbprop}

  This gives a partial affirmative answer to the main conjecture in \cite{FN}.

  There is another proof of Proposition \ref{p:basic} of more classical type as described in \ref{clpf}. 
  But the proof we give  here is in the spirit of \ref{123} and of the forthcoming paper \cite{KNU7}. In our approach, the positivity of the Hermitian form of $H^1(X/Y)$ giving the Hodge metric is reduced to the positivity of the local integral.

\medskip

For the proof of Proposition \ref{p:basic}, by \ref{sp(t)}, we may assume that the base $Y$ is an fs log point. We assume this in the proof of Proposition \ref{p:basic} which is given in \ref{H1a}--\ref{S5HA2}.

\begin{sbpara}\label{H1a} 
 Denote the pre-$\Z$-PLH $H^1(X/Y)$ of weight $1$ on $Y$ by $H$.
Let 
$$V=\Gamma(X, \omega^1_{X/Y}).$$

 Let $\cT$ be the canonical toric variety associated to the fs log point $Y$.

Consider the homomorphism $H_{\cO} \to \cO^{\log}_Y\otimes_{\Z} H_{\Z}$ on $Y^{\log}$. 
For $s\in \cT_{\triv}$, this induces a homomorphism $$a_s: V \subset H^1_{\dR}(X/Y)=H_{\cO}\overset{\cong}\to  H_{\cT,\Z,s}\otimes_{\Z} \C$$ 
at $s$. Let  $$\bar a_s: V \to H_{\cT,\Z,s}\otimes_\Z \C$$ be the composition of $a_s$ and  the complex conjugation on $H_{\cT,\Z, s}\otimes_\Z \C$. 

Using the cup product $\cup: H_{\Z}\times H_{\Z}=R^1f^{\log}_*\Z \times R^1f^{\log}_*\Z\to R^2f^{\log}_*\Z\cong \Z$, define the Hermitian form  
$$(\cdot,\cdot)_s : V \times V\to \C$$
by  $$(u, v)_s:= ia_s(u) \cup {\bar a}_s(v)\quad \text{for}\; u,v\in V.$$

\medskip

 By \ref{sp(t)}, to prove that $H$ is an LVPH, it is enough to prove that  for each $p\in Y_{[:]}$,  the conditions (i) and (ii) below are satisfied if $s\in \cT_{\triv}$ is sufficiently near to $p$ in $\cT_{[:]}$.

\medskip

(i) $(\cdot, \cdot)_s$ is positive definite.

(ii) The map $(a_s, {\bar a}_s): V\oplus V \to H_{\cT,\Z,s}\otimes_{\Z}\C$ is an isomorphism. 

\medskip

Actually, (ii) follows from (i). First the injectivity of the map in (ii) is shown as follows. If $(u,v)\in V \oplus V$ is in the kernel, ${\bar a}_s(v)=-a_s(u)$ and hence $(v,v)_s= ia_s(v) \cup {\bar a}_s(v)=-i a_ s(v) \cup a_s(u)= 0$ by $F^1\cup F^1=0$ (\ref{PLH2}). Hence  $v=0$ by (i)  and hence $u=0$. Since the annihilator of $V$ for the bilinear form $\cup$ on $H_{\cT, \Z, s}\otimes_{\Z} \C$  is $V$ (\ref{PLH2}), we have $\dim_\C (H_{\cT, \Z, s}\otimes_{\Z}\C) \leq  2\dim_{\C}(V)$. Hence the injectivity of the map $(a_s, {\bar a}_s)$ implies its bijectivity. 

\end{sbpara}

\begin{sbpara}\label{H1b} We denote the composite map  $V\to \cO_Y^{\log}\otimes_{\Z} H_\Z\to A^{\an,\log}_{3, Y, \C} \otimes_{\Z} H_\Z$ given on $Y^{\log}_{[:]}$  by $a$. 
Let $\bar a: V\to 
A^{\an,\log}_{3,Y, \C} \otimes_\Z H_\Z$ 
be the complex conjugate of $a$. 
Let $$(\cdot, \cdot): V \times V \to A^{\an,\log}_{3,Y,\C}$$ be the Hermitian form
$(u, v)\mapsto i a(u) \cup {\bar a}(v)$. The image of this map is in $\tau_*A^{\an,\log}_{3,Y,\C}= A^{\an}_{3,Y, \C}$. 

If $s\in \cT_{\triv}$ is sufficiently near to $p$ in $\cT_{[:]}$, we have 
$$(u, v)_s=  (u, v)(s),$$
where the right-hand-side is the value at $s$ of the extension of $(u, v)\in A^{\an,\log}_{3, Y,\C}$ to $A^{\an, \log}_{3, Y, \cT,\C}$ on an open neighborhood of $p$. 

\end{sbpara}

Now we consider the integration.

\begin{sblem}\label{H1c} 
We denote $\int_{X_{[:]}/Y_{[:]}}$ (Section $\ref{ss:2nint}$) by $\int$ for simplicity.   
For $u, v\in V$, let $$I(u,v):=i\int u \wedge \bar v\in A_{3,Y,\C}.$$
Then we have
$$(u,v)=I(u,v) \quad\text{in}\;\;A^{\log}_{3,Y,\C}.$$

\end{sblem}

\begin{pf} The composite map $V \overset{a}\to A^{\an, \log}_{3,Y,\C} \otimes_{\Z} R^1f_{[:]*}^{\log}\Z\to A^{\log}_{3,Y,\C} \otimes_{\Z} R^1f_{[:]*}^{\log}\Z$  
coincides with the composite map $$V \to  
f^{\log}_{[:]*}(A^{1,\log}_{3,X/Y,d=0,\C})\to A^{\log}_{3,Y,\C}\otimes_{\Z} R^1f^{\log}_{[:]*}\Z,$$
where the last arrow is that in Section \ref{ss:intcoh2}.

Because  the quasi-isomorphism $A^{\log}_{3,Y}\to A^{\bullet,\log}_{3,X/Y}$ is compatible with the product structures,  
we have a commutative diagram 
$$\begin{matrix}  f^{\log}_{[:]*}(A^{1,\log}_{3,X/Y,d=0})\times  f^{\log}_{[:]*}(A^{1,\log}_{3,X/Y,d=0})& \to&  f^{\log}_{[:]*}(A^{2,\log}_{3,X/Y})\\
\downarrow && \downarrow \\
A^{\log}_{3,Y}\otimes_{\Z} R^1f^{\log}_{[:]*}\Z \times A^{\log}_{3,Y}\otimes_{\Z} R^1f^{\log}_{[:]*}\Z  &\to& A^{\log}_{3,Y}\end{matrix}$$
in which the left vertical arrow is the integration in Section \ref{ss:intcoh2}, the right vertical arrow  is the integration  in Section \ref{ss:2nint}, the upper horizontal arrow is the wedge product, and the lower horizontal arrow is the cup product. 

These prove Lemma \ref{H1c}.
\end{pf}

\begin{sbpara}\label{511} Fix $p\in Y_{[:]}$.

We define a decreasing filtration $(V^j)_{1\leq j\leq m+2}$ of $V$ as follows.

Let the filtration $M_{Y,p}=M^{(0)}\supsetneq M^{(1)}\supsetneq \dots \supsetneq M^{(m)}=\cO_{Y,p}^\times$ on $M_{Y,p}$ and $f_{j,k}\in M_{Y,p}$ be as in \ref{m,f_jk} defined replacing $X$ there by the present $Y$.

 Let $y_{j,k}=-(2\pi)^{-1}\log |f_{j,k}|\in A_{3,Y,p}^{\an}$.  Put $y_{m+1,1}=1$, $y_{m+2,1}=0$. 

Let $p'\in X_{[:]}$ lying over $p$. Let $q$ be a local coordinate function of $X$ over $Y$ at $p'$  (Theorem $\ref{SEisS}$) and let $y=-(2\pi)^{-1}\log |q|\in A_{3,X,p'}^{\an}$. A local coordinate function $q$ 
 is not unique,  but  for each $1\leq j\leq m+1$, whether the value of $y/y_{j,1}$ at $p'$ is $0$ or $\infty$ or in $(0, \infty)$ is independent of the choice of $q$.

Let $V^j$ be the set of all elements $w\in V$ such that if $p'\in X_{[:]}$ lies over $p$ and if $w$ is a base of $\omega^1_{X/Y,p'}$ at $p'$, then the local coordinate function of $X$ over $Y$ at $p'$ satisfies the condition that the value of $y/y_{j,1}$ at $p'$ is not $\infty$.

Then $V^j$ is a $\C$-subspace of $V$. 
We have $V^1=V$ and  $V^{m+2}=0$.

\end{sbpara}

\begin{sbpara}

For example, assume that $Y$ is the standard  log point, $X$ is semi-stable over $Y$, and $p$ is the unique point of $Y_{[:]}=Y$. Then $$V=V^1\supset V^2\supset V^3=0,$$ and $V^2$ is the subset of $\Gamma(X, \omega^1_{X/Y})$ consisting of all elements $w$ such that at every singular point of $X$,  $w$ is not a local base of $\omega^1_{X/Y}$. That is, $V^2$ is the image of $\Gamma(X, \Omega^1_{X/Y})$.  

\end{sbpara}

\begin{sbprop}\label{I1} $(1)$ Let $u, v\in V$. If $1\leq j\leq m+2$ and if either one of $u$, $v$ belongs to $V^j$,  then $I(u,v)\in y_{j,1}A_{2,Y, p,\C}$. 

$(2)$ Let $1\leq j\leq m+1$ and consider the Hermitian form $V^j \times V^j \to \C$  which sends $(u,v)\in V \times V$ to the image of $y_{j,1}^{-1}I(u,v)\in A_{2,Y, p,\C}$  
under the evaluation $A_{2, Y, p, \C}\to \C$ at $p$. Then this Hermitian form induces a positive definite Hermitian form $V^j/V^{j+1}\times V^j/V^{j+1}\to \C$. 

$(3)$ Take a splitting $V= \bigoplus_j V_{(j)}$ of the filtration $(V^j)_j$ ($V_{(j)}\subset V^j$ and $V_{(j)}\overset{\cong}\to V^j/V^{j+1}$). Then we have$:$

$(3.1)$ The image of the Hermitian form 
$$I(\cdot,\cdot)': V \times V\to A_{3, Y, p}\;;\; (u, v) \mapsto {\textstyle\sum}_{j,k} y_{j,1}^{-1/2}y_{k,1}^{-1/2}I(u_{(j)}, v_{(k)})\in A_{3,Y, p,\C}$$ 
is contained in $A_{2, Y,p, \C}$.
Here $u_{(j)}$ denotes the $V_{(j)}$-component of $u$.

$(3.2)$ Let the Hermitian form $I(\cdot, \cdot)'(p): V \times V \to \C$ be the composition of the Hermitian form in $(3.1)$ and the evaluation 
 $A_{2, Y, p, \C}\to \C$ at $p$.
 Then this Hermitian form is positive definite. We have $I(V_{(j)}, V_{(k)})'(p)=0$ if $j\neq k$. 
\end{sbprop}

\begin{pf}  
Let $p'\in X_{[:]}$ be a point lying over $p$, and let $q\in M_{X,p'}$ be a local coordinate function  (Theorem \ref{SEisS}) of $X$ over $Y$ at $p'$. Let $u\in \omega^1_{X/Y, p'}$. 
Let $q=e^{2\pi i(x+iy)}$ with $x$ and $y$ being real. 
  Since $\omega^1_{X/Y}$ has the base $d\log(q)=2\pi i (dx+idy)$ at $p'$, we have the following (i), (ii), and (iii).

(i) Let $u, v\in \omega^1_{X/Y, p'}$. Then $u\wedge \bar v\in A_{2,X, p',\C} dx\wedge dy$.

(ii) Let $u,v\in \omega^1_{X/Y,p'}$. Write $q-q(p')$ by $q_0$. (Hence $q_0\in \cO_{X,p'}$ and  has value $0$ at $p'$.) 
  Unless both $u$ and $v$ are bases of $\omega^1_{X/Y, p'}$, 
$u\wedge \bar v$ belongs to $(A_{2,X, p', \C}q_0+ A_{2, X, p', \C}\bar q_0)dx\wedge dy$.

(iii) Let $u\in \omega^1_{X/Y,p'}$. Then $iu \wedge \bar u=gdx\wedge dy$ for some $g\in A_{2,X, p'}$ satisfying  (P) (\ref{pos}). If $u$ is a base of  $\omega^1_{X/Y,p'}$, we have $g(p')>0$.

\medskip

We prove (1).  Take an open covering $(U_{\lam})_{\lam}$ of $X_{[:]}$ such that for any $\lam$ there is a local coordinate function $q_{\lam}\in M(U_{\lam})$  for any $p' \in U_{\lam}$ (Theorem \ref{SEisS}). 
  We may assume that the index set is finite, since $X \to Y$ is projective.
  Then there is a partition of unity $1=\sum s_{\lam}$, $s_{\lam} \in A_{1,X}$, %
  by Theorem \ref{softhm}.
  By Proposition \ref{tildeJ}, it is sufficient to prove that  %
$y_{j,1}^{-1}s_{\lam}w$ belongs to $\tilde \cE$ there, where $w=u\wedge \bar v$. 

  Fix a $\lam$ and write $q_{\lam}$ as $q$. 
  Let $q=e^{2\pi i(x+iy)}$ with $x$ and $y$ being real. 
  Let $p'\in U_{\lam}$. 

Assume first that $y/y_{j,1}$ does not have the value $\infty$ at $p$. Then around $p$, by the above (i), we have
$y_{j,1}^{-1}s_{\la}w = 
(y_{j,1}^{-1}y)\cdot y\cdot s_{\la}fdx/y \wedge dy/y$ for some $f\in A_{2,\C}$. 
  Since $y_{j,1}^{-1}y \in A_2$, this presentation shows that $y_{j,1}^{-1}s_{\la}w$ belongs to $\tilde \cE$ around $p'$ (we take $ y_{j,1}^{-1}$ as $h_1$ in the definition of $\tilde \cE$ in Proposition \ref{tildeJ}). 

Next assume that $y/y_{j,1}$ has value $\infty$ at $p$. Then either $u$ or $v$ is not a base of $\omega^1_{X/Y,p'}$. Take a $c<0$. Then  by the above (ii), around $p'$,
 we have $y_{j,1}^{-1}s_{\la}w=(y_{j,1}^{-1}y^c)\cdot y\cdot s_{\la}y^{1-c}(fq_0+g{\bar q_0})dx/y\wedge dy/y$ for some $f, g \in A_{2,\C}$. 
   Since $y_{j,1}^{-1}y^c \in A_2$ and $y^{1-c}q, y^{1-c}\bar q\in A_{2,\bC}$, 
this presentation shows that $y_{j,1}^{-1}s_{\la}w$ belongs to $\tilde \cE$ around $p'$ (again we take $y_{j,1}^{-1}$ as $h_1$ in the definition of $\tilde \cE$ in Proposition \ref{tildeJ}). 
 
We prove (2). Let the notation be the same as in the proof of (1). 
  We choose the partition of unity to be with (P) by Theorem \ref{softhm}. 
  Let $u \in V^j\smallsetminus V^{j+1}$. 
  We apply Proposition \ref{tildeJ2}  to $w=iy_{j,1}^{-1}s_{\lam}u \wedge \bar u$. 
By  Proposition \ref{tildeJ2} (1) and by the first part of the above (iii), $\int w$ satisfies (P) (\ref{pos}). Since $u\notin V^{j+1}$, Proposition \ref{3.1fiber} shows that there exist $\la$ and $p'\in U_{\lam}$ lying over $p$ such that $u$ is a base of $\omega^1_{X/Y,p'}$,  $y/y_{j,1}\in A_2$ is invertible at $p'$, and  $s_{\lam}(p')>0$. Hence by Proposition \ref{tildeJ2} (2) and by the second part of the above (iii), $y_{j,1}^{-1}\int iu\wedge \bar u$ has value $>0$ at $p$. This proves (2). 
 
(3) follows from (1) and (2).
\end{pf}

 \begin{sbpara}\label{S5HA2} We now complete the proof of Proposition \ref{p:basic}. 
It is enough to prove that the condition (i) in \ref{H1a} is satisfied. 

Let $p\in Y_{[:]}$. Let the filtration $(V^j)_j$ on $V$ be as in \ref{511}.  As in Proposition \ref{I1} (3),  fix a splitting 
$V=\bigoplus_j V_{(j)}$ of the filtration $(V^j)_j$. Let  $y_{j,1}$ be as in \ref{511}.

Define
$$(\cdot,\cdot)'_s:  V \times V \to \C\;;\; (u,v)\mapsto {\textstyle\sum}_{j,k}\;  y_{j,1}(s)^{-1/2}y_{k,1}(s)^{-1/2}(u_{(j)}, v_{(k)})_s.$$

If the image of $s\in \cT_{\triv}$ is sufficiently near to $p$ in $\cT_{[:]}$, then by Lemma \ref{H1c}, 
$(u,v)'_s$ ($u,v\in V$) is the value at $s$ of the extension of $I(u,v)'\in A^{\an}_{3, Y}$ to $A^{\an}_{3,Y,\cT}$. By Proposition \ref{I1} (3), this local section of $A^{\an}_{3,Y}$ belongs to $A_{2,Y}$ and hence belongs to $A^{\an}_{2,Y}$ by Proposition \ref{pt13}. 
\

By this and by Proposition \ref{A2cont}, when $s\in \cT_{\triv}$ converges to $p$ in $\cT_{[:]}$, the Hermitian form $(\cdot, \cdot)'_s$ on $V$ converges to the Hermitian form $I(\cdot,\cdot)'(p)$ on $V$ in Proposition \ref{I1} (3). Since $I(\cdot,\cdot)'(p)$ is positive definite (Proposition \ref{I1} (3)), 
 if  $s\in \cT_{\triv}$ is sufficiently near to $p$, the Hermitian form $(\cdot, \cdot)'_s$ is positive definite and hence so is $(\cdot, \cdot)_s$.

\end{sbpara}

\begin{sbpara} 
Assume that $Y$ is an fs log point and let $p\in Y_{[:]}$. Then the filtration $(V^j)_j$ ($V=\Gamma(X, \omega^1_{X/Y})$) is understood as follows in terms of the LVPH $H=H^1(X/Y)$. The relative monodromy filtrations $W^{(j)}$ ($1\leq j\leq m$) on $H_\Q$ at $p$ induce relative monodromy filtrations $W^{(j)}$ on $H^1_{\dR}(X/Y)=\tau_{Y*}(\cO_Y^{\log}\otimes_\Q H_\Q)$ ($\tau_Y: Y^{\log}_{[:]}\to Y_{[:]}$). They satisfy $W^{(j)}_{-1}=0$, $W^{(j)}_2=H^1_{\dR}(X/Y)$, $W^{(j)}_0$ and $W^{(j)}_1$ are the annihilators of each other for the polarization, and 
$$W^{(1)}_0\subset W^{(2)}_0 \subset \dots \subset W^{(m)}_0\subset W^{(m)}_1\subset\dots \subset W^{(2)}_1\subset W^{(1)}_1.$$ We have
$$V^j=V \cap W^{(j)}_1.$$ 
This can be proved by using the norm estimate   \cite{CKS1} Theorem (5.21)
(the relation of the relative monodromy filtrations and the growth of Hodge metric in degeneration) and the estimate of the Hodge metric obtained in Proposition \ref{I1}. We omit the details of the proof. 

\end{sbpara}

\begin{sbrem}

In the classical Hodge theory, harmonic forms and their integrals play important roles to have  the Hodge structure associated to a non-singular  space. In the forthcoming paper \cite{KNU7}, to generalize the result Proposition \ref{p:basic} to the result on higher direct images for the higher relative dimension, we will discuss log harmonic forms and their integrals in log geometry.

\end{sbrem}

\begin{sbpara}\label{clpf}

Another proof of the main result Proposition \ref{p:basic} of this Section \ref{ss:curve} by using the reduction to the case where 

\medskip

(+) $X$ is semi-stable over a standard log point $Y$

\medskip
\noindent
is as follows. By \cite{IKN}, we are reduced to the case where $Y$ is an fs log point. 
For each element $N$ of the interior of the monodromy cone of the base, 
we can show that on $H^1$, we have $N^2=0$ and $H^1_\Q\to \Q\; ;\; a \mapsto  Na\cup a $ is semi-positive definite.
  This is reduced to the case (+) by using a suitable morphism $Y'\to Y$ from a standard log point and 
by using a variant of the semi-stable reduction theorem (cf.\ \cite{FN} 5.7).
  For two elements $N_1$ and $N_2$ of the monodromy cone of the base, by the semi-positivity, we have $\Ker(N_1+N_2)= \Ker(N_1)\cap \Ker(N_2)$. Since $tN_2=N_1+N_3$ for some $t>0$ and for some element $N_3$ of the interior of the monodromy cone, we have that $\Ker(N)$ and $\text{Im}(N)$ ($=$ the annihilator of $\Ker(N)$ for the cup product pairing) are independent of $N$. By the reduction to the case (+), we can show that $\Ker(N)/\text{Im}(N)$ is a Hodge structure of weight $1$.  From this, we can see that $H^1$ is a PLH of weight $1$. 

Our proof given above using our integration is essentially different from this. 
\end{sbpara}

\noindent {\rm Kazuya KATO
\\
Department of mathematics
\\
University of Chicago
\\
Chicago, Illinois, 60637, USA}
\\
{\tt kkato@uchicago.edu}

\bigskip

\noindent {\rm Chikara NAKAYAMA
\\
Department of Economics 
\\
Hitotsubashi University 
\\
2-1 Naka, Kunitachi, Tokyo 186-8601, Japan}
\\
{\tt c.nakayama@r.hit-u.ac.jp}

\bigskip

\noindent
{\rm Sampei USUI
\\
Graduate School of Science
\\
Osaka University
\\
Toyonaka, Osaka, 560-0043, Japan}
\\
{\tt usui@math.sci.osaka-u.ac.jp}

\end{document}